# A note on the counterfeit coins problem

Li An-Ping

Beijing 100085, P.R.China
apli0001@sina.com

Abstract

In this paper, we will present an algorithm to resolve the counterfeit coins problem in the case that the number of false coins is unknown in advance.

The counterfeit coins problem is a well-known combinatorial search problem, namely, search for the counterfeit coins in a number of coins by a balance, where the fakes and the normals have same semblance but different weights. The problem has several versions, for instance, the fakes are known lighter or heavier than the normals, and the number of the fakes is assumed known. For the detail results and researches in this subject refer to see the papers [1]~[7].

In this paper, we will consider a more general case that the number of the fakes is unknown beforehand. In general, it is assumed that the coins will be permitted to be remarked by numbers in order to be distinguished each other. Our main result is as following

**Proposition 1.** Suppose that $\mathcal{S}$ is a set of $n$ coins with same semblance, in which possibly there are some counterfeit coins, which are heavier (or lighter) than the normals. Denoted by $g(n)$ the least number of weighings need to sort $\mathcal{S}$ into the normals and fakes by a balance, assumed that additional normal coins will be available if needed, then it has

$$\lceil n \cdot \log_3 2 \rceil \leq g(n) \leq \lceil 7n/11 \rceil, \tag{1}$$

but with a exception that $g(3) = 3$.

The proof of Proposition 1 will largely rely on the following lemma, which is one of our unpublished results obtained in 1996.

**Lemma 1.** Let $\mathcal{A}$ be a set of $11$ coins with two possible distinct weights, then $\mathcal{A}$ can be sorted by a balance with at most $7$ weighings, and in the sixth weighing it will be found if $\mathcal{A}$ contains only one kind weight of coins.

Proof. It will be suffice to provide a feasible algorithm for prove the lemma. We have succeeded such an algorithm, which sketch is provided as an appendix in the end of the paper, for to describe the whole algorithm require much more space. The readers are encouraged yet to design own algorithms for it as a puzzle. □

**Lemma 2.** Let $\mathcal{S}$ be same as in Proposition 1, $\mathcal{S} = \mathcal{A} \cup \mathcal{B}, |\mathcal{B}| = 3$. If there is a algorithm $F$ with at most $k$ weighings, $k \geq 1$, to find all the fakes in $\mathcal{A}$, then there is also a algorithm $\tilde{F}$ to find all the fakes in $\mathcal{S}$ with at most $k + 2$ weighings.

Proof. We apply the algorithm $F$ on set $\mathcal{A}$ up to $k - 1$ weighings. Let $\Gamma$ be a non-empty output after $k - 1$ weighings, clearly, $|\Gamma| = 1, 2,$ or $3$.

If $|\Gamma| = 1$, that is, $\Gamma$ contains only one objective, and so which is the set of all the fakes in $\mathcal{A}$, the rest we need to do is to find the fakes in set $\mathcal{B}$, thus in this case,

$$g(|\mathcal{S}|) \leq k - 1 + g(3) \leq k + 2.$$

If $|\Gamma|=2$, suppose that $\Gamma=\{X,Y\}$, then there is one coin $x$  $x\in X$ (or $Y$) such that $x\notin Y$ (or $X$), let $\mathcal{B}'=x\cup\mathcal{B}$, it is easy to know we are provided to search for all the fakes in set $\mathcal{B}'$. On the other hand, it is not difficult to know that $g(4)=3$, hence,

$$g(|\mathcal{S}|)\leq k-1+g(4)\leq k+2.$$

If $|\Gamma|=3$, suppose that $\Gamma=\{X,Y,Z\}$. If there are two of $\{X,Y,Z\}$, say, $X$ and $Y$, such that $X\not\subseteq Y\cup Z$ and $Y\not\subseteq X\cup Z$, then take $x\in X\setminus(Y\cup Z)$ and $y\in Y\setminus(X\cup Z)$. It is easy to know that which one of three objectives in the set $\Delta=\{x,y,\varnothing\}$ is fakes will tell which of $\{X,Y,Z\}$ is the set of fakes in $\mathcal{A}$. So it will be suffice to search all the fakes in sets $\Delta$ and $\mathcal{B}$ in three weighings, such an algorithm is described in Page 21.

Thereby, assume that $X\subseteq Y\cup Z$ and $Y\subseteq X\cup Z$. If there is no one of $Z\cap X$ and $Z\cap Y$ is contained in the other one, then take $x\in(Z\cap X)\setminus(Z\cap Y)$, $y\in(Z\cap Y)\setminus(Z\cap X)$. Similar to the above, which one of three objectives in the set $\Delta'=\{x,y,x\cdot y\}$ is fakes will tell which of $\{X,Y,Z\}$ is the set of fakes in $\mathcal{A}$. We have to find all the fakes in sets $\Delta'$ and $\mathcal{B}$ in three weighings. The algorithm is presented in Page 22.

If $(Z\cap X)\subset(Z\cap Y)$, this implies that $X\subset Y$. Moreover, from $Y=X\cup(Z\cap Y)$, it has $(Z\cap Y)\not\subseteq X$. If $X\not\subseteq(Z\cap Y)$, then take $x\in X\setminus(Z\cap Y)$ and $z\in(Z\cap Y)\setminus X$, then $\Delta''=\{x,z,x\cdot z\}$, obviously, the case is similar to $\Delta'$. So, the rest to be discussed is $X\subseteq(Z\cap Y)$. From $Y\subseteq X\cup Z$, it has $Y\subset Z$. Take $y\in Y\setminus X$ and $z\in Z\setminus Y$, then which one of three objectives in the set $\Delta'''=\{z\cdot y,y,\varnothing\}$ is fakes will tell which of $\{X,Y,Z\}$ is the set of fakes in $\mathcal{A}$. We have to find all the fakes in sets $\Delta'''$ and $\mathcal{B}$ in three weighings. The sketch of such an algorithm has been shown in the page 23 of the appendices. $\square$

Proof of Proposition 1.  The left-hand of (1) is trivial for the information theoretical bound of $g(n)$ is equal to $\lceil n\cdot\log_3 2\rceil$, so we are provided to prove the right-hand of (1).

It is not difficult to verified that for $n\leq 11, n\neq 3$, it has $g(n)\leq\lceil 7n/11\rceil$, and $g(3)=3$. For

$n > 11$, let $n = 11 \cdot m + r$, $0 \le r < 11$, and $\mathcal{S} = \bigcup_{0 \le i \le m} \mathcal{A}_i$, $|\mathcal{A}_i| = 11, 0 \le i < m, |\mathcal{A}_m| = r$. If $r \ne 3$, by Lemma 1, there is

$$g(n) \le 7m + \lceil 7r/11 \rceil \le \lceil 7n/11 \rceil.$$

If $r = 3$, by Lemma 1 and Lemma 2, it has

$$g(n) \le 7m + 2 \le 7m + \lceil 7 \times 3/11 \rceil \le \lceil 7n/11 \rceil. \qquad \square$$

In the following, we will investigate the sorting case, that is, without provided of additional normal coins.

**Proposition 2.** Suppose that $\mathcal{S}$ is a set of $n$ coins with same semblance and with two possible distinct weights. Denoted by $\bar{g}(n)$ the least number of weighings needed to sort $\mathcal{S}$ in the weights by a balance, then it has

$$\lceil \log_3(2^n - 1) \rceil \le \bar{g}(n) \le \lceil 7n/11 \rceil. \qquad (2)$$

To prove Proposition 2, it needs following initial result

**Lemma 3.** For $n \le 11$, there is an algorithm with at most $\lceil 7n/11 \rceil$ weighings to sort $\mathcal{S}$, and it will be found if all the coins in $\mathcal{S}$ are of same weight in the $(\lceil 7n/11 \rceil - 1)$-th weighing except $n = 3$.

Proof. It is easy to verify the conclusions in the Lemma for $n < 6$ directly. Moreover, the required algorithm for $n = 7$ may be derived from the one for $n = 6$, and the algorithm for $n = 8$ may be derived from the one for $n = 7$. The required algorithm for $n = 9$ may be obtained from the algorithms for $n = 6$ and Lemma 2, the one for $n = 10$ may be derived from the one for $n = 9$. The case $n = 11$ is as in Lemma 1. So, the rest algorithm to be designed is the one for $n = 6$, such an algorithm is provided in the page 24. $\qquad \square$

The Proof of Proposition 2. The lower bound of (2) is the information theoretic bound of $\bar{g}(n)$, so, to be proved is the upper bound of (2). With Lemma 3, we may assume that $n > 11$. Let $n = 11 \cdot k + r, 0 \le r < 11,$ and

$$\mathcal{S} = \bigcup_{1 \le i \le k} \mathcal{A}_i \cup \mathcal{B}, \quad |\mathcal{A}_i| = 11, 1 \le i \le k, |\mathcal{B}| = r.$$

Denoted by $\mathcal{Q} = \bigcup_{1 \le i \le k} \mathcal{A}_i$, then by Lemma 1, we can spend at most $7k$ weighings to complete a sorting of the set $\mathcal{Q}$, and if $r \ne 3$, by Lemma 3, at most $7k + \lceil 7r/11 \rceil$ weighings will be

enough to complete a sorting of whole set $\mathcal{S}$, hence

$$\bar{g}(n) \leq 7k + \lceil 7r/11 \rceil = \lceil 7n/11 \rceil.$$

So, the rest to be discussed is the case $r = 3$. With Lemma 2, it spends at most $7k + 2$ weighings to complete a sorting of the set $\mathcal{S}$, hence

$$\bar{g}(n) \leq 7k + 2 = 7k + \lceil 7 \times 3/11 \rceil = \lceil 7n/11 \rceil.$$

So, the proof is finished. □

It maybe should be mentioned that in Lemma 2 the third algorithm we presented, that is the last one to deal with the case $X \subset Y \subset Z$, will need a normal coin, clearly which may be provided unless $Z = \mathcal{Q}$, however, this case will not occur in applying the algorithm described in Lemma 1, for we will find whether the coins in set $\mathcal{Q}$ are of same weight by $7k - 1 (= 6k + k - 1)$ weighings if it is, so, the objective set $\Gamma = \{\mathcal{Q}\}$, rather than $\Gamma = \{X, Y, Z\}$.

In addition, as a derivative result, there is a corollary of Lemma 1

**Corollary 1.** Let $\mathcal{A}_i, 1 \leq i \leq 11,$ be $11$ sets of coins, each $\mathcal{A}_i$ consists of two coins, one is normal and one is fake, the normals and the fakes have distinct weights, then at most spend $7$ weighings by a balance to find all the fakes in $\mathcal{A}_i$'s, $1 \leq i \leq 11$.

Proof. From each $\mathcal{A}_i$ take a coin, $1 \leq i \leq 11$, to form a coin set $\mathcal{A}$, applying Lemma 1 to $\mathcal{A}$, if there are two kinds of coins in $\mathcal{A}$, then each $\mathcal{A}_i, 1 \leq i \leq 11,$ will be sorted, otherwise, that is, if $\mathcal{A}$ contains only one kind of coins, then a additional weighing will work well. □

Likely, there may be the algorithms independent of Lemma 1 for the result of Corollary 1, which will be left to the readers.

**Acknowledgement:**

The author greatly appreciates Professor L. Lovasz for his endorsement for this paper, and many thanks to arXiv for publish it.
Finally, the author should acknowledge Hagen von Eitzen for he had pointed out the flaws in the statement and the proof of Proposition 1 in the earlier versions, and told me the counter-example $n = 19$ for the conjecture $\bar{g}(n) = \lceil n \cdot \log_3 2 \rceil$ once proposed in the original version of this paper.

**Appendices**     Algorithm to Sort 11 Coins

The expression $L : R$ means a comparison or a weighing with $L$ and $R$ are the left-hand and the right-hand of the balance respectively, and for simplicity, the expression {1,2,3} means the total weight of the coins with marks 1,2 and 3, and so on.

Fig.22
 ↖(<)
Fig.21 ←(=) {3,4,7,10} :{1,2,6,11} ←(<) {5,7,9} :{6,8,10} ←(<) {1,7,8} :{2,9,10} ←(<) {1,2,3} :{4,5,6}
 ↙(>)                                    ↓(=)
Fig.20

                          {2,9,11} :{ 3,4,5}
                        ↙(<)   ↓(=)   ↘(>)
                       Fig.19  Fig.18  Fig.17                     (=)

Fig.16 ↖(<)
Fig.15 ←(=) {2,3,6} :{5,8,9} ←(<) {1,5,7,9} :{2,6,8,10}
       ↙(>)                        ↓(=)
Fig.14

                                   {1,2,4} :{ 3,7,10}
                                 ↙(<)   ↓(=)   ↘(>)
                                Fig.13  Fig.12  Fig.11          (=)

Fig.10 ↖(<)
Fig.9 ←(=) {1,5,6} :{2,3,7} ←(<) {7,9} :{8,10} ←(<) {1,7,8} :{ 4,9,10}
       ↙(>)
Fig.8                              ↓(=)
                     Fig.7 ←(<) {5,7,10} :{ 6,8,9}                (=)
                                   ↓(=)
                                  Fig..6

Fig.5 ↖(<)
Fig.4 ←(=) {1,3,4}:{8,9,11} ←(<) {1,2,7} :{ 3,4,8}
       ↙(>)
Fig.3                                                            (=)

                                         Fig.2 ←(<) {1,3} :{ 2,4}
                                                              ↓(=)
                                                             Fig.1

Fig. 0

## Fig.1

- {1,2} :{10,11}
  - (<) → {10} :{11}  {1} :{2} (=)
    - (>) → (7,8,9,10)
    - (<) → (11)
  - (=) → {1} :{2}
    - (<) → {10} :{11}
    - (=) → (~)
    - (>) → {10} :{11}
  - (>) → {1} :{2} (=) {10} :{11}
    - (>) → (1,4)
    - (<) → (2,3,5,6)

Fig.1

## Fig.2

- {1} :{7}
  - (<) → {3,5} :{1,11} (=) {9} :{10}
    - (>) → {6} :{11}
    - (<) → {9} :{10}
  - (=) → {5,11} :{6,9}
    - (<) → {10} :{11}
    - (=) → {9} :{10}
    - (>) → {10} :{11}
  - (>) → {5,11} :{6,9} (=) {5} :{6}
    - (>) → {10} :{11}
    - (<) → {10} :{11}

Fig.2

## Fig.3

- {5} :{6}
  - (<) → {1,9} :{8,11} (=) {10} :{11}
    - (>) → {5,9} :{7,10}
    - (<) → {1} :{7}
  - (=) → {2} :{7}
    - (<) → {10} :{11}
    - (=) → {1} :{3}
    - (>) → {9} :{11}
  - (>) → {1,9} :{8,11} (=) {10} :{11}
    - (>) → {6,9} :{7,10}
    - (<) → {1} :{7}

Fig.3

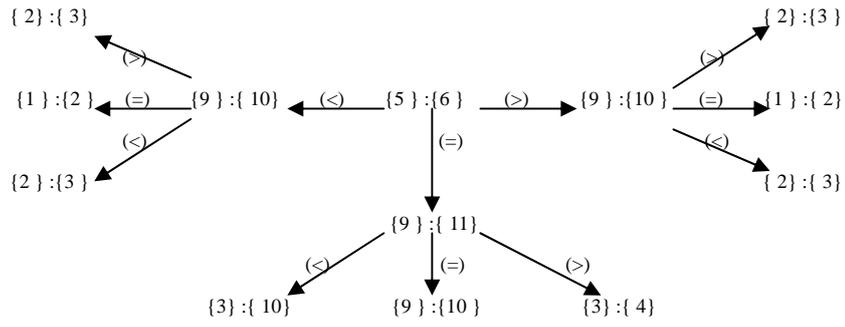

Fig.4

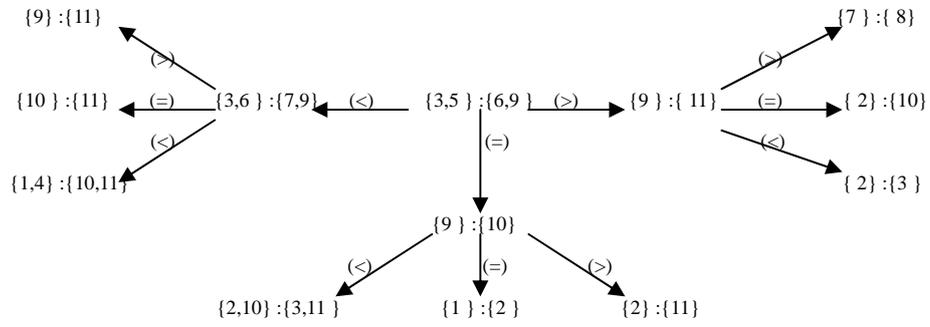

Fig.5

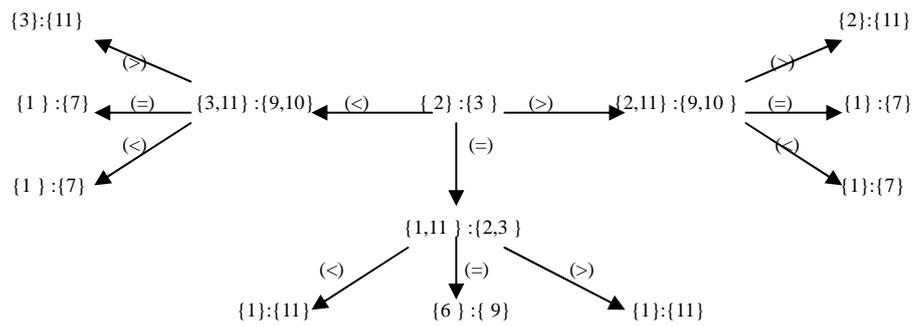

Fig.6

## Fig.7

- {1}:{11}
- {6}:{8}  (=) {1,11}:{2,8} (<) {2,5}:{3,7} (>) {1,6}:{5,11} (=) {6}:{8}
- {4,5}:{10,11}
- {1,4}:{10,11}
- {5}:{11}
- {1,10}:{8,11}
  - (<) {8}:{11}
  - (=) {6}:{9}
  - (>) {2}:{11}

Fig.7

## Fig.8

- {2,3}:{9,11}
- {9}:{11} (=) {5,9}:{6,11} (<) {2}:{3} (>) {5,9}:{6,11} (=) {9}:{11}
- {2,3}:{9,11}
- {2,3}:{9,11}
- {2,3}:{9,11}
- {5,6}:{9,11}
  - (<) {9}:{11}
  - (=) {9}:{10}
  - (>) {9}:{11}

Fig.8

## Fig.9

- {2,3}:{9,11}
- {9}:{11} (=) {5,9}:{6,11} (<) {2}:{3} (>) {5,9}:{6,11} (=) {9}:{11}
- {2,3}:{9,11}
- {2,3}:{9,11}
- {2,3}:{9,11}
- {5,6}:{9,11}
  - (<) {9}:{11}
  - (=) {9}:{10}
  - (>) {9}:{11}

Fig.9

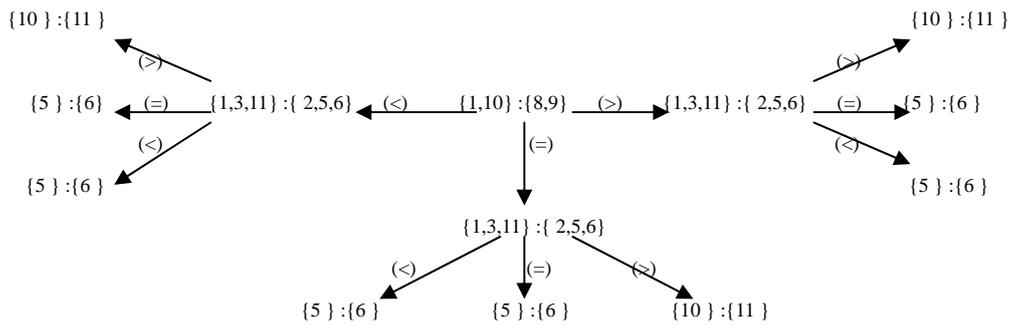

Fig.10

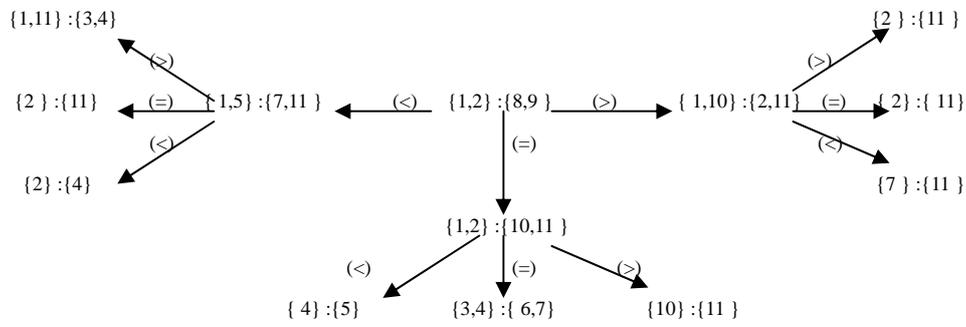

Fig.11

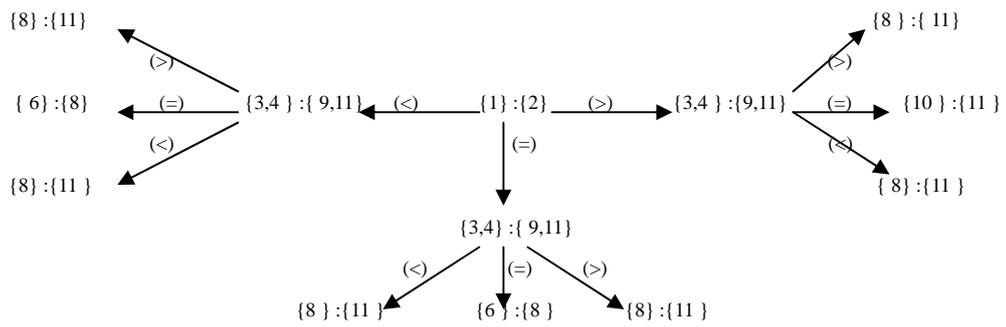

Fig.12

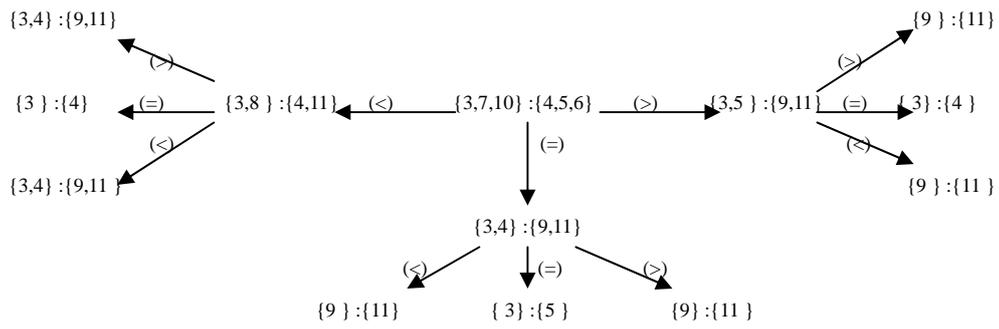

Fig.13

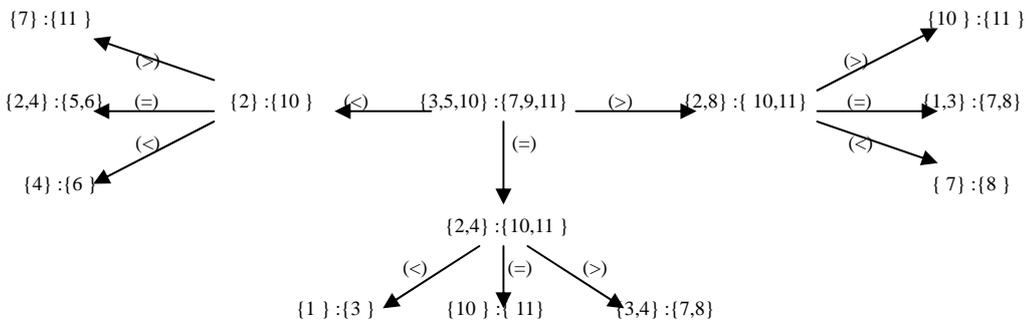

Fig.14

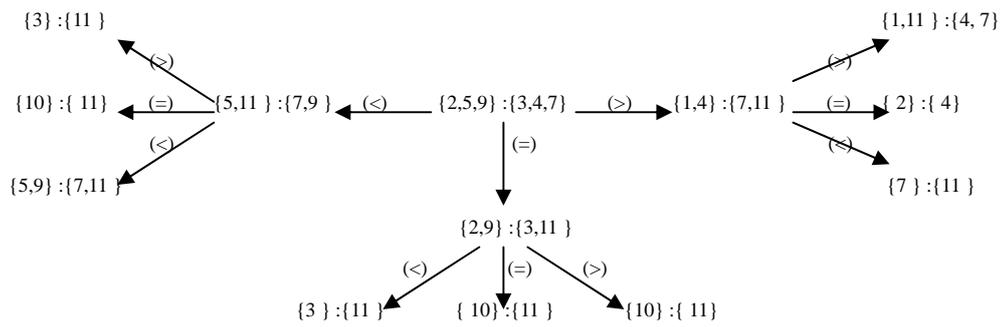

Fig.15

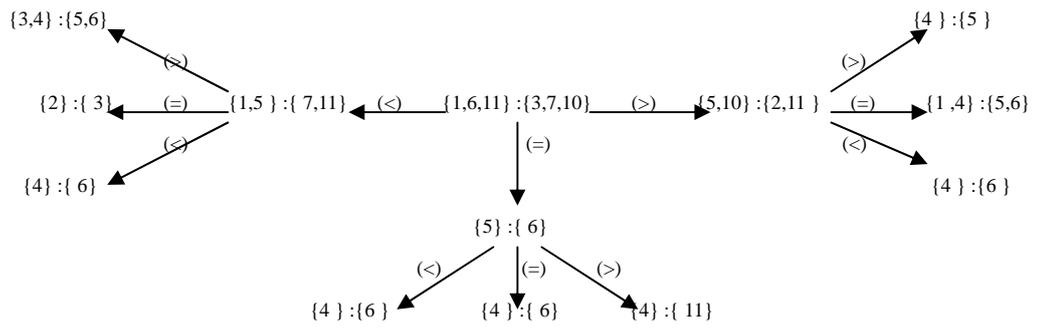

Fig.16

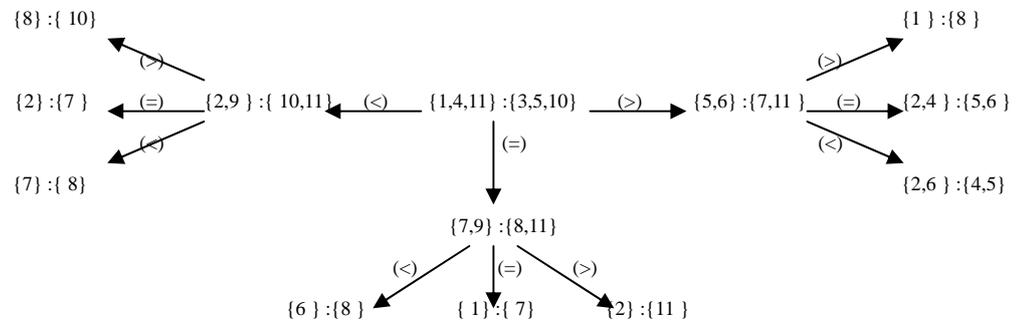

Fig.17

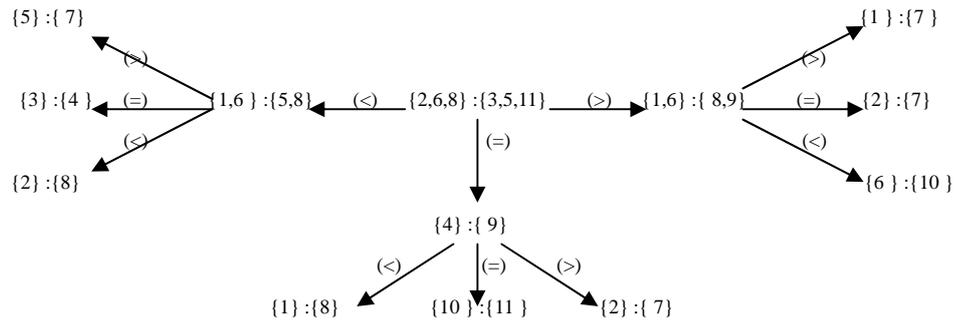

Fig.18

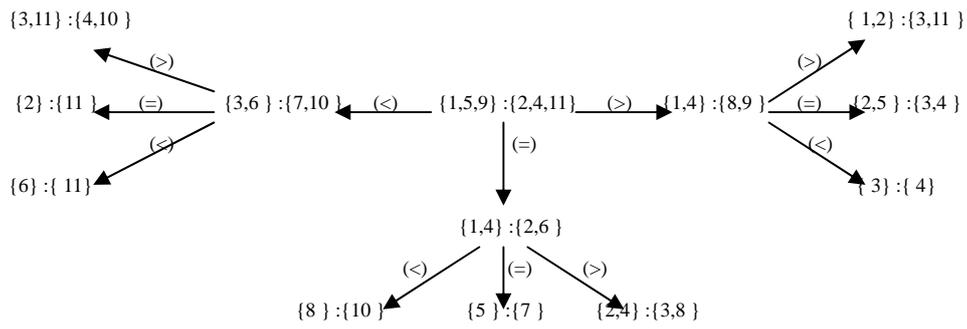

Fig.19

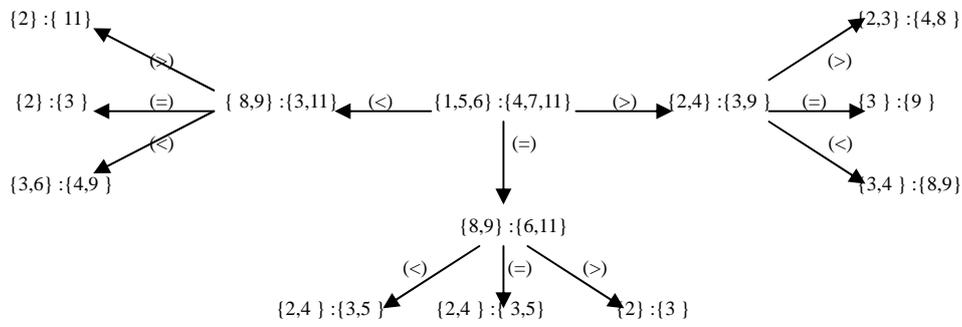

Fig.20

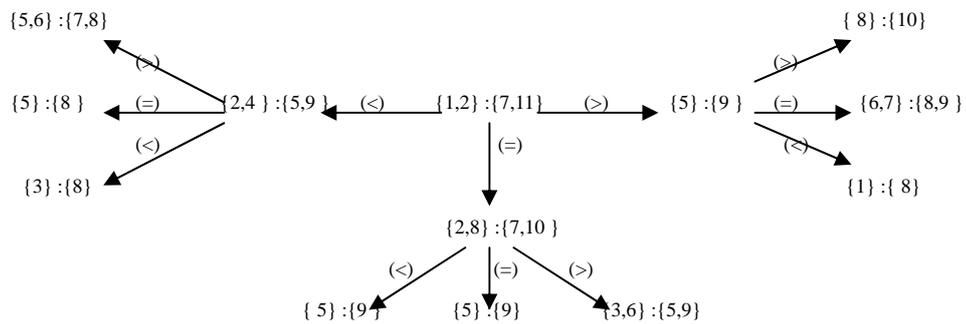

Fig.21

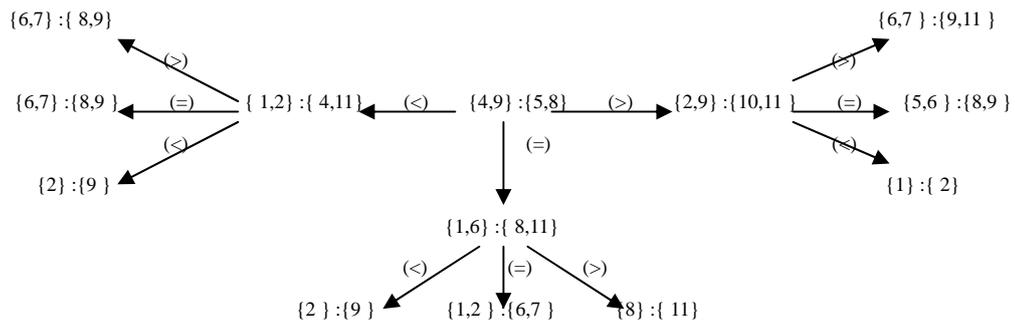

Fig.22

Where $L > R$, or $L < R$, or $L = R$ means $L$ heavier than, or lighter than, or equal to $R$ respectively, and the outputs (x), (a,b,c) and (~) means null, the coins with marks $a, b, c$ are heavier ones, and all coins have same weights respectively. The symmetric part of the algorithm has been omitted.

# The second algorithm to sort 11 coins

The difference between the new and the previous one is mainly the Fig.1 to Fig.5 and the third weighing on the direction (=, =) in Fig.0, and the other keep unchanged.

Fig.22
(<)
Fig.21 ←(=)— {3,4,7,10} :{1,2,6,11} ←(<)— {5,7,9} :{6,8,10} ←(<)— {1,7,8} :{2,9,10} ←(<)— {1,2,3} :{4,5,6}
(>)
Fig.20
(=)
{2,9,11} :{ 3,4,5}                                    (=)
(<)  (=)  (>)
Fig.19  Fig.18  Fig.17

Fig.16  (<)
Fig.15 ←(=)— {2,3,6} :{5,8,9} ←(<)— {1,5,7,9} :{2,6,8,10}                    (=)
(>)
Fig.14                                       (=)

{1,2,4} :{ 3,7,10}
(<)  (=)  (>)
Fig.13  Fig.12  Fig.11

Fig.10 (<)
Fig.9 ←(=)— {1,5,6} :{2,3,7} ←(<)— {7,9} :{8,10} ←(<)— {1,7,8} :{ 4,9,10}
(>)
Fig.8                                       (=)

Fig.7 ←(<)— {5,7,10} :{ 6,8,9}
(=)
Fig..6                                      (=)

Fig.5'  (<)
Fig.4' ←(=)— {1,9}:{4,10} ←(<)— {2, 7} :{ 3, 8}
(>)
Fig.3'                                      (=)

Fig.2' ←(<)— {5} :{ 6}
(=)
Fig.1'

Fig. 0'

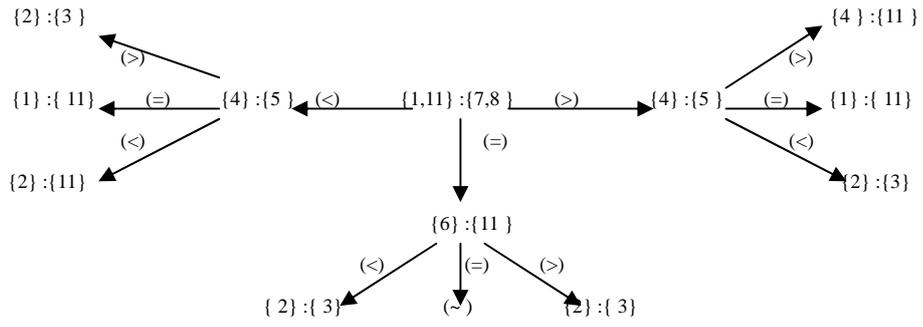

Fig.1'

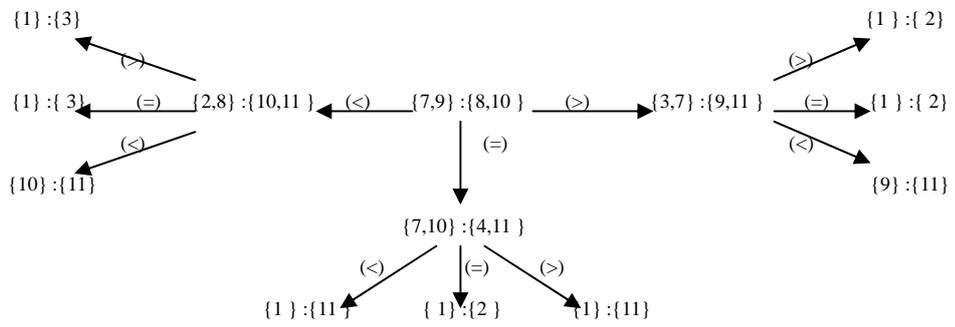

Fig.2'

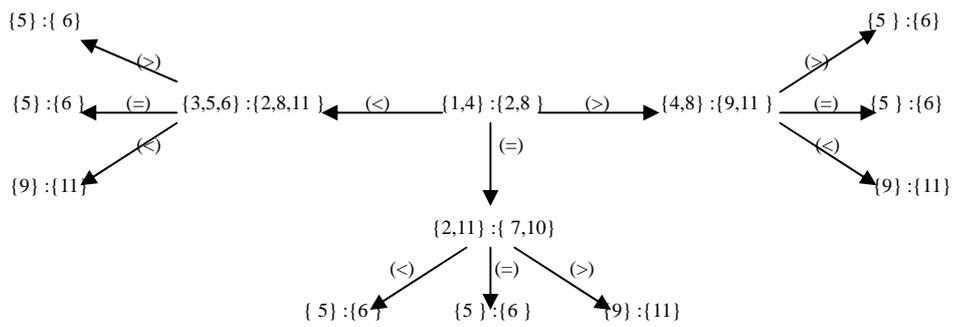

Fig.3'

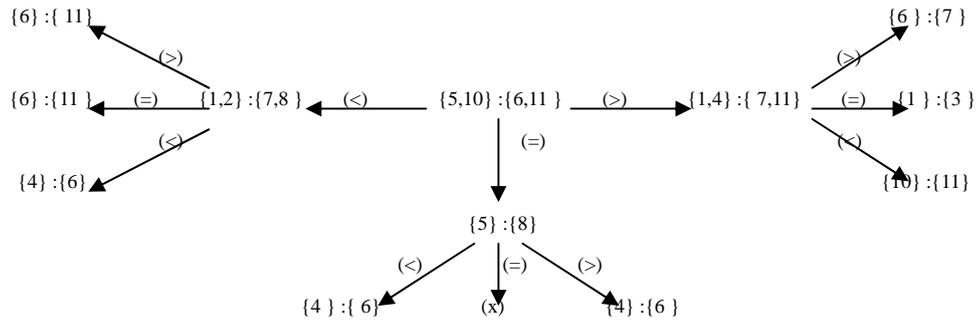

Fig.4'

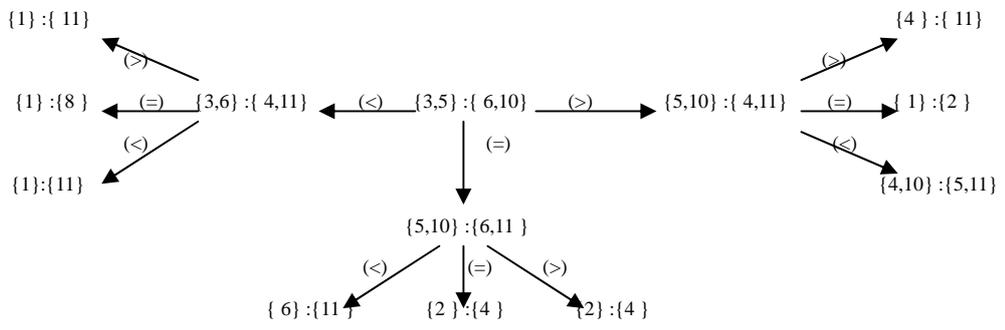

Fig.5'

# The third algorithm to sort 11 coins

The difference between the third and the first one are in the Figs.3,4,5 and the weighing on the direction (=, =, <) in Fig.0, and the other keep unchanged.

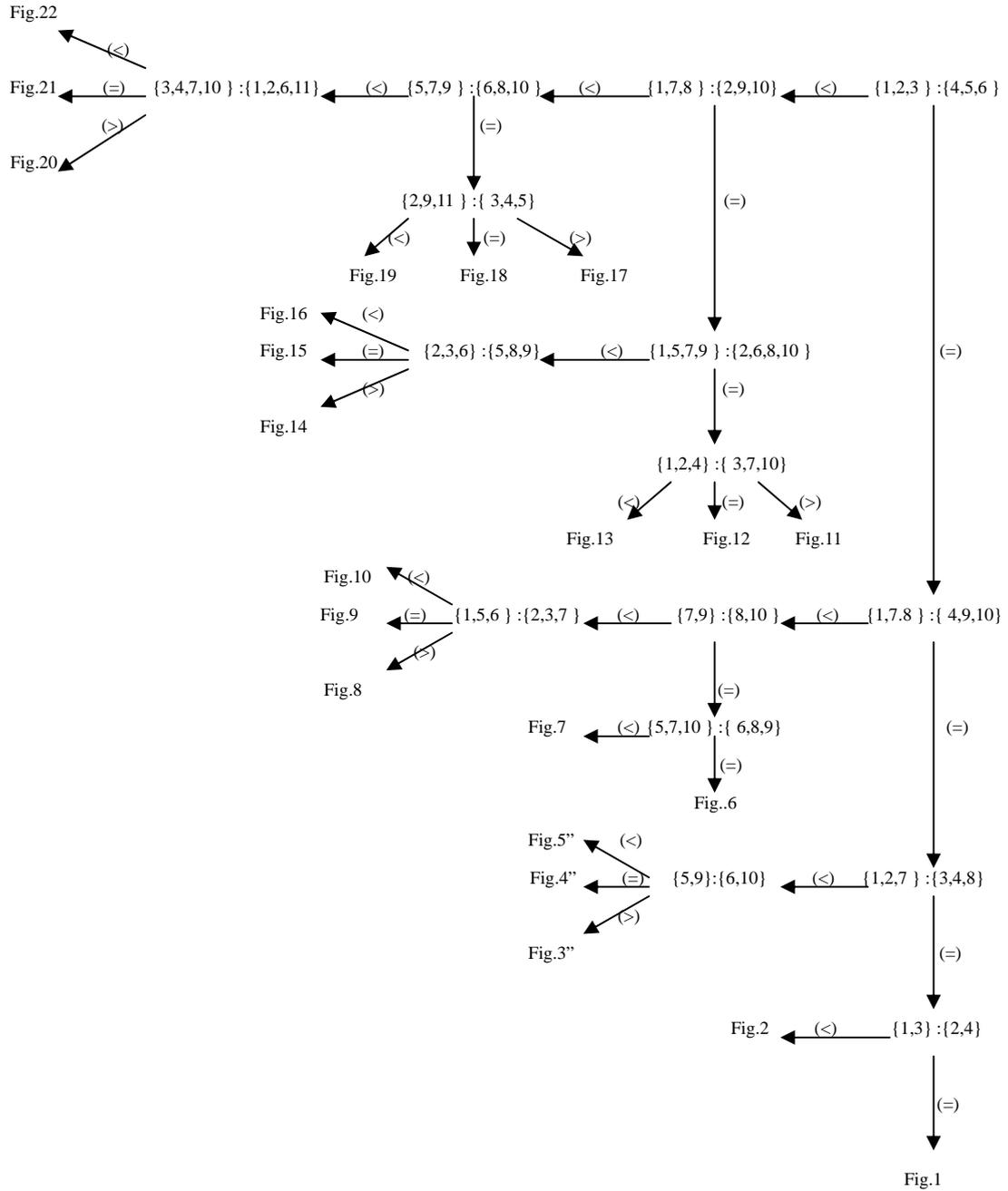

Fig. 0"

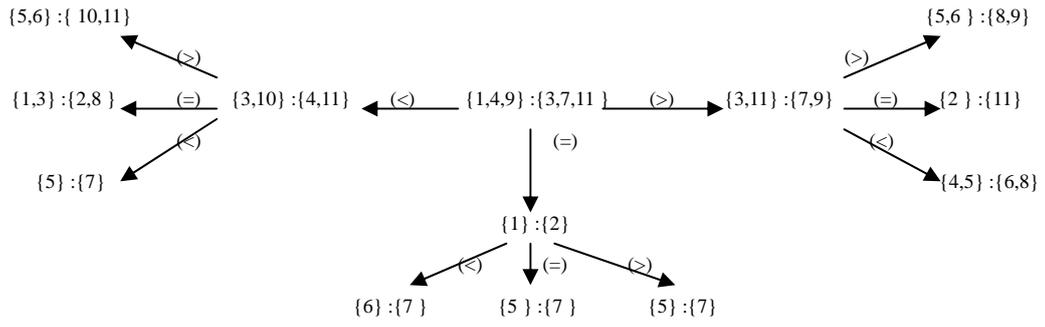

Fig.3"

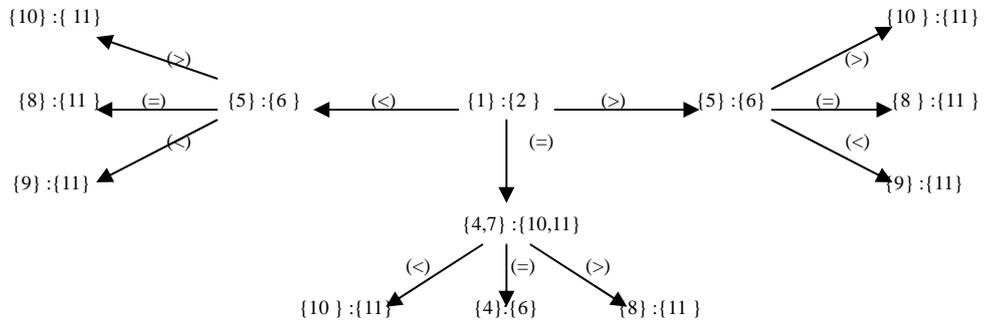

Fig.4"

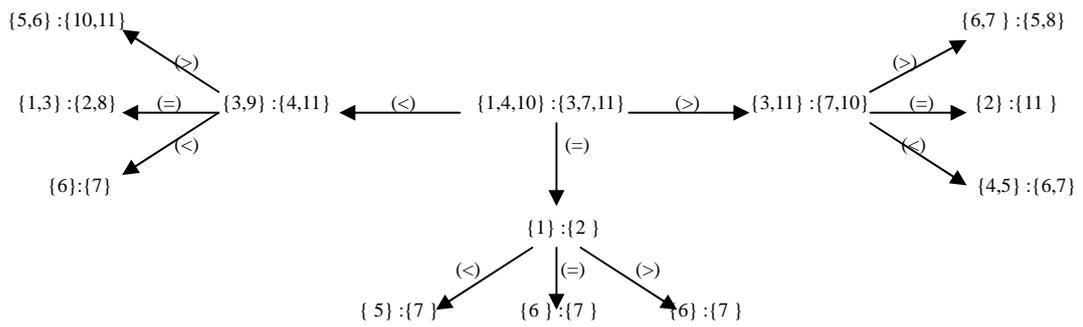

Fig.5"

Algorithm for Lemma 2 (1)

$b_2.b_3.y$

$(<)$

$b_2.b_3 \xleftarrow{(=)} b_1b_2 : b_3 y \xleftarrow{(<)} b_1 y : b_2 b_3 \xleftarrow{(<)} xb_1 : yb_2$

$(>)$

$b_2$

$(=)$ $(>)$

$b_3.y \xleftarrow{(<)} b_2 : b_3 \quad\quad b_1 : y$

$(=)$ $(>)$ $(<)\;(=)\;(>)$ $(=)$

$b_1.b_2.b_3.y \quad b_2.y \quad y \quad b_1.b_2.y \quad \sim$

$b_1.b_2.b_3$

$(<)$

$b_1.b_2 \xleftarrow{(=)} y : b_3 \xleftarrow{(<)} xb_3 : b_1 b_2$

$(>)$

$b_1.y$

$(>)$ $(=)$

$b_3 \xleftarrow{(<)} b_2 : b_3 \quad\quad y : b_2$

$(=)$ $(>)$ $(>)\;(=)\;(<)$

$b_2.b_3.x \quad \sim \quad b_1.b_3.y \quad \varnothing \quad x.b_2$

Fig. A

$\Delta = \{\varnothing, x, y\}, \;\; \mathcal{B} = \{b_1, b_2, b_3\}$.

$\sim\; :=$ Null

Algorithm for Lemma 2    (2)

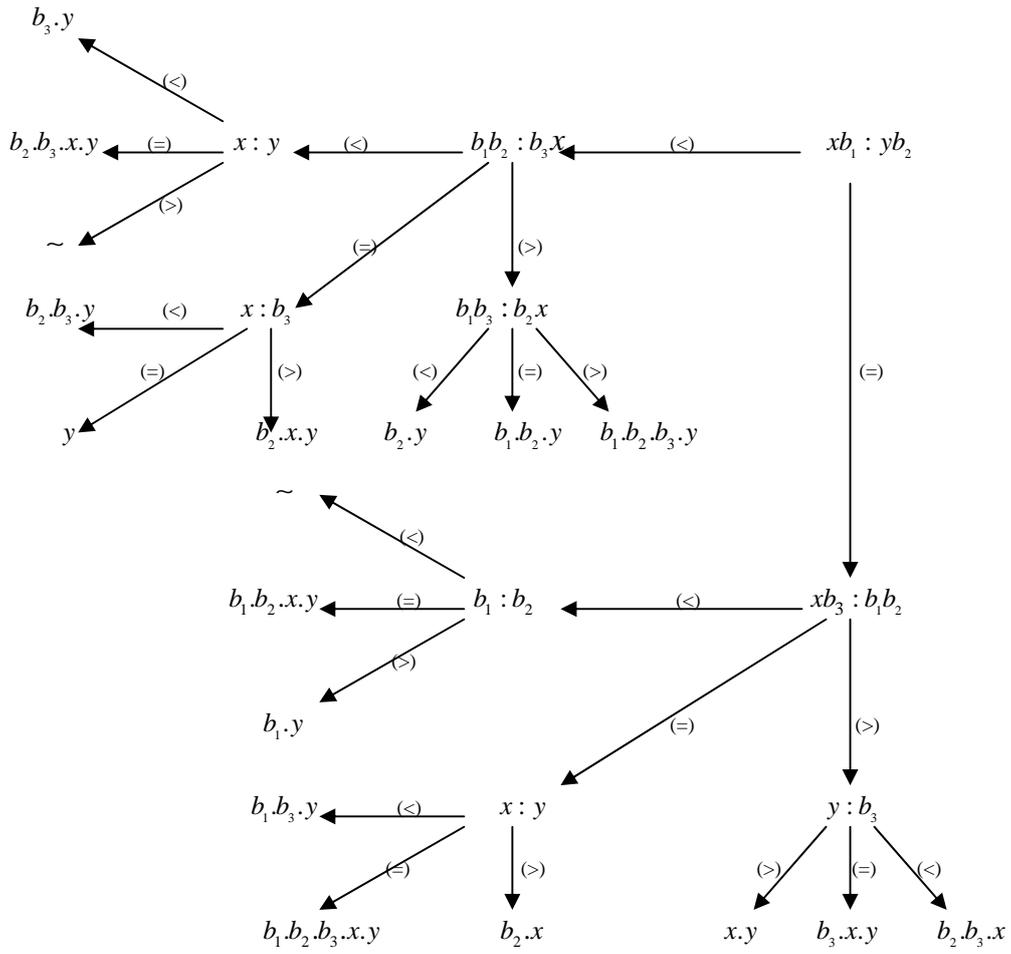

Fig. B

$\Delta = \{x, y, x.y\}$, $\mathcal{B} = \{b_1, b_2, b_3\}$.

Algorithm for Lemma 2     (3)

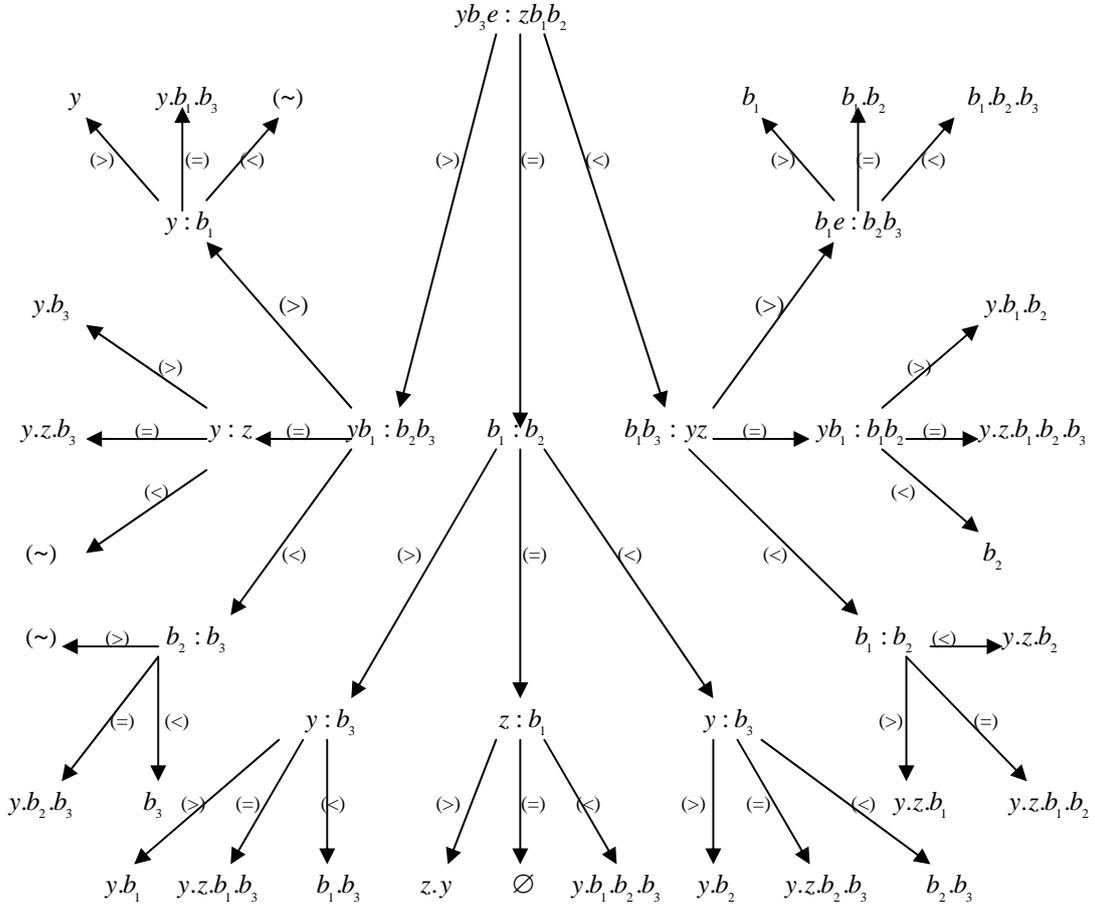

Fig. C

$\Delta = \{\varnothing, y, y.z\}$, $\mathcal{B} = \{b_1, b_2, b_3\}$, $e$ is a normal coin.

Algorithm $g(6) = 4$ in Lemma 3

3.5.6

5.6 ←(=)— {2}:{3} ←(<)— {4}:{6} ←(<)— {1,4}:{2,5} ←(<)— {1,2,3}:{4,5,6}

(<) above, (>) below to 2.5.6

{4}:{6} (=) → {1,5}:{2,3}; (>) → 2.4.5

{1,5}:{2,3}: (<) 2.3.4.5.6, (=) 2.4.5.6, (>) 5

{1,4}:{2,5} (=) ↓

3.4.5.6 ←(<)— {2}:{3} ←(<)— {1,5}:{3,6}
6 ←(=)—
2.4.6 ←(>)—

{1,5}:{3,6} (=) → {3}:{6};  (>) → {1}:{4}

{3}:{6}: (<) 4.5.6, (=) ~, (>) 3.4.5
{1}:{4}: (<) 4.5, (=) 1.2.4.5.6, (>) 1.5.6

{1,2,3}:{4,5,6} (=) ↓

2.6 ←(<)— {1}:{2} ←(<)— {3}:{6} ←(<)— {1,4}:{2,5}
1.2.5.6 ←(=)—
~ ←(>)—

{3}:{6} (=) → {2}:{3};  (>) → {2}:{3}

left {2}:{3}: (<) ~, (=) 2.3.5.6, (>) 2.5
right {2}:{3}: (<) 3.5, (=) 2.3.4.5, (>) ~

{1,4}:{2,5} (=) ↓

2.3.4.6
3.6 ←(=)— {1}:{2} ←(<)— {1,2,}:{3,6}
1.3.5.6 ←(>)—

{1,2,}:{3,6} (=) → 1.2.3.4.5.6;  (>) → {1}:{2}

{1}:{2}: (<) 2.4, (=) 1.2.4.5, (>) 1.5

Fig. D

The completed list of weighings of the first algorithm to sort 11 coins

*The 1-th weighing*
w( ) = {1,2,3}:{4,5,6}
*The 2-th weighing*
w(0) = {1,7,8}:{4,9,10}          w(1) = {1,7,8}:{2,9,10}          w(2) = {1,7,8}:{2,9,10}
*The 3-th weighing*
w(0,0) = {1,2,7}:{3,4,8}         w(1,0) = {1,5,7,9}:{2,6,8,10}    w(2,0) = {1,5,7,9}:{2,6,8,10}
w(0,1) = {7,9}:{8,10}            w(1,1) = {5,7,9}:{6,8,10}        w(2,1) = {5,7,9}:{6,8,10}
w(0,2) = {7,9}:{8,10}            w(1,2) = {5,7,9}:{6,8,10}        w(2,2) = {5,7,9}:{6,8,10}
*The 4-th weighing*
w(0,0,0) = {1,3}:{2,4}           w(1,0,0) = {1,2,4}:{3,7,10}      w(2,0,0) = {1,2,4}:{3,7,10}
w(0,1,0) = {5,7,10}:{6,8,9}      w(1,1,0) = {2,9,11}:{3,4,5}      w(2,1,0) = {1,7,11}:{3,4,5}
w(0,2,0) = {5,7,10}:{6,8,9}      w(1,2,0) = {1,7,11}:{3,4,5}      w(2,2,0) = {2,9,11}:{3,4,5}
w(0,0,1) = {1,3,4}:{8,9,11}      w(1,0,1) = {2,3,6}:{5,8,9}       w(2,0,1) = {1,3,5}:{6,8,9}
w(0,1,1) = {1,5,6}:{2,3,7}       w(1,1,1) = {3,4,7,10}:{1,2,6,11} w(2,1,1) = {3,4,7,10}:{1,2,5,11}
w(0,2,1) = {1,5,6}:{2,3,8}       w(1,2,1) = {3,4,8,9}:{1,2,6,11}  w(2,2,1) = {3,4,8,9}:{1,2,5,11}
w(0,0,2) = {1,3,4}:{8,9,11}      w(1,0,2) = {1,3,5}:{6,8,9}       w(2,0,2) = {2,3,6}:{5,8,9}
w(0,1,2) = {1,5,6}:{2,3,8}       w(1,1,2) = {3,4,8,9}:{1,2,5,11}  w(2,1,2) = {3,4,8,9}:{1,2,6,11}
w(0,2,2) = {1,5,6}:{2,3,7}       w(1,2,2) = {3,4,7,10}:{1,2,5,11} w(2,2,2) = {3,4,7,10}:{1,2,6,11}
*The 5-th weighing*
w(0,0,0,0) = {1,2}:{10,11}       w(1,0,0,0) = {1}:{2}             w(2,0,0,0) = {1}:{2}
w(0,1,0,0) = {2}:{3}             w(1,1,0,0) = {2,6,8}:{3,5,11}    w(2,1,0,0) = {1,6,10}:{3,5,11}
w(0,2,0,0) = {2}:{3}             w(1,2,0,0) = {1,6,10}:{3,5,11}   w(2,2,0,0) = {2,6,8}:{3,5,11}
w(0,0,1,0) = {5}:{6}             w(1,0,1,0) = {2,5,9}:{3,4,7}     w(2,0,1,0) = {1,6,8}:{3,4,10}
w(0,1,1,0) = {2}:{3}             w(1,1,1,0) = {1,2}:{7,11}        w(2,1,1,0) = {1,2}:{10,11}
w(0,2,1,0) = {2}:{3}             w(1,2,1,0) = {1,2}:{9,11}        w(2,2,1,0) = {1,2}:{8,11}
w(0,0,2,0) = {5}:{6}             w(1,0,2,0) = {1,6,8}:{3,4,10}    w(2,0,2,0) = {2,5,9}:{3,4,7}
w(0,1,2,0) = {2}:{3}             w(1,1,2,0) = {1,2}:{8,11}        w(2,1,2,0) = {1,2}:{9,11}
w(0,2,2,0) = {2}:{3}             w(1,2,2,0) = {1,2}:{10,11}       w(2,2,2,0) = {1,2}:{7,11}
w(0,0,0,1) = {1}:{7}             w(1,0,0,1) = {3,7,10}:{4,5,6}    w(2,0,0,1) = {1,2}:{8,9}
w(0,1,0,1) = {2,5}:{3,7}         w(1,1,0,1) = {1,5,9}:{2,4,11}    w(2,1,0,1) = {2,4,11}:{3,5,8}
w(0,2,0,1) = {2,6}:{3,8}         w(1,2,0,1) = {2,5,7}:{1,4,11}    w(2,2,0,1) = {1,4,11}:{3,5,10}
w(0,0,1,1) = {3,5}:{6,9}         w(1,0,1,1) = {1,6,11}:{3,7,10}   w(2,0,1,1) = {3,6,7}:{8,10,11}
w(0,1,1,1) = {1,10}:{8,9}        w(1,1,1,1) = {4,9}:{5,8}         w(2,1,1,1) = {2,5,6}:{4,10,11}
w(0,2,1,1) = {2}:{3}             w(1,2,1,1) = {4,7}:{5,10}        w(2,2,1,1) = {1,5,6}:{4,8,11}
w(0,0,2,1) = {5}:{6}             w(1,0,2,1) = {2,5,11}:{3,7,10}   w(2,0,2,1) = {3,5,10}:{7,9,11}
w(0,1,2,1) = {1,9}:{7,10}        w(1,1,2,1) = {4,10}:{6,7}        w(2,1,2,1) = {2,5,6}:{4,9,11}
w(0,2,2,1) = {2}:{3}             w(1,2,2,1) = {4,8}:{6,9}         w(2,2,2,1) = {1,5,6}:{4,7,11}
w(0,0,0,2) = {1}:{7}             w(1,0,0,2) = {1,2}:{8,9}         w(2,0,0,2) = {3,7,10}:{4,5,6}
w(0,1,0,2) = {2,6}:{3,8}         w(1,1,0,2) = {1,4,11}:{3,5,10}   w(2,1,0,2) = {2,5,7}:{1,4,11}
w(0,2,0,2) = {2,5}:{3,7}         w(1,2,0,2) = {2,4,11}:{3,5,8}    w(2,2,0,2) = {1,5,9}:{2,4,11}
w(0,0,1,2) = {5}:{6}             w(1,0,1,2) = {3,5,10}:{7,9,11}   w(2,0,1,2) = {2,5,11}:{3,7,10}

| | | |
|---|---|---|
| w(0,1,1,2) = {2}:{3} | w(1,1,1,2) = {1,5,6}:{4,7,11} | w(2,1,1,2) = {4,8}:{6,9} |
| w(0,2,1,2) = {1,9}:{7,10} | w(1,2,1,2) = {2,5,6}:{4,9,11} | w(2,2,1,2) = {4,10}:{6,7} |
| w(0,0,2,2) = {3,5}:{6,9} | w(1,0,2,2) = {3,6,7}:{8,10,11} | w(2,0,2,2) = {1,6,11}:{3,7,10} |
| w(0,1,2,2) = {2}:{3} | w(1,1,2,2) = {1,5,6}:{4,8,11} | w(2,1,2,2) = {4,7}:{5,10} |
| w(0,2,2,2) = {1,10}:{8,9} | w(1,2,2,2) = {2,5,6}:{4,10,11} | w(2,2,2,2) = {4,9}:{5,8} |

*The 6-th weighing*

| | | |
|---|---|---|
| w(0,0,0,0,0) = {1}:{2} | w(1,0,0,0,0) = {3,4}:{9,11} | w(2,0,0,0,0) = {3,4}:{9,11} |
| w(0,1,0,0,0) = {1,11}:{2,3} | w(1,1,0,0,0) = {4}:{9} | w(2,1,0,0,0) = {4}:{7} |
| w(0,2,0,0,0) = {1,11}:{2,3} | w(1,2,0,0,0) = {4}:{7} | w(2,2,0,0,0) = {4}:{9} |
| w(0,0,1,0,0) = {9}:{11} | w(1,0,1,0,0) = {2,9}:{3,11} | w(2,0,1,0,0) = {1,8}:{3,11} |
| w(0,1,1,0,0) = {5,6}:{9,11} | w(1,1,1,0,0) = {2,8}:{7,10} | w(2,1,1,0,0) = {1,9}:{7,10} |
| w(0,2,1,0,0) = {5,6}:{10,11} | w(1,2,1,0,0) = {1,10}:{8,9} | w(2,2,1,0,0) = {2,7}:{8,9} |
| w(0,0,2,0,0) = {9}:{11} | w(1,0,2,0,0) = {1,8}:{3,11} | w(2,0,2,0,0) = {2,9}:{3,11} |
| w(0,1,2,0,0) = {5,6}:{10,11} | w(1,1,2,0,0) = {2,7}:{8,9} | w(2,1,2,0,0) = {1,10}:{8,9} |
| w(0,2,2,0,0) = {5,6}:{9,11} | w(1,2,2,0,0) = {1,9}:{7,10} | w(2,2,2,0,0) = {2,8}:{7,10} |
| w(0,0,0,1,0) = {5,11}:{6,9} | w(1,0,0,1,0) = {3,4}:{9,11} | w(2,0,0,1,0) = {1,2}:{10,11} |
| w(0,1,0,1,0) = {1,10}:{8,11} | w(1,1,0,1,0) = {1,4}:{2,6} | w(2,1,0,1,0) = {7,9}:{10,11} |
| w(0,2,0,1,0) = {1,9}:{7,11} | w(1,2,0,1,0) = {2,4}:{1,6} | w(2,2,0,1,0) = {7,9}:{8,11} |
| w(0,0,1,1,0) = {9}:{10} | w(1,0,1,1,0) = {5}:{6} | w(2,0,1,1,0) = {1,4}:{7,11} |
| w(0,1,1,1,0) = {1,3,11}:{2,5,6} | w(1,1,1,1,0) = {1,6}:{8,11} | w(2,1,1,1,0) = {8,9}:{5,11} |
| w(0,2,1,1,0) = {5,6}:{10,11} | w(1,2,1,1,0) = {2,6}:{10,11} | w(2,2,1,1,0) = {7,10}:{5,11} |
| w(0,0,2,1,0) = {2}:{7} | w(1,0,2,1,0) = {6}:{5} | w(2,0,2,1,0) = {2,4}:{10,11} |
| w(0,1,2,1,0) = {1,3,11}:{2,5,6} | w(1,1,2,1,0) = {1,5}:{7,11} | w(2,1,2,1,0) = {7,10}:{6,11} |
| w(0,2,2,1,0) = {5,6}:{9,11} | w(1,2,2,1,0) = {2,5}:{9,11} | w(2,2,2,1,0) = {8,9}:{6,11} |
| w(0,0,0,2,0) = {5,11}:{6,9} | w(1,0,0,2,0) = {1,2}:{10,11} | w(2,0,0,2,0) = {3,4}:{9,11} |
| w(0,1,0,2,0) = {1,9}:{7,11} | w(1,1,0,2,0) = {7,9}:{8,11} | w(2,1,0,2,0) = {2,4}:{1,6} |
| w(0,2,0,2,0) = {1,10}:{8,11} | w(1,2,0,2,0) = {7,9}:{10,11} | w(2,2,0,2,0) = {1,4}:{2,6} |
| w(0,0,1,2,0) = {2}:{7} | w(1,0,1,2,0) = {2,4}:{10,11} | w(2,0,1,2,0) = {6}:{5} |
| w(0,1,1,2,0) = {5,6}:{9,11} | w(1,1,1,2,0) = {8,9}:{6,11} | w(2,1,1,2,0) = {2,5}:{9,11} |
| w(0,2,1,2,0) = {1,3,11}:{2,5,6} | w(1,2,1,2,0) = {7,10}:{6,11} | w(2,2,1,2,0) = {1,5}:{7,11} |
| w(0,0,2,2,0) = {9}:{10} | w(1,0,2,2,0) = {1,4}:{7,11} | w(2,0,2,2,0) = {5}:{6} |
| w(0,1,2,2,0) = {5,6}:{10,11} | w(1,1,2,2,0) = {7,10}:{5,11} | w(2,1,2,2,0) = {2,6}:{10,11} |
| w(0,2,2,2,0) = {1,3,11}:{2,5,6} | w(1,2,2,2,0) = {8,9}:{5,11} | w(2,2,2,2,0) = {1,6}:{8,11} |
| w(0,0,0,0,1) = {10}:{11} | w(1,0,0,0,1) = {3,4}:{9,11} | w(2,0,0,0,1) = {3,4}:{9,11} |
| w(0,1,0,0,1) = {3,11}:{9,10} | w(1,1,0,0,1) = {1,6}:{5,8} | w(2,1,0,0,1) = {2,6}:{7,10} |
| w(0,2,0,0,1) = {2,11}:{9,10} | w(1,2,0,0,1) = {2,6}:{5,10} | w(2,2,0,0,1) = {1,6}:{8,9} |
| w(0,0,1,0,1) = {9}:{10} | w(1,0,1,0,1) = {5,11}:{7,9} | w(2,0,1,0,1) = {2,4}:{10,11} |
| w(0,1,1,0,1) = {5,9}:{6,11} | w(1,1,1,0,1) = {2,4}:{5,9} | w(2,1,1,0,1) = {6}:{8} |
| w(0,2,1,0,1) = {5,10}:{6,11} | w(1,2,1,0,1) = {1,4}:{5,7} | w(2,2,1,0,1) = {6}:{10} |
| w(0,0,2,0,1) = {9}:{10} | w(1,0,2,0,1) = {6,11}:{8,10} | w(2,0,2,0,1) = {1,4}:{7,11} |
| w(0,1,2,0,1) = {5,10}:{6,11} | w(1,1,2,0,1) = {2,4}:{6,10} | w(2,1,2,0,1) = {5}:{7} |
| w(0,2,2,0,1) = {5,9}:{6,11} | w(1,2,2,0,1) = {1,4}:{6,8} | w(2,2,2,0,1) = {5}:{9} |
| w(0,0,0,1,1) = {3,5}:{1,11} | w(1,0,0,1,1) = {3,8}:{4,11} | w(2,0,0,1,1) = {1,10}:{2,11} |
| w(0,1,0,1,1) = {1,11}:{2,8} | w(1,1,0,1,1) = {3,6}:{7,10} | w(2,1,0,1,1) = {5,6}:{9,11} |
| w(0,2,0,1,1) = {1,5}:{6,11} | w(1,2,0,1,1) = {3,6}:{8,9} | w(2,2,0,1,1) = {5,6}:{7,11} |
| w(0,0,1,1,1) = {3,6}:{7,9} | w(1,0,1,1,1) = {1,5}:{7,11} | w(2,0,1,1,1) = {1,9}:{7,11} |
| w(0,1,1,1,1) = {1,3,11}:{2,5,6} | w(1,1,1,1,1) = {1,2}:{4,11} | w(2,1,1,1,1) = {1,4}:{3,8} |
| w(0,2,1,1,1) = {5,10}:{6,11} | w(1,2,1,1,1) = {1,2}:{4,11} | w(2,2,1,1,1) = {2,4}:{3,10} |
| w(0,0,2,1,1) = {1,9}:{8,11} | w(1,0,2,1,1) = {2,6}:{10,11} | w(2,0,2,1,1) = {2,8}:{10,11} |
| w(0,1,2,1,1) = {1,3,11}:{2,5,6} | w(1,1,2,1,1) = {1,2}:{4,11} | w(2,1,2,1,1) = {1,4}:{3,7} |
| w(0,2,2,1,1) = {5,9}:{6,11} | w(1,2,2,1,1) = {1,2}:{4,11} | w(2,2,2,1,1) = {2,4}:{3,9} |
| w(0,0,0,2,1) = {5,11}:{6,9} | w(1,0,0,2,1) = {1,5}:{7,11} | w(2,0,0,2,1) = {3,5}:{9,11} |
| w(0,1,0,2,1) = {1,11}:{2,7} | w(1,1,0,2,1) = {2,9}:{10,11} | w(2,1,0,2,1) = {2,4}:{7,10} |
| w(0,2,0,2,1) = {1,6}:{5,11} | w(1,2,0,2,1) = {1,7}:{8,11} | w(2,2,0,2,1) = {1,4}:{8,9} |
| w(0,0,1,2,1) = {1,9}:{8,11} | w(1,0,1,2,1) = {2}:{10} | w(2,0,1,2,1) = {6,7}:{1,11} |
| w(0,1,1,2,1) = {5,9}:{6,11} | w(1,1,1,2,1) = {8,9}:{3,11} | w(2,1,1,2,1) = {1,8}:{7,11} |
| w(0,2,1,2,1) = {1,3,11}:{2,5,6} | w(1,2,1,2,1) = {7,10}:{3,11} | w(2,2,1,2,1) = {2,10}:{9,11} |
| w(0,0,2,2,1) = {9}:{11} | w(1,0,2,2,1) = {1}:{7} | w(2,0,2,2,1) = {5,10}:{2,11} |
| w(0,1,2,2,1) = {5,10}:{6,11} | w(1,1,2,2,1) = {7,10}:{3,11} | w(2,1,2,2,1) = {1,7}:{8,11} |
| w(0,2,2,2,1) = {1,3,11}:{2,5,6} | w(1,2,2,2,1) = {8,9}:{3,11} | w(2,2,2,2,1) = {2,9}:{10,11} |
| w(0,0,0,0,2) = {10}:{11} | w(1,0,0,0,2) = {3,4}:{9,11} | w(2,0,0,0,2) = {3,4}:{9,11} |
| w(0,1,0,0,2) = {2,11}:{9,10} | w(1,1,0,0,2) = {1,6}:{8,9} | w(2,1,0,0,2) = {2,6}:{5,10} |
| w(0,2,0,0,2) = {3,11}:{9,10} | w(1,2,0,0,2) = {2,6}:{7,10} | w(2,2,0,0,2) = {1,6}:{5,8} |
| w(0,0,1,0,2) = {9}:{10} | w(1,0,1,0,2) = {1,4}:{7,11} | w(2,0,1,0,2) = {6,11}:{8,10} |
| w(0,1,1,0,2) = {5,9}:{6,11} | w(1,1,1,0,2) = {5}:{9} | w(2,1,1,0,2) = {1,4}:{6,8} |
| w(0,2,1,0,2) = {5,10}:{6,11} | w(1,2,1,0,2) = {5}:{7} | w(2,2,1,0,2) = {2,4}:{6,10} |
| w(0,0,2,0,2) = {9}:{10} | w(1,0,2,0,2) = {2,4}:{10,11} | w(2,0,2,0,2) = {5,11}:{7,9} |
| w(0,1,2,0,2) = {5,10}:{6,11} | w(1,1,2,0,2) = {6}:{10} | w(2,1,2,0,2) = {1,4}:{5,7} |
| w(0,2,2,0,2) = {5,9}:{6,11} | w(1,2,2,0,2) = {6}:{8} | w(2,2,2,0,2) = {2,4}:{5,9} |
| w(0,0,0,1,2) = {5,11}:{6,9} | w(1,0,0,1,2) = {3,5}:{9,11} | w(2,0,0,1,2) = {1,5}:{7,11} |
| w(0,1,0,1,2) = {1,6}:{5,11} | w(1,1,0,1,2) = {1,4}:{8,9} | w(2,1,0,1,2) = {1,7}:{8,11} |
| w(0,2,0,1,2) = {1,11}:{2,7} | w(1,2,0,1,2) = {2,4}:{7,10} | w(2,2,0,1,2) = {2,9}:{10,11} |
| w(0,0,1,1,2) = {9}:{11} | w(1,0,1,1,2) = {5,10}:{2,11} | w(2,0,1,1,2) = {1}:{7} |
| w(0,1,1,1,2) = {1,3,11}:{2,5,6} | w(1,1,1,1,2) = {2,9}:{10,11} | w(2,1,1,1,2) = {8,9}:{3,11} |
| w(0,2,1,1,2) = {5,10}:{6,11} | w(1,2,1,1,2) = {1,7}:{8,11} | w(2,2,1,1,2) = {7,10}:{3,11} |
| w(0,0,2,1,2) = {1,9}:{8,11} | w(1,0,2,1,2) = {6,7}:{1,11} | w(2,0,2,1,2) = {2}:{10} |
| w(0,1,2,1,2) = {1,3,11}:{2,5,6} | w(1,1,2,1,2) = {2,10}:{9,11} | w(2,1,2,1,2) = {7,10}:{3,11} |
| w(0,2,2,1,2) = {5,9}:{6,11} | w(1,2,2,1,2) = {1,8}:{7,11} | w(2,2,2,1,2) = {8,9}:{3,11} |
| w(0,0,0,2,2) = {3,5}:{1,11} | w(1,0,0,2,2) = {1,10}:{2,11} | w(2,0,0,2,2) = {3,8}:{4,11} |
| w(0,1,0,2,2) = {1,5}:{6,11} | w(1,1,0,2,2) = {5,6}:{7,11} | w(2,1,0,2,2) = {3,6}:{8,9} |
| w(0,2,0,2,2) = {1,11}:{2,8} | w(1,2,0,2,2) = {5,6}:{9,11} | w(2,2,0,2,2) = {3,6}:{7,10} |
| w(0,0,1,2,2) = {1,9}:{8,11} | w(1,0,1,2,2) = {2,8}:{10,11} | w(2,0,1,2,2) = {2,6}:{10,11} |
| w(0,1,1,2,2) = {5,9}:{6,11} | w(1,1,1,2,2) = {2,4}:{3,9} | w(2,1,1,2,2) = {1,2}:{4,11} |
| w(0,2,1,2,2) = {1,3,11}:{2,5,6} | w(1,2,1,2,2) = {1,4}:{3,7} | w(2,2,1,2,2) = {1,2}:{4,11} |
| w(0,0,2,2,2) = {3,6}:{7,9} | w(1,0,2,2,2) = {1,9}:{7,11} | w(2,0,2,2,2) = {1,5}:{7,11} |
| w(0,1,2,2,2) = {5,10}:{6,11} | w(1,1,2,2,2) = {2,4}:{3,10} | w(2,1,2,2,2) = {1,2}:{4,11} |
| w(0,2,2,2,2) = {1,3,11}:{2,5,6} | w(1,2,2,2,2) = {1,4}:{3,8} | w(2,2,2,2,2) = {1,2}:{4,11} |

*The 7-th weighing*

| | | |
|---|---|---|
| w(0,0,0,0,0) = { }:{ } | w(1,0,0,0,0) = {6}:{8} | w(2,0,0,0,0) = {6}:{8} |
| w(0,1,0,0,0) = {6}:{9} | w(1,1,0,0,0) = {10}:{11} | w(2,1,0,0,0) = {8}:{11} |
| w(0,2,0,0,0) = {6}:{9} | w(1,2,0,0,0) = {8}:{11} | w(2,2,0,0,0) = {10}:{11} |
| w(0,0,1,0,0) = {9}:{10} | w(1,0,1,0,0) = {10}:{11} | w(2,0,1,0,0) = {7}:{11} |
| w(0,1,1,0,0) = {9}:{10} | w(1,1,1,0,0) = {5}:{9} | w(2,1,1,0,0) = {6}:{8} |
| w(0,2,1,0,0) = {10}:{9} | w(1,2,1,0,0) = {5}:{7} | w(2,2,1,0,0) = {6}:{10} |
| w(0,0,2,0,0) = {9}:{10} | w(1,0,2,0,0) = {7}:{11} | w(2,0,2,0,0) = {10}:{11} |
| w(0,1,2,0,0) = {10}:{9} | w(1,1,2,0,0) = {6}:{10} | w(2,1,2,0,0) = {5}:{7} |
| w(0,2,2,0,0) = {9}:{10} | w(1,2,2,0,0) = {6}:{8} | w(2,2,2,0,0) = {5}:{9} |
| w(0,0,0,1,0) = {9}:{10} | w(1,0,0,1,0) = {3}:{5} | w(2,0,0,1,0) = {3,4}:{6,7} |
| w(0,1,0,1,0) = {6}:{9} | w(1,1,0,1,0) = {5}:{7} | w(2,1,0,1,0) = {2}:{9} |
| w(0,2,0,1,0) = {5}:{10} | w(1,2,0,1,0) = {5}:{9} | w(2,2,0,1,0) = {1}:{7} |
| w(0,0,1,1,0) = {1}:{2} | w(1,0,1,1,0) = {4}:{6} | w(2,0,1,1,0) = {7}:{11} |
| w(0,1,1,1,0) = {5}:{6} | w(1,1,1,1,0) = {1,2}:{6,7} | w(2,1,1,1,0) = {1,4}:{3,6} |
| w(0,2,1,1,0) = {10}:{9} | w(1,2,1,1,0) = {1,2}:{6,9} | w(2,2,1,1,0) = {2,4}:{3,6} |
| w(0,0,2,1,0) = {1}:{3} | w(1,0,2,1,0) = {4}:{5} | w(2,0,2,1,0) = {10}:{11} |
| w(0,1,2,1,0) = {5}:{6} | w(1,1,2,1,0) = {1,2}:{5,8} | w(2,1,2,1,0) = {1,4}:{3,5} |
| w(0,2,2,1,0) = {9}:{10} | w(1,2,2,1,0) = {1,2}:{5,10} | w(2,2,2,1,0) = {2,4}:{3,5} |
| w(0,0,0,2,0) = {9}:{10} | w(1,0,0,2,0) = {3,4}:{6,7} | w(2,0,0,2,0) = {3}:{5} |
| w(0,1,0,2,0) = {5}:{10} | w(1,1,0,2,0) = {1}:{7} | w(2,1,0,2,0) = {5}:{9} |
| w(0,2,0,2,0) = {6}:{9} | w(1,2,0,2,0) = {2}:{9} | w(2,2,0,2,0) = {5}:{7} |
| w(0,0,1,2,0) = {1}:{3} | w(1,0,1,2,0) = {10}:{11} | w(2,0,1,2,0) = {4}:{5} |
| w(0,1,1,2,0) = {9}:{10} | w(1,1,1,2,0) = {2,4}:{3,5} | w(2,1,1,2,0) = {1,2}:{5,10} |
| w(0,2,1,2,0) = {5}:{6} | w(1,2,1,2,0) = {1,4}:{3,5} | w(2,2,1,2,0) = {1,2}:{5,8} |
| w(0,0,2,2,0) = {1}:{2} | w(1,0,2,2,0) = {7}:{11} | w(2,0,2,2,0) = {4}:{6} |
| w(0,1,2,2,0) = {10}:{9} | w(1,1,2,2,0) = {2,4}:{3,6} | w(2,1,2,2,0) = {1,2}:{6,9} |
| w(0,2,2,2,0) = {5}:{6} | w(1,2,2,2,0) = {1,4}:{3,6} | w(2,2,2,2,0) = {1,2}:{6,7} |
| w(0,0,0,0,1) = {1}:{2} | w(1,0,0,0,1) = {6}:{8} | w(2,0,0,0,1) = {10}:{11} |
| w(0,1,0,0,1) = {1}:{7} | w(1,1,0,0,1) = {3}:{4} | w(2,1,0,0,1) = {1}:{9} |
| w(0,2,0,0,1) = {1}:{7} | w(1,2,0,0,1) = {3}:{4} | w(2,2,0,0,1) = {2}:{7} |
| w(0,0,1,0,1) = {1}:{2} | w(1,0,1,0,1) = {10}:{11} | w(2,0,1,0,1) = {1}:{4} |
| w(0,1,1,0,1) = {9}:{11} | w(1,1,1,0,1) = {5}:{8} | w(2,1,1,0,1) = {5,10}:{8,9} |
| w(0,2,1,0,1) = {10}:{11} | w(1,2,1,0,1) = {5}:{10} | w(2,2,1,0,1) = {5,8}:{7,10} |
| w(0,0,2,0,1) = {1}:{2} | w(1,0,2,0,1) = {7}:{11} | w(2,0,2,0,1) = {2}:{4} |
| w(0,1,2,0,1) = {10}:{11} | w(1,1,2,0,1) = {6}:{7} | w(2,1,2,0,1) = {6,9}:{7,10} |
| w(0,2,2,0,1) = {9}:{11} | w(1,2,2,0,1) = {6}:{9} | w(2,2,2,0,1) = {6,7}:{8,9} |
| w(0,0,0,1,1) = {9}:{10} | w(1,0,0,1,1) = {3}:{4} | w(2,0,0,1,1) = {2}:{11} |
| w(0,1,0,1,1) = {6}:{8} | w(1,1,0,1,1) = {2}:{11} | w(2,1,0,1,1) = {1,4}:{5,6} |
| w(0,2,0,1,1) = {5}:{7} | w(1,2,0,1,1) = {1}:{11} | w(2,2,0,1,1) = {2,4}:{5,6} |
| w(0,0,1,1,1) = {10}:{11} | w(1,0,1,1,1) = {2}:{3} | w(2,0,1,1,1) = {2,3}:{9,10} |
| w(0,1,1,1,1) = {5}:{6} | w(1,1,1,1,1) = {6,7}:{8,9} | w(2,1,1,1,1) = {3}:{8} |
| w(0,2,1,1,1) = {10}:{11} | w(1,2,1,1,1) = {6,9}:{7,10} | w(2,2,1,1,1) = {3}:{10} |
| w(0,0,2,1,1) = {10}:{11} | w(1,0,2,1,1) = {1}:{3} | w(2,0,2,1,1) = {1,3}:{7,8} |
| w(0,1,2,1,1) = {5}:{6} | w(1,1,2,1,1) = {5,8}:{7,10} | w(2,1,2,1,1) = {3}:{7} |

w(0,2,2,1,1,0) = {9}:{11}       w(1,2,2,1,1,0) = {5,10}:{8,9}   w(2,2,2,1,1,0) = {3}:{9}
w(0,0,0,2,1,0) = {5}:{6}        w(1,0,0,2,1,0) = {2}:{11}       w(2,0,0,2,1,0) = {3}:{4}
w(0,1,0,2,1,0) = {5}:{7}        w(1,1,0,2,1,0) = {2}:{7}        w(2,1,0,2,1,0) = {1,5}:{3,4}
w(0,2,0,2,1,0) = {6}:{8}        w(1,2,0,2,1,0) = {1}:{9}        w(2,2,0,2,1,0) = {2,5}:{3,4}
w(0,0,1,2,1,0) = {10}:{11}      w(1,0,1,2,1,0) = {2,4}:{5,6}    w(2,0,1,2,1,0) = {2,4}:{5,6}
w(0,1,1,2,1,0) = {9}:{11}       w(1,1,1,2,1,0) = {2}:{3}        w(2,1,1,2,1,0) = {5,6}:{8,9}
w(0,2,1,2,1,0) = {5}:{6}        w(1,2,1,2,1,0) = {1}:{3}        w(2,2,1,2,1,0) = {5,6}:{7,10}
w(0,0,2,2,1,0) = {2}:{10}       w(1,0,2,2,1,0) = {1,4}:{5,6}    w(2,0,2,2,1,0) = {1,4}:{5,6}
w(0,1,2,2,1,0) = {10}:{11}      w(1,1,2,2,1,0) = {2}:{3}        w(2,1,2,2,1,0) = {5,6}:{7,10}
w(0,2,2,2,1,0) = {5}:{6}        w(1,2,2,2,1,0) = {1}:{3}        w(2,2,2,2,1,0) = {5,6}:{8,9}
w(0,0,0,0,2,0) = {1}:{2}        w(1,0,0,0,2,0) = {10}:{11}      w(2,0,0,0,2,0) = {6}:{8}
w(0,1,0,0,2,0) = {1}:{7}        w(1,1,0,0,2,0) = {2}:{7}        w(2,1,0,0,2,0) = {3}:{4}
w(0,2,0,0,2,0) = {1}:{7}        w(1,2,0,0,2,0) = {1}:{9}        w(2,2,0,0,2,0) = {3}:{4}
w(0,0,1,0,2,0) = {1}:{2}        w(1,0,1,0,2,0) = {2}:{4}        w(2,0,1,0,2,0) = {7}:{11}
w(0,1,1,0,2,0) = {9}:{11}       w(1,1,1,0,2,0) = {6,7}:{8,9}    w(2,1,1,0,2,0) = {6}:{9}
w(0,2,1,0,2,0) = {10}:{11}      w(1,2,1,0,2,0) = {6,9}:{7,10}   w(2,2,1,0,2,0) = {6}:{7}
w(0,0,2,0,2,0) = {1}:{2}        w(1,0,2,0,2,0) = {1}:{4}        w(2,0,2,0,2,0) = {10}:{11}
w(0,1,2,0,2,0) = {10}:{11}      w(1,1,2,0,2,0) = {5,8}:{7,10}   w(2,1,2,0,2,0) = {5}:{10}
w(0,2,2,0,2,0) = {9}:{11}       w(1,2,2,0,2,0) = {5,10}:{8,9}   w(2,2,2,0,2,0) = {5}:{8}
w(0,0,0,1,2,0) = {5}:{6}        w(1,0,0,1,2,0) = {3}:{4}        w(2,0,0,1,2,0) = {2}:{11}
w(0,1,0,1,2,0) = {6}:{8}        w(1,1,0,1,2,0) = {2,5}:{3,4}    w(2,1,0,1,2,0) = {1}:{9}
w(0,2,0,1,2,0) = {5}:{7}        w(1,2,0,1,2,0) = {1,5}:{3,4}    w(2,2,0,1,2,0) = {2}:{7}
w(0,0,1,1,2,0) = {2}:{10}       w(1,0,1,1,2,0) = {1,4}:{5,6}    w(2,0,1,1,2,0) = {1,4}:{5,6}
w(0,1,1,1,2,0) = {5}:{6}        w(1,1,1,1,2,0) = {5,6}:{8,9}    w(2,1,1,1,2,0) = {1}:{3}
w(0,2,1,1,2,0) = {10}:{11}      w(1,2,1,1,2,0) = {5,6}:{7,10}   w(2,2,1,1,2,0) = {2}:{3}
w(0,0,2,1,2,0) = {10}:{11}      w(1,0,2,1,2,0) = {2,4}:{5,6}    w(2,0,2,1,2,0) = {2,4}:{5,6}
w(0,1,2,1,2,0) = {5}:{6}        w(1,1,2,1,2,0) = {5,6}:{7,10}   w(2,1,2,1,2,0) = {1}:{3}
w(0,2,2,1,2,0) = {9}:{11}       w(1,2,2,1,2,0) = {5,6}:{8,9}    w(2,2,2,1,2,0) = {2}:{3}
w(0,0,0,2,2,0) = {9}:{10}       w(1,0,0,2,2,0) = {2}:{11}       w(2,0,0,2,2,0) = {3}:{4}
w(0,1,0,2,2,0) = {5}:{7}        w(1,1,0,2,2,0) = {2,4}:{5,6}    w(2,1,0,2,2,0) = {1}:{11}
w(0,2,0,2,2,0) = {6}:{8}        w(1,2,0,2,2,0) = {1,4}:{5,6}    w(2,2,0,2,2,0) = {2}:{11}
w(0,0,1,2,2,0) = {10}:{11}      w(1,0,1,2,2,0) = {1,3}:{7,8}    w(2,0,1,2,2,0) = {1}:{3}
w(0,1,1,2,2,0) = {9}:{11}       w(1,1,1,2,2,0) = {3}:{9}        w(2,1,1,2,2,0) = {5,10}:{8,9}
w(0,2,1,2,2,0) = {5}:{6}        w(1,2,1,2,2,0) = {3}:{7}        w(2,2,1,2,2,0) = {5,8}:{7,10}
w(0,0,2,2,2,0) = {10}:{11}      w(1,0,2,2,2,0) = {2,3}:{9,10}   w(2,0,2,2,2,0) = {2}:{3}
w(0,1,2,2,2,0) = {10}:{11}      w(1,1,2,2,2,0) = {3}:{10}       w(2,1,2,2,2,0) = {6,9}:{7,10}
w(0,2,2,2,2,0) = {5}:{6}        w(1,2,2,2,2,0) = {3}:{8}        w(2,2,2,2,2,0) = {6,7}:{8,9}
w(0,0,0,0,0,1) = {10}:{11}      w(1,0,0,0,0,1) = {8}:{11}       w(2,0,0,0,0,1) = {8}:{11}
w(0,1,0,0,0,1) = {1}:{11}       w(1,1,0,0,0,1) = {1}:{8}        w(2,1,0,0,0,1) = {1}:{9}
w(0,2,0,0,0,1) = {1}:{11}       w(1,2,0,0,0,1) = {2}:{10}       w(2,2,0,0,0,1) = {2}:{7}
w(0,0,1,0,0,1) = {3}:{10}       w(1,0,1,0,0,1) = {3}:{11}       w(2,0,1,0,0,1) = {7}:{11}
w(0,1,1,0,0,1) = {9}:{11}       w(1,1,1,0,0,1) = {5}:{9}        w(2,1,1,0,0,1) = {3,5}:{6,8}
w(0,2,1,0,0,1) = {10}:{11}      w(1,2,1,0,0,1) = {5}:{7}        w(2,2,1,0,0,1) = {3,5}:{6,10}
w(0,0,2,0,0,1) = {3}:{4}        w(1,0,2,0,0,1) = {3}:{11}       w(2,0,2,0,0,1) = {10}:{11}

w(0,1,2,0,0,1) = {10}:{11}   w(1,1,2,0,0,1) = {6}:{10}        w(2,1,2,0,0,1) = {3,6}:{5,7}
w(0,2,2,0,0,1) = {9}:{11}    w(1,2,2,0,0,1) = {6}:{8}         w(2,2,2,0,0,1) = {3,6}:{5,9}
w(0,0,0,1,0,1) = {10}:{11}   w(1,0,0,1,0,1) = {9}:{11}        w(2,0,0,1,0,1) = {10}:{11}
w(0,1,0,1,0,1) = {8}:{11}    w(1,1,0,1,0,1) = {8}:{10}        w(2,1,0,1,0,1) = {1}:{11}
w(0,2,0,1,0,1) = {2}:{11}    w(1,2,0,1,0,1) = {10}:{8}        w(2,2,0,1,0,1) = {2}:{11}
w(0,0,1,1,0,1) = {2,10}:{3,11}  w(1,0,1,1,0,1) = {4}:{6}      w(2,0,1,1,0,1) = {3,4}:{9,10}
w(0,1,1,1,0,1) = {5}:{6}     w(1,1,1,1,0,1) = {2}:{9}         w(2,1,1,1,0,1) = {1}:{3}
w(0,2,1,1,0,1) = {10}:{11}   w(1,2,1,1,0,1) = {1}:{7}         w(2,2,1,1,0,1) = {2}:{3}
w(0,0,2,1,0,1) = {9}:{11}    w(1,0,2,1,0,1) = {4}:{5}         w(2,0,2,1,0,1) = {3,4}:{7,8}
w(0,1,2,1,0,1) = {5}:{6}     w(1,1,2,1,0,1) = {2}:{10}        w(2,1,2,1,0,1) = {1}:{3}
w(0,2,2,1,0,1) = {9}:{11}    w(1,2,2,1,0,1) = {1}:{8}         w(2,2,2,1,0,1) = {2}:{3}
w(0,0,0,2,0,1) = {10}:{11}   w(1,0,0,2,0,1) = {4}:{5}         w(2,0,0,2,0,1) = {9}:{11}
w(0,1,0,2,0,1) = {7}:{11}    w(1,1,0,2,0,1) = {6}:{8}         w(2,1,0,2,0,1) = {1,4}:{3,10}
w(0,2,0,2,0,1) = {2}:{11}    w(1,2,0,2,0,1) = {6}:{10}        w(2,2,0,2,0,1) = {2,4}:{3,8}
w(0,0,1,2,0,1) = {10}:{11}   w(1,0,1,2,0,1) = {1}:{3}         w(2,0,1,2,0,1) = {4}:{11}
w(0,1,1,2,0,1) = {9}:{11}    w(1,1,1,2,0,1) = {2,4}:{3,5}     w(2,1,1,2,0,1) = {9}:{11}
w(0,2,1,2,0,1) = {9}:{11}    w(1,2,1,2,0,1) = {1,4}:{3,5}     w(2,2,1,2,0,1) = {7}:{11}
w(0,0,2,2,0,1) = {2}:{11}    w(1,0,2,2,0,1) = {2}:{3}         w(2,0,2,2,0,1) = {4}:{11}
w(0,1,2,2,0,1) = {10}:{11}   w(1,1,2,2,0,1) = {2,4}:{3,6}     w(2,1,2,2,0,1) = {10}:{11}
w(0,2,2,2,0,1) = {10}:{11}   w(1,2,2,2,0,1) = {1,4}:{3,6}     w(2,2,2,2,0,1) = {8}:{11}
w(0,0,0,0,1,1) = {}:{}       w(1,0,0,0,1,1) = {8}:{11}        w(2,0,0,0,1,1) = {8}:{11}
w(0,1,0,0,1,1) = {1}:{7}     w(1,1,0,0,1,1) = {2}:{8}         w(2,1,0,0,1,1) = {2}:{9}
w(0,2,0,0,1,1) = {2}:{11}    w(1,2,0,0,1,1) = {1}:{10}        w(2,2,0,0,1,1) = {1}:{7}
w(0,0,1,0,1,1) = {2}:{3}     w(1,0,1,0,1,1) = {5,9}:{7,11}    w(2,0,1,0,1,1) = {2,11}:{4,10}
w(0,1,1,0,1,1) = {2,3}:{9,11}  w(1,1,1,0,1,1) = {3}:{8}       w(2,1,1,0,1,1) = {9}:{7}
w(0,2,1,0,1,1) = {2,3}:{10,11}  w(1,2,1,0,1,1) = {3}:{10}     w(2,2,1,0,1,1) = {7}:{9}
w(0,0,2,0,1,1) = {2}:{3}     w(1,0,2,0,1,1) = {6,8}:{10,11}   w(2,0,2,0,1,1) = {1,11}:{4,7}
w(0,1,2,0,1,1) = {2,3}:{10,11}  w(1,1,2,0,1,1) = {3}:{7}      w(2,1,2,0,1,1) = {10}:{8}
w(0,2,2,0,1,1) = {2,3}:{9,11}  w(1,2,2,0,1,1) = {3}:{9}       w(2,2,2,0,1,1) = {8}:{10}
w(0,0,0,1,1,1) = {9}:{10}    w(1,0,0,1,1,1) = {3,4}:{9,11}    w(2,0,0,1,1,1) = {2}:{11}
w(0,1,0,1,1,1) = {4,5}:{10,11}  w(1,1,0,1,1,1) = {6}:{11}     w(2,1,0,1,1,1) = {2}:{10}
w(0,2,0,1,1,1) = {1,4}:{9,11}  w(1,2,0,1,1,1) = {6}:{11}      w(2,2,0,1,1,1) = {1}:{8}
w(0,0,1,1,1,1) = {1,4}:{10,11}  w(1,0,1,1,1,1) = {4}:{6}      w(2,0,1,1,1,1) = {7}:{11}
w(0,1,1,1,1,1) = {5}:{6}     w(1,1,1,1,1,1) = {2}:{9}         w(2,1,1,1,1,1) = {1,3}:{4,9}
w(0,2,1,1,1,1) = {2,3}:{10,11}  w(1,2,1,1,1,1) = {1}:{7}      w(2,2,1,1,1,1) = {2,3}:{4,7}
w(0,0,2,1,1,1) = {6,9}:{7,10}  w(1,0,2,1,1,1) = {4}:{5}       w(2,0,2,1,1,1) = {10}:{11}
w(0,1,2,1,1,1) = {5}:{6}     w(1,1,2,1,1,1) = {2}:{10}        w(2,1,2,1,1,1) = {1,3}:{4,10}
w(0,2,2,1,1,1) = {2,3}:{9,11}  w(1,2,2,1,1,1) = {1}:{8}       w(2,2,2,1,1,1) = {2,3}:{4,8}
w(0,0,0,2,1,1) = {10}:{11}   w(1,0,0,2,1,1) = {2}:{4}         w(2,0,0,2,1,1) = {9}:{11}
w(0,1,0,2,1,1) = {4,6}:{9,11}  w(1,1,0,2,1,1) = {7}:{8}       w(2,1,0,2,1,1) = {1,2}:{3,11}
w(0,2,0,2,1,1) = {1,4}:{10,11}  w(1,2,0,2,1,1) = {9}:{10}     w(2,2,0,2,1,1) = {1,2}:{3,11}
w(0,0,1,2,1,1) = {1}:{7}     w(1,0,1,2,1,1) = {4}:{6}         w(2,0,1,2,1,1) = {4}:{6}
w(0,1,1,2,1,1) = {2,3}:{9,11}  w(1,1,1,2,1,1) = {3,6}:{4,9}   w(2,1,1,2,1,1) = {5,10}:{8,11}
w(0,2,1,2,1,1) = {9}:{11}    w(1,2,1,2,1,1) = {3,6}:{4,7}     w(2,2,1,2,1,1) = {5,8}:{10,11}

| | | |
|---|---|---|
| w(0,0,2,2,1,1) = {7}:{8} | w(1,0,2,2,1,1) = {4}:{5} | w(2,0,2,2,1,1) = {4}:{5} |
| w(0,1,2,2,1,1) = {2,3}:{10,11} | w(1,1,2,2,1,1) = {3,5}:{4,10} | w(2,1,2,2,1,1) = {6,9}:{7,11} |
| w(0,2,2,2,1,1) = {10}:{11} | w(1,2,2,2,1,1) = {3,5}:{4,8} | w(2,2,2,2,1,1) = {6,7}:{9,11} |
| w(0,0,0,0,2,1) = { }:{ } | w(1,0,0,0,2,1) = {8}:{11} | w(2,0,0,0,2,1) = {8}:{11} |
| w(0,1,0,0,2,1) = {1}:{7} | w(1,1,0,0,2,1) = {6}:{10} | w(2,1,0,0,2,1) = {5}:{9} |
| w(0,2,0,0,2,1) = {3}:{11} | w(1,2,0,0,2,1) = {6}:{8} | w(2,2,0,0,2,1) = {5}:{7} |
| w(0,0,1,0,2,1) = {2}:{3} | w(1,0,1,0,2,1) = {7}:{11} | w(2,0,1,0,2,1) = {3}:{11} |
| w(0,1,1,0,2,1) = {2,3}:{9,11} | w(1,1,1,0,2,1) = {1}:{8} | w(2,1,1,0,2,1) = {5,6}:{9,10} |
| w(0,2,1,0,2,1) = {2,3}:{10,11} | w(1,2,1,0,2,1) = {2}:{10} | w(2,2,1,0,2,1) = {5,6}:{7,8} |
| w(0,0,2,0,2,1) = {2}:{3} | w(1,0,2,0,2,1) = {10}:{11} | w(2,0,2,0,2,1) = {3}:{11} |
| w(0,1,2,0,2,1) = {2,3}:{10,11} | w(1,1,2,0,2,1) = {1}:{7} | w(2,1,2,0,2,1) = {5,6}:{9,10} |
| w(0,2,2,0,2,1) = {2,3}:{9,11} | w(1,2,2,0,2,1) = {2}:{9} | w(2,2,2,0,2,1) = {5,6}:{7,8} |
| w(0,0,0,1,2,1) = {10}:{11} | w(1,0,0,1,2,1) = {9}:{11} | w(2,0,0,1,2,1) = {1,11}:{3,4} |
| w(0,1,0,1,2,1) = {5}:{11} | w(1,1,0,1,2,1) = {3}:{4} | w(2,1,0,1,2,1) = {10}:{8} |
| w(0,2,0,1,2,1) = {1}:{11} | w(1,2,0,1,2,1) = {3}:{4} | w(2,2,0,1,2,1) = {8}:{10} |
| w(0,0,1,1,2,1) = {2}:{3} | w(1,0,1,1,2,1) = {4}:{6} | w(2,0,1,1,2,1) = {10}:{11} |
| w(0,1,1,1,2,1) = {5}:{6} | w(1,1,1,1,2,1) = {1}:{2} | w(2,1,1,1,2,1) = {1}:{11} |
| w(0,2,1,1,2,1) = {2,3}:{10,11} | w(1,2,1,1,2,1) = {2}:{1} | w(2,2,1,1,2,1) = {2}:{11} |
| w(0,0,2,1,2,1) = {5,9}:{7,10} | w(1,0,2,1,2,1) = {4}:{5} | w(2,0,2,1,2,1) = {7}:{11} |
| w(0,1,2,1,2,1) = {5}:{6} | w(1,1,2,1,2,1) = {1}:{2} | w(2,1,2,1,2,1) = {1}:{11} |
| w(0,2,2,1,2,1) = {2,3}:{9,11} | w(1,2,2,1,2,1) = {2}:{1} | w(2,2,2,1,2,1) = {2}:{11} |
| w(0,0,0,2,2,1) = {6}:{11} | w(1,0,0,2,2,1) = {7}:{11} | w(2,0,0,2,2,1) = {3,4}:{9,11} |
| w(0,1,0,2,2,1) = {6}:{11} | w(1,1,0,2,2,1) = {2,6}:{4,5} | w(2,1,0,2,2,1) = {3,11}:{4,8} |
| w(0,2,0,2,2,1) = {1}:{11} | w(1,2,0,2,2,1) = {1,6}:{4,5} | w(2,2,0,2,2,1) = {3,11}:{4,10} |
| w(0,0,1,2,2,1) = {1}:{7} | w(1,0,1,2,2,1) = {7}:{8} | w(2,0,1,2,2,1) = {3,4}:{5,6} |
| w(0,1,1,2,2,1) = {2,3}:{9,11} | w(1,1,1,2,2,1) = {3,4}:{8,9} | w(2,1,1,2,2,1) = {5,10}:{8,9} |
| w(0,2,1,2,2,1) = {9}:{11} | w(1,2,1,2,2,1) = {3,4}:{7,10} | w(2,2,1,2,2,1) = {5,8}:{7,10} |
| w(0,0,2,2,2,1) = {9}:{11} | w(1,0,2,2,2,1) = {10}:{9} | w(2,0,2,2,2,1) = {3,4}:{5,6} |
| w(0,1,2,2,2,1) = {2,3}:{10,11} | w(1,1,2,2,2,1) = {3,4}:{7,10} | w(2,1,2,2,2,1) = {6,9}:{7,10} |
| w(0,2,2,2,2,1) = {10}:{11} | w(1,2,2,2,2,1) = {3,4}:{8,9} | w(2,2,2,2,2,1) = {6,7}:{8,9} |
| w(0,0,0,0,0,2) = {10}:{11} | w(1,0,0,0,0,2) = {8}:{11} | w(2,0,0,0,0,2) = {8}:{11} |
| w(0,1,0,0,0,2) = {1}:{11} | w(1,1,0,0,0,2) = {2}:{7} | w(2,1,0,0,0,2) = {2}:{10} |
| w(0,2,0,0,0,2) = {1}:{11} | w(1,2,0,0,0,2) = {1}:{9} | w(2,2,0,0,0,2) = {1}:{8} |
| w(0,0,1,0,0,2) = {3}:{4} | w(1,0,1,0,0,2) = {10}:{11} | w(2,0,1,0,0,2) = {3}:{11} |
| w(0,1,1,0,0,2) = {9}:{11} | w(1,1,1,0,0,2) = {3,6}:{5,9} | w(2,1,1,0,0,2) = {6}:{8} |
| w(0,2,1,0,0,2) = {10}:{11} | w(1,2,1,0,0,2) = {3,6}:{5,7} | w(2,2,1,0,0,2) = {6}:{10} |
| w(0,0,2,0,0,2) = {3}:{10} | w(1,0,2,0,0,2) = {7}:{11} | w(2,0,2,0,0,2) = {3}:{11} |
| w(0,1,2,0,0,2) = {10}:{11} | w(1,1,2,0,0,2) = {3,5}:{6,10} | w(2,1,2,0,0,2) = {5}:{7} |
| w(0,2,2,0,0,2) = {9}:{11} | w(1,2,2,0,0,2) = {3,5}:{6,8} | w(2,2,2,0,0,2) = {5}:{9} |
| w(0,0,0,1,0,2) = {10}:{11} | w(1,0,0,1,0,2) = {9}:{11} | w(2,0,0,1,0,2) = {4}:{5} |
| w(0,1,0,1,0,2) = {2}:{11} | w(1,1,0,1,0,2) = {2,4}:{3,8} | w(2,1,0,1,0,2) = {6}:{10} |
| w(0,2,0,1,0,2) = {7}:{11} | w(1,2,0,1,0,2) = {1,4}:{3,10} | w(2,2,0,1,0,2) = {6}:{8} |
| w(0,0,1,1,0,2) = {2}:{11} | w(1,0,1,1,0,2) = {4}:{11} | w(2,0,1,1,0,2) = {2}:{3} |
| w(0,1,1,1,0,2) = {10}:{11} | w(1,1,1,1,0,2) = {8}:{11} | w(2,1,1,1,0,2) = {1,4}:{3,6} |

w(0,2,1,1,0,2) = {10}:{11}        w(1,2,1,1,0,2) = {10}:{11}        w(2,2,1,1,0,2) = {2,4}:{3,6}
w(0,0,2,1,0,2) = {10}:{11}        w(1,0,2,1,0,2) = {4}:{11}         w(2,0,2,1,0,2) = {1}:{3}
w(0,1,2,1,0,2) = {9}:{11}         w(1,1,2,1,0,2) = {7}:{11}         w(2,1,2,1,0,2) = {1,4}:{3,5}
w(0,2,2,1,0,2) = {9}:{11}         w(1,2,2,1,0,2) = {9}:{11}         w(2,2,2,1,0,2) = {2,4}:{3,5}
w(0,0,0,2,0,2) = {10}:{11}        w(1,0,0,2,0,2) = {10}:{11}        w(2,0,0,2,0,2) = {9}:{11}
w(0,1,0,2,0,2) = {2}:{11}         w(1,1,0,2,0,2) = {2}:{11}         w(2,1,0,2,0,2) = {10}:{8}
w(0,2,0,2,0,2) = {8}:{11}         w(1,2,0,2,0,2) = {1}:{11}         w(2,2,0,2,0,2) = {8}:{10}
w(0,0,1,2,0,2) = {9}:{11}         w(1,0,1,2,0,2) = {3,4}:{7,8}      w(2,0,1,2,0,2) = {4}:{5}
w(0,1,1,2,0,2) = {9}:{11}         w(1,1,1,2,0,2) = {2}:{3}          w(2,1,1,2,0,2) = {1}:{8}
w(0,2,1,2,0,2) = {5}:{6}          w(1,2,1,2,0,2) = {1}:{3}          w(2,2,1,2,0,2) = {2}:{10}
w(0,0,2,2,0,2) = {2,10}:{3,11}    w(1,0,2,2,0,2) = {3,4}:{9,10}     w(2,0,2,2,0,2) = {4}:{6}
w(0,1,2,2,0,2) = {10}:{11}        w(1,1,2,2,0,2) = {2}:{3}          w(2,1,2,2,0,2) = {1}:{7}
w(0,2,2,2,0,2) = {5}:{6}          w(1,2,2,2,0,2) = {1}:{3}          w(2,2,2,2,0,2) = {2}:{9}
w(0,0,0,0,1,2) = { }:{ }          w(1,0,0,0,1,2) = {8}:{11}         w(2,0,0,0,1,2) = {8}:{11}
w(0,1,0,0,1,2) = {3}:{11}         w(1,1,0,0,1,2) = {5}:{7}          w(2,1,0,0,1,2) = {6}:{8}
w(0,2,0,0,1,2) = {1}:{7}          w(1,2,0,0,1,2) = {5}:{9}          w(2,2,0,0,1,2) = {6}:{10}
w(0,0,1,0,1,2) = {2}:{3}          w(1,0,1,0,1,2) = {3}:{11}         w(2,0,1,0,1,2) = {10}:{11}
w(0,1,1,0,1,2) = {2,3}:{9,11}     w(1,1,1,0,1,2) = {5,6}:{7,8}      w(2,1,1,0,1,2) = {2}:{9}
w(0,2,1,0,1,2) = {2,3}:{10,11}    w(1,2,1,0,1,2) = {5,6}:{9,10}     w(2,2,1,0,1,2) = {1}:{7}
w(0,0,2,0,1,2) = {2}:{3}          w(1,0,2,0,1,2) = {3}:{11}         w(2,0,2,0,1,2) = {7}:{11}
w(0,1,2,0,1,2) = {2,3}:{10,11}    w(1,1,2,0,1,2) = {5,6}:{7,8}      w(2,1,2,0,1,2) = {2}:{10}
w(0,2,2,0,1,2) = {2,3}:{9,11}     w(1,2,2,0,1,2) = {5,6}:{9,10}     w(2,2,2,0,1,2) = {1}:{8}
w(0,0,0,1,1,2) = {6}:{11}         w(1,0,0,1,1,2) = {3,4}:{9,11}     w(2,0,0,1,1,2) = {7}:{11}
w(0,1,0,1,1,2) = {1}:{11}         w(1,1,0,1,1,2) = {3,11}:{4,10}    w(2,1,0,1,1,2) = {1,6}:{4,5}
w(0,2,0,1,1,2) = {6}:{11}         w(1,2,0,1,1,2) = {3,11}:{4,8}     w(2,2,0,1,1,2) = {2,6}:{4,5}
w(0,0,1,1,1,2) = {9}:{11}         w(1,0,1,1,1,2) = {3,4}:{5,6}      w(2,0,1,1,1,2) = {10}:{9}
w(0,1,1,1,1,2) = {10}:{11}        w(1,1,1,1,1,2) = {6,7}:{8,9}      w(2,1,1,1,1,2) = {3,4}:{8,9}
w(0,2,1,1,1,2) = {2,3}:{10,11}    w(1,2,1,1,1,2) = {6,9}:{7,10}     w(2,2,1,1,1,2) = {3,4}:{7,10}
w(0,0,2,1,1,2) = {1}:{7}          w(1,0,2,1,1,2) = {3,4}:{5,6}      w(2,0,2,1,1,2) = {7}:{8}
w(0,1,2,1,1,2) = {9}:{11}         w(1,1,2,1,1,2) = {5,8}:{7,10}     w(2,1,2,1,1,2) = {3,4}:{7,10}
w(0,2,2,1,1,2) = {2,3}:{9,11}     w(1,2,2,1,1,2) = {5,10}:{8,9}     w(2,2,2,1,1,2) = {3,4}:{8,9}
w(0,0,0,2,1,2) = {10}:{11}        w(1,0,0,2,1,2) = {1,11}:{3,4}     w(2,0,0,2,1,2) = {9}:{11}
w(0,1,0,2,1,2) = {1}:{11}         w(1,1,0,2,1,2) = {8}:{10}         w(2,1,0,2,1,2) = {3}:{4}
w(0,2,0,2,1,2) = {5}:{11}         w(1,2,0,2,1,2) = {10}:{8}         w(2,2,0,2,1,2) = {3}:{4}
w(0,0,1,2,1,2) = {5,9}:{7,10}     w(1,0,1,2,1,2) = {7}:{11}         w(2,0,1,2,1,2) = {4}:{5}
w(0,1,1,2,1,2) = {2,3}:{9,11}     w(1,1,1,2,1,2) = {2}:{11}         w(2,1,1,2,1,2) = {2}:{1}
w(0,2,1,2,1,2) = {5}:{6}          w(1,2,1,2,1,2) = {1}:{11}         w(2,2,1,2,1,2) = {1}:{2}
w(0,0,2,2,1,2) = {2}:{3}          w(1,0,2,2,1,2) = {10}:{11}        w(2,0,2,2,1,2) = {4}:{6}
w(0,1,2,2,1,2) = {2,3}:{10,11}    w(1,1,2,2,1,2) = {2}:{11}         w(2,1,2,2,1,2) = {2}:{1}
w(0,2,2,2,1,2) = {5}:{6}          w(1,2,2,2,1,2) = {1}:{11}         w(2,2,2,2,1,2) = {1}:{2}
w(0,0,0,0,2,2) = { }:{ }          w(1,0,0,0,2,2) = {8}:{11}         w(2,0,0,0,2,2) = {8}:{11}
w(0,1,0,0,2,2) = {2}:{11}         w(1,1,0,0,2,2) = {1}:{7}          w(2,1,0,0,2,2) = {1}:{10}
w(0,2,0,0,2,2) = {1}:{7}          w(1,2,0,0,2,2) = {2}:{9}          w(2,2,0,0,2,2) = {2}:{8}
w(0,0,1,0,2,2) = {2}:{3}          w(1,0,1,0,2,2) = {1,11}:{4,7}     w(2,0,1,0,2,2) = {6,8}:{10,11}

| | | |
|---|---|---|
| w(0,1,1,0,2,2) = {2,3}:{9,11} | w(1,1,1,0,2,2) = {8}:{10} | w(2,1,1,0,2,2) = {3}:{9} |
| w(0,2,1,0,2,2) = {2,3}:{10,11} | w(1,2,1,0,2,2) = {10}:{8} | w(2,2,1,0,2,2) = {3}:{7} |
| w(0,0,2,0,2,2) = {2}:{3} | w(1,0,2,0,2,2) = {2,11}:{4,10} | w(2,0,2,0,2,2) = {5,9}:{7,11} |
| w(0,1,2,0,2,2) = {2,3}:{10,11} | w(1,1,2,0,2,2) = {7}:{9} | w(2,1,2,0,2,2) = {3}:{10} |
| w(0,2,2,0,2,2) = {2,3}:{9,11} | w(1,2,2,0,2,2) = {9}:{7} | w(2,2,2,0,2,2) = {3}:{8} |
| w(0,0,0,1,2,2) = {10}:{11} | w(1,0,0,1,2,2) = {9}:{11} | w(2,0,0,1,2,2) = {2}:{4} |
| w(0,1,0,1,2,2) = {1,4}:{10,11} | w(1,1,0,1,2,2) = {1,2}:{3,11} | w(2,1,0,1,2,2) = {9}:{10} |
| w(0,2,0,1,2,2) = {4,6}:{9,11} | w(1,2,0,1,2,2) = {1,2}:{3,11} | w(2,2,0,1,2,2) = {7}:{8} |
| w(0,0,1,1,2,2) = {7}:{8} | w(1,0,1,1,2,2) = {4}:{5} | w(2,0,1,1,2,2) = {4}:{5} |
| w(0,1,1,1,2,2) = {10}:{11} | w(1,1,1,1,2,2) = {6,7}:{9,11} | w(2,1,1,1,2,2) = {3,5}:{4,8} |
| w(0,2,1,1,2,2) = {2,3}:{10,11} | w(1,2,1,1,2,2) = {6,9}:{7,11} | w(2,2,1,1,2,2) = {3,5}:{4,10} |
| w(0,0,2,1,2,2) = {1}:{7} | w(1,0,2,1,2,2) = {4}:{6} | w(2,0,2,1,2,2) = {4}:{6} |
| w(0,1,2,1,2,2) = {9}:{11} | w(1,1,2,1,2,2) = {5,8}:{10,11} | w(2,1,2,1,2,2) = {3,6}:{4,7} |
| w(0,2,2,1,2,2) = {2,3}:{9,11} | w(1,2,2,1,2,2) = {5,10}:{8,11} | w(2,2,2,1,2,2) = {3,6}:{4,9} |
| w(0,0,0,2,2,2) = {9}:{10} | w(1,0,0,2,2,2) = {2}:{11} | w(2,0,0,2,2,2) = {3,4}:{9,11} |
| w(0,1,0,2,2,2) = {1,4}:{9,11} | w(1,1,0,2,2,2) = {1}:{8} | w(2,1,0,2,2,2) = {6}:{11} |
| w(0,2,0,2,2,2) = {4,5}:{10,11} | w(1,2,0,2,2,2) = {2}:{10} | w(2,2,0,2,2,2) = {6}:{11} |
| w(0,0,1,2,2,2) = {6,9}:{7,10} | w(1,0,1,2,2,2) = {10}:{11} | w(2,0,1,2,2,2) = {4}:{5} |
| w(0,1,1,2,2,2) = {2,3}:{9,11} | w(1,1,1,2,2,2) = {2,3}:{4,8} | w(2,1,1,2,2,2) = {1}:{8} |
| w(0,2,1,2,2,2) = {5}:{6} | w(1,2,1,2,2,2) = {1,3}:{4,10} | w(2,2,1,2,2,2) = {2}:{10} |
| w(0,0,2,2,2,2) = {1,4}:{10,11} | w(1,0,2,2,2,2) = {7}:{11} | w(2,0,2,2,2,2) = {4}:{6} |
| w(0,1,2,2,2,2) = {2,3}:{10,11} | w(1,1,2,2,2,2) = {2,3}:{4,7} | w(2,1,2,2,2,2) = {1}:{7} |
| w(0,2,2,2,2,2) = {5}:{6} | w(1,2,2,2,2,2) = {1,3}:{4,9} | w(2,2,2,2,2,2) = {2}:{9} |

List 1

1:= "<", 0:= "=", 2:= ">", and w(x) is the weighing in the direction x.

The completed map of the first algorithm to sort 11 coins

| | | |
|---|---|---|
| f(0,0,0,0,0,0)    = 2047 | f(1,0,0,0,0,0) = 1468 | f(2,0,0,0,0,0) = 579 |
| f(0,1,0,0,0,0) = 1855 | f(1,1,0,0,0,0) = 2046 | f(2,1,0,0,0,0) = 2 |
| f(0,2,0,0,0,0) = 192 | f(1,2,0,0,0,0) = 2045 | f(2,2,0,0,0,0) = 1 |
| f(0,0,1,0,0,0) = 138 | f(1,0,1,0,0,0) = 1786 | f(2,0,1,0,0,0) = 134 |
| f(0,1,1,0,0,0) = 1974 | f(1,1,1,0,0,0) = 944 | f(2,1,1,0,0,0) = 1551 |
| f(0,2,1,0,0,0) = 137 | f(1,2,1,0,0,0) = 752 | f(2,2,1,0,0,0) = 1167 |
| f(0,0,2,0,0,0) = 1909 | f(1,0,2,0,0,0) = 1913 | f(2,0,2,0,0,0) = 261 |
| f(0,1,2,0,0,0) = 1910 | f(1,1,2,0,0,0) = 880 | f(2,1,2,0,0,0) = 1295 |
| f(0,2,2,0,0,0) = 73 | f(1,2,2,0,0,0) = 496 | f(2,2,2,0,0,0) = 1103 |
| f(1,0,0,1,0,0) = 2044 | f(2,0,0,1,0,0) = 1991 | f(0,1,0,1,0,0) = 1838 |
| f(1,1,0,1,0,0) = 300 | f(2,1,0,1,0,0) = 774 | f(0,2,0,1,0,0) = 225 |
| f(1,2,0,1,0,0) = 108 | f(2,2,0,1,0,0) = 197 | f(1,0,1,1,0,0) = 696 |
| f(2,0,1,1,0,0) = 2031 | f(0,1,1,1,0,0) = 1674 | f(1,1,1,1,0,0) = 1594 |
| f(2,1,1,1,0,0) = 1827 | f(0,2,1,1,0,0) = 128 | f(1,2,1,1,0,0) = 1209 |
| f(2,2,1,1,0,0) = 1251 | f(1,0,2,1,0,0) = 376 | f(2,0,2,1,0,0) = 2015 |
| f(0,1,2,1,0,0) = 1354 | f(1,1,2,1,0,0) = 1338 | f(2,1,2,1,0,0) = 1811 |
| f(0,2,2,1,0,0) = 64 | f(1,2,2,1,0,0) = 1145 | f(2,2,2,1,0,0) = 1235 |
| f(1,0,0,2,0,0) = 56 | f(2,0,0,2,0,0) = 3 | f(0,1,0,2,0,0) = 1822 |
| f(1,1,0,2,0,0) = 1850 | f(2,1,0,2,0,0) = 1939 | f(0,2,0,2,0,0) = 209 |
| f(1,2,0,2,0,0) = 1273 | f(2,2,0,2,0,0) = 1747 | f(1,0,1,2,0,0) = 32 |
| f(2,0,1,2,0,0) = 1671 | f(0,1,1,2,0,0) = 1983 | f(1,1,1,2,0,0) = 812 |
| f(2,1,1,2,0,0) = 902 | f(0,2,1,2,0,0) = 693 | f(1,2,1,2,0,0) = 236 |
| f(2,2,1,2,0,0) = 709 | f(1,0,2,2,0,0) = 16 | f(2,0,2,2,0,0) = 1351 |
| f(0,1,2,2,0,0) = 1919 | f(1,1,2,2,0,0) = 796 | f(2,1,2,2,0,0) = 838 |
| f(0,2,2,2,0,0) = 373 | f(1,2,2,2,0,0) = 220 | f(2,2,2,2,0,0) = 453 |
| f(0,0,0,0,1,0) = 1984 | f(1,0,0,0,1,0) = 1534 | f(2,0,0,0,1,0) = 66 |
| f(0,1,0,0,1,0) = 1804 | f(1,1,0,0,1,0) = 1854 | f(2,1,0,0,1,0) = 1943 |
| f(0,2,0,0,1,0) = 245 | f(1,2,0,0,1,0) = 1277 | f(2,2,0,0,1,0) = 1751 |
| f(0,0,1,0,1,0) = 1060 | f(2,0,1,0,1,0) = 710 | f(1,1,1,0,1,0) = 1832 |
| f(2,1,1,0,1,0) = 1702 | f(1,2,1,0,1,0) = 1256 | f(2,2,1,0,1,0) = 1701 |
| f(0,0,2,0,1,0) = 1003 | f(2,0,2,0,1,0) = 837 | f(1,1,2,0,1,0) = 1816 |
| f(2,1,2,0,1,0) = 1366 | f(1,2,2,0,1,0) = 1240 | f(2,2,2,0,1,0) = 1365 |
| f(0,0,0,1,1,0) = 1134 | f(2,0,0,1,1,0) = 388 | f(1,1,0,1,1,0) = 1662 |
| f(2,1,0,1,1,0) = 919 | f(1,2,0,1,1,0) = 1469 | f(2,2,0,1,1,0) = 727 |
| f(0,0,1,1,1,0) = 2020 | f(1,0,1,1,1,0) = 648 | f(2,0,1,1,1,0) = 1699 |
| f(0,1,1,1,1,0) = 1930 | f(1,1,1,1,1,0) = 1714 | f(2,1,1,1,1,0) = 1538 |
| f(1,2,1,1,1,0) = 1713 | f(2,2,1,1,1,0) = 1153 | f(0,0,2,1,1,0) = 2027 |
| f(1,0,2,1,1,0) = 328 | f(2,0,2,1,1,0) = 1363 | f(0,1,2,1,1,0) = 1866 |
| f(1,1,2,1,1,0) = 1394 | f(2,1,2,1,1,0) = 1282 | f(1,2,2,1,1,0) = 1393 |
| f(2,2,2,1,1,0) = 1089 | f(1,0,0,2,1,0) = 392 | f(1,1,0,2,1,0) = 2034 |
| f(2,1,0,2,1,0) = 1799 | f(1,2,0,2,1,0) = 2033 | f(2,2,0,2,1,0) = 1223 |

f(0,0,1,2,1,0,0) = 36
f(1,1,1,2,1,0,0) = 520
f(1,2,1,2,1,0,0) = 136
f(1,0,2,2,1,0,0) = 1048
f(2,1,2,2,1,0,0) = 790
f(2,2,2,2,1,0,0) = 213
f(2,0,0,0,2,0,0) = 513
f(2,1,0,0,2,0,0) = 770
f(2,2,0,0,2,0,0) = 193
f(1,1,1,0,2,0,0) = 682
f(2,2,1,0,2,0,0) = 231
f(1,1,2,0,2,0,0) = 346
f(2,2,2,0,2,0,0) = 215
f(2,1,0,1,2,0,0) = 14
f(0,0,1,1,2,0,0) = 1428
f(0,1,1,1,2,0,0) = 1546
f(1,2,1,1,2,0,0) = 1257
f(1,0,2,1,2,0,0) = 408
f(1,1,2,1,2,0,0) = 1818
f(2,2,2,1,2,0,0) = 1527
f(1,1,0,2,2,0,0) = 1320
f(2,2,0,2,2,0,0) = 385
f(2,0,1,2,2,0,0) = 1719
f(0,2,1,2,2,0,0) = 181
f(0,0,2,2,2,0,0) = 27
f(1,1,2,2,2,0,0) = 894
f(1,2,2,2,2,0,0) = 509
f(2,0,0,0,0,1,0) = 1987
f(2,1,0,0,0,1,0) = 582
f(2,2,0,0,0,1,0) = 389
f(2,0,1,0,0,1,0) = 1126
f(2,1,1,0,0,1,0) = 514
f(2,2,1,0,0,1,0) = 129
f(2,0,2,0,0,1,0) = 1557
f(2,1,2,0,0,1,0) = 258
f(2,2,2,0,0,1,0) = 65
f(2,0,0,1,0,1,0) = 1604
f(2,1,0,1,0,1,0) = 1951
f(2,2,0,1,0,1,0) = 1759
f(2,0,1,1,0,1,0) = 1766
f(2,1,1,1,0,1,0) = 1575
f(2,2,1,1,0,1,0) = 1191
f(2,0,2,1,0,1,0) = 1877
f(2,1,2,1,0,1,0) = 1303

f(1,0,1,2,1,0,0) = 1064
f(2,1,1,2,1,0,0) = 806
f(2,2,1,2,1,0,0) = 229
f(2,0,2,2,1,0,0) = 1623
f(0,2,2,2,1,0,0) = 501
f(0,0,0,0,2,0,0) = 63
f(0,1,0,0,2,0,0) = 1802
f(0,2,0,0,2,0,0) = 243
f(0,0,1,0,2,0,0) = 1044
f(2,1,1,0,2,0,0) = 807
f(0,0,2,0,2,0,0) = 987
f(2,1,2,0,2,0,0) = 791
f(2,0,0,1,2,0,0) = 1655
f(1,2,0,1,2,0,0) = 248
f(1,0,1,1,2,0,0) = 424
f(1,1,1,1,2,0,0) = 1834
f(2,2,1,1,2,0,0) = 1783
f(2,0,2,1,2,0,0) = 983
f(2,1,2,1,2,0,0) = 1911
f(0,0,0,2,2,0,0) = 913
f(2,1,0,2,2,0,0) = 578
f(0,0,1,2,2,0,0) = 20
f(1,1,1,2,2,0,0) = 958
f(1,2,1,2,2,0,0) = 765
f(1,0,2,2,2,0,0) = 348
f(2,1,2,2,2,0,0) = 334
f(2,2,2,2,2,0,0) = 333
f(0,1,0,0,0,1,0) = 822
f(0,2,0,0,0,1,0) = 246
f(0,0,1,0,0,1,0) = 1740
f(0,1,1,0,0,1,0) = 1920
f(0,2,1,0,0,1,0) = 1737
f(0,0,2,0,0,1,0) = 1587
f(0,1,2,0,0,1,0) = 1856
f(0,2,2,0,0,1,0) = 1481
f(0,0,0,1,0,1,0) = 418
f(0,1,0,1,0,1,0) = 1454
f(0,2,0,1,0,1,0) = 1262
f(0,0,1,1,0,1,0) = 1664
f(0,1,1,1,0,1,0) = 650
f(0,2,1,1,0,1,0) = 1728
f(0,0,2,1,0,1,0) = 1344
f(0,1,2,1,0,1,0) = 330
f(0,2,2,1,0,1,0) = 1472

f(2,0,1,2,1,0,0) = 1639
f(0,2,1,2,1,0,0) = 757
f(0,0,2,2,1,0,0) = 619
f(1,1,2,2,1,0,0) = 264
f(1,2,2,2,1,0,0) = 72
f(1,0,0,0,2,0,0) = 1981
f(1,1,0,0,2,0,0) = 296
f(1,2,0,0,2,0,0) = 104
f(1,0,1,0,2,0,0) = 1210
f(1,2,1,0,2,0,0) = 681
f(1,0,2,0,2,0,0) = 1337
f(1,2,2,0,2,0,0) = 345
f(1,1,0,1,2,0,0) = 824
f(2,2,0,1,2,0,0) = 13
f(2,0,1,1,2,0,0) = 999
f(2,1,1,1,2,0,0) = 1975
f(0,0,2,1,2,0,0) = 2011
f(0,1,2,1,2,0,0) = 1290
f(1,2,2,1,2,0,0) = 1241
f(1,0,0,2,2,0,0) = 1659
f(1,2,0,2,2,0,0) = 1128
f(1,0,1,2,2,0,0) = 684
f(2,1,1,2,2,0,0) = 654
f(2,2,1,2,2,0,0) = 653
f(2,0,2,2,2,0,0) = 1399
f(0,2,2,2,2,0,0) = 117
f(1,0,0,0,1,0) = 1456
f(1,1,0,0,0,1,0) = 816
f(1,2,0,0,0,1,0) = 240
f(1,0,1,0,0,1,0) = 1716
f(1,1,1,0,0,1,0) = 544
f(1,2,1,0,0,1,0) = 160
f(1,0,2,0,0,1,0) = 1396
f(1,1,2,0,0,1,0) = 272
f(1,2,2,0,0,1,0) = 80
f(1,0,0,1,0,1,0) = 2032
f(1,1,0,1,0,1,0) = 1022
f(1,2,0,1,0,1,0) = 1021
f(1,0,1,1,0,1,0) = 2024
f(1,1,1,1,0,1,0) = 1978
f(1,2,1,1,0,1,0) = 1785
f(1,0,2,1,0,1,0) = 2008
f(1,1,2,1,0,1,0) = 1914
f(1,2,2,1,0,1,0) = 1529

| | | |
|---|---|---|
| f(2,2,2,1,0,1,0) = 1111 | f(0,0,0,2,0,1,0) = 365 | f(1,0,0,2,0,1,0) = 1080 |
| f(2,0,0,2,0,1,0) = 1411 | f(0,1,0,2,0,1,0) = 1630 | f(1,1,0,2,0,1,0) = 1458 |
| f(2,1,0,2,0,1,0) = 1891 | f(0,2,0,2,0,1,0) = 1246 | f(1,2,0,2,0,1,0) = 1649 |
| f(2,2,0,2,0,1,0) = 1507 | f(0,0,1,2,0,1,0) = 76 | f(1,0,1,2,0,1,0) = 608 |
| f(2,0,1,2,0,1,0) = 663 | f(0,1,1,2,0,1,0) = 1929 | f(1,1,1,2,0,1,0) = 556 |
| f(2,1,1,2,0,1,0) = 1924 | f(0,2,1,2,0,1,0) = 691 | f(1,2,1,2,0,1,0) = 172 |
| f(2,2,1,2,0,1,0) = 1732 | f(0,0,2,2,0,1,0) = 585 | f(1,0,2,2,0,1,0) = 592 |
| f(2,0,2,2,0,1,0) = 359 | f(0,1,2,2,0,1,0) = 1865 | f(1,1,2,2,0,1,0) = 284 |
| f(2,1,2,2,0,1,0) = 1860 | f(0,2,2,2,0,1,0) = 371 | f(1,2,2,2,0,1,0) = 92 |
| f(2,2,2,2,0,1,0) = 1476 | f(0,0,0,0,1,1) = 1024 | f(1,0,0,0,1,1,0) = 1522 |
| f(2,0,0,0,1,1,0) = 1474 | f(0,1,0,0,1,1,0) = 780 | f(1,1,0,0,1,1,0) = 1552 |
| f(2,1,0,0,1,1,0) = 1998 | f(0,2,0,0,1,1,0) = 1013 | f(1,2,0,0,1,1,0) = 1168 |
| f(2,2,0,0,1,1,0) = 1997 | f(0,0,1,0,1,1,0) = 1774 | f(1,0,1,0,1,1,0) = 1004 |
| f(2,0,1,0,1,1,0) = 1668 | f(0,1,1,0,1,1,0) = 1572 | f(1,1,1,0,1,1,0) = 1972 |
| f(2,1,1,0,1,1,0) = 1670 | f(0,2,1,0,1,1,0) = 1197 | f(1,2,1,0,1,1,0) = 1780 |
| f(2,2,1,0,1,1,0) = 1669 | f(0,0,2,0,1,1,0) = 1569 | f(1,0,2,0,1,1,0) = 988 |
| f(2,0,2,0,1,1,0) = 1348 | f(0,1,2,0,1,1,0) = 1316 | f(1,1,2,0,1,1,0) = 1908 |
| f(2,1,2,0,1,1,0) = 1350 | f(0,2,2,0,1,1,0) = 1133 | f(1,2,2,0,1,1,0) = 1524 |
| f(2,2,2,0,1,1,0) = 1349 | f(1,0,0,1,1,1,0) = 1656 | f(2,0,0,1,1,1,0) = 1478 |
| f(0,1,0,1,1,1,0) = 1006 | f(1,1,0,1,1,1,0) = 538 | f(2,1,0,1,1,1,0) = 1942 |
| f(0,2,0,1,1,1,0) = 228 | f(1,2,0,1,1,1,0) = 153 | f(2,2,0,1,1,1,0) = 1749 |
| f(0,0,1,1,1,1,0) = 384 | f(1,0,1,1,1,1,0) = 1000 | f(2,0,1,1,1,1,0) = 1730 |
| f(0,1,1,1,1,1,0) = 906 | f(1,1,1,1,1,1,0) = 1584 | f(2,1,1,1,1,1,0) = 1666 |
| f(0,2,1,1,1,1,0) = 1188 | f(1,2,1,1,1,1,0) = 1200 | f(2,2,1,1,1,1,0) = 1665 |
| f(0,0,2,1,1,1,0) = 2018 | f(1,0,2,1,1,1,0) = 984 | f(2,0,2,1,1,1,0) = 1857 |
| f(0,1,2,1,1,1,0) = 842 | f(1,1,2,1,1,1,0) = 1328 | f(2,1,2,1,1,1,0) = 1346 |
| f(0,2,2,1,1,1,0) = 1124 | f(1,2,2,1,1,1,0) = 1136 | f(2,2,2,1,1,1,0) = 1345 |
| f(0,0,0,2,1,1,0) = 356 | f(1,0,0,2,1,1,0) = 1530 | f(2,0,0,2,1,1,0) = 1419 |
| f(0,1,0,2,1,1,0) = 990 | f(1,1,0,2,1,1,0) = 1840 | f(2,1,0,2,1,1,0) = 1793 |
| f(0,2,0,2,1,1,0) = 212 | f(1,2,0,2,1,1,0) = 1264 | f(2,2,0,2,1,1,0) = 1218 |
| f(0,0,1,2,1,1,0) = 174 | f(1,0,1,2,1,1,0) = 1640 | f(2,0,1,2,1,1,0) = 647 |
| f(0,1,1,2,1,1,0) = 1581 | f(1,1,1,2,1,1,0) = 1836 | f(2,1,1,2,1,1,0) = 1638 |
| f(0,2,1,2,1,1,0) = 755 | f(1,2,1,2,1,1,0) = 1260 | f(2,2,1,2,1,1,0) = 1445 |
| f(0,0,2,2,1,1,0) = 1067 | f(1,0,2,2,1,1,0) = 1624 | f(2,0,2,2,1,1,0) = 327 |
| f(0,1,2,2,1,1,0) = 1325 | f(1,1,2,2,1,1,0) = 1820 | f(2,1,2,2,1,1,0) = 1430 |
| f(0,2,2,2,1,1,0) = 499 | f(1,2,2,2,1,1,0) = 1244 | f(2,2,2,2,1,1,0) = 1621 |
| f(0,0,0,0,2,1) = 1087 | f(1,0,0,0,2,1,0) = 1969 | f(2,0,0,0,2,1,0) = 1921 |
| f(0,1,0,0,2,1,0) = 778 | f(1,1,0,0,2,1,0) = 1018 | f(2,1,0,0,2,1,0) = 915 |
| f(0,2,0,0,2,1,0) = 1011 | f(1,2,0,0,2,1,0) = 1017 | f(2,2,0,0,2,1,0) = 723 |
| f(0,0,1,0,2,1,0) = 1758 | f(1,0,1,0,2,1,0) = 1778 | f(2,0,1,0,2,1,0) = 643 |
| f(0,1,1,0,2,1,0) = 1570 | f(1,1,1,0,2,1,0) = 810 | f(2,1,1,0,2,1,0) = 422 |
| f(0,2,1,0,2,1,0) = 1195 | f(1,2,1,0,2,1,0) = 233 | f(2,2,1,0,2,1,0) = 613 |
| f(0,0,2,0,2,1,0) = 1553 | f(1,0,2,0,2,1,0) = 1905 | f(2,0,2,0,2,1,0) = 323 |
| f(0,1,2,0,2,1,0) = 1314 | f(1,1,2,0,2,1,0) = 794 | f(2,1,2,0,2,1,0) = 598 |

f(0,2,2,0,2,1,0) = 1131
f(0,0,0,1,2,1,0) = 427
f(0,1,0,1,2,1,0) = 1438
f(0,2,0,1,2,1,0) = 226
f(0,0,1,1,2,1,0) = 1718
f(0,1,1,1,2,1,0) = 522
f(0,2,1,1,2,1,0) = 1186
f(0,0,2,1,2,1,0) = 2002
f(0,1,2,1,2,1,0) = 266
f(0,2,2,1,2,1,0) = 1122
f(0,0,0,2,2,1,0) = 1953
f(0,1,0,2,2,1,0) = 1646
f(0,2,0,2,2,1,0) = 210
f(0,0,1,2,2,1,0) = 158
f(0,1,1,2,2,1,0) = 1579
f(0,2,1,2,2,1,0) = 179
f(0,0,2,2,2,1,0) = 603
f(0,1,2,2,2,1,0) = 1323
f(0,2,2,2,2,1,0) = 115
f(1,0,0,0,0,2,0) = 60
f(1,1,0,0,0,2,0) = 1658
f(1,2,0,0,0,2,0) = 1465
f(1,0,1,0,0,2,0) = 490
f(1,1,1,0,0,2,0) = 1982
f(1,2,1,0,0,2,0) = 1789
f(1,0,2,0,0,2,0) = 921
f(1,1,2,0,0,2,0) = 1918
f(1,2,2,0,0,2,0) = 1533
f(1,0,0,1,0,2,0) = 636
f(1,1,0,1,0,2,0) = 540
f(1,2,0,1,0,2,0) = 156
f(1,0,1,1,0,2,0) = 1688
f(1,1,1,1,0,2,0) = 571
f(1,2,1,1,0,2,0) = 187
f(1,0,2,1,0,2,0) = 1384
f(1,1,2,1,0,2,0) = 315
f(1,2,2,1,0,2,0) = 123
f(1,0,0,2,0,2,0) = 443
f(1,1,0,2,0,2,0) = 288
f(1,2,0,2,0,2,0) = 96
f(1,0,1,2,0,2,0) = 170
f(1,1,1,2,0,2,0) = 936
f(1,2,1,2,0,2,0) = 744
f(1,0,2,2,0,2,0) = 281

f(1,2,2,0,2,1,0) = 217
f(1,0,0,1,2,1,0) = 1992
f(1,1,0,1,2,1,0) = 924
f(1,2,0,1,2,1,0) = 732
f(1,0,1,1,2,1,0) = 1448
f(1,1,1,1,2,1,0) = 1824
f(1,2,1,1,2,1,0) = 1248
f(1,0,2,1,2,1,0) = 1432
f(1,1,2,1,2,1,0) = 1808
f(1,2,2,1,2,1,0) = 1232
f(1,0,0,2,2,1,0) = 1146
f(1,1,0,2,2,1,0) = 1896
f(1,2,0,2,2,1,0) = 1512
f(1,0,1,2,2,1,0) = 553
f(1,1,1,2,2,1,0) = 956
f(1,2,1,2,2,1,0) = 764
f(1,0,2,2,2,1,0) = 90
f(1,1,2,2,2,1,0) = 892
f(1,2,2,2,2,1,0) = 508
f(2,0,0,0,0,2,0) = 591
f(2,1,0,0,0,2,0) = 1807
f(2,2,0,0,0,2,0) = 1231
f(2,0,1,0,0,2,0) = 651
f(2,1,1,0,0,2,0) = 1967
f(2,2,1,0,0,2,0) = 1775
f(2,0,2,0,0,2,0) = 331
f(2,1,2,0,0,2,0) = 1887
f(2,2,2,0,0,2,0) = 1503
f(2,0,0,1,0,2,0) = 967
f(2,1,0,1,0,2,0) = 398
f(2,2,0,1,0,2,0) = 589
f(2,0,1,1,0,2,0) = 1455
f(2,1,1,1,0,2,0) = 1955
f(2,2,1,1,0,2,0) = 1763
f(2,0,2,1,0,2,0) = 1439
f(2,1,2,1,0,2,0) = 1875
f(2,2,2,1,0,2,0) = 1491
f(2,0,0,2,0,2,0) = 15
f(2,1,0,2,0,2,0) = 1026
f(2,2,0,2,0,2,0) = 1025
f(2,0,1,2,0,2,0) = 39
f(2,1,1,2,0,2,0) = 518
f(2,2,1,2,0,2,0) = 133
f(2,0,2,2,0,2,0) = 23

f(2,2,2,0,2,1,0) = 405
f(2,0,0,1,2,1,0) = 1094
f(2,1,0,1,2,1,0) = 1038
f(2,2,0,1,2,1,0) = 1037
f(2,0,1,1,2,1,0) = 486
f(2,1,1,1,2,1,0) = 1591
f(2,2,1,1,2,1,0) = 1207
f(2,0,2,1,2,1,0) = 917
f(2,1,2,1,2,1,0) = 1335
f(2,2,2,1,2,1,0) = 1143
f(2,0,0,2,2,1,0) = 1611
f(2,1,0,2,2,1,0) = 1986
f(2,2,0,2,2,1,0) = 1985
f(2,0,1,2,2,1,0) = 1687
f(2,1,1,2,2,1,0) = 1678
f(2,2,1,2,2,1,0) = 1677
f(2,0,2,2,2,1,0) = 1383
f(2,1,2,2,2,1,0) = 1358
f(2,2,2,2,2,1,0) = 1357
f(0,1,0,0,0,2,0) = 1801
f(0,2,0,0,0,2,0) = 1225
f(0,0,1,0,0,2,0) = 460
f(0,1,1,0,0,2,0) = 566
f(0,2,1,0,0,2,0) = 191
f(0,0,2,0,0,2,0) = 307
f(0,1,2,0,0,2,0) = 310
f(0,2,2,0,0,2,0) = 127
f(0,0,0,1,0,2,0) = 1682
f(0,1,0,1,0,2,0) = 801
f(0,2,0,1,0,2,0) = 417
f(0,0,1,1,0,2,0) = 1462
f(0,1,1,1,0,2,0) = 1676
f(0,2,1,1,0,2,0) = 182
f(0,0,2,1,0,2,0) = 1971
f(0,1,2,1,0,2,0) = 1356
f(0,2,2,1,0,2,0) = 118
f(0,0,0,2,0,2,0) = 1629
f(0,1,0,2,0,2,0) = 785
f(0,2,0,2,0,2,0) = 593
f(0,0,1,2,0,2,0) = 703
f(0,1,1,2,0,2,0) = 575
f(0,2,1,2,0,2,0) = 1717
f(0,0,2,2,0,2,0) = 383
f(0,1,2,2,0,2,0) = 319

| | | |
|---|---|---|
| f(1,1,2,2,0,2,0) = 856 | f(2,1,2,2,0,2,0) = 262 | f(0,2,2,2,0,2,0) = 1397 |
| f(1,2,2,2,0,2,0) = 472 | f(2,2,2,2,0,2,0) = 69 | f(0,0,0,0,1,2) = 960 |
| f(1,0,0,0,1,2,0) = 126 | f(2,0,0,0,1,2,0) = 78 | f(0,1,0,0,1,2,0) = 1036 |
| f(1,1,0,0,1,2,0) = 1324 | f(2,1,0,0,1,2,0) = 1030 | f(0,2,0,0,1,2,0) = 1269 |
| f(1,2,0,0,1,2,0) = 1132 | f(2,2,0,0,1,2,0) = 1029 | f(0,0,1,0,1,2,0) = 494 |
| f(1,0,1,0,1,2,0) = 1724 | f(2,0,1,0,1,2,0) = 142 | f(0,1,1,0,1,2,0) = 916 |
| f(1,1,1,0,1,2,0) = 1642 | f(2,1,1,0,1,2,0) = 1830 | f(0,2,1,0,1,2,0) = 733 |
| f(1,2,1,0,1,2,0) = 1449 | f(2,2,1,0,1,2,0) = 1253 | f(0,0,2,0,1,2,0) = 289 |
| f(1,0,2,0,1,2,0) = 1404 | f(2,0,2,0,1,2,0) = 269 | f(0,1,2,0,1,2,0) = 852 |
| f(1,1,2,0,1,2,0) = 1434 | f(2,1,2,0,1,2,0) = 1814 | f(0,2,2,0,1,2,0) = 477 |
| f(1,2,2,0,1,2,0) = 1625 | f(2,2,2,0,1,2,0) = 1237 | f(0,0,0,1,1,2,0) = 94 |
| f(1,0,0,1,1,2,0) = 436 | f(2,0,0,1,1,2,0) = 901 | f(0,1,0,1,1,2,0) = 1837 |
| f(1,1,0,1,1,2,0) = 62 | f(2,1,0,1,1,2,0) = 535 | f(0,2,0,1,1,2,0) = 401 |
| f(1,2,0,1,1,2,0) = 61 | f(2,2,0,1,1,2,0) = 151 | f(0,0,1,1,1,2,0) = 1444 |
| f(1,0,1,1,1,2,0) = 664 | f(2,0,1,1,1,2,0) = 1957 | f(0,1,1,1,1,2,0) = 1932 |
| f(1,1,1,1,1,2,0) = 690 | f(2,1,1,1,1,2,0) = 1539 | f(0,2,1,1,1,2,0) = 724 |
| f(1,2,1,1,1,2,0) = 689 | f(2,2,1,1,1,2,0) = 1155 | f(0,0,2,1,1,2,0) = 1889 |
| f(1,0,2,1,1,2,0) = 360 | f(2,0,2,1,1,2,0) = 1494 | f(0,1,2,1,1,2,0) = 1868 |
| f(1,1,2,1,1,2,0) = 370 | f(2,1,2,1,1,2,0) = 1283 | f(0,2,2,1,1,2,0) = 468 |
| f(1,2,2,1,1,2,0) = 369 | f(2,2,2,1,1,2,0) = 1091 | f(0,0,0,2,1,2,0) = 1620 |
| f(1,0,0,2,1,2,0) = 953 | f(2,0,0,2,1,2,0) = 55 | f(0,1,0,2,1,2,0) = 1821 |
| f(1,1,0,2,1,2,0) = 1010 | f(2,1,0,2,1,2,0) = 1315 | f(0,2,0,2,1,2,0) = 609 |
| f(1,2,0,2,1,2,0) = 1009 | f(2,2,0,2,1,2,0) = 1123 | f(0,0,1,2,1,2,0) = 45 |
| f(1,0,1,2,1,2,0) = 1130 | f(2,0,1,2,1,2,0) = 615 | f(0,1,1,2,1,2,0) = 925 |
| f(1,1,1,2,1,2,0) = 904 | f(2,1,1,2,1,2,0) = 815 | f(0,2,1,2,1,2,0) = 1781 |
| f(1,2,1,2,1,2,0) = 712 | f(2,2,1,2,1,2,0) = 239 | f(0,0,2,2,1,2,0) = 329 |
| f(1,0,2,2,1,2,0) = 1561 | f(2,0,2,2,1,2,0) = 599 | f(0,1,2,2,1,2,0) = 861 |
| f(1,1,2,2,1,2,0) = 840 | f(2,1,2,2,1,2,0) = 799 | f(0,2,2,2,1,2,0) = 1525 |
| f(1,2,2,2,1,2,0) = 456 | f(2,2,2,2,1,2,0) = 223 | f(0,0,0,0,2,2) = 1023 |
| f(1,0,0,0,2,2,0) = 573 | f(2,0,0,0,2,2,0) = 525 | f(0,1,0,0,2,2,0) = 1034 |
| f(1,1,0,0,2,2,0) = 50 | f(2,1,0,0,2,2,0) = 879 | f(0,2,0,0,2,2,0) = 1267 |
| f(1,2,0,0,2,2,0) = 49 | f(2,2,0,0,2,2,0) = 495 | f(0,0,1,0,2,2,0) = 478 |
| f(1,0,1,0,2,2,0) = 699 | f(2,0,1,0,2,2,0) = 1059 | f(0,1,1,0,2,2,0) = 914 |
| f(1,1,1,0,2,2,0) = 698 | f(2,1,1,0,2,2,0) = 523 | f(0,2,1,0,2,2,0) = 731 |
| f(1,2,1,0,2,2,0) = 697 | f(2,2,1,0,2,2,0) = 139 | f(0,0,2,0,2,2,0) = 273 |
| f(1,0,2,0,2,2,0) = 379 | f(2,0,2,0,2,2,0) = 1043 | f(0,1,2,0,2,2,0) = 850 |
| f(1,1,2,0,2,2,0) = 378 | f(2,1,2,0,2,2,0) = 267 | f(0,2,2,0,2,2,0) = 475 |
| f(1,2,2,0,2,2,0) = 377 | f(2,2,2,0,2,2,0) = 75 | f(0,0,0,1,2,2,0) = 1691 |
| f(1,0,0,1,2,2,0) = 628 | f(2,0,0,1,2,2,0) = 517 | f(0,1,0,1,2,2,0) = 1835 |
| f(1,1,0,1,2,2,0) = 829 | f(2,1,0,1,2,2,0) = 783 | f(0,2,0,1,2,2,0) = 1057 |
| f(1,2,0,1,2,2,0) = 254 | f(2,2,0,1,2,2,0) = 207 | f(0,0,1,1,2,2,0) = 980 |
| f(1,0,1,1,2,2,0) = 1720 | f(2,0,1,1,2,2,0) = 423 | f(0,1,1,1,2,2,0) = 1548 |
| f(1,1,1,1,2,2,0) = 426 | f(2,1,1,1,2,2,0) = 803 | f(0,2,1,1,2,2,0) = 722 |
| f(1,2,1,1,2,2,0) = 617 | f(2,2,1,1,2,2,0) = 227 | f(0,0,2,1,2,2,0) = 1873 |
| f(1,0,2,1,2,2,0) = 1400 | f(2,0,2,1,2,2,0) = 407 | f(0,1,2,1,2,2,0) = 1292 |
| f(1,1,2,1,2,2,0) = 602 | f(2,1,2,1,2,2,0) = 787 | f(0,2,2,1,2,2,0) = 466 |
| f(1,2,2,1,2,2,0) = 409 | f(2,2,2,1,2,2,0) = 211 | f(1,0,0,2,2,2,0) = 569 |
| f(2,0,0,2,2,2,0) = 391 | f(0,1,0,2,2,2,0) = 1819 | f(1,1,0,2,2,2,0) = 298 |
| f(2,1,0,2,2,2,0) = 1894 | f(0,2,0,2,2,2,0) = 1041 | f(1,2,0,2,2,2,0) = 105 |
| f(2,2,0,2,2,2,0) = 1509 | f(0,0,1,2,2,2,0) = 29 | f(1,0,1,2,2,2,0) = 190 |
| f(2,0,1,2,2,2,0) = 1063 | f(0,1,1,2,2,2,0) = 923 | f(1,1,1,2,2,2,0) = 702 |
| f(2,1,1,2,2,2,0) = 911 | f(0,2,1,2,2,2,0) = 1205 | f(1,2,1,2,2,2,0) = 701 |
| f(2,2,1,2,2,2,0) = 719 | f(0,0,2,2,2,2,0) = 1663 | f(1,0,2,2,2,2,0) = 317 |
| f(2,0,2,2,2,2,0) = 1047 | f(0,1,2,2,2,2,0) = 859 | f(1,1,2,2,2,2,0) = 382 |
| f(2,1,2,2,2,2,0) = 847 | f(0,2,2,2,2,2,0) = 1141 | f(1,2,2,2,2,2,0) = 381 |
| f(2,2,2,2,2,2,0) = 463 | f(2,0,0,0,0,0,1) = 1999 | f(0,1,0,0,0,0,1) = 768 |
| f(1,1,0,0,0,0,1) = 1470 | f(2,1,0,0,0,0,1) = 1847 | f(1,2,0,0,0,0,1) = 1661 |
| f(2,2,0,0,0,0,1) = 1271 | f(0,0,1,0,0,0,1) = 694 | f(2,0,1,0,0,0,1) = 1711 |
| f(0,1,1,0,0,0,1) = 512 | f(1,1,1,0,0,0,1) = 928 | f(2,1,1,0,0,0,1) = 898 |
| f(1,2,1,0,0,0,1) = 736 | f(2,2,1,0,0,0,1) = 705 | f(0,0,2,0,0,0,1) = 576 |
| f(2,0,2,0,0,0,1) = 1375 | f(0,1,2,0,0,0,1) = 256 | f(1,1,2,0,0,0,1) = 848 |
| f(2,1,2,0,0,0,1) = 834 | f(1,2,2,0,0,0,1) = 464 | f(2,2,2,0,0,0,1) = 449 |
| f(0,0,0,1,0,0,1) = 1698 | f(1,0,0,1,0,0,1) = 624 | f(2,0,0,1,0,0,1) = 2039 |
| f(1,1,0,1,0,0,1) = 876 | f(2,1,0,1,0,0,1) = 781 | f(0,2,0,1,0,0,1) = 2017 |
| f(1,2,0,1,0,0,1) = 492 | f(2,2,0,1,0,0,1) = 206 | f(0,0,1,1,0,0,1) = 1162 |
| f(1,0,1,1,0,0,1) = 688 | f(2,0,1,1,0,0,1) = 1447 | f(0,1,1,1,0,0,1) = 1710 |
| f(1,1,1,1,0,0,1) = 1568 | f(2,1,1,1,0,0,1) = 1831 | f(1,2,1,1,0,0,1) = 1184 |
| f(2,2,1,1,0,0,1) = 1255 | f(0,0,2,1,0,0,1) = 1398 | f(1,0,2,1,0,0,1) = 368 |
| f(2,0,2,1,0,0,1) = 1431 | f(0,1,2,1,0,0,1) = 1390 | f(1,1,2,1,0,0,1) = 1296 |
| f(2,1,2,1,0,0,1) = 1815 | f(1,2,2,1,0,0,1) = 1104 | f(2,2,2,1,0,0,1) = 1239 |
| f(0,0,0,2,0,0,1) = 1645 | f(1,0,0,2,0,0,1) = 2043 | f(1,1,0,2,0,0,1) = 2042 |
| f(2,1,0,2,0,0,1) = 771 | f(0,2,0,2,0,0,1) = 2001 | f(1,2,0,2,0,0,1) = 2041 |
| f(2,2,0,2,0,0,1) = 195 | f(0,0,1,2,0,0,1) = 140 | f(1,0,1,2,0,0,1) = 1068 |
| f(2,0,1,2,0,0,1) = 695 | f(0,1,1,2,0,0,1) = 521 | f(1,1,1,2,0,0,1) = 1980 |
| f(2,1,1,2,0,0,1) = 516 | f(0,2,1,2,0,0,1) = 673 | f(1,2,1,2,0,0,1) = 1788 |
| f(2,2,1,2,0,0,1) = 132 | f(0,0,2,2,0,0,1) = 74 | f(1,0,2,2,0,0,1) = 1052 |
| f(2,0,2,2,0,0,1) = 375 | f(0,1,2,2,0,0,1) = 265 | f(1,1,2,2,0,0,1) = 1916 |
| f(2,1,2,2,0,0,1) = 260 | f(0,2,2,2,0,0,1) = 353 | f(1,2,2,2,0,0,1) = 1532 |
| f(2,2,2,2,0,0,1) = 68 | f(0,0,0,0,1,0,1) = 2038 | f(2,0,0,0,1,0,1) = 1486 |
| f(0,1,0,0,1,0,1) = 1996 | f(1,1,0,0,1,0,1) = 1848 | f(2,1,0,0,1,0,1) = 1798 |
| f(0,2,0,0,1,0,1) = 204 | f(1,2,0,0,1,0,1) = 1272 | f(2,2,0,0,1,0,1) = 1221 |
| f(0,0,1,0,1,0,1) = 1198 | f(1,0,1,0,1,0,1) = 1452 | f(2,0,1,0,1,0,1) = 1742 |
| f(0,1,1,0,1,0,1) = 1556 | f(1,1,1,0,1,0,1) = 1960 | f(2,1,1,0,1,0,1) = 1958 |
| f(0,2,1,0,1,0,1) = 1181 | f(1,2,1,0,1,0,1) = 1768 | f(2,2,1,0,1,0,1) = 1765 |
| f(0,0,2,0,1,0,1) = 34 | f(1,0,2,0,1,0,1) = 1436 | f(2,0,2,0,1,0,1) = 1869 |
| f(0,1,2,0,1,0,1) = 1300 | f(1,1,2,0,1,0,1) = 1880 | f(2,1,2,0,1,0,1) = 1878 |
| f(0,2,2,0,1,0,1) = 1117 | f(1,2,2,0,1,0,1) = 1496 | f(2,2,2,0,1,0,1) = 1493 |
| f(0,0,0,1,1,0,1) = 714 | f(1,0,0,1,1,0,1) = 1016 | f(2,0,0,1,1,0,1) = 1988 |
| f(0,1,0,1,1,0,1) = 1420 | f(1,1,0,1,1,0,1) = 1564 | f(2,1,0,1,1,0,1) = 918 |
| f(0,2,0,1,1,0,1) = 237 | f(1,2,0,1,1,0,1) = 1180 | f(2,2,0,1,1,0,1) = 725 |
| f(0,0,1,1,1,0,1) = 1518 | f(1,0,1,1,1,0,1) = 1692 | f(2,0,1,1,1,0,1) = 2019 |
| f(0,1,1,1,1,0,1) = 1966 | f(1,1,1,1,1,0,1) = 1970 | f(2,1,1,1,1,0,1) = 1923 |
| f(0,2,1,1,1,0,1) = 1172 | f(1,2,1,1,1,0,1) = 1777 | f(2,2,1,1,1,0,1) = 1731 |
| f(0,0,2,1,1,0,1) = 1378 | f(1,0,2,1,1,0,1) = 1388 | f(2,0,2,1,1,0,1) = 2003 |
| f(0,1,2,1,1,0,1) = 1902 | f(1,1,2,1,1,0,1) = 1906 | f(2,1,2,1,1,0,1) = 1859 |
| f(0,2,2,1,1,0,1) = 1108 | f(1,2,2,1,1,0,1) = 1521 | f(2,2,2,1,1,0,1) = 1475 |
| f(0,0,0,2,1,0,1) = 1636 | f(1,0,0,2,1,0,1) = 1464 | f(2,0,0,2,1,0,1) = 11 |
| f(0,1,0,2,1,0,1) = 1612 | f(1,1,0,2,1,0,1) = 864 | f(2,1,0,2,1,0,1) = 1903 |
| f(0,2,0,2,1,0,1) = 221 | f(1,2,0,2,1,0,1) = 480 | f(2,2,0,2,1,0,1) = 1519 |
| f(0,0,1,2,1,0,1) = 1069 | f(1,0,1,2,1,0,1) = 1056 | f(2,0,1,2,1,0,1) = 549 |
| f(0,1,1,2,1,0,1) = 1565 | f(1,1,1,2,1,0,1) = 1964 | f(2,1,1,2,1,0,1) = 998 |
| f(0,2,1,2,1,0,1) = 737 | f(1,2,1,2,1,0,1) = 1772 | f(2,2,1,2,1,0,1) = 997 |
| f(0,0,2,2,1,0,1) = 545 | f(1,0,2,2,1,0,1) = 1040 | f(2,0,2,2,1,0,1) = 86 |
| f(0,1,2,2,1,0,1) = 1309 | f(1,1,2,2,1,0,1) = 1884 | f(2,1,2,2,1,0,1) = 982 |
| f(0,2,2,2,1,0,1) = 481 | f(1,2,2,2,1,0,1) = 1500 | f(2,2,2,2,1,0,1) = 981 |
| f(0,0,0,0,2,0,1) = 54 | f(2,0,0,0,2,0,1) = 1933 | f(0,1,0,0,2,0,1) = 1994 |
| f(1,1,0,0,2,0,1) = 872 | f(2,1,0,0,2,0,1) = 779 | f(0,2,0,0,2,0,1) = 202 |
| f(1,2,0,0,2,0,1) = 488 | f(2,2,0,0,2,0,1) = 203 | f(0,0,1,0,2,0,1) = 1182 |
| f(2,0,1,0,2,0,1) = 1703 | f(0,1,1,0,2,0,1) = 1554 | f(1,1,1,0,2,0,1) = 954 |
| f(2,1,1,0,2,0,1) = 903 | f(0,2,1,0,2,0,1) = 1179 | f(1,2,1,0,2,0,1) = 761 |
| f(2,2,1,0,2,0,1) = 711 | f(0,0,2,0,2,0,1) = 18 | f(2,0,2,0,2,0,1) = 1367 |
| f(0,1,2,0,2,0,1) = 1298 | f(1,1,2,0,2,0,1) = 890 | f(2,1,2,0,2,0,1) = 839 |
| f(0,2,2,0,2,0,1) = 1115 | f(1,2,2,0,2,0,1) = 505 | f(2,2,2,0,2,0,1) = 455 |
| f(0,0,0,1,2,0,1) = 1707 | f(1,0,0,1,2,0,1) = 584 | f(2,0,0,1,2,0,1) = 1541 |
| f(0,1,0,1,2,0,1) = 394 | f(1,1,0,1,2,0,1) = 828 | f(2,1,0,1,2,0,1) = 782 |
| f(0,2,0,1,2,0,1) = 1259 | f(1,2,0,1,2,0,1) = 252 | f(2,2,0,1,2,0,1) = 205 |
| f(0,0,1,1,2,0,1) = 2004 | f(1,0,1,1,2,0,1) = 416 | f(2,0,1,1,2,0,1) = 166 |
| f(0,1,1,1,2,0,1) = 1582 | f(1,1,1,1,2,0,1) = 1962 | f(2,1,1,1,2,0,1) = 1446 |
| f(0,2,1,1,2,0,1) = 1170 | f(1,2,1,1,2,0,1) = 1769 | f(2,2,1,1,2,0,1) = 1637 |
| f(0,0,2,1,2,0,1) = 1362 | f(1,0,2,1,2,0,1) = 400 | f(2,0,2,1,2,0,1) = 277 |
| f(0,1,2,1,2,0,1) = 1326 | f(1,1,2,1,2,0,1) = 1882 | f(2,1,2,1,2,0,1) = 1622 |
| f(0,2,2,1,2,0,1) = 1106 | f(1,2,2,1,2,0,1) = 1497 | f(2,2,2,1,2,0,1) = 1429 |
| f(0,0,0,2,2,0,1) = 1589 | f(2,0,0,2,2,0,1) = 971 | f(0,1,0,2,2,0,1) = 586 |
| f(1,1,0,2,2,0,1) = 1312 | f(2,1,0,2,2,0,1) = 1602 | f(0,2,0,2,2,0,1) = 1243 |
| f(1,2,0,2,2,0,1) = 1120 | f(2,2,0,2,2,0,1) = 1409 | f(0,0,1,2,2,0,1) = 1053 |
| f(1,0,1,2,2,0,1) = 1790 | f(2,0,1,2,2,0,1) = 1174 | f(0,1,1,2,2,0,1) = 1563 |
| f(1,1,1,2,2,0,1) = 952 | f(2,1,1,2,2,0,1) = 910 | f(0,2,1,2,2,0,1) = 161 |
| f(1,2,1,2,2,0,1) = 760 | f(2,2,1,2,2,0,1) = 717 | f(0,0,2,2,2,0,1) = 1051 |
| f(1,0,2,2,2,0,1) = 1917 | f(2,0,2,2,2,0,1) = 1317 | f(0,1,2,2,2,0,1) = 1307 |
| f(1,1,2,2,2,0,1) = 888 | f(2,1,2,2,2,0,1) = 846 | f(0,2,2,2,2,0,1) = 97 |
| f(1,2,2,2,2,0,1) = 504 | f(2,2,2,2,2,0,1) = 461 | f(0,0,0,0,0,1,1) = 1078 |
| f(1,0,0,0,0,1,1) = 1072 | f(2,0,0,0,0,1,1) = 1603 | f(0,1,0,0,0,1,1) = 1846 |

f(1,1,0,0,0,1,1) = 912
f(1,2,0,0,0,1,1) = 720
f(1,0,1,0,0,1,1) = 1696
f(1,1,1,0,0,1,1) = 800
f(1,2,1,0,0,1,1) = 224
f(1,0,2,0,0,1,1) = 1360
f(1,1,2,0,0,1,1) = 784
f(1,2,2,0,0,1,1) = 208
f(1,0,0,1,0,1,1) = 1648
f(1,1,0,1,0,1,1) = 830
f(1,2,0,1,0,1,1) = 253
f(1,0,1,1,0,1,1) = 2016
f(1,1,1,1,0,1,1) = 1952
f(1,2,1,1,0,1,1) = 1760
f(1,0,2,1,0,1,1) = 2000
f(1,1,2,1,0,1,1) = 1872
f(1,2,2,1,0,1,1) = 1488
f(2,0,0,2,0,1,1) = 1027
f(2,1,0,2,0,1,1) = 1895
f(2,2,0,2,0,1,1) = 1511
f(2,0,1,2,0,1,1) = 1750
f(2,1,1,2,0,1,1) = 1540
f(2,2,1,2,0,1,1) = 1156
f(2,0,2,2,0,1,1) = 1893
f(2,1,2,2,0,1,1) = 1284
f(2,2,2,2,0,1,1) = 1092
f(0,1,0,0,1,1,1) = 972
f(0,2,0,0,1,1,1) = 2037
f(0,0,1,0,1,1,1) = 676
f(0,1,1,0,1,1,1) = 1956
f(1,2,1,0,1,1,1) = 1784
f(0,1,2,0,1,1,1) = 1892
f(1,2,2,0,1,1,1) = 1528
f(2,0,0,1,1,1,1) = 1412
f(2,1,0,1,1,1,1) = 772
f(2,2,0,1,1,1,1) = 196
f(2,0,1,1,1,1,1) = 1154
f(2,1,1,1,1,1,1) = 1922
f(2,2,1,1,1,1,1) = 1729
f(2,0,2,1,1,1,1) = 1281
f(2,1,2,1,1,1,1) = 1858
f(2,2,2,1,1,1,1) = 1473
f(2,0,0,2,1,1,1) = 1035
f(2,1,0,2,1,1,1) = 1797

f(2,1,0,0,0,1,1) = 966
f(2,2,0,0,0,1,1) = 965
f(2,0,1,0,0,1,1) = 1158
f(2,1,1,0,0,1,1) = 642
f(2,2,1,0,0,1,1) = 641
f(2,0,2,0,0,1,1) = 1285
f(2,1,2,0,0,1,1) = 322
f(2,2,2,0,0,1,1) = 321
f(2,0,0,1,0,1,1) = 1028
f(2,1,0,1,0,1,1) = 1558
f(2,2,0,1,0,1,1) = 1173
f(2,0,1,1,0,1,1) = 2023
f(2,1,1,1,0,1,1) = 1686
f(2,2,1,1,0,1,1) = 1685
f(2,0,2,1,0,1,1) = 2007
f(2,1,2,1,0,1,1) = 1382
f(2,2,2,1,0,1,1) = 1381
f(0,1,0,2,0,1,1) = 1054
f(0,2,0,2,0,1,1) = 1233
f(0,0,1,2,0,1,1) = 1100
f(0,1,1,2,0,1,1) = 1545
f(0,2,1,2,0,1,1) = 1715
f(0,0,2,2,0,1,1) = 1609
f(0,1,2,2,0,1,1) = 1289
f(0,2,2,2,0,1,1) = 1395
f(1,0,0,0,1,1,1) = 1138
f(1,1,0,0,1,1,1) = 1944
f(1,2,0,0,1,1,1) = 1752
f(1,0,1,0,1,1,1) = 2028
f(1,1,1,0,1,1,1) = 1976
f(1,0,2,0,1,1,1) = 2012
f(1,1,2,0,1,1,1) = 1912
f(0,0,0,1,1,1,1) = 1738
f(0,1,0,1,1,1,1) = 2030
f(0,2,0,1,1,1,1) = 1252
f(0,0,1,1,1,1,1) = 1408
f(0,1,1,1,1,1,1) = 942
f(0,2,1,1,1,1,1) = 1764
f(0,0,2,1,1,1,1) = 1634
f(0,1,2,1,1,1,1) = 878
f(0,2,2,1,1,1,1) = 1508
f(0,0,0,2,1,1,1) = 1380
f(0,1,0,2,1,1,1) = 2014
f(0,2,0,2,1,1,1) = 1236

f(0,2,0,0,0,1,1) = 1270
f(0,0,1,0,0,1,1) = 1673
f(0,1,1,0,0,1,1) = 1536
f(0,2,1,0,0,1,1) = 1161
f(0,0,2,0,0,1,1) = 1098
f(0,1,2,0,0,1,1) = 1280
f(0,2,2,0,0,1,1) = 1097
f(0,0,0,1,0,1,1) = 1442
f(0,1,0,1,0,1,1) = 1070
f(0,2,0,1,0,1,1) = 1249
f(0,0,1,1,0,1,1) = 1700
f(0,1,1,1,0,1,1) = 686
f(0,2,1,1,0,1,1) = 1152
f(0,0,2,1,0,1,1) = 1600
f(0,1,2,1,0,1,1) = 366
f(0,2,2,1,0,1,1) = 1088
f(0,0,0,2,0,1,1) = 1389
f(1,1,0,2,0,1,1) = 1946
f(1,2,0,2,0,1,1) = 1753
f(1,0,1,2,0,1,1) = 1644
f(1,1,1,2,0,1,1) = 1596
f(1,2,1,2,0,1,1) = 1212
f(1,0,2,2,0,1,1) = 1628
f(1,1,2,2,0,1,1) = 1340
f(1,2,2,2,0,1,1) = 1148
f(2,0,0,0,1,1,1) = 1090
f(2,1,0,0,1,1,1) = 1796
f(2,2,0,0,1,1,1) = 1220
f(2,0,1,0,1,1,1) = 644
f(0,2,1,0,1,1,1) = 1773
f(2,0,2,0,1,1,1) = 324
f(0,2,2,0,1,1,1) = 1517
f(1,0,0,1,1,1,1) = 2040
f(1,1,0,1,1,1,1) = 1560
f(1,2,0,1,1,1,1) = 1176
f(1,0,1,1,1,1,1) = 992
f(1,1,1,1,1,1,1) = 1968
f(1,2,1,1,1,1,1) = 1776
f(1,0,2,1,1,1,1) = 976
f(1,1,2,1,1,1,1) = 1904
f(1,2,2,1,1,1,1) = 1520
f(1,0,0,2,1,1,1) = 1416
f(1,1,0,2,1,1,1) = 1936
f(1,2,0,2,1,1,1) = 1744

| | | |
|---|---|---|
| f(2,2,0,2,1,1,1) = 1222 | f(0,0,1,2,1,1,1) = 750 | f(1,0,1,2,1,1,1) = 1632 |
| f(2,0,1,2,1,1,1) = 1573 | f(0,1,1,2,1,1,1) = 1965 | f(1,1,1,2,1,1,1) = 1544 |
| f(2,1,1,2,1,1,1) = 2022 | f(0,2,1,2,1,1,1) = 1779 | f(1,2,1,2,1,1,1) = 1160 |
| f(2,2,1,2,1,1,1) = 2021 | f(1,0,2,2,1,1,1) = 1616 | f(2,0,2,2,1,1,1) = 1110 |
| f(0,1,2,2,1,1,1) = 1901 | f(1,1,2,2,1,1,1) = 1288 | f(2,1,2,2,1,1,1) = 2006 |
| f(0,2,2,2,1,1,1) = 1523 | f(1,2,2,2,1,1,1) = 1096 | f(2,2,2,2,1,1,1) = 2005 |
| f(1,0,0,0,2,1,1) = 1585 | f(2,0,0,0,2,1,1) = 1537 | f(0,1,0,0,2,1,1) = 970 |
| f(1,1,0,0,2,1,1) = 922 | f(2,1,0,0,2,1,1) = 773 | f(0,2,0,0,2,1,1) = 2035 |
| f(1,2,0,0,2,1,1) = 729 | f(2,2,0,0,2,1,1) = 198 | f(0,0,1,0,2,1,1) = 660 |
| f(1,0,1,0,2,1,1) = 1202 | f(2,0,1,0,2,1,1) = 1667 | f(0,1,1,0,2,1,1) = 1954 |
| f(1,1,1,0,2,1,1) = 938 | f(2,1,1,0,2,1,1) = 935 | f(0,2,1,0,2,1,1) = 1771 |
| f(1,2,1,0,2,1,1) = 745 | f(2,2,1,0,2,1,1) = 743 | f(1,0,2,0,2,1,1) = 1329 |
| f(2,0,2,0,2,1,1) = 1347 | f(0,1,2,0,2,1,1) = 1890 | f(1,1,2,0,2,1,1) = 858 |
| f(2,1,2,0,2,1,1) = 855 | f(0,2,2,0,2,1,1) = 1515 | f(1,2,2,0,2,1,1) = 473 |
| f(2,2,2,0,2,1,1) = 471 | f(0,0,0,1,2,1,1) = 1451 | f(1,0,0,1,2,1,1) = 1608 |
| f(2,0,0,1,2,1,1) = 70 | f(0,1,0,1,2,1,1) = 1418 | f(1,1,0,1,2,1,1) = 920 |
| f(2,1,0,1,2,1,1) = 1422 | f(0,2,0,1,2,1,1) = 1250 | f(1,2,0,1,2,1,1) = 728 |
| f(2,2,0,1,2,1,1) = 1613 | f(0,0,1,1,2,1,1) = 1684 | f(1,0,1,1,2,1,1) = 1440 |
| f(2,0,1,1,2,1,1) = 1510 | f(0,1,1,1,2,1,1) = 558 | f(1,1,1,1,2,1,1) = 1578 |
| f(2,1,1,1,2,1,1) = 1062 | f(0,2,1,1,2,1,1) = 1762 | f(1,2,1,1,2,1,1) = 1193 |
| f(2,2,1,1,2,1,1) = 1061 | f(0,0,2,1,2,1,1) = 1618 | f(1,0,2,1,2,1,1) = 1424 |
| f(2,0,2,1,2,1,1) = 1941 | f(0,1,2,1,2,1,1) = 302 | f(1,1,2,1,2,1,1) = 1306 |
| f(2,1,2,1,2,1,1) = 1046 | f(0,2,2,1,2,1,1) = 1506 | f(1,2,2,1,2,1,1) = 1113 |
| f(2,2,2,1,2,1,1) = 1045 | f(0,0,0,2,2,1,1) = 1937 | f(1,0,0,2,2,1,1) = 1083 |
| f(2,0,0,2,2,1,1) = 1995 | f(0,1,0,2,2,1,1) = 1610 | f(1,1,0,2,2,1,1) = 1800 |
| f(2,1,0,2,2,1,1) = 962 | f(0,2,0,2,2,1,1) = 1234 | f(1,2,0,2,2,1,1) = 1224 |
| f(2,2,0,2,2,1,1) = 961 | f(0,0,1,2,2,1,1) = 734 | f(1,0,1,2,2,1,1) = 1708 |
| f(2,0,1,2,2,1,1) = 1683 | f(0,1,1,2,2,1,1) = 1963 | f(1,1,1,2,2,1,1) = 948 |
| f(2,1,1,2,2,1,1) = 1934 | f(0,2,1,2,2,1,1) = 1203 | f(1,2,1,2,2,1,1) = 756 |
| f(2,2,1,2,2,1,1) = 1741 | f(0,0,2,2,2,1,1) = 1627 | f(1,0,2,2,2,1,1) = 1372 |
| f(2,0,2,2,2,1,1) = 1379 | f(0,1,2,2,2,1,1) = 1899 | f(1,1,2,2,2,1,1) = 884 |
| f(2,1,2,2,2,1,1) = 1870 | f(0,2,2,2,2,1,1) = 1139 | f(1,2,2,2,2,1,1) = 500 |
| f(2,2,2,2,2,1,1) = 1485 | f(0,0,0,0,0,2,1) = 1033 | f(1,0,0,0,0,2,1) = 1084 |
| f(2,0,0,0,0,2,1) = 1615 | f(0,1,0,0,0,2,1) = 1792 | f(2,1,0,0,0,2,1) = 1805 |
| f(0,2,0,0,0,2,1) = 1216 | f(2,2,0,0,0,2,1) = 1230 | f(0,0,1,0,0,2,1) = 393 |
| f(1,0,1,0,0,2,1) = 1514 | f(2,0,1,0,0,2,1) = 1675 | f(0,1,1,0,0,2,1) = 1590 |
| f(1,1,1,0,0,2,1) = 1690 | f(2,1,1,0,0,2,1) = 1935 | f(0,2,1,0,0,2,1) = 1215 |
| f(1,2,1,0,0,2,1) = 1689 | f(2,2,1,0,0,2,1) = 1743 | f(0,0,2,0,0,2,1) = 883 |
| f(1,0,2,0,0,2,1) = 1945 | f(2,0,2,0,0,2,1) = 1355 | f(0,1,2,0,0,2,1) = 1334 |
| f(1,1,2,0,0,2,1) = 1386 | f(2,1,2,0,0,2,1) = 1871 | f(0,2,2,0,0,2,1) = 1151 |
| f(1,2,2,0,0,2,1) = 1385 | f(2,2,2,0,0,2,1) = 1487 | f(0,0,0,1,0,2,1) = 1426 |
| f(1,0,0,1,0,2,1) = 1660 | f(2,0,0,1,0,2,1) = 1015 | f(0,1,0,1,0,2,1) = 1825 |
| f(1,1,0,1,0,2,1) = 1948 | f(2,1,0,1,0,2,1) = 974 | f(0,2,0,1,0,2,1) = 1441 |
| f(1,2,0,1,0,2,1) = 1756 | f(2,2,0,1,0,2,1) = 973 | f(0,0,1,1,0,2,1) = 1484 |

f(1,0,1,1,0,2,1) = 1680     f(2,0,1,1,0,2,1) = 933     f(1,1,1,1,0,2,1) = 1595
f(2,1,1,1,0,2,1) = 1959     f(0,2,1,1,0,2,1) = 1206     f(1,2,1,1,0,2,1) = 1211
f(2,2,1,1,0,2,1) = 1767     f(0,0,2,1,0,2,1) = 1331     f(1,0,2,1,0,2,1) = 1376
f(2,0,2,1,0,2,1) = 470     f(1,1,2,1,0,2,1) = 1339     f(2,1,2,1,0,2,1) = 1879
f(0,2,2,1,0,2,1) = 1142     f(1,2,2,1,0,2,1) = 1147     f(2,2,2,1,0,2,1) = 1495
f(0,0,0,2,0,2,1) = 1373     f(1,0,0,2,0,2,1) = 1467     f(2,0,0,2,0,2,1) = 1039
f(0,1,0,2,0,2,1) = 1809     f(1,1,0,2,0,2,1) = 1888     f(2,1,0,2,0,2,1) = 1410
f(0,2,0,2,0,2,1) = 1617     f(1,2,0,2,0,2,1) = 1504     f(2,2,0,2,0,2,1) = 1601
f(0,0,1,2,0,2,1) = 1727     f(1,0,1,2,0,2,1) = 746     f(0,1,1,2,0,2,1) = 1599
f(1,1,1,2,0,2,1) = 940     f(2,1,1,2,0,2,1) = 678     f(0,2,1,2,0,2,1) = 1697
f(1,2,1,2,0,2,1) = 748     f(2,2,1,2,0,2,1) = 677     f(0,0,2,2,0,2,1) = 1407
f(1,0,2,2,0,2,1) = 857     f(0,1,2,2,0,2,1) = 1343     f(1,1,2,2,0,2,1) = 860
f(2,1,2,2,0,2,1) = 342     f(0,2,2,2,0,2,1) = 1377     f(1,2,2,2,0,2,1) = 476
f(2,2,2,2,0,2,1) = 341     f(1,0,0,0,1,2,1) = 1150     f(2,0,0,0,1,2,1) = 1102
f(1,1,0,0,1,2,1) = 1900     f(2,1,0,0,1,2,1) = 1414     f(0,2,0,0,1,2,1) = 1228
f(1,2,0,0,1,2,1) = 1516     f(2,2,0,0,1,2,1) = 1605     f(0,0,1,0,1,2,1) = 1453
f(1,0,1,0,1,2,1) = 1704     f(2,0,1,0,1,2,1) = 1166     f(0,1,1,0,1,2,1) = 1940
f(1,1,1,0,1,2,1) = 2026     f(2,1,1,0,1,2,1) = 1829     f(0,2,1,0,1,2,1) = 1757
f(1,2,1,0,1,2,1) = 2025     f(2,2,1,0,1,2,1) = 1254     f(1,0,2,0,1,2,1) = 1368
f(2,0,2,0,1,2,1) = 1293     f(0,1,2,0,1,2,1) = 1876     f(1,1,2,0,1,2,1) = 2010
f(2,1,2,0,1,2,1) = 1813     f(0,2,2,0,1,2,1) = 1501     f(1,2,2,0,1,2,1) = 2009
f(2,2,2,0,1,2,1) = 1238     f(0,0,0,1,1,2,1) = 1118     f(1,0,0,1,1,2,1) = 1460
f(2,0,0,1,1,2,1) = 1925     f(0,1,0,1,1,2,1) = 1828     f(1,1,0,1,1,2,1) = 58
f(2,1,0,1,1,2,1) = 534     f(0,2,0,1,1,2,1) = 1425     f(1,2,0,1,1,2,1) = 57
f(2,2,0,1,1,2,1) = 149     f(1,0,1,1,1,2,1) = 656     f(2,0,1,1,1,2,1) = 1443
f(1,1,1,1,1,2,1) = 946     f(2,1,1,1,1,2,1) = 1931     f(0,2,1,1,1,2,1) = 1748
f(1,2,1,1,1,2,1) = 753     f(2,2,1,1,1,2,1) = 1739     f(0,0,2,1,1,2,1) = 354
f(1,0,2,1,1,2,1) = 352     f(2,0,2,1,1,2,1) = 1427     f(1,1,2,1,1,2,1) = 882
f(2,1,2,1,1,2,1) = 1867     f(0,2,2,1,1,2,1) = 1492     f(1,2,2,1,1,2,1) = 497
f(2,2,2,1,1,2,1) = 1483     f(0,0,0,2,1,2,1) = 1364     f(1,0,0,2,1,2,1) = 440
f(2,0,0,2,1,2,1) = 1079     f(0,1,0,2,1,2,1) = 1812     f(1,1,0,2,1,2,1) = 818
f(2,1,0,2,1,2,1) = 1803     f(0,2,0,2,1,2,1) = 1633     f(1,2,0,2,1,2,1) = 241
f(2,2,0,2,1,2,1) = 1227     f(0,0,1,2,1,2,1) = 1005     f(1,0,1,2,1,2,1) = 1194
f(2,0,1,2,1,2,1) = 726     f(0,1,1,2,1,2,1) = 1949     f(1,1,1,2,1,2,1) = 1928
f(2,1,1,2,1,2,1) = 805     f(0,2,1,2,1,2,1) = 1761     f(1,2,1,2,1,2,1) = 1736
f(2,2,1,2,1,2,1) = 230     f(1,0,2,2,1,2,1) = 1305     f(2,0,2,2,1,2,1) = 869
f(0,1,2,2,1,2,1) = 1885     f(1,1,2,2,1,2,1) = 1864     f(2,1,2,2,1,2,1) = 789
f(0,2,2,2,1,2,1) = 1505     f(1,2,2,2,1,2,1) = 1480     f(2,2,2,2,1,2,1) = 214
f(1,0,0,0,2,2,1) = 1597     f(2,0,0,0,2,2,1) = 1549     f(1,1,0,0,2,2,1) = 626
f(2,1,0,0,2,2,1) = 870     f(0,2,0,0,2,2,1) = 1226     f(1,2,0,0,2,2,1) = 433
f(2,2,0,0,2,2,1) = 485     f(0,0,1,0,2,2,1) = 1437     f(1,0,1,0,2,2,1) = 186
f(2,0,1,0,2,2,1) = 1635     f(0,1,1,0,2,2,1) = 1938     f(1,1,1,0,2,2,1) = 570
f(2,1,1,0,2,2,1) = 907     f(0,2,1,0,2,2,1) = 1755     f(1,2,1,0,2,2,1) = 185
f(2,2,1,0,2,2,1) = 715     f(1,0,2,0,2,2,1) = 313     f(2,0,2,0,2,2,1) = 1619

f(0,1,2,0,2,2,1) = 1874
f(0,2,2,0,2,2,1) = 1499
f(0,0,0,1,2,2,1) = 1435
f(1,1,0,1,2,2,1) = 1853
f(1,2,0,1,2,2,1) = 1278
f(1,0,1,1,2,2,1) = 1712
f(0,2,1,1,2,2,1) = 1746
f(0,0,2,1,2,2,1) = 338
f(2,1,2,1,2,2,1) = 851
f(2,2,2,1,2,2,1) = 467
f(2,0,0,2,2,2,1) = 1415
f(2,1,0,2,2,2,1) = 1606
f(2,2,0,2,2,2,1) = 1413
f(2,0,1,2,2,2,1) = 150
f(2,1,1,2,2,2,1) = 646
f(2,2,1,2,2,2,1) = 645
f(2,0,2,2,2,2,1) = 293
f(2,1,2,2,2,2,1) = 326
f(2,2,2,2,2,2,1) = 325
f(2,1,0,0,0,0,2) = 386
f(2,2,0,0,0,0,2) = 577
f(1,1,1,0,0,0,2) = 1598
f(1,2,1,0,0,0,2) = 1213
f(1,0,2,0,0,0,2) = 336
f(0,2,2,0,0,0,2) = 1535
f(0,0,0,1,0,0,2) = 402
f(1,1,0,1,0,0,2) = 1852
f(2,2,0,1,0,0,2) = 5
f(2,0,1,1,0,0,2) = 995
f(2,1,1,1,0,0,2) = 515
f(2,2,1,1,0,0,2) = 131
f(2,0,2,1,0,0,2) = 979
f(2,1,2,1,0,0,2) = 259
f(2,2,2,1,0,0,2) = 67
f(2,0,0,2,0,0,2) = 1423
f(2,1,0,2,0,0,2) = 1555
f(0,0,1,2,0,0,2) = 649
f(1,1,1,2,0,0,2) = 808
f(1,2,1,2,0,0,2) = 232
f(1,0,2,2,0,0,2) = 600
f(2,1,2,2,0,0,2) = 863
f(2,2,2,2,0,0,2) = 479
f(0,1,0,0,1,0,2) = 1845
f(0,2,0,0,1,0,2) = 53

f(1,1,2,0,2,2,1) = 314
f(1,2,2,0,2,2,1) = 121
f(1,0,0,1,2,2,1) = 1652
f(2,1,0,1,2,2,1) = 590
f(2,2,0,1,2,2,1) = 397
f(1,1,1,1,2,2,1) = 1450
f(1,2,1,1,2,2,1) = 1641
f(1,0,2,1,2,2,1) = 1392
f(0,2,2,1,2,2,1) = 1490
f(0,0,0,2,2,2,1) = 565
f(0,1,0,2,2,2,1) = 786
f(0,2,0,2,2,2,1) = 1651
f(0,0,1,2,2,2,1) = 989
f(0,1,1,2,2,2,1) = 1947
f(0,2,1,2,2,2,1) = 1185
f(0,0,2,2,2,2,1) = 630
f(0,1,2,2,2,2,1) = 1883
f(0,2,2,2,2,2,1) = 1121
f(1,0,0,0,0,0,2) = 48
f(0,2,0,0,0,0,2) = 1279
f(0,0,1,0,0,0,2) = 1471
f(2,1,1,0,0,0,2) = 1583
f(2,2,1,0,0,0,2) = 1199
f(1,1,2,0,0,0,2) = 1342
f(1,2,2,0,0,0,2) = 1149
f(2,0,0,1,0,0,2) = 4
f(2,1,0,1,0,0,2) = 6
f(0,0,1,1,0,0,2) = 1973
f(0,1,1,1,0,0,2) = 1694
f(0,2,1,1,0,0,2) = 1782
f(0,0,2,1,0,0,2) = 1907
f(0,1,2,1,0,0,2) = 1374
f(0,2,2,1,0,0,2) = 1526
f(0,0,0,2,0,0,2) = 349
f(0,1,0,2,0,0,2) = 30
f(1,2,0,2,0,0,2) = 1266
f(1,0,1,2,0,0,2) = 616
f(2,1,1,2,0,0,2) = 943
f(2,2,1,2,0,0,2) = 751
f(2,0,2,2,0,0,2) = 1359
f(0,2,2,2,0,0,2) = 337
f(0,0,0,0,1,0,2) = 1993
f(1,1,0,0,1,0,2) = 1844
f(1,2,0,0,1,0,2) = 1268

f(2,1,2,0,2,2,1) = 843
f(2,2,2,0,2,2,1) = 459
f(0,1,0,1,2,2,1) = 802
f(0,2,0,1,2,2,1) = 1459
f(0,0,1,1,2,2,1) = 404
f(2,1,1,1,2,2,1) = 931
f(2,2,1,1,2,2,1) = 739
f(1,1,2,1,2,2,1) = 1626
f(1,2,2,1,2,2,1) = 1433
f(1,0,0,2,2,2,1) = 1593
f(1,1,0,2,2,2,1) = 1466
f(1,2,0,2,2,2,1) = 1657
f(1,0,1,2,2,2,1) = 1214
f(1,1,1,2,2,2,1) = 568
f(1,2,1,2,2,2,1) = 184
f(1,0,2,2,2,2,1) = 1341
f(1,1,2,2,2,2,1) = 312
f(1,2,2,2,2,2,1) = 120
f(1,1,0,0,0,0,2) = 776
f(1,2,0,0,0,0,2) = 200
f(1,0,1,0,0,0,2) = 672
f(0,2,1,0,0,0,2) = 1791
f(0,0,2,0,0,0,2) = 1353
f(2,1,2,0,0,0,2) = 1311
f(2,2,2,0,0,0,2) = 1119
f(0,1,0,1,0,0,2) = 46
f(1,2,0,1,0,0,2) = 1276
f(1,0,1,1,0,0,2) = 1672
f(1,1,1,1,0,0,2) = 1979
f(1,2,1,1,0,0,2) = 1787
f(1,0,2,1,0,0,2) = 1352
f(1,1,2,1,0,0,2) = 1915
f(1,2,2,1,0,0,2) = 1531
f(1,0,0,2,0,0,2) = 8
f(1,1,0,2,0,0,2) = 1841
f(2,2,0,2,0,0,2) = 1171
f(2,0,1,2,0,0,2) = 1679
f(0,2,1,2,0,0,2) = 657
f(0,0,2,2,0,0,2) = 885
f(1,1,2,2,0,0,2) = 792
f(1,2,2,2,0,0,2) = 216
f(1,0,0,0,1,0,2) = 114
f(2,1,0,0,1,0,2) = 1559
f(2,2,0,0,1,0,2) = 1175

f(0,0,1,0,1,0,2) = 2029   f(1,0,1,0,1,0,2) = 680    f(0,1,1,0,1,0,2) = 932
f(1,1,1,0,1,0,2) = 1592   f(2,1,1,0,1,0,2) = 1542   f(0,2,1,0,1,0,2) = 749
f(1,2,1,0,1,0,2) = 1208   f(2,2,1,0,1,0,2) = 1157   f(0,0,2,0,1,0,2) = 865
f(1,0,2,0,1,0,2) = 344    f(0,1,2,0,1,0,2) = 868    f(1,1,2,0,1,0,2) = 1336
f(2,1,2,0,1,0,2) = 1286   f(0,2,2,0,1,0,2) = 493    f(1,2,2,0,1,0,2) = 1144
f(2,2,2,0,1,0,2) = 1093   f(0,0,0,1,1,0,2) = 458    f(1,0,0,1,1,0,2) = 1076
f(0,1,0,1,1,0,2) = 804    f(1,1,0,1,1,0,2) = 638    f(2,1,0,1,1,0,2) = 927
f(0,2,0,1,1,0,2) = 1461   f(1,2,0,1,1,0,2) = 445    f(2,2,0,1,1,0,2) = 735
f(0,0,1,1,1,0,2) = 996    f(1,0,1,1,1,0,2) = 730    f(2,0,1,1,1,0,2) = 130
f(0,1,1,1,1,0,2) = 1950   f(1,1,1,1,1,0,2) = 1586   f(2,1,1,1,1,0,2) = 1543
f(0,2,1,1,1,0,2) = 740    f(1,2,1,1,1,0,2) = 1201   f(2,2,1,1,1,0,2) = 1159
f(0,0,2,1,1,0,2) = 994    f(1,0,2,1,1,0,2) = 873    f(2,0,2,1,1,0,2) = 257
f(0,1,2,1,1,0,2) = 1886   f(1,1,2,1,1,0,2) = 1330   f(2,1,2,1,1,0,2) = 1287
f(0,2,2,1,1,0,2) = 484    f(1,2,2,1,1,0,2) = 1137   f(2,2,2,1,1,0,2) = 1095
f(0,0,0,2,1,0,2) = 340    f(1,0,0,2,1,0,2) = 506    f(2,0,0,2,1,0,2) = 1463
f(0,1,0,2,1,0,2) = 788    f(1,1,0,2,1,0,2) = 1842   f(2,1,0,2,1,0,2) = 1795
f(0,2,0,2,1,0,2) = 1653   f(1,2,0,2,1,0,2) = 1265   f(2,2,0,2,1,0,2) = 1219
f(0,0,1,2,1,0,2) = 685    f(1,0,1,2,1,0,2) = 1770   f(2,0,1,2,1,0,2) = 1647
f(0,1,1,2,1,0,2) = 941    f(1,1,1,2,1,0,2) = 618    f(2,1,1,2,1,0,2) = 550
f(0,2,1,2,1,0,2) = 721    f(1,2,1,2,1,0,2) = 425    f(2,2,1,2,1,0,2) = 165
f(0,0,2,2,1,0,2) = 43     f(1,0,2,2,1,0,2) = 1881   f(2,0,2,2,1,0,2) = 1631
f(0,1,2,2,1,0,2) = 877    f(1,1,2,2,1,0,2) = 410    f(2,1,2,2,1,0,2) = 278
f(0,2,2,2,1,0,2) = 465    f(1,2,2,2,1,0,2) = 601    f(2,2,2,2,1,0,2) = 85
f(0,0,0,0,2,0,2) = 9      f(1,0,0,0,2,0,2) = 561    f(0,1,0,0,2,0,2) = 1843
f(1,1,0,0,2,0,2) = 826    f(2,1,0,0,2,0,2) = 775    f(0,2,0,0,2,0,2) = 51
f(1,2,0,0,2,0,2) = 249    f(2,2,0,0,2,0,2) = 199    f(0,0,1,0,2,0,2) = 2013
f(1,0,1,0,2,0,2) = 178    f(2,0,1,0,2,0,2) = 611    f(0,1,1,0,2,0,2) = 930
f(1,1,1,0,2,0,2) = 554    f(2,1,1,0,2,0,2) = 551    f(0,2,1,0,2,0,2) = 747
f(1,2,1,0,2,0,2) = 169    f(2,2,1,0,2,0,2) = 167    f(0,0,2,0,2,0,2) = 849
f(1,0,2,0,2,0,2) = 305    f(2,0,2,0,2,0,2) = 595    f(0,1,2,0,2,0,2) = 866
f(1,1,2,0,2,0,2) = 282    f(2,1,2,0,2,0,2) = 279    f(0,2,2,0,2,0,2) = 491
f(1,2,2,0,2,0,2) = 89     f(2,2,2,0,2,0,2) = 87     f(0,0,0,1,2,0,2) = 411
f(1,0,0,1,2,0,2) = 2036   f(2,0,0,1,2,0,2) = 583    f(0,1,0,1,2,0,2) = 1826
f(1,1,0,1,2,0,2) = 528    f(2,1,0,1,2,0,2) = 1567   f(0,2,0,1,2,0,2) = 435
f(1,2,0,1,2,0,2) = 144    f(2,2,0,1,2,0,2) = 1183   f(0,0,1,1,2,0,2) = 1502
f(1,0,1,1,2,0,2) = 1961   f(2,0,1,1,2,0,2) = 1007   f(0,1,1,1,2,0,2) = 1566
f(1,1,1,1,2,0,2) = 1066   f(2,1,1,1,2,0,2) = 547    f(0,2,1,1,2,0,2) = 738
f(1,2,1,1,2,0,2) = 1065   f(2,2,1,1,2,0,2) = 163    f(0,0,2,1,2,0,2) = 978
f(1,0,2,1,2,0,2) = 1498   f(2,0,2,1,2,0,2) = 991    f(0,1,2,1,2,0,2) = 1310
f(1,1,2,1,2,0,2) = 1050   f(2,1,2,1,2,0,2) = 275    f(0,2,2,1,2,0,2) = 482
f(1,2,2,1,2,0,2) = 1049   f(2,2,2,1,2,0,2) = 83     f(0,0,0,2,2,0,2) = 1333
f(1,0,0,2,2,0,2) = 59     f(2,0,0,2,2,0,2) = 1031   f(0,1,0,2,2,0,2) = 1810
f(1,1,0,2,2,0,2) = 1322   f(2,1,0,2,2,0,2) = 867    f(0,2,0,2,2,0,2) = 627
f(1,2,0,2,2,0,2) = 1129   f(2,2,0,2,2,0,2) = 483    f(0,0,1,2,2,0,2) = 669

f(1,0,1,2,2,0,2) = 44
f(1,1,1,2,2,0,2) = 572
f(1,2,1,2,2,0,2) = 188
f(1,0,2,2,2,0,2) = 28
f(1,1,2,2,2,0,2) = 316
f(1,2,2,2,2,0,2) = 124
f(1,0,0,0,0,1,2) = 432
f(1,1,0,0,0,1,2) = 817
f(0,0,1,0,0,1,2) = 1164
f(0,1,1,0,0,1,2) = 896
f(0,2,1,0,0,1,2) = 713
f(0,0,2,0,0,1,2) = 1654
f(0,1,2,0,0,1,2) = 832
f(0,2,2,0,0,1,2) = 457
f(0,0,0,1,0,1,2) = 674
f(0,1,0,1,0,1,2) = 430
f(0,2,0,1,0,1,2) = 238
f(0,0,1,1,0,1,2) = 640
f(1,1,1,1,0,1,2) = 1706
f(1,2,1,1,0,1,2) = 1705
f(2,0,2,1,0,1,2) = 1301
f(2,1,2,1,0,1,2) = 1299
f(2,2,2,1,0,1,2) = 1107
f(2,0,0,2,0,1,2) = 387
f(2,1,0,2,0,1,2) = 291
f(2,2,0,2,0,1,2) = 99
f(2,0,1,2,0,1,2) = 671
f(2,1,1,2,0,1,2) = 900
f(0,0,2,2,0,1,2) = 563
f(0,1,2,2,0,1,2) = 841
f(1,2,2,2,0,1,2) = 88
f(2,0,0,0,1,1,2) = 450
f(2,1,0,0,1,1,2) = 1614
f(1,0,1,0,1,1,2) = 428
f(1,1,1,0,1,1,2) = 1588
f(1,2,1,0,1,1,2) = 1204
f(1,0,2,0,1,1,2) = 412
f(1,1,2,0,1,1,2) = 1332
f(1,2,2,0,1,1,2) = 1140
f(1,0,0,1,1,1,2) = 632
f(1,1,0,1,1,1,2) = 634
f(1,2,0,1,1,1,2) = 441
f(1,0,1,1,1,1,2) = 1754
f(1,1,1,1,1,1,2) = 1722

f(2,0,1,2,2,0,2) = 659
f(2,1,1,2,2,0,2) = 526
f(2,2,1,2,2,0,2) = 141
f(2,0,2,2,2,0,2) = 355
f(2,1,2,2,2,0,2) = 270
f(2,2,2,2,2,0,2) = 77
f(2,0,0,0,0,1,2) = 963
f(0,2,0,0,0,1,2) = 255
f(1,0,1,0,0,1,2) = 692
f(1,1,1,0,0,1,2) = 560
f(1,2,1,0,0,1,2) = 176
f(1,0,2,0,0,1,2) = 372
f(1,1,2,0,0,1,2) = 304
f(1,2,2,0,0,1,2) = 112
f(1,0,0,1,0,1,2) = 1008
f(1,1,0,1,0,1,2) = 446
f(1,2,0,1,0,1,2) = 637
f(2,0,1,1,0,1,2) = 1190
f(2,1,1,1,0,1,2) = 1571
f(2,2,1,1,0,1,2) = 1187
f(0,1,2,1,0,1,2) = 350
f(0,2,2,1,0,1,2) = 448
f(0,0,0,2,0,1,2) = 621
f(0,1,0,2,0,1,2) = 606
f(0,2,0,2,0,1,2) = 222
f(0,0,1,2,0,1,2) = 716
f(0,1,1,2,0,1,2) = 905
f(1,2,1,2,0,1,2) = 168
f(1,0,2,2,0,1,2) = 1114
f(1,1,2,2,0,1,2) = 280
f(2,2,2,2,0,1,2) = 452
f(0,1,0,0,1,1,2) = 821
f(1,2,0,0,1,1,2) = 1177
f(2,0,1,0,1,1,2) = 1734
f(2,1,1,0,1,1,2) = 1926
f(2,2,1,0,1,1,2) = 1733
f(2,0,2,0,1,1,2) = 1861
f(2,1,2,0,1,1,2) = 1862
f(2,2,2,0,1,1,2) = 1477
f(2,0,0,1,1,1,2) = 454
f(2,1,0,1,1,1,2) = 390
f(2,2,0,1,1,1,2) = 581
f(2,0,1,1,1,1,2) = 706
f(2,1,1,1,1,1,2) = 1927

f(0,1,1,2,2,0,2) = 939
f(0,2,1,2,2,0,2) = 145
f(0,0,2,2,2,0,2) = 529
f(0,1,2,2,2,0,2) = 875
f(0,2,2,2,2,0,2) = 81
f(0,0,0,0,1,2) = 1014
f(0,1,0,0,0,1,2) = 831
f(1,2,0,0,0,1,2) = 242
f(2,0,1,0,0,1,2) = 102
f(2,1,1,0,0,1,2) = 662
f(2,2,1,0,0,1,2) = 661
f(2,0,2,0,0,1,2) = 533
f(2,1,2,0,0,1,2) = 358
f(2,2,2,0,0,1,2) = 357
f(2,0,0,1,0,1,2) = 580
f(2,1,0,1,0,1,2) = 543
f(2,2,0,1,0,1,2) = 159
f(0,1,1,1,0,1,2) = 670
f(0,2,1,1,0,1,2) = 704
f(0,0,2,1,0,1,2) = 320
f(1,1,2,1,0,1,2) = 1370
f(1,2,2,1,0,1,2) = 1369
f(1,0,0,2,0,1,2) = 1032
f(1,1,0,2,0,1,2) = 1074
f(1,2,0,2,0,1,2) = 1073
f(1,0,1,2,0,1,2) = 1577
f(1,1,1,2,0,1,2) = 552
f(2,2,1,2,0,1,2) = 708
f(2,0,2,2,0,1,2) = 367
f(2,1,2,2,0,1,2) = 836
f(1,0,0,0,1,1,2) = 498
f(1,1,0,0,1,1,2) = 1562
f(2,2,0,0,1,1,2) = 1421
f(0,1,1,0,1,1,2) = 548
f(0,2,1,0,1,1,2) = 173
f(0,0,2,0,1,1,2) = 610
f(0,1,2,0,1,1,2) = 292
f(0,2,2,0,1,1,2) = 109
f(0,0,0,1,1,1,2) = 1482
f(0,1,0,1,1,1,2) = 396
f(0,2,0,1,1,1,2) = 1261
f(0,0,1,1,1,1,2) = 1417
f(0,1,1,1,1,1,2) = 926
f(0,2,1,1,1,1,2) = 164

f(1,2,1,1,1,2) = 1721   f(2,2,1,1,1,2) = 1735   f(0,0,2,1,1,2) = 1058
f(1,0,2,1,1,2) = 1897   f(2,0,2,1,1,2) = 833    f(0,1,2,1,1,2) = 862
f(1,1,2,1,1,2) = 1402   f(2,1,2,1,1,2) = 1863   f(0,2,2,1,1,2) = 100
f(1,2,2,1,1,2) = 1401   f(2,2,2,1,1,2) = 1479   f(0,0,0,2,1,2) = 612
f(2,0,0,2,1,2) = 395    f(0,1,0,2,1,2) = 588    f(1,0,0,2,1,2) = 1650
f(2,1,0,2,1,2) = 769    f(0,2,0,2,1,2) = 1245   f(1,2,0,2,1,2) = 1457
f(2,2,0,2,1,2) = 194    f(0,0,1,2,1,2) = 1709   f(2,0,1,2,1,2) = 655
f(0,1,1,2,1,2) = 557    f(1,1,1,2,1,2) = 1580   f(2,1,1,2,1,2) = 614
f(1,2,1,2,1,2) = 1196   f(2,2,1,2,1,2) = 421    f(0,0,2,2,1,2) = 1643
f(2,0,2,2,1,2) = 335    f(0,1,2,2,1,2) = 301    f(1,1,2,2,1,2) = 1308
f(2,1,2,2,1,2) = 406    f(1,2,2,2,1,2) = 1116   f(2,2,2,2,1,2) = 597
f(1,0,0,0,2,1,2) = 945  f(2,0,0,0,2,1,2) = 897  f(0,1,0,0,2,1,2) = 819
f(1,1,0,0,2,1,2) = 442  f(2,1,0,0,2,1,2) = 531  f(1,2,0,0,2,1,2) = 633
f(2,2,0,0,2,1,2) = 147  f(1,0,1,0,2,1,2) = 754  f(2,0,1,0,2,1,2) = 679
f(0,1,1,0,2,1,2) = 546  f(1,1,1,0,2,1,2) = 809  f(2,1,1,0,2,1,2) = 38
f(0,2,1,0,2,1,2) = 171  f(1,2,1,0,2,1,2) = 234  f(2,2,1,0,2,1,2) = 37
f(0,0,2,0,2,1,2) = 594  f(1,0,2,0,2,1,2) = 881  f(2,0,2,0,2,1,2) = 343
f(0,1,2,0,2,1,2) = 290  f(1,1,2,0,2,1,2) = 793  f(2,1,2,0,2,1,2) = 22
f(0,2,2,0,2,1,2) = 107  f(1,2,2,0,2,1,2) = 218  f(2,2,2,0,2,1,2) = 21
f(0,0,0,1,2,1,2) = 683  f(1,0,0,1,2,1,2) = 968  f(2,0,0,1,2,1,2) = 1607
f(0,1,0,1,2,1,2) = 414  f(1,1,0,1,2,1,2) = 820  f(2,1,0,1,2,1,2) = 1806
f(0,2,0,1,2,1,2) = 235  f(1,2,0,1,2,1,2) = 244  f(2,2,0,1,2,1,2) = 1229
f(1,0,1,1,2,1,2) = 1178 f(2,0,1,1,2,1,2) = 742  f(0,1,1,1,2,1,2) = 542
f(1,1,1,1,2,1,2) = 1833 f(2,1,1,1,2,1,2) = 567  f(0,2,1,1,2,1,2) = 162
f(1,2,1,1,2,1,2) = 1258 f(2,2,1,1,2,1,2) = 183  f(0,0,2,1,2,1,2) = 1042
f(1,0,2,1,2,1,2) = 1321 f(2,0,2,1,2,1,2) = 853  f(0,1,2,1,2,1,2) = 286
f(1,1,2,1,2,1,2) = 1817 f(2,1,2,1,2,1,2) = 311  f(0,2,2,1,2,1,2) = 98
f(1,2,2,1,2,1,2) = 1242 f(2,2,2,1,2,1,2) = 119  f(0,0,0,2,2,1,2) = 929
f(1,0,0,2,2,1,2) = 122  f(2,0,0,2,2,1,2) = 587  f(0,1,0,2,2,1,2) = 622
f(1,1,0,2,2,1,2) = 1898 f(2,1,0,2,2,1,2) = 1990 f(0,2,0,2,2,1,2) = 219
f(1,2,0,2,2,1,2) = 1513 f(2,2,0,2,2,1,2) = 1989 f(0,0,1,2,2,1,2) = 1693
f(1,0,1,2,2,1,2) = 620  f(2,0,1,2,2,1,2) = 1695 f(0,1,1,2,2,1,2) = 555
f(1,1,1,2,2,1,2) = 564  f(2,1,1,2,2,1,2) = 1550 f(1,2,1,2,2,1,2) = 180
f(2,2,1,2,2,1,2) = 1165 f(1,0,2,2,2,1,2) = 604  f(2,0,2,2,2,1,2) = 1391
f(0,1,2,2,2,1,2) = 299  f(1,1,2,2,2,1,2) = 308  f(2,1,2,2,2,1,2) = 1294
f(1,2,2,2,2,1,2) = 116  f(2,2,2,2,2,1,2) = 1101 f(0,0,0,0,2,2) = 969
f(1,0,0,0,0,2,2) = 444  f(2,0,0,0,0,2,2) = 975  f(0,1,0,0,0,2,2) = 777
f(1,1,0,0,0,2,2) = 1082 f(2,1,0,0,0,2,2) = 1327 f(0,2,0,0,0,2,2) = 201
f(1,2,0,0,0,2,2) = 1081 f(2,2,0,0,0,2,2) = 1135 f(0,0,1,0,0,2,2) = 949
f(1,0,1,0,0,2,2) = 762  f(2,0,1,0,0,2,2) = 687  f(0,1,1,0,0,2,2) = 950
f(1,1,1,0,0,2,2) = 1726 f(2,1,1,0,0,2,2) = 1839 f(0,2,1,0,0,2,2) = 767
f(1,2,1,0,0,2,2) = 1725 f(2,2,1,0,0,2,2) = 1263 f(0,0,2,0,0,2,2) = 374
f(1,0,2,0,0,2,2) = 889  f(2,0,2,0,0,2,2) = 351  f(0,1,2,0,0,2,2) = 886
f(1,1,2,0,0,2,2) = 1406 f(2,1,2,0,0,2,2) = 1823 f(0,2,2,0,0,2,2) = 511

f(1,2,2,0,0,2,2) = 1405
f(2,2,2,0,0,2,2) = 1247
f(0,0,0,1,0,2,2) = 658
f(1,0,0,1,0,2,2) = 1020
f(0,1,0,1,0,2,2) = 814
f(1,1,0,1,0,2,2) = 536
f(2,1,0,1,0,2,2) = 294
f(0,2,0,1,0,2,2) = 993
f(1,2,0,1,0,2,2) = 152
f(2,2,0,1,0,2,2) = 101
f(0,0,1,1,0,2,2) = 438
f(1,0,1,1,0,2,2) = 154
f(2,0,1,1,0,2,2) = 419
f(0,1,1,1,0,2,2) = 652
f(1,1,1,1,0,2,2) = 955
f(2,1,1,1,0,2,2) = 899
f(0,2,1,1,0,2,2) = 758
f(1,2,1,1,0,2,2) = 763
f(2,2,1,1,0,2,2) = 707
f(0,0,2,1,0,2,2) = 947
f(1,0,2,1,0,2,2) = 297
f(2,0,2,1,0,2,2) = 403
f(0,1,2,1,0,2,2) = 332
f(1,1,2,1,0,2,2) = 891
f(2,1,2,1,0,2,2) = 835
f(0,2,2,1,0,2,2) = 502
f(1,2,2,1,0,2,2) = 507
f(2,2,2,1,0,2,2) = 451
f(0,0,0,2,0,2,2) = 605
f(1,0,0,2,0,2,2) = 1019
f(2,0,0,2,0,2,2) = 399
f(0,1,0,2,0,2,2) = 798
f(1,1,0,2,0,2,2) = 874
f(2,1,0,2,0,2,2) = 1794
f(0,2,0,2,0,2,2) = 977
f(1,2,0,2,0,2,2) = 489
f(2,2,0,2,0,2,2) = 1217
f(0,0,1,2,0,2,2) = 447
f(1,0,1,2,0,2,2) = 40
f(2,0,1,2,0,2,2) = 47
f(0,1,1,2,0,2,2) = 959
f(1,1,1,2,0,2,2) = 666
f(2,1,1,2,0,2,2) = 559
f(0,2,1,2,0,2,2) = 1681
f(1,2,1,2,0,2,2) = 665
f(2,2,1,2,0,2,2) = 175
f(0,0,2,2,0,2,2) = 347
f(1,0,2,2,0,2,2) = 24
f(2,0,2,2,0,2,2) = 31
f(0,1,2,2,0,2,2) = 895
f(1,1,2,2,0,2,2) = 362
f(2,1,2,2,0,2,2) = 287
f(0,2,2,2,0,2,2) = 1361
f(1,2,2,2,0,2,2) = 361
f(2,2,2,2,0,2,2) = 95
f(1,0,0,0,1,2,2) = 510
f(2,0,0,0,1,2,2) = 462
f(0,1,0,0,1,2,2) = 12
f(1,1,0,0,1,2,2) = 1849
f(2,1,0,0,1,2,2) = 1318
f(0,2,0,0,1,2,2) = 1077
f(1,2,0,0,1,2,2) = 1274
f(2,2,0,0,1,2,2) = 1125
f(1,0,1,0,1,2,2) = 700
f(2,0,1,0,1,2,2) = 718
f(0,1,1,0,1,2,2) = 532
f(1,1,1,0,1,2,2) = 1576
f(2,1,1,0,1,2,2) = 1574
f(0,2,1,0,1,2,2) = 157
f(1,2,1,0,1,2,2) = 1192
f(2,2,1,0,1,2,2) = 1189
f(0,0,2,0,1,2,2) = 1387
f(1,0,2,0,1,2,2) = 380
f(2,0,2,0,1,2,2) = 845
f(0,1,2,0,1,2,2) = 276
f(1,1,2,0,1,2,2) = 1304
f(2,1,2,0,1,2,2) = 1302
f(0,2,2,0,1,2,2) = 93
f(1,2,2,0,1,2,2) = 1112
f(2,2,2,0,1,2,2) = 1109
f(0,0,0,1,1,2,2) = 110
f(1,0,0,1,1,2,2) = 52
f(2,0,0,1,1,2,2) = 964
f(0,1,0,1,1,2,2) = 813
f(1,1,0,1,1,2,2) = 1086
f(2,1,0,1,1,2,2) = 823
f(0,2,0,1,1,2,2) = 437
f(1,2,0,1,1,2,2) = 1085
f(2,2,0,1,1,2,2) = 247
f(0,0,1,1,1,2,2) = 420
f(1,0,1,1,1,2,2) = 668
f(2,0,1,1,1,2,2) = 675
f(0,1,1,1,1,2,2) = 908
f(1,1,1,1,1,2,2) = 562
f(2,1,1,1,1,2,2) = 1547
f(0,2,1,1,1,2,2) = 148
f(1,2,1,1,1,2,2) = 177
f(2,2,1,1,1,2,2) = 1163
f(0,0,2,1,1,2,2) = 1313
f(1,0,2,1,1,2,2) = 364
f(2,0,2,1,1,2,2) = 339
f(0,1,2,1,1,2,2) = 844
f(1,1,2,1,1,2,2) = 306
f(2,1,2,1,1,2,2) = 1291
f(0,2,2,1,1,2,2) = 84
f(1,2,2,1,1,2,2) = 113
f(2,2,2,1,1,2,2) = 1099
f(0,0,0,2,1,2,2) = 596
f(1,0,0,2,1,2,2) = 1977
f(2,0,0,2,1,2,2) = 439
f(0,1,0,2,1,2,2) = 797
f(1,1,0,2,1,2,2) = 434
f(2,1,0,2,1,2,2) = 1319
f(0,2,0,2,1,2,2) = 629
f(1,2,0,2,1,2,2) = 625
f(2,2,0,2,1,2,2) = 1127
f(0,0,1,2,1,2,2) = 429
f(1,0,1,2,1,2,2) = 106
f(2,0,1,2,1,2,2) = 623
f(0,1,1,2,1,2,2) = 541
f(1,1,1,2,1,2,2) = 1002
f(2,1,1,2,1,2,2) = 934
f(0,2,1,2,1,2,2) = 1745
f(1,2,1,2,1,2,2) = 1001
f(2,2,1,2,1,2,2) = 741
f(0,0,2,2,1,2,2) = 363
f(1,0,2,2,1,2,2) = 537
f(2,0,2,2,1,2,2) = 607
f(0,1,2,2,1,2,2) = 285
f(1,1,2,2,1,2,2) = 986
f(2,1,2,2,1,2,2) = 854
f(0,2,2,2,1,2,2) = 1489

f(1,2,2,2,1,2,2) = 985         f(2,2,2,2,1,2,2) = 469         f(1,0,0,0,2,2,2) = 957
f(2,0,0,0,2,2,2) = 909         f(0,1,0,0,2,2,2) = 10          f(1,1,0,0,2,2,2) = 827
f(2,1,0,0,2,2,2) = 295         f(0,2,0,0,2,2,2) = 1075        f(1,2,0,0,2,2,2) = 251
f(2,2,0,0,2,2,2) = 103         f(1,0,1,0,2,2,2) = 1723        f(2,0,1,0,2,2,2) = 35
f(0,1,1,0,2,2,2) = 530         f(2,1,1,0,2,2,2) = 519         f(0,2,1,0,2,2,2) = 155
f(2,2,1,0,2,2,2) = 135         f(0,0,2,0,2,2,2) = 1371        f(1,0,2,0,2,2,2) = 1403
f(2,0,2,0,2,2,2) = 19          f(0,1,2,0,2,2,2) = 274         f(2,1,2,0,2,2,2) = 263
f(0,2,2,0,2,2,2) = 91          f(2,2,2,0,2,2,2) = 71          f(0,0,0,1,2,2,2) = 667
f(1,0,0,1,2,2,2) = 1012        f(2,0,0,1,2,2,2) = 631         f(0,1,0,1,2,2,2) = 811
f(1,1,0,1,2,2,2) = 825         f(2,1,0,1,2,2,2) = 303         f(0,2,0,1,2,2,2) = 33
f(1,2,0,1,2,2,2) = 250         f(2,2,0,1,2,2,2) = 111         f(1,0,1,1,2,2,2) = 937
f(2,0,1,1,2,2,2) = 431         f(0,1,1,1,2,2,2) = 524         f(1,1,1,1,2,2,2) = 42
f(2,1,1,1,2,2,2) = 951         f(0,2,1,1,2,2,2) = 146         f(1,2,1,1,2,2,2) = 41
f(2,2,1,1,2,2,2) = 759         f(0,0,2,1,2,2,2) = 1297        f(1,0,2,1,2,2,2) = 474
f(2,0,2,1,2,2,2) = 415         f(0,1,2,1,2,2,2) = 268         f(1,1,2,1,2,2,2) = 26
f(2,1,2,1,2,2,2) = 887         f(0,2,2,1,2,2,2) = 82          f(1,2,2,1,2,2,2) = 25
f(2,2,2,1,2,2,2) = 503         f(0,0,0,2,2,2,2) = 309         f(1,0,0,2,2,2,2) = 635
f(2,0,0,2,2,2,2) = 7           f(0,1,0,2,2,2,2) = 795         f(1,1,0,2,2,2,2) = 1851
f(2,1,0,2,2,2,2) = 871         f(0,2,0,2,2,2,2) = 17          f(1,2,0,2,2,2,2) = 1275
f(2,2,0,2,2,2,2) = 487         f(0,0,1,2,2,2,2) = 413         f(1,0,1,2,2,2,2) = 766
f(2,0,1,2,2,2,2) = 1071        f(0,1,1,2,2,2,2) = 539         f(1,1,1,2,2,2,2) = 574
f(2,1,1,2,2,2,2) = 527         f(0,2,1,2,2,2,2) = 1169        f(1,2,1,2,2,2,2) = 189
f(2,2,1,2,2,2,2) = 143         f(0,0,2,2,2,2,2) = 639         f(1,0,2,2,2,2,2) = 893
f(2,0,2,2,2,2,2) = 1055        f(0,1,2,2,2,2,2) = 283         f(1,1,2,2,2,2,2) = 318
f(2,1,2,2,2,2,2) = 271         f(0,2,2,2,2,2,2) = 1105        f(1,2,2,2,2,2,2) = 125
f(2,2,2,2,2,2,2) = 79

List 2

The outcome f(x) is the objective set in the direction x.

We know that there is a 1-1 correlation between the subsets of the integers in $[1,n]$ and the positive integers less than $2^n$ in binary forms: $\{a_i\} \to \sum_i 2^{a_i-1}$. So we substitute subset $\{a_i\}$ with integer $\sum_i 2^{a_i-1}$ for simplicity. For instance, f(2,2,0,2,0,1,2) = 99 indicates that in the direction (2,2,0,2,0,1,2) the fake coins set is {1,2,6,7}.

The completed list of weighings of the second algorithm to sort 11 coins

*The 1-th weighing*
w( ) = {1,2,3}:{4,5,6}

*The 2-th weighing*
w(0) = {1,7,8}:{4,9,10}     w(1) = {1,7,8}:{2,9,10}     w(2) = {1,7,8}:{2,9,10}

*The 3-th weighing*
w(0,0) = {2,7}:{3,8}         w(1,0) = {1,5,7,9}:{2,6,8,10}   w(2,0) = {1,5,7,9}:{2,6,8,10}
w(0,1) = {7,9}:{8,10}        w(1,1) = {5,7,9}:{6,8,10}       w(2,1) = {5,7,9}:{6,8,10}
w(0,2) = {7,9}:{8,10}        w(1,2) = {5,7,9}:{6,8,10}       w(2,2) = {5,7,9}:{6,8,10}

*The 4-th weighing*
w(0,0,0) = {5}:{6}           w(1,0,0) = {1,2,4}:{3,7,10}     w(2,0,0) = {1,2,4}:{3,7,10}
w(0,1,0) = {5,7,10}:{6,8,9}  w(1,1,0) = {2,9,11}:{3,4,5}     w(2,1,0) = {1,7,11}:{3,4,5}
w(0,2,0) = {5,7,10}:{6,8,9}  w(1,2,0) = {1,7,11}:{3,4,5}     w(2,2,0) = {2,9,11}:{3,4,5}
w(0,0,1) = {1,9}:{4,10}      w(1,0,1) = {2,3,6}:{5,8,9}      w(2,0,1) = {1,3,5}:{6,8,9}
w(0,1,1) = {1,5,6}:{2,3,7}   w(1,1,1) = {3,4,7,10}:{1,2,6,11} w(2,1,1) = {3,4,7,10}:{1,2,5,11}
w(0,2,1) = {1,5,6}:{2,3,8}   w(1,2,1) = {3,4,8,9}:{1,2,6,11}  w(2,2,1) = {3,4,8,9}:{1,2,5,11}
w(0,0,2) = {1,9}:{4,10}      w(1,0,2) = {1,3,5}:{6,8,9}       w(2,0,2) = {2,3,6}:{5,8,9}
w(0,1,2) = {1,5,6}:{2,3,8}   w(1,1,2) = {3,4,8,9}:{1,2,5,11}  w(2,1,2) = {3,4,8,9}:{1,2,6,11}
w(0,2,2) = {1,5,6}:{2,3,7}   w(1,2,2) = {3,4,7,10}:{1,2,5,11} w(2,2,2) = {3,4,7,10}:{1,2,6,11}

*The 5-th weighing*
w(0,0,0,0) = {1,11}:{7,8}    w(1,0,0,0) = {1}:{2}            w(2,0,0,0) = {1}:{2}
w(0,1,0,0) = {2}:{3}         w(1,1,0,0) = {2,6,8}:{3,5,11}   w(2,1,0,0) = {1,6,10}:{3,5,11}
w(0,2,0,0) = {2}:{3}         w(1,2,0,0) = {1,6,10}:{3,5,11}  w(2,2,0,0) = {2,6,8}:{3,5,11}
w(0,0,1,0) = {5,10}:{6,11}   w(1,0,1,0) = {2,5,9}:{3,4,7}    w(2,0,1,0) = {1,6,8}:{3,4,10}
w(0,1,1,0) = {2}:{3}         w(1,1,1,0) = {1,2}:{7,11}       w(2,1,1,0) = {1,2}:{10,11}
w(0,2,1,0) = {2}:{3}         w(1,2,1,0) = {1,2}:{9,11}       w(2,2,1,0) = {1,2}:{8,11}
w(0,0,2,0) = {5,10}:{6,11}   w(1,0,2,0) = {1,6,8}:{3,4,10}   w(2,0,2,0) = {2,5,9}:{3,4,7}
w(0,1,2,0) = {2}:{3}         w(1,1,2,0) = {1,2}:{8,11}       w(2,1,2,0) = {1,2}:{9,11}
w(0,2,2,0) = {2}:{3}         w(1,2,2,0) = {1,2}:{10,11}      w(2,2,2,0) = {1,2}:{7,11}
w(0,0,0,1) = {7,9}:{8,10}    w(1,0,0,1) = {3,7,10}:{4,5,6}   w(2,0,0,1) = {1,2}:{8,9}
w(0,1,0,1) = {2,5}:{3,7}     w(1,1,0,1) = {1,5,9}:{2,4,11}   w(2,1,0,1) = {2,4,11}:{3,5,8}
w(0,2,0,1) = {2,6}:{3,8}     w(1,2,0,1) = {2,5,7}:{1,4,11}   w(2,2,0,1) = {1,4,11}:{3,5,10}
w(0,0,1,1) = {3,5}:{6,10}    w(1,0,1,1) = {1,6,11}:{3,7,10}  w(2,0,1,1) = {3,6,7}:{8,10,11}
w(0,1,1,1) = {1,10}:{8,9}    w(1,1,1,1) = {4,9}:{5,8}        w(2,1,1,1) = {2,5,6}:{4,10,11}
w(0,2,1,1) = {2}:{3}         w(1,2,1,1) = {4,7}:{5,10}       w(2,2,1,1) = {1,5,6}:{4,8,11}
w(0,0,2,1) = {1,4}:{2,8}     w(1,0,2,1) = {2,5,11}:{3,7,10}  w(2,0,2,1) = {3,5,10}:{7,9,11}
w(0,1,2,1) = {1,9}:{7,10}    w(1,1,2,1) = {4,10}:{6,7}       w(2,1,2,1) = {2,5,6}:{4,9,11}
w(0,2,2,1) = {2}:{3}         w(1,2,2,1) = {4,8}:{6,9}        w(2,2,2,1) = {1,5,6}:{4,7,11}
w(0,0,0,2) = {7,9}:{8,10}    w(1,0,0,2) = {1,2}:{8,9}        w(2,0,0,2) = {3,7,10}:{4,5,6}
w(0,1,0,2) = {2,6}:{3,8}     w(1,1,0,2) = {1,4,11}:{3,5,10}  w(2,1,0,2) = {2,5,7}:{1,4,11}
w(0,2,0,2) = {2,5}:{3,7}     w(1,2,0,2) = {2,4,11}:{3,5,8}   w(2,2,0,2) = {1,5,9}:{2,4,11}
w(0,0,1,2) = {1,4}:{2,8}     w(1,0,1,2) = {3,5,10}:{7,9,11}  w(2,0,1,2) = {2,5,11}:{3,7,10}

| | | |
|---|---|---|
| w(0,1,1,2) = {2}:{3} | w(1,1,1,2) = {1,5,6}:{4,7,11} | w(2,1,1,2) = {4,8}:{6,9} |
| w(0,2,1,2) = {1,9}:{7,10} | w(1,2,1,2) = {2,5,6}:{4,9,11} | w(2,2,1,2) = {4,10}:{6,7} |
| w(0,0,2,2) = {3,5}:{6,10} | w(1,0,2,2) = {3,6,7}:{8,10,11} | w(2,0,2,2) = {1,6,11}:{3,7,10} |
| w(0,1,2,2) = {2}:{3} | w(1,1,2,2) = {1,5,6}:{4,8,11} | w(2,1,2,2) = {4,7}:{5,10} |
| w(0,2,2,2) = {1,10}:{8,9} | w(1,2,2,2) = {2,5,6}:{4,10,11} | w(2,2,2,2) = {4,9}:{5,8} |

*The 6-th weighing*

| | | |
|---|---|---|
| w(0,0,0,0,0) = {6}:{11} | w(1,0,0,0,0) = {3,4}:{9,11} | w(2,0,0,0,0) = {3,4}:{9,11} |
| w(0,1,0,0,0) = {1,11}:{2,3} | w(1,1,0,0,0) = {4}:{9} | w(2,1,0,0,0) = {4}:{7} |
| w(0,2,0,0,0) = {1,11}:{2,3} | w(1,2,0,0,0) = {4}:{7} | w(2,2,0,0,0) = {4}:{9} |
| w(0,0,1,0,0) = {5}:{8} | w(1,0,1,0,0) = {2,9}:{3,11} | w(2,0,1,0,0) = {1,8}:{3,11} |
| w(0,1,1,0,0) = {5,6}:{9,11} | w(1,1,1,0,0) = {2,8}:{7,10} | w(2,1,1,0,0) = {1,9}:{7,10} |
| w(0,2,1,0,0) = {5,6}:{10,11} | w(1,2,1,0,0) = {1,10}:{8,9} | w(2,2,1,0,0) = {2,7}:{8,9} |
| w(0,0,2,0,0) = {5}:{8} | w(1,0,2,0,0) = {1,8}:{3,11} | w(2,0,2,0,0) = {2,9}:{3,11} |
| w(0,1,2,0,0) = {5,6}:{10,11} | w(1,1,2,0,0) = {2,7}:{8,9} | w(2,1,2,0,0) = {1,10}:{8,9} |
| w(0,2,2,0,0) = {5,6}:{9,11} | w(1,2,2,0,0) = {1,9}:{7,10} | w(2,2,2,0,0) = {2,8}:{7,10} |
| w(0,0,0,1,0) = {7,10}:{4,11} | w(1,0,0,1,0) = {3,4}:{9,11} | w(2,0,0,1,0) = {1,2}:{10,11} |
| w(0,1,0,1,0) = {1,10}:{8,11} | w(1,1,0,1,0) = {1,4}:{2,6} | w(2,1,0,1,0) = {7,9}:{10,11} |
| w(0,2,0,1,0) = {1,9}:{7,11} | w(1,2,0,1,0) = {2,4}:{1,6} | w(2,2,0,1,0) = {7,9}:{8,11} |
| w(0,0,1,1,0) = {5,10}:{6,11} | w(1,0,1,1,0) = {5}:{6} | w(2,0,1,1,0) = {1,4}:{7,11} |
| w(0,1,1,1,0) = {1,3,11}:{2,5,6} | w(1,1,1,1,0) = {1,6}:{8,11} | w(2,1,1,1,0) = {8,9}:{5,11} |
| w(0,2,1,1,0) = {5,6}:{10,11} | w(1,2,1,1,0) = {2,6}:{10,11} | w(2,2,1,1,0) = {7,10}:{5,11} |
| w(0,0,2,1,0) = {2,11}:{7,10} | w(1,0,2,1,0) = {6}:{5} | w(2,0,2,1,0) = {2,4}:{10,11} |
| w(0,1,2,1,0) = {1,3,11}:{2,5,6} | w(1,1,2,1,0) = {1,5}:{7,11} | w(2,1,2,1,0) = {7,10}:{6,11} |
| w(0,2,2,1,0) = {5,6}:{9,11} | w(1,2,2,1,0) = {2,5}:{9,11} | w(2,2,2,1,0) = {8,9}:{6,11} |
| w(0,0,0,2,0) = {7,10}:{4,11} | w(1,0,0,2,0) = {1,2}:{10,11} | w(2,0,0,2,0) = {3,4}:{9,11} |
| w(0,1,0,2,0) = {1,9}:{7,11} | w(1,1,0,2,0) = {7,9}:{8,11} | w(2,1,0,2,0) = {2,4}:{1,6} |
| w(0,2,0,2,0) = {1,10}:{8,11} | w(1,2,0,2,0) = {7,9}:{10,11} | w(2,2,0,2,0) = {1,4}:{2,6} |
| w(0,0,1,2,0) = {2,11}:{7,10} | w(1,0,1,2,0) = {2,4}:{10,11} | w(2,0,1,2,0) = {6}:{5} |
| w(0,1,1,2,0) = {5,6}:{9,11} | w(1,1,1,2,0) = {8,9}:{6,11} | w(2,1,1,2,0) = {2,5}:{9,11} |
| w(0,2,1,2,0) = {1,3,11}:{2,5,6} | w(1,2,1,2,0) = {7,10}:{6,11} | w(2,2,1,2,0) = {1,5}:{7,11} |
| w(0,0,2,2,0) = {5,10}:{6,11} | w(1,0,2,2,0) = {1,4}:{7,11} | w(2,0,2,2,0) = {5}:{6} |
| w(0,1,2,2,0) = {5,6}:{10,11} | w(1,1,2,2,0) = {7,10}:{5,11} | w(2,1,2,2,0) = {2,6}:{10,11} |
| w(0,2,2,2,0) = {1,3,11}:{2,5,6} | w(1,2,2,2,0) = {8,9}:{5,11} | w(2,2,2,2,0) = {1,6}:{8,11} |
| w(0,0,0,0,1) = {4}:{5} | w(1,0,0,0,1) = {3,4}:{9,11} | w(2,0,0,0,1) = {3,4}:{9,11} |
| w(0,1,0,0,1) = {3,11}:{9,10} | w(1,1,0,0,1) = {1,6}:{5,8} | w(2,1,0,0,1) = {2,6}:{7,10} |
| w(0,2,0,0,1) = {2,11}:{9,10} | w(1,2,0,0,1) = {2,6}:{5,10} | w(2,2,0,0,1) = {1,6}:{8,9} |
| w(0,0,1,0,1) = {1,2}:{7,8} | w(1,0,1,0,1) = {5,11}:{7,9} | w(2,0,1,0,1) = {2,4}:{10,11} |
| w(0,1,1,0,1) = {5,9}:{6,11} | w(1,1,1,0,1) = {2,4}:{5,9} | w(2,1,1,0,1) = {6}:{8} |
| w(0,2,1,0,1) = {5,10}:{6,11} | w(1,2,1,0,1) = {1,4}:{5,7} | w(2,2,1,0,1) = {6}:{10} |
| w(0,0,2,0,1) = {1,4}:{7,11} | w(1,0,2,0,1) = {6,11}:{8,10} | w(2,0,2,0,1) = {1,4}:{7,11} |
| w(0,1,2,0,1) = {5,10}:{6,11} | w(1,1,2,0,1) = {2,4}:{6,10} | w(2,1,2,0,1) = {5}:{7} |
| w(0,2,2,0,1) = {5,9}:{6,11} | w(1,2,2,0,1) = {1,4}:{6,8} | w(2,2,2,0,1) = {5}:{9} |
| w(0,0,0,1,1) = {2,8}:{10,11} | w(1,0,0,1,1) = {3,8}:{4,11} | w(2,0,0,1,1) = {1,10}:{2,11} |
| w(0,1,0,1,1) = {1,11}:{2,8} | w(1,1,0,1,1) = {3,6}:{7,10} | w(2,1,0,1,1) = {5,6}:{9,11} |

w(0,2,0,1,1) = {1,5}:{6,11}     w(1,2,0,1,1) = {3,6}:{8,9}      w(2,2,0,1,1) = {5,6}:{7,11}
w(0,0,1,1,1) = {3,6}:{4,11}     w(1,0,1,1,1) = {1,5}:{7,11}     w(2,0,1,1,1) = {1,9}:{7,11}
w(0,1,1,1,1) = {1,3,11}:{2,5,6} w(1,1,1,1,1) = {1,2}:{4,11}     w(2,1,1,1,1) = {1,4}:{3,8}
w(0,2,1,1,1) = {5,10}:{6,11}    w(1,2,1,1,1) = {1,2}:{4,11}     w(2,2,1,1,1) = {2,4}:{3,10}
w(0,0,2,1,1) = {4,8}:{9,11}     w(1,0,2,1,1) = {2,6}:{10,11}    w(2,0,2,1,1) = {2,8}:{10,11}
w(0,1,2,1,1) = {1,3,11}:{2,5,6} w(1,1,2,1,1) = {1,2}:{4,11}     w(2,1,2,1,1) = {1,4}:{3,7}
w(0,2,2,1,1) = {5,9}:{6,11}     w(1,2,2,1,1) = {1,2}:{4,11}     w(2,2,2,1,1) = {2,4}:{3,9}
w(0,0,0,2,1) = {3,7}:{9,11}     w(1,0,0,2,1) = {1,5}:{7,11}     w(2,0,0,2,1) = {3,5}:{9,11}
w(0,1,0,2,1) = {1,11}:{2,7}     w(1,1,0,2,1) = {2,9}:{10,11}    w(2,1,0,2,1) = {2,4}:{7,10}
w(0,2,0,2,1) = {1,6}:{5,11}     w(1,2,0,2,1) = {1,7}:{8,11}     w(2,2,0,2,1) = {1,4}:{8,9}
w(0,0,1,2,1) = {3,5,6}:{2,8,11} w(1,0,1,2,1) = {2}:{10}         w(2,0,1,2,1) = {6,7}:{1,11}
w(0,1,1,2,1) = {5,9}:{6,11}     w(1,1,1,2,1) = {8,9}:{3,11}     w(2,1,1,2,1) = {1,8}:{7,11}
w(0,2,1,2,1) = {1,3,11}:{2,5,6} w(1,2,1,2,1) = {7,10}:{3,11}    w(2,2,1,2,1) = {2,10}:{9,11}
w(0,0,2,2,1) = {5,10}:{4,11}    w(1,0,2,2,1) = {1}:{7}          w(2,0,2,2,1) = {5,10}:{2,11}
w(0,1,2,2,1) = {5,10}:{6,11}    w(1,1,2,2,1) = {7,10}:{3,11}    w(2,1,2,2,1) = {1,7}:{8,11}
w(0,2,2,2,1) = {1,3,11}:{2,5,6} w(1,2,2,2,1) = {8,9}:{3,11}     w(2,2,2,2,1) = {2,9}:{10,11}
w(0,0,0,0,2) = {4}:{5}          w(1,0,0,0,2) = {3,4}:{9,11}     w(2,0,0,0,2) = {3,4}:{9,11}
w(0,1,0,0,2) = {2,11}:{9,10}    w(1,1,0,0,2) = {1,6}:{8,9}      w(2,1,0,0,2) = {2,6}:{5,10}
w(0,2,0,0,2) = {3,11}:{9,10}    w(1,2,0,0,2) = {2,6}:{7,10}     w(2,2,0,0,2) = {1,6}:{5,8}
w(0,0,1,0,2) = {1,4}:{7,11}     w(1,0,1,0,2) = {1,4}:{7,11}     w(2,0,1,0,2) = {6,11}:{8,10}
w(0,1,1,0,2) = {5,9}:{6,11}     w(1,1,1,0,2) = {5}:{9}          w(2,1,1,0,2) = {1,4}:{6,8}
w(0,2,1,0,2) = {5,10}:{6,11}    w(1,2,1,0,2) = {5}:{7}          w(2,2,1,0,2) = {2,4}:{6,10}
w(0,0,2,0,2) = {1,2}:{7,8}      w(1,0,2,0,2) = {2,4}:{10,11}    w(2,0,2,0,2) = {5,11}:{7,9}
w(0,1,2,0,2) = {5,10}:{6,11}    w(1,1,2,0,2) = {6}:{10}         w(2,1,2,0,2) = {1,4}:{5,7}
w(0,2,2,0,2) = {5,9}:{6,11}     w(1,2,2,0,2) = {6}:{8}          w(2,2,2,0,2) = {2,4}:{5,9}
w(0,0,0,1,2) = {3,7}:{9,11}     w(1,0,0,1,2) = {3,5}:{9,11}     w(2,0,0,1,2) = {1,5}:{7,11}
w(0,1,0,1,2) = {1,6}:{5,11}     w(1,1,0,1,2) = {1,4}:{8,9}      w(2,1,0,1,2) = {1,7}:{8,11}
w(0,2,0,1,2) = {1,11}:{2,7}     w(1,2,0,1,2) = {2,4}:{7,10}     w(2,2,0,1,2) = {2,9}:{10,11}
w(0,0,1,1,2) = {5,10}:{4,11}    w(1,0,1,1,2) = {5,10}:{2,11}    w(2,0,1,1,2) = {1}:{7}
w(0,1,1,1,2) = {1,3,11}:{2,5,6} w(1,1,1,1,2) = {2,9}:{10,11}    w(2,1,1,1,2) = {8,9}:{3,11}
w(0,2,1,1,2) = {5,10}:{6,11}    w(1,2,1,1,2) = {1,7}:{8,11}     w(2,2,1,1,2) = {7,10}:{3,11}
w(0,0,2,1,2) = {3,5,6}:{2,8,11} w(1,0,2,1,2) = {6,7}:{1,11}     w(2,0,2,1,2) = {2}:{10}
w(0,1,2,1,2) = {1,3,11}:{2,5,6} w(1,1,2,1,2) = {2,10}:{9,11}    w(2,1,2,1,2) = {7,10}:{3,11}
w(0,2,2,1,2) = {5,9}:{6,11}     w(1,2,2,1,2) = {1,8}:{7,11}     w(2,2,2,1,2) = {8,9}:{3,11}
w(0,0,0,2,2) = {2,8}:{10,11}    w(1,0,0,2,2) = {1,10}:{2,11}    w(2,0,0,2,2) = {3,8}:{4,11}
w(0,1,0,2,2) = {1,5}:{6,11}     w(1,1,0,2,2) = {5,6}:{7,11}     w(2,1,0,2,2) = {3,6}:{8,9}
w(0,2,0,2,2) = {1,11}:{2,8}     w(1,2,0,2,2) = {5,6}:{9,11}     w(2,2,0,2,2) = {3,6}:{7,10}
w(0,0,1,2,2) = {4,8}:{9,11}     w(1,0,1,2,2) = {2,8}:{10,11}    w(2,0,1,2,2) = {2,6}:{10,11}
w(0,1,1,2,2) = {5,9}:{6,11}     w(1,1,1,2,2) = {2,4}:{3,9}      w(2,1,1,2,2) = {1,2}:{4,11}
w(0,2,1,2,2) = {1,3,11}:{2,5,6} w(1,2,1,2,2) = {1,4}:{3,7}      w(2,2,1,2,2) = {1,2}:{4,11}
w(0,0,2,2,2) = {3,6}:{4,11}     w(1,0,2,2,2) = {1,9}:{7,11}     w(2,0,2,2,2) = {1,5}:{7,11}
w(0,1,2,2,2) = {5,10}:{6,11}    w(1,1,2,2,2) = {2,4}:{3,10}     w(2,1,2,2,2) = {1,2}:{4,11}
w(0,2,2,2,2) = {1,3,11}:{2,5,6} w(1,2,2,2,2) = {1,4}:{3,8}      w(2,2,2,2,2) = {1,2}:{4,11}

*The 7-th weighing*

| | | |
|---|---|---|
| w(0,0,0,0,0,0) = { }:{ } | w(1,0,0,0,0,0) = {6}:{8} | w(2,0,0,0,0,0) = {6}:{8} |
| w(0,1,0,0,0,0) = {6}:{9} | w(1,1,0,0,0,0) = {10}:{11} | w(2,1,0,0,0,0) = {8}:{11} |
| w(0,2,0,0,0,0) = {6}:{9} | w(1,2,0,0,0,0) = {8}:{11} | w(2,2,0,0,0,0) = {10}:{11} |
| w(0,0,1,0,0,0) = { }:{ } | w(1,0,1,0,0,0) = {10}:{11} | w(2,0,1,0,0,0) = {7}:{11} |
| w(0,1,1,0,0,0) = {9}:{10} | w(1,1,1,0,0,0) = {5}:{9} | w(2,1,1,0,0,0) = {6}:{8} |
| w(0,2,1,0,0,0) = {10}:{9} | w(1,2,1,0,0,0) = {5}:{7} | w(2,2,1,0,0,0) = {6}:{10} |
| w(0,0,2,0,0,0) = { }:{ } | w(1,0,2,0,0,0) = {7}:{11} | w(2,0,2,0,0,0) = {10}:{11} |
| w(0,1,2,0,0,0) = {10}:{9} | w(1,1,2,0,0,0) = {6}:{10} | w(2,1,2,0,0,0) = {5}:{7} |
| w(0,2,2,0,0,0) = {9}:{10} | w(1,2,2,0,0,0) = {6}:{8} | w(2,2,2,0,0,0) = {5}:{9} |
| w(0,0,0,1,0,0) = {1}:{2} | w(1,0,0,1,0,0) = {3}:{5} | w(2,0,0,1,0,0) = {3,4}:{6,7} |
| w(0,1,0,1,0,0) = {6}:{9} | w(1,1,0,1,0,0) = {5}:{7} | w(2,1,0,1,0,0) = {2}:{9} |
| w(0,2,0,1,0,0) = {5}:{10} | w(1,2,0,1,0,0) = {5}:{9} | w(2,2,0,1,0,0) = {1}:{7} |
| w(0,0,1,1,0,0) = {2}:{4} | w(1,0,1,1,0,0) = {4}:{6} | w(2,0,1,1,0,0) = {7}:{11} |
| w(0,1,1,1,0,0) = {5}:{6} | w(1,1,1,1,0,0) = {1,2}:{6,7} | w(2,1,1,1,0,0) = {1,4}:{3,6} |
| w(0,2,1,1,0,0) = {10}:{9} | w(1,2,1,1,0,0) = {1,2}:{6,9} | w(2,2,1,1,0,0) = {2,4}:{3,6} |
| w(0,0,2,1,0,0) = {5}:{6} | w(1,0,2,1,0,0) = {4}:{5} | w(2,0,2,1,0,0) = {10}:{11} |
| w(0,1,2,1,0,0) = {5}:{6} | w(1,1,2,1,0,0) = {1,2}:{5,8} | w(2,1,2,1,0,0) = {1,4}:{3,5} |
| w(0,2,2,1,0,0) = {9}:{10} | w(1,2,2,1,0,0) = {1,2}:{5,10} | w(2,2,2,1,0,0) = {2,4}:{3,5} |
| w(0,0,0,2,0,0) = {1}:{2} | w(1,0,0,2,0,0) = {3,4}:{6,7} | w(2,0,0,2,0,0) = {3}:{5} |
| w(0,1,0,2,0,0) = {5}:{10} | w(1,1,0,2,0,0) = {1}:{7} | w(2,1,0,2,0,0) = {5}:{9} |
| w(0,2,0,2,0,0) = {6}:{9} | w(1,2,0,2,0,0) = {2}:{9} | w(2,2,0,2,0,0) = {5}:{7} |
| w(0,0,1,2,0,0) = {5}:{6} | w(1,0,1,2,0,0) = {10}:{11} | w(2,0,1,2,0,0) = {4}:{5} |
| w(0,1,1,2,0,0) = {9}:{10} | w(1,1,1,2,0,0) = {2,4}:{3,5} | w(2,1,1,2,0,0) = {1,2}:{5,10} |
| w(0,2,1,2,0,0) = {5}:{6} | w(1,2,1,2,0,0) = {1,4}:{3,5} | w(2,2,1,2,0,0) = {1,2}:{5,8} |
| w(0,0,2,2,0,0) = {2}:{4} | w(1,0,2,2,0,0) = {7}:{11} | w(2,0,2,2,0,0) = {4}:{6} |
| w(0,1,2,2,0,0) = {10}:{9} | w(1,1,2,2,0,0) = {2,4}:{3,6} | w(2,1,2,2,0,0) = {1,2}:{6,9} |
| w(0,2,2,2,0,0) = {5}:{6} | w(1,2,2,2,0,0) = {1,4}:{3,6} | w(2,2,2,2,0,0) = {1,2}:{6,7} |
| w(0,0,0,0,1,0) = {1}:{11} | w(1,0,0,0,1,0) = {6}:{8} | w(2,0,0,0,1,0) = {10}:{11} |
| w(0,1,0,0,1,0) = {1}:{7} | w(1,1,0,0,1,0) = {3}:{4} | w(2,1,0,0,1,0) = {1}:{9} |
| w(0,2,0,0,1,0) = {1}:{7} | w(1,2,0,0,1,0) = {3}:{4} | w(2,2,0,0,1,0) = {2}:{7} |
| w(0,0,1,0,1,0) = {6}:{11} | w(1,0,1,0,1,0) = {10}:{11} | w(2,0,1,0,1,0) = {1}:{4} |
| w(0,1,1,0,1,0) = {9}:{11} | w(1,1,1,0,1,0) = {5}:{8} | w(2,1,1,0,1,0) = {5,10}:{8,9} |
| w(0,2,1,0,1,0) = {10}:{11} | w(1,2,1,0,1,0) = {5}:{10} | w(2,2,1,0,1,0) = {5,8}:{7,10} |
| w(0,0,2,0,1,0) = {1}:{3} | w(1,0,2,0,1,0) = {7}:{11} | w(2,0,2,0,1,0) = {2}:{4} |
| w(0,1,2,0,1,0) = {10}:{11} | w(1,1,2,0,1,0) = {6}:{7} | w(2,1,2,0,1,0) = {6,9}:{7,10} |
| w(0,2,2,0,1,0) = {9}:{11} | w(1,2,2,0,1,0) = {6}:{9} | w(2,2,2,0,1,0) = {6,7}:{8,9} |
| w(0,0,0,1,1,0) = {1}:{3} | w(1,0,0,1,1,0) = {3}:{4} | w(2,0,0,1,1,0) = {2}:{11} |
| w(0,1,0,1,1,0) = {6}:{8} | w(1,1,0,1,1,0) = {2}:{11} | w(2,1,0,1,1,0) = {1,4}:{5,6} |
| w(0,2,0,1,1,0) = {5}:{7} | w(1,2,0,1,1,0) = {1}:{11} | w(2,2,0,1,1,0) = {2,4}:{5,6} |
| w(0,0,1,1,1,0) = {1}:{8} | w(1,0,1,1,1,0) = {2}:{3} | w(2,0,1,1,1,0) = {2,3}:{9,10} |
| w(0,1,1,1,1,0) = {5}:{6} | w(1,1,1,1,1,0) = {6,7}:{8,9} | w(2,1,1,1,1,0) = {3}:{8} |
| w(0,2,1,1,1,0) = {10}:{11} | w(1,2,1,1,1,0) = {6,9}:{7,10} | w(2,2,1,1,1,0) = {3}:{10} |
| w(0,0,2,1,1,0) = {5}:{6} | w(1,0,2,1,1,0) = {1}:{3} | w(2,0,2,1,1,0) = {1,3}:{7,8} |
| w(0,1,2,1,1,0) = {5}:{6} | w(1,1,2,1,1,0) = {5,8}:{7,10} | w(2,1,2,1,1,0) = {3}:{7} |

w(0,2,2,1,1,0) = {9}:{11}          w(1,2,2,1,1,0) = {5,10}:{8,9}      w(2,2,2,1,1,0) = {3}:{9}
w(0,0,0,2,1,0) = {1}:{2}           w(1,0,0,2,1,0) = {2}:{11}          w(2,0,0,2,1,0) = {3}:{4}
w(0,1,0,2,1,0) = {5}:{7}           w(1,1,0,2,1,0) = {2}:{7}           w(2,1,0,2,1,0) = {1,5}:{3,4}
w(0,2,0,2,1,0) = {6}:{8}           w(1,2,0,2,1,0) = {1}:{9}           w(2,2,0,2,1,0) = {2,5}:{3,4}
w(0,0,1,2,1,0) = {5}:{6}           w(1,0,1,2,1,0) = {2,4}:{5,6}       w(2,0,1,2,1,0) = {2,4}:{5,6}
w(0,1,1,2,1,0) = {9}:{11}          w(1,1,1,2,1,0) = {2}:{3}           w(2,1,1,2,1,0) = {5,6}:{8,9}
w(0,2,1,2,1,0) = {5}:{6}           w(1,2,1,2,1,0) = {1}:{3}           w(2,2,1,2,1,0) = {5,6}:{7,10}
w(0,0,2,2,1,0) = {1}:{2}           w(1,0,2,2,1,0) = {1,4}:{5,6}       w(2,0,2,2,1,0) = {1,4}:{5,6}
w(0,1,2,2,1,0) = {10}:{11}         w(1,1,2,2,1,0) = {2}:{3}           w(2,1,2,2,1,0) = {5,6}:{7,10}
w(0,2,2,2,1,0) = {5}:{6}           w(1,2,2,2,1,0) = {1}:{3}           w(2,2,2,2,1,0) = {5,6}:{8,9}
w(0,0,0,0,2,0) = {1}:{11}          w(1,0,0,0,2,0) = {10}:{11}         w(2,0,0,0,2,0) = {6}:{8}
w(0,1,0,0,2,0) = {1}:{7}           w(1,1,0,0,2,0) = {2}:{7}           w(2,1,0,0,2,0) = {3}:{4}
w(0,2,0,0,2,0) = {1}:{7}           w(1,2,0,0,2,0) = {1}:{9}           w(2,2,0,0,2,0) = {3}:{4}
w(0,0,1,0,2,0) = {1}:{3}           w(1,0,1,0,2,0) = {2}:{4}           w(2,0,1,0,2,0) = {7}:{11}
w(0,1,1,0,2,0) = {9}:{11}          w(1,1,1,0,2,0) = {6,7}:{8,9}       w(2,1,1,0,2,0) = {6}:{9}
w(0,2,1,0,2,0) = {10}:{11}         w(1,2,1,0,2,0) = {6,9}:{7,10}      w(2,2,1,0,2,0) = {6}:{7}
w(0,0,2,0,2,0) = {6}:{11}          w(1,0,2,0,2,0) = {1}:{4}           w(2,0,2,0,2,0) = {10}:{11}
w(0,1,2,0,2,0) = {10}:{11}         w(1,1,2,0,2,0) = {5,8}:{7,10}      w(2,1,2,0,2,0) = {5}:{10}
w(0,2,2,0,2,0) = {9}:{11}          w(1,2,2,0,2,0) = {5,10}:{8,9}      w(2,2,2,0,2,0) = {5}:{8}
w(0,0,0,1,2,0) = {1}:{2}           w(1,0,0,1,2,0) = {3}:{4}           w(2,0,0,1,2,0) = {2}:{11}
w(0,1,0,1,2,0) = {6}:{8}           w(1,1,0,1,2,0) = {2,5}:{3,4}       w(2,1,0,1,2,0) = {1}:{9}
w(0,2,0,1,2,0) = {5}:{7}           w(1,2,0,1,2,0) = {1,5}:{3,4}       w(2,2,0,1,2,0) = {2}:{7}
w(0,0,1,1,2,0) = {1}:{2}           w(1,0,1,1,2,0) = {1,4}:{5,6}       w(2,0,1,1,2,0) = {1,4}:{5,6}
w(0,1,1,1,2,0) = {5}:{6}           w(1,1,1,1,2,0) = {5,6}:{8,9}       w(2,1,1,1,2,0) = {1}:{3}
w(0,2,1,1,2,0) = {10}:{11}         w(1,2,1,1,2,0) = {5,6}:{7,10}      w(2,2,1,1,2,0) = {2}:{3}
w(0,0,2,1,2,0) = {5}:{6}           w(1,0,2,1,2,0) = {2,4}:{5,6}       w(2,0,2,1,2,0) = {2,4}:{5,6}
w(0,1,2,1,2,0) = {5}:{6}           w(1,1,2,1,2,0) = {5,6}:{7,10}      w(2,1,2,1,2,0) = {1}:{3}
w(0,2,2,1,2,0) = {9}:{11}          w(1,2,2,1,2,0) = {5,6}:{8,9}       w(2,2,2,1,2,0) = {2}:{3}
w(0,0,0,2,2,0) = {1}:{3}           w(1,0,0,2,2,0) = {2}:{11}          w(2,0,0,2,2,0) = {3}:{4}
w(0,1,0,2,2,0) = {5}:{7}           w(1,1,0,2,2,0) = {2,4}:{5,6}       w(2,1,0,2,2,0) = {1}:{11}
w(0,2,0,2,2,0) = {6}:{8}           w(1,2,0,2,2,0) = {1,4}:{5,6}       w(2,2,0,2,2,0) = {2}:{11}
w(0,0,1,2,2,0) = {5}:{6}           w(1,0,1,2,2,0) = {1,3}:{7,8}       w(2,0,1,2,2,0) = {1}:{3}
w(0,1,1,2,2,0) = {9}:{11}          w(1,1,1,2,2,0) = {3}:{9}           w(2,1,1,2,2,0) = {5,10}:{8,9}
w(0,2,1,2,2,0) = {5}:{6}           w(1,2,1,2,2,0) = {3}:{7}           w(2,2,1,2,2,0) = {5,8}:{7,10}
w(0,0,2,2,2,0) = {1}:{8}           w(1,0,2,2,2,0) = {2,3}:{9,10}      w(2,0,2,2,2,0) = {2}:{3}
w(0,1,2,2,2,0) = {10}:{11}         w(1,1,2,2,2,0) = {3}:{10}          w(2,1,2,2,2,0) = {6,9}:{7,10}
w(0,2,2,2,2,0) = {5}:{6}           w(1,2,2,2,2,0) = {3}:{8}           w(2,2,2,2,2,0) = {6,7}:{8,9}
w(0,0,0,0,0,1) = {2}:{3}           w(1,0,0,0,0,1) = {8}:{11}          w(2,0,0,0,0,1) = {8}:{11}
w(0,1,0,0,0,1) = {1}:{11}          w(1,1,0,0,0,1) = {1}:{8}           w(2,1,0,0,0,1) = {1}:{9}
w(0,2,0,0,0,1) = {1}:{11}          w(1,2,0,0,0,1) = {2}:{10}          w(2,2,0,0,0,1) = {2}:{7}
w(0,0,1,0,0,1) = {4}:{6}           w(1,0,1,0,0,1) = {3}:{11}          w(2,0,1,0,0,1) = {7}:{11}
w(0,1,1,0,0,1) = {9}:{11}          w(1,1,1,0,0,1) = {5}:{9}           w(2,1,1,0,0,1) = {3,5}:{6,8}
w(0,2,1,0,0,1) = {10}:{11}         w(1,2,1,0,0,1) = {5}:{7}           w(2,2,1,0,0,1) = {3,5}:{6,10}
w(0,0,2,0,0,1) = {4}:{6}           w(1,0,2,0,0,1) = {3}:{11}          w(2,0,2,0,0,1) = {10}:{11}

| | | |
|---|---|---|
| w(0,1,2,0,0,1) = {10}:{11} | w(1,1,2,0,0,1) = {6}:{10} | w(2,1,2,0,0,1) = {3,6}:{5,7} |
| w(0,2,2,0,0,1) = {9}:{11} | w(1,2,2,0,0,1) = {6}:{8} | w(2,2,2,0,0,1) = {3,6}:{5,9} |
| w(0,0,0,1,0,1) = {1}:{11} | w(1,0,0,1,0,1) = {9}:{11} | w(2,0,0,1,0,1) = {10}:{11} |
| w(0,1,0,1,0,1) = {8}:{11} | w(1,1,0,1,0,1) = {8}:{10} | w(2,1,0,1,0,1) = {1}:{11} |
| w(0,2,0,1,0,1) = {2}:{11} | w(1,2,0,1,0,1) = {10}:{8} | w(2,2,0,1,0,1) = {2}:{11} |
| w(0,0,1,1,0,1) = {6}:{11} | w(1,0,1,1,0,1) = {4}:{6} | w(2,0,1,1,0,1) = {3,4}:{9,10} |
| w(0,1,1,1,0,1) = {5}:{6} | w(1,1,1,1,0,1) = {2}:{9} | w(2,1,1,1,0,1) = {1}:{3} |
| w(0,2,1,1,0,1) = {10}:{11} | w(1,2,1,1,0,1) = {1}:{7} | w(2,2,1,1,0,1) = {2}:{3} |
| w(0,0,2,1,0,1) = {9}:{11} | w(1,0,2,1,0,1) = {4}:{5} | w(2,0,2,1,0,1) = {3,4}:{7,8} |
| w(0,1,2,1,0,1) = {5}:{6} | w(1,1,2,1,0,1) = {2}:{10} | w(2,1,2,1,0,1) = {1}:{3} |
| w(0,2,2,1,0,1) = {9}:{11} | w(1,2,2,1,0,1) = {1}:{8} | w(2,2,2,1,0,1) = {2}:{3} |
| w(0,0,0,2,0,1) = {1}:{11} | w(1,0,0,2,0,1) = {4}:{5} | w(2,0,0,2,0,1) = {9}:{11} |
| w(0,1,0,2,0,1) = {7}:{11} | w(1,1,0,2,0,1) = {6}:{8} | w(2,1,0,2,0,1) = {1,4}:{3,10} |
| w(0,2,0,2,0,1) = {2}:{11} | w(1,2,0,2,0,1) = {6}:{10} | w(2,2,0,2,0,1) = {2,4}:{3,8} |
| w(0,0,1,2,0,1) = {5}:{6} | w(1,0,1,2,0,1) = {1}:{3} | w(2,0,1,2,0,1) = {4}:{11} |
| w(0,1,1,2,0,1) = {9}:{11} | w(1,1,1,2,0,1) = {2,4}:{3,5} | w(2,1,1,2,0,1) = {9}:{11} |
| w(0,2,1,2,0,1) = {9}:{11} | w(1,2,1,2,0,1) = {1,4}:{3,5} | w(2,2,1,2,0,1) = {7}:{11} |
| w(0,0,2,2,0,1) = {2}:{4} | w(1,0,2,2,0,1) = {2}:{3} | w(2,0,2,2,0,1) = {4}:{11} |
| w(0,1,2,2,0,1) = {10}:{11} | w(1,1,2,2,0,1) = {2,4}:{3,6} | w(2,1,2,2,0,1) = {10}:{11} |
| w(0,2,2,2,0,1) = {10}:{11} | w(1,2,2,2,0,1) = {1,4}:{3,6} | w(2,2,2,2,0,1) = {8}:{11} |
| w(0,0,0,0,1,1) = {2}:{11} | w(1,0,0,0,1,1) = {8}:{11} | w(2,0,0,0,1,1) = {8}:{11} |
| w(0,1,0,0,1,1) = {1}:{7} | w(1,1,0,0,1,1) = {2}:{8} | w(2,1,0,0,1,1) = {2}:{9} |
| w(0,2,0,0,1,1) = {2}:{11} | w(1,2,0,0,1,1) = {1}:{10} | w(2,2,0,0,1,1) = {1}:{7} |
| w(0,0,1,0,1,1) = {4}:{6} | w(1,0,1,0,1,1) = {5,9}:{7,11} | w(2,0,1,0,1,1) = {2,11}:{4,10} |
| w(0,1,1,0,1,1) = {2,3}:{9,11} | w(1,1,1,0,1,1) = {3}:{8} | w(2,1,1,0,1,1) = {9}:{7} |
| w(0,2,1,0,1,1) = {2,3}:{10,11} | w(1,2,1,0,1,1) = {3}:{10} | w(2,2,1,0,1,1) = {7}:{9} |
| w(0,0,2,0,1,1) = {6}:{7} | w(1,0,2,0,1,1) = {6,8}:{10,1 | w(2,0,2,0,1,1) = {1,11}:{4,7} |
| w(0,1,2,0,1,1) = {2,3}:{10,11} | w(1,1,2,0,1,1) = {3}:{7} | w(2,1,2,0,1,1) = {10}:{8} |
| w(0,2,2,0,1,1) = {2,3}:{9,11} | w(1,2,2,0,1,1) = {3}:{9} | w(2,2,2,0,1,1) = {8}:{10} |
| w(0,0,0,1,1,1) = {10}:{11} | w(1,0,0,1,1,1) = {3,4}:{9,11} | w(2,0,0,1,1,1) = {2}:{11} |
| w(0,1,0,1,1,1) = {4,5}:{10,11} | w(1,1,0,1,1,1) = {6}:{11} | w(2,1,0,1,1,1) = {2}:{10} |
| w(0,2,0,1,1,1) = {1,4}:{9,11} | w(1,2,0,1,1,1) = {6}:{11} | w(2,2,0,1,1,1) = {1}:{8} |
| w(0,0,1,1,1,1) = {1}:{11} | w(1,0,1,1,1,1) = {4}:{6} | w(2,0,1,1,1,1) = {7}:{11} |
| w(0,1,1,1,1,1) = {5}:{6} | w(1,1,1,1,1,1) = {2}:{9} | w(2,1,1,1,1,1) = {1,3}:{4,9} |
| w(0,2,1,1,1,1) = {2,3}:{10,11} | w(1,2,1,1,1,1) = {1}:{7} | w(2,2,1,1,1,1) = {2,3}:{4,7} |
| w(0,0,2,1,1,1) = {5}:{6} | w(1,0,2,1,1,1) = {4}:{5} | w(2,0,2,1,1,1) = {10}:{11} |
| w(0,1,2,1,1,1) = {5}:{6} | w(1,1,2,1,1,1) = {2}:{10} | w(2,1,2,1,1,1) = {1,3}:{4,10} |
| w(0,2,2,1,1,1) = {2,3}:{9,11} | w(1,2,2,1,1,1) = {1}:{8} | w(2,2,2,1,1,1) = {2,3}:{4,8} |
| w(0,0,0,2,1,1) = {1}:{2} | w(1,0,0,2,1,1) = {2}:{4} | w(2,0,0,2,1,1) = {9}:{11} |
| w(0,1,0,2,1,1) = {4,6}:{9,11} | w(1,1,0,2,1,1) = {7}:{8} | w(2,1,0,2,1,1) = {1,2}:{3,11} |
| w(0,2,0,2,1,1) = {1,4}:{10,11} | w(1,2,0,2,1,1) = {9}:{10} | w(2,2,0,2,1,1) = {1,2}:{3,11} |
| w(0,0,1,2,1,1) = {9}:{11} | w(1,0,1,2,1,1) = {4}:{6} | w(2,0,1,2,1,1) = {4}:{6} |
| w(0,1,1,2,1,1) = {2,3}:{9,11} | w(1,1,1,2,1,1) = {3,6}:{4,9} | w(2,1,1,2,1,1) = {5,10}:{8,11} |
| w(0,2,1,2,1,1) = {9}:{11} | w(1,2,1,2,1,1) = {3,6}:{4,7} | w(2,2,1,2,1,1) = {5,8}:{10,11} |
| w(0,0,2,2,1,1) = {4}:{11} | w(1,0,2,2,1,1) = {4}:{5} | w(2,0,2,2,1,1) = {4}:{5} |
| w(0,1,2,2,1,1) = {2,3}:{10,11} | w(1,1,2,2,1,1) = {3,5}:{4,10} | w(2,1,2,2,1,1) = {6,9}:{7,11} |
| w(0,2,2,2,1,1) = {10}:{11} | w(1,2,2,2,1,1) = {3,5}:{4,8} | w(2,2,2,2,1,1) = {6,7}:{9,11} |
| w(0,0,0,0,2,1) = {2}:{3} | w(1,0,0,0,2,1) = {8}:{11} | w(2,0,0,0,2,1) = {8}:{11} |
| w(0,1,0,0,2,1) = {1}:{7} | w(1,1,0,0,2,1) = {6}:{10} | w(2,1,0,0,2,1) = {5}:{9} |
| w(0,2,0,0,2,1) = {3}:{11} | w(1,2,0,0,2,1) = {6}:{8} | w(2,2,0,0,2,1) = {5}:{7} |
| w(0,0,1,0,2,1) = {10}:{11} | w(1,0,1,0,2,1) = {7}:{11} | w(2,0,1,0,2,1) = {3}:{11} |
| w(0,1,1,0,2,1) = {2,3}:{9,11} | w(1,1,1,0,2,1) = {1}:{8} | w(2,1,1,0,2,1) = {5,6}:{9,10} |
| w(0,2,1,0,2,1) = {2,3}:{10,11} | w(1,2,1,0,2,1) = {2}:{10} | w(2,2,1,0,2,1) = {5,6}:{7,8} |
| w(0,0,2,0,2,1) = {6}:{11} | w(1,0,2,0,2,1) = {10}:{11} | w(2,0,2,0,2,1) = {3}:{11} |
| w(0,1,2,0,2,1) = {2,3}:{10,11} | w(1,1,2,0,2,1) = {1}:{7} | w(2,1,2,0,2,1) = {5,6}:{9,10} |
| w(0,2,2,0,2,1) = {2,3}:{9,11} | w(1,2,2,0,2,1) = {2}:{9} | w(2,2,2,0,2,1) = {5,6}:{7,8} |
| w(0,0,0,1,2,1) = {9}:{11} | w(1,0,0,1,2,1) = {9}:{11} | w(2,0,0,1,2,1) = {1,11}:{3,4} |
| w(0,1,0,1,2,1) = {5}:{11} | w(1,1,0,1,2,1) = {3}:{4} | w(2,1,0,1,2,1) = {10}:{8} |
| w(0,2,0,1,2,1) = {1}:{11} | w(1,2,0,1,2,1) = {3}:{4} | w(2,2,0,1,2,1) = {8}:{10} |
| w(0,0,1,1,2,1) = {4,10}:{5,11} | w(1,0,1,1,2,1) = {4}:{6} | w(2,0,1,1,2,1) = {10}:{11} |
| w(0,1,1,1,2,1) = {5}:{6} | w(1,1,1,1,2,1) = {1}:{2} | w(2,1,1,1,2,1) = {1}:{11} |
| w(0,2,1,1,2,1) = {2,3}:{10,11} | w(1,2,1,1,2,1) = {2}:{1} | w(2,2,1,1,2,1) = {2}:{11} |
| w(0,0,2,1,2,1) = {5}:{6} | w(1,0,2,1,2,1) = {4}:{5} | w(2,0,2,1,2,1) = {7}:{11} |
| w(0,1,2,1,2,1) = {5}:{6} | w(1,1,2,1,2,1) = {1}:{2} | w(2,1,2,1,2,1) = {1}:{11} |
| w(0,2,2,1,2,1) = {2,3}:{9,11} | w(1,2,2,1,2,1) = {2}:{1} | w(2,2,2,1,2,1) = {2}:{11} |
| w(0,0,0,2,2,1) = {1}:{3} | w(1,0,0,2,2,1) = {7}:{11} | w(2,0,0,2,2,1) = {3,4}:{9,11} |
| w(0,1,0,2,2,1) = {6}:{11} | w(1,1,0,2,2,1) = {2,6}:{4,5} | w(2,1,0,2,2,1) = {3,11}:{4,8} |
| w(0,2,0,2,2,1) = {1}:{11} | w(1,2,0,2,2,1) = {1,6}:{4,5} | w(2,2,0,2,2,1) = {3,11}:{4,10} |
| w(0,0,1,2,2,1) = {9}:{11} | w(1,0,1,2,2,1) = {7}:{8} | w(2,0,1,2,2,1) = {3,4}:{5,6} |
| w(0,1,1,2,2,1) = {2,3}:{9,11} | w(1,1,1,2,2,1) = {3,4}:{8,9} | w(2,1,1,2,2,1) = {5,10}:{8,9} |
| w(0,2,1,2,2,1) = {9}:{11} | w(1,2,1,2,2,1) = {3,4}:{7,10} | w(2,2,1,2,2,1) = {5,8}:{7,10} |
| w(0,0,2,2,2,1) = {1}:{11} | w(1,0,2,2,2,1) = {10}:{9} | w(2,0,2,2,2,1) = {3,4}:{5,6} |
| w(0,1,2,2,2,1) = {2,3}:{10,11} | w(1,1,2,2,2,1) = {3,4}:{7,10} | w(2,1,2,2,2,1) = {6,9}:{7,10} |
| w(0,2,2,2,2,1) = {10}:{11} | w(1,2,2,2,2,1) = {3,4}:{8,9} | w(2,2,2,2,2,1) = {6,7}:{8,9} |
| w(0,0,0,0,0,2) = {2}:{3} | w(1,0,0,0,0,2) = {8}:{11} | w(2,0,0,0,0,2) = {8}:{11} |
| w(0,1,0,0,0,2) = {1}:{11} | w(1,1,0,0,0,2) = {2}:{7} | w(2,1,0,0,0,2) = {2}:{10} |
| w(0,2,0,0,0,2) = {1}:{11} | w(1,2,0,0,0,2) = {1}:{9} | w(2,2,0,0,0,2) = {1}:{8} |
| w(0,0,1,0,0,2) = {4}:{6} | w(1,0,1,0,0,2) = {10}:{11} | w(2,0,1,0,0,2) = {3}:{11} |
| w(0,1,1,0,0,2) = {9}:{11} | w(1,1,1,0,0,2) = {3,6}:{5,9} | w(2,1,1,0,0,2) = {6}:{8} |
| w(0,2,1,0,0,2) = {10}:{11} | w(1,2,1,0,0,2) = {3,6}:{5,7} | w(2,2,1,0,0,2) = {6}:{10} |
| w(0,0,2,0,0,2) = {4}:{6} | w(1,0,2,0,0,2) = {7}:{11} | w(2,0,2,0,0,2) = {3}:{11} |
| w(0,1,2,0,0,2) = {10}:{11} | w(1,1,2,0,0,2) = {3,5}:{6,10} | w(2,1,2,0,0,2) = {5}:{7} |
| w(0,2,2,0,0,2) = {9}:{11} | w(1,2,2,0,0,2) = {3,5}:{6,8} | w(2,2,2,0,0,2) = {5}:{9} |
| w(0,0,0,1,0,2) = {1}:{11} | w(1,0,0,1,0,2) = {9}:{11} | w(2,0,0,1,0,2) = {4}:{5} |
| w(0,1,0,1,0,2) = {2}:{11} | w(1,1,0,1,0,2) = {2,4}:{3,8} | w(2,1,0,1,0,2) = {6}:{10} |
| w(0,2,0,1,0,2) = {7}:{11} | w(1,2,0,1,0,2) = {1,4}:{3,10} | w(2,2,0,1,0,2) = {6}:{8} |
| w(0,0,1,1,0,2) = {2}:{4} | w(1,0,1,1,0,2) = {4}:{11} | w(2,0,1,1,0,2) = {2}:{3} |
| w(0,1,1,1,0,2) = {10}:{11} | w(1,1,1,1,0,2) = {8}:{11} | w(2,1,1,1,0,2) = {1,4}:{3,6} |
| w(0,2,1,1,0,2) = {10}:{11} | w(1,2,1,1,0,2) = {10}:{11} | w(2,2,1,1,0,2) = {2,4}:{3,6} |
| w(0,0,2,1,0,2) = {5}:{6} | w(1,0,2,1,0,2) = {4}:{11} | w(2,0,2,1,0,2) = {1}:{3} |
| w(0,1,2,1,0,2) = {9}:{11} | w(1,1,2,1,0,2) = {7}:{11} | w(2,1,2,1,0,2) = {1,4}:{3,5} |
| w(0,2,2,1,0,2) = {9}:{11} | w(1,2,2,1,0,2) = {9}:{11} | w(2,2,2,1,0,2) = {2,4}:{3,5} |
| w(0,0,0,2,0,2) = {1}:{11} | w(1,0,0,2,0,2) = {10}:{11} | w(2,0,0,2,0,2) = {9}:{11} |
| w(0,1,0,2,0,2) = {2}:{11} | w(1,1,0,2,0,2) = {2}:{11} | w(2,1,0,2,0,2) = {10}:{8} |
| w(0,2,0,2,0,2) = {8}:{11} | w(1,2,0,2,0,2) = {1}:{11} | w(2,2,0,2,0,2) = {8}:{10} |
| w(0,0,1,2,0,2) = {9}:{11} | w(1,0,1,2,0,2) = {3,4}:{7,8} | w(2,0,1,2,0,2) = {4}:{5} |
| w(0,1,1,2,0,2) = {9}:{11} | w(1,1,1,2,0,2) = {2}:{3} | w(2,1,1,2,0,2) = {1}:{8} |
| w(0,2,1,2,0,2) = {5}:{6} | w(1,2,1,2,0,2) = {1}:{3} | w(2,2,1,2,0,2) = {2}:{10} |
| w(0,0,2,2,0,2) = {6}:{11} | w(1,0,2,2,0,2) = {3,4}:{9,10} | w(2,0,2,2,0,2) = {4}:{6} |
| w(0,1,2,2,0,2) = {10}:{11} | w(1,1,2,2,0,2) = {2}:{3} | w(2,1,2,2,0,2) = {1}:{7} |
| w(0,2,2,2,0,2) = {5}:{6} | w(1,2,2,2,0,2) = {1}:{3} | w(2,2,2,2,0,2) = {2}:{9} |
| w(0,0,0,0,1,2) = {2}:{3} | w(1,0,0,0,1,2) = {8}:{11} | w(2,0,0,0,1,2) = {8}:{11} |
| w(0,1,0,0,1,2) = {3}:{11} | w(1,1,0,0,1,2) = {5}:{7} | w(2,1,0,0,1,2) = {6}:{8} |
| w(0,2,0,0,1,2) = {1}:{7} | w(1,2,0,0,1,2) = {5}:{9} | w(2,2,0,0,1,2) = {6}:{10} |
| w(0,0,1,0,1,2) = {6}:{11} | w(1,0,1,0,1,2) = {3}:{11} | w(2,0,1,0,1,2) = {10}:{11} |
| w(0,1,1,0,1,2) = {2,3}:{9,11} | w(1,1,1,0,1,2) = {5,6}:{7,8} | w(2,1,1,0,1,2) = {2}:{9} |
| w(0,2,1,0,1,2) = {2,3}:{10,11} | w(1,2,1,0,1,2) = {5,6}:{9,10} | w(2,2,1,0,1,2) = {1}:{7} |
| w(0,0,2,0,1,2) = {10}:{11} | w(1,0,2,0,1,2) = {3}:{11} | w(2,0,2,0,1,2) = {7}:{11} |
| w(0,1,2,0,1,2) = {2,3}:{10,11} | w(1,1,2,0,1,2) = {5,6}:{7,8} | w(2,1,2,0,1,2) = {2}:{10} |
| w(0,2,2,0,1,2) = {2,3}:{9,11} | w(1,2,2,0,1,2) = {5,6}:{9,10} | w(2,2,2,0,1,2) = {1}:{8} |
| w(0,0,0,1,1,2) = {1}:{3} | w(1,0,0,1,1,2) = {3,4}:{9,11} | w(2,0,0,1,1,2) = {7}:{11} |
| w(0,1,0,1,1,2) = {1}:{11} | w(1,1,0,1,1,2) = {3,11}:{4,10} | w(2,1,0,1,1,2) = {1,6}:{4,5} |
| w(0,2,0,1,1,2) = {6}:{11} | w(1,2,0,1,1,2) = {3,11}:{4,8} | w(2,2,0,1,1,2) = {2,6}:{4,5} |
| w(0,0,1,1,1,2) = {1}:{11} | w(1,0,1,1,1,2) = {3,4}:{5,6} | w(2,0,1,1,1,2) = {10}:{9} |
| w(0,1,1,1,1,2) = {10}:{11} | w(1,1,1,1,1,2) = {6,7}:{8,9} | w(2,1,1,1,1,2) = {3,4}:{8,9} |
| w(0,2,1,1,1,2) = {2,3}:{10,11} | w(1,2,1,1,1,2) = {6,9}:{7,10} | w(2,2,1,1,1,2) = {3,4}:{7,10} |
| w(0,0,2,1,1,2) = {9}:{11} | w(1,0,2,1,1,2) = {3,4}:{5,6} | w(2,0,2,1,1,2) = {7}:{8} |
| w(0,1,2,1,1,2) = {9}:{11} | w(1,1,2,1,1,2) = {5,8}:{7,10} | w(2,1,2,1,1,2) = {3,4}:{7,10} |
| w(0,2,2,1,1,2) = {2,3}:{9,11} | w(1,2,2,1,1,2) = {5,10}:{8,9} | w(2,2,2,1,1,2) = {3,4}:{8,9} |
| w(0,0,0,2,1,2) = {9}:{11} | w(1,0,0,2,1,2) = {1,11}:{3,4} | w(2,0,0,2,1,2) = {9}:{11} |
| w(0,1,0,2,1,2) = {1}:{11} | w(1,1,0,2,1,2) = {8}:{10} | w(2,1,0,2,1,2) = {3}:{4} |
| w(0,2,0,2,1,2) = {5}:{11} | w(1,2,0,2,1,2) = {10}:{8} | w(2,2,0,2,1,2) = {3}:{4} |
| w(0,0,1,2,1,2) = {5}:{6} | w(1,0,1,2,1,2) = {7}:{11} | w(2,0,1,2,1,2) = {4}:{5} |
| w(0,1,1,2,1,2) = {2,3}:{9,11} | w(1,1,1,2,1,2) = {2}:{11} | w(2,1,1,2,1,2) = {2}:{1} |
| w(0,2,1,2,1,2) = {5}:{6} | w(1,2,1,2,1,2) = {1}:{11} | w(2,2,1,2,1,2) = {1}:{2} |
| w(0,0,2,2,1,2) = {4,10}:{5,11} | w(1,0,2,2,1,2) = {10}:{11} | w(2,0,2,2,1,2) = {4}:{6} |
| w(0,1,2,2,1,2) = {2,3}:{10,11} | w(1,1,2,2,1,2) = {2}:{11} | w(2,1,2,2,1,2) = {2}:{1} |
| w(0,2,2,2,1,2) = {5}:{6} | w(1,2,2,2,1,2) = {1}:{11} | w(2,2,2,2,1,2) = {1}:{2} |
| w(0,0,0,0,2,2) = {2}:{11} | w(1,0,0,0,2,2) = {8}:{11} | w(2,0,0,0,2,2) = {8}:{11} |
| w(0,1,0,0,2,2) = {2}:{11} | w(1,1,0,0,2,2) = {1}:{7} | w(2,1,0,0,2,2) = {1}:{10} |
| w(0,2,0,0,2,2) = {1}:{7} | w(1,2,0,0,2,2) = {2}:{9} | w(2,2,0,0,2,2) = {2}:{8} |
| w(0,0,1,0,2,2) = {6}:{7} | w(1,0,1,0,2,2) = {1,11}:{4,7} | w(2,0,1,0,2,2) = {6,8}:{10,11} |
| w(0,1,1,0,2,2) = {2,3}:{9,11} | w(1,1,1,0,2,2) = {8}:{10} | w(2,1,1,0,2,2) = {3}:{9} |
| w(0,2,1,0,2,2) = {2,3}:{10,11} | w(1,2,1,0,2,2) = {10}:{8} | w(2,2,1,0,2,2) = {3}:{7} |
| w(0,0,2,0,2,2) = {4}:{6} | w(1,0,2,0,2,2) = {2,11}:{4,10} | w(2,0,2,0,2,2) = {5,9}:{7,11} |
| w(0,1,2,0,2,2) = {2,3}:{10,11} | w(1,1,2,0,2,2) = {7}:{9} | w(2,1,2,0,2,2) = {3}:{10} |
| w(0,2,2,0,2,2) = {2,3}:{9,11} | w(1,2,2,0,2,2) = {9}:{7} | w(2,2,2,0,2,2) = {3}:{8} |
| w(0,0,0,1,2,2) = {1}:{2} | w(1,0,0,1,2,2) = {9}:{11} | w(2,0,0,1,2,2) = {2}:{4} |
| w(0,1,0,1,2,2) = {1,4}:{10,11} | w(1,1,0,1,2,2) = {1,2}:{3,11} | w(2,1,0,1,2,2) = {9}:{10} |
| w(0,2,0,1,2,2) = {4,6}:{9,11} | w(1,2,0,1,2,2) = {1,2}:{3,11} | w(2,2,0,1,2,2) = {7}:{8} |
| w(0,0,1,1,2,2) = {4}:{11} | w(1,0,1,1,2,2) = {4}:{5} | w(2,0,1,1,2,2) = {4}:{5} |
| w(0,1,1,1,2,2) = {10}:{11} | w(1,1,1,1,2,2) = {6,7}:{9,11} | w(2,1,1,1,2,2) = {3,5}:{4,8} |
| w(0,2,1,1,2,2) = {2,3}:{10,11} | w(1,2,1,1,2,2) = {6,9}:{7,11} | w(2,2,1,1,2,2) = {3,5}:{4,10} |
| w(0,0,2,1,2,2) = {9}:{11} | w(1,0,2,1,2,2) = {4}:{6} | w(2,0,2,1,2,2) = {4}:{6} |
| w(0,1,2,1,2,2) = {9}:{11} | w(1,1,2,1,2,2) = {5,8}:{10,11} | w(2,1,2,1,2,2) = {3,6}:{4,7} |
| w(0,2,2,1,2,2) = {2,3}:{9,11} | w(1,2,2,1,2,2) = {5,10}:{8,11} | w(2,2,2,1,2,2) = {3,6}:{4,9} |
| w(0,0,0,2,2,2) = {10}:{11} | w(1,0,0,2,2,2) = {2}:{11} | w(2,0,0,2,2,2) = {3,4}:{9,11} |
| w(0,1,0,2,2,2) = {1,4}:{9,11} | w(1,1,0,2,2,2) = {1}:{8} | w(2,1,0,2,2,2) = {6}:{11} |
| w(0,2,0,2,2,2) = {4,5}:{10,11} | w(1,2,0,2,2,2) = {2}:{10} | w(2,2,0,2,2,2) = {6}:{11} |
| w(0,0,1,2,2,2) = {5}:{6} | w(1,0,1,2,2,2) = {10}:{11} | w(2,0,1,2,2,2) = {4}:{5} |
| w(0,1,1,2,2,2) = {2,3}:{9,11} | w(1,1,1,2,2,2) = {2,3}:{4,8} | w(2,1,1,2,2,2) = {1}:{8} |
| w(0,2,1,2,2,2) = {5}:{6} | w(1,2,1,2,2,2) = {1,3}:{4,10} | w(2,2,1,2,2,2) = {2}:{10} |
| w(0,0,2,2,2,2) = {1}:{11} | w(1,0,2,2,2,2) = {7}:{11} | w(2,0,2,2,2,2) = {4}:{6} |
| w(0,1,2,2,2,2) = {2,3}:{10,11} | w(1,1,2,2,2,2) = {2,3}:{4,7} | w(2,1,2,2,2,2) = {1}:{7} |
| w(0,2,2,2,2,2) = {5}:{6} | w(1,2,2,2,2,2) = {1,3}:{4,9} | w(2,2,2,2,2,2) = {2}:{9} |

List 3

The completed map of the second algorithm to sort 11 coins

f(0,0,0,0,0,0) = 2047          f(1,0,0,0,0,0) = 1468          f(2,0,0,0,0,0) = 579
f(0,1,0,0,0,0) = 1855          f(1,1,0,0,0,0) = 2046          f(2,1,0,0,0,0) = 2
f(0,2,0,0,0,0) = 192           f(1,2,0,0,0,0) = 2045          f(2,2,0,0,0,0) = 1
f(1,0,1,0,0,0) = 1786          f(2,0,1,0,0,0) = 134           f(0,1,1,0,0,0) = 1974
f(1,1,1,0,0,0) = 944           f(2,1,1,0,0,0) = 1551          f(0,2,1,0,0,0) = 137
f(1,2,1,0,0,0) = 752           f(2,2,1,0,0,0) = 1167          f(1,0,2,0,0,0) = 1913
f(2,0,2,0,0,0) = 261           f(0,1,2,0,0,0) = 1910          f(1,1,2,0,0,0) = 880
f(2,1,2,0,0,0) = 1295          f(0,2,2,0,0,0) = 73            f(1,2,2,0,0,0) = 496
f(2,2,2,0,0,0) = 1103          f(1,0,0,1,0,0) = 2044          f(2,0,0,1,0,0) = 1991
f(0,1,0,1,0,0) = 1838          f(1,1,0,1,0,0) = 300           f(2,1,0,1,0,0) = 774
f(0,2,0,1,0,0) = 225           f(1,2,0,1,0,0) = 108           f(2,2,0,1,0,0) = 197
f(0,0,1,1,0,0) = 1727          f(1,0,1,1,0,0) = 696           f(2,0,1,1,0,0) = 2031
f(0,1,1,1,0,0) = 1674          f(1,1,1,1,0,0) = 1594          f(2,1,1,1,0,0) = 1827
f(0,2,1,1,0,0) = 128           f(1,2,1,1,0,0) = 1209          f(2,2,1,1,0,0) = 1251
f(0,0,2,1,0,0) = 74            f(1,0,2,1,0,0) = 376           f(2,0,2,1,0,0) = 2015
f(0,1,2,1,0,0) = 1354          f(1,1,2,1,0,0) = 1338          f(2,1,2,1,0,0) = 1811
f(0,2,2,1,0,0) = 64            f(1,2,2,1,0,0) = 1145          f(2,2,2,1,0,0) = 1235
f(1,0,0,2,0,0) = 56            f(2,0,0,2,0,0) = 3             f(0,1,0,2,0,0) = 1822
f(1,1,0,2,0,0) = 1850          f(2,1,0,2,0,0) = 1939          f(0,2,0,2,0,0) = 209
f(1,2,0,2,0,0) = 1273          f(2,2,0,2,0,0) = 1747          f(0,0,1,2,0,0) = 1973
f(1,0,1,2,0,0) = 32            f(2,0,1,2,0,0) = 1671          f(0,1,1,2,0,0) = 1983
f(1,1,1,2,0,0) = 812           f(2,1,1,2,0,0) = 902           f(0,2,1,2,0,0) = 693
f(1,2,1,2,0,0) = 236           f(2,2,1,2,0,0) = 709           f(0,0,2,2,0,0) = 320
f(1,0,2,2,0,0) = 16            f(2,0,2,2,0,0) = 1351          f(0,1,2,2,0,0) = 1919
f(1,1,2,2,0,0) = 796           f(2,1,2,2,0,0) = 838           f(0,2,2,2,0,0) = 373
f(1,2,2,2,0,0) = 220           f(2,2,2,2,0,0) = 453           f(0,0,0,0,1,0) = 960
f(1,0,0,0,1,0) = 1534          f(2,0,0,0,1,0) = 66            f(0,1,0,0,1,0) = 1804
f(1,1,0,0,1,0) = 1854          f(2,1,0,0,1,0) = 1943          f(0,2,0,0,1,0) = 245
f(1,2,0,0,1,0) = 1277          f(2,2,0,0,1,0) = 1751          f(0,0,1,0,1,0) = 1060
f(2,0,1,0,1,0) = 710           f(1,1,1,0,1,0) = 1832          f(2,1,1,0,1,0) = 1702
f(1,2,1,0,1,0) = 1256          f(2,2,1,0,1,0) = 1701          f(0,0,2,0,1,0) = 34
f(2,0,2,0,1,0) = 837           f(1,1,2,0,1,0) = 1816          f(2,1,2,0,1,0) = 1366
f(1,2,2,0,1,0) = 1240          f(2,2,2,0,1,0) = 1365          f(0,0,0,1,1,0) = 1698
f(2,0,0,1,1,0) = 388           f(1,1,0,1,1,0) = 1662          f(2,1,0,1,1,0) = 919
f(1,2,0,1,1,0) = 1469          f(2,2,0,1,1,0) = 727           f(0,0,1,1,1,0) = 1709
f(1,0,1,1,1,0) = 648           f(2,0,1,1,1,0) = 1699          f(0,1,1,1,1,0) = 1930
f(1,1,1,1,1,0) = 1714          f(2,1,1,1,1,0) = 1538          f(1,2,1,1,1,0) = 1713
f(2,2,1,1,1,0) = 1153          f(0,0,2,1,1,0) = 630           f(1,0,2,1,1,0) = 328
f(2,0,2,1,1,0) = 1363          f(0,1,2,1,1,0) = 1866          f(1,1,2,1,1,0) = 1394
f(2,1,2,1,1,0) = 1282          f(1,2,2,1,1,0) = 1393          f(2,2,2,1,1,0) = 1089
f(0,0,0,2,1,0) = 667           f(1,0,0,2,1,0) = 392           f(1,1,0,2,1,0) = 2034

| | | |
|---|---|---|
| f(2,1,0,2,1,0,0) = 1799 | f(1,2,0,2,1,0,0) = 2033 | f(2,2,0,2,1,0,0) = 1223 |
| f(0,0,1,2,1,0,0) = 1462 | f(1,0,1,2,1,0,0) = 1064 | f(2,0,1,2,1,0,0) = 1639 |
| f(1,1,1,2,1,0,0) = 520 | f(2,1,1,2,1,0,0) = 806 | f(0,2,1,2,1,0,0) = 757 |
| f(1,2,1,2,1,0,0) = 136 | f(2,2,1,2,1,0,0) = 229 | f(1,0,2,2,1,0,0) = 1048 |
| f(2,0,2,2,1,0,0) = 1623 | f(1,1,2,2,1,0,0) = 264 | f(2,1,2,2,1,0,0) = 790 |
| f(0,2,2,2,1,0,0) = 501 | f(1,2,2,2,1,0,0) = 72 | f(2,2,2,2,1,0,0) = 213 |
| f(0,0,0,0,2,0,0) = 1087 | f(1,0,0,0,2,0,0) = 1981 | f(2,0,0,0,2,0,0) = 513 |
| f(0,1,0,0,2,0,0) = 1802 | f(1,1,0,0,2,0,0) = 296 | f(2,1,0,0,2,0,0) = 770 |
| f(0,2,0,0,2,0,0) = 243 | f(1,2,0,0,2,0,0) = 104 | f(2,2,0,0,2,0,0) = 193 |
| f(0,0,1,0,2,0,0) = 2013 | f(1,0,1,0,2,0,0) = 1210 | f(1,1,1,0,2,0,0) = 682 |
| f(2,1,1,0,2,0,0) = 807 | f(1,2,1,0,2,0,0) = 681 | f(2,2,1,0,2,0,0) = 231 |
| f(0,0,2,0,2,0,0) = 987 | f(1,0,2,0,2,0,0) = 1337 | f(1,1,2,0,2,0,0) = 346 |
| f(2,1,2,0,2,0,0) = 791 | f(1,2,2,0,2,0,0) = 345 | f(2,2,2,0,2,0,0) = 215 |
| f(0,0,0,1,2,0,0) = 1380 | f(2,0,0,1,2,0,0) = 1655 | f(1,1,0,1,2,0,0) = 824 |
| f(2,1,0,1,2,0,0) = 14 | f(1,2,0,1,2,0,0) = 248 | f(2,2,0,1,2,0,0) = 13 |
| f(1,0,1,1,2,0,0) = 424 | f(2,0,1,1,2,0,0) = 999 | f(0,1,1,1,2,0,0) = 1546 |
| f(1,1,1,1,2,0,0) = 1834 | f(2,1,1,1,2,0,0) = 1975 | f(1,2,1,1,2,0,0) = 1257 |
| f(2,2,1,1,2,0,0) = 1783 | f(0,0,2,1,2,0,0) = 585 | f(1,0,2,1,2,0,0) = 408 |
| f(2,0,2,1,2,0,0) = 983 | f(0,1,2,1,2,0,0) = 1290 | f(1,1,2,1,2,0,0) = 1818 |
| f(2,1,2,1,2,0,0) = 1911 | f(1,2,2,1,2,0,0) = 1241 | f(2,2,2,1,2,0,0) = 1527 |
| f(0,0,0,2,2,0,0) = 349 | f(1,0,0,2,2,0,0) = 1659 | f(1,1,0,2,2,0,0) = 1320 |
| f(2,1,0,2,2,0,0) = 578 | f(1,2,0,2,2,0,0) = 1128 | f(2,2,0,2,2,0,0) = 385 |
| f(0,0,1,2,2,0,0) = 1417 | f(1,0,1,2,2,0,0) = 684 | f(2,0,1,2,2,0,0) = 1719 |
| f(1,1,1,2,2,0,0) = 958 | f(2,1,1,2,2,0,0) = 654 | f(0,2,1,2,2,0,0) = 181 |
| f(1,2,1,2,2,0,0) = 765 | f(2,2,1,2,2,0,0) = 653 | f(0,0,2,2,2,0,0) = 338 |
| f(1,0,2,2,2,0,0) = 348 | f(2,0,2,2,2,0,0) = 1399 | f(1,1,2,2,2,0,0) = 894 |
| f(2,1,2,2,2,0,0) = 334 | f(0,2,2,2,2,0,0) = 117 | f(1,2,2,2,2,0,0) = 509 |
| f(2,2,2,2,2,0,0) = 333 | f(0,0,0,0,0,1,0) = 1993 | f(1,0,0,0,0,1,0) = 1456 |
| f(2,0,0,0,0,1,0) = 1987 | f(0,1,0,0,0,1,0) = 822 | f(1,1,0,0,0,1,0) = 816 |
| f(2,1,0,0,0,1,0) = 582 | f(0,2,0,0,0,1,0) = 246 | f(1,2,0,0,0,1,0) = 240 |
| f(2,2,0,0,0,1,0) = 389 | f(0,0,1,0,0,1,0) = 1005 | f(1,0,1,0,0,1,0) = 1716 |
| f(2,0,1,0,0,1,0) = 1126 | f(0,1,1,0,0,1,0) = 1920 | f(1,1,1,0,0,1,0) = 544 |
| f(2,1,1,0,0,1,0) = 514 | f(0,2,1,0,0,1,0) = 1737 | f(1,2,1,0,0,1,0) = 160 |
| f(2,2,1,0,0,1,0) = 129 | f(0,0,2,0,0,1,0) = 1003 | f(1,0,2,0,0,1,0) = 1396 |
| f(2,0,2,0,0,1,0) = 1557 | f(0,1,2,0,0,1,0) = 1856 | f(1,1,2,0,0,1,0) = 272 |
| f(2,1,2,0,0,1,0) = 258 | f(0,2,2,0,0,1,0) = 1481 | f(1,2,2,0,0,1,0) = 80 |
| f(2,2,2,0,0,1,0) = 65 | f(0,0,0,1,0,1,0) = 1451 | f(1,0,0,1,0,1,0) = 2032 |
| f(2,0,0,1,0,1,0) = 1604 | f(0,1,0,1,0,1,0) = 1454 | f(1,1,0,1,0,1,0) = 1022 |
| f(2,1,0,1,0,1,0) = 1951 | f(0,2,0,1,0,1,0) = 1262 | f(1,2,0,1,0,1,0) = 1021 |
| f(2,2,0,1,0,1,0) = 1759 | f(0,0,1,1,0,1,0) = 1198 | f(1,0,1,1,0,1,0) = 2024 |
| f(2,0,1,1,0,1,0) = 1766 | f(0,1,1,1,0,1,0) = 650 | f(1,1,1,1,0,1,0) = 1978 |
| f(2,1,1,1,0,1,0) = 1575 | f(0,2,1,1,0,1,0) = 1728 | f(1,2,1,1,0,1,0) = 1785 |
| f(2,2,1,1,0,1,0) = 1191 | f(0,0,2,1,0,1,0) = 576 | f(1,0,2,1,0,1,0) = 2008 |
| f(2,0,2,1,0,1,0) = 1877 | f(0,1,2,1,0,1,0) = 330 | f(1,1,2,1,0,1,0) = 1914 |

f(2,1,2,1,0,1,0) = 1303
f(2,2,2,1,0,1,0) = 1111
f(2,0,0,2,0,1,0) = 1411
f(2,1,0,2,0,1,0) = 1891
f(2,2,0,2,0,1,0) = 1507
f(2,0,1,2,0,1,0) = 663
f(2,1,1,2,0,1,0) = 1924
f(2,2,1,2,0,1,0) = 1732
f(2,0,2,2,0,1,0) = 359
f(2,1,2,2,0,1,0) = 1860
f(2,2,2,2,0,1,0) = 1476
f(2,0,0,0,1,1,0) = 1474
f(2,1,0,0,1,1,0) = 1998
f(2,2,0,0,1,1,0) = 1997
f(2,0,1,0,1,1,0) = 1668
f(2,1,1,0,1,1,0) = 1670
f(2,2,1,0,1,1,0) = 1669
f(2,0,2,0,1,1,0) = 1348
f(2,1,2,0,1,1,0) = 1350
f(2,2,2,0,1,1,0) = 1349
f(2,0,0,1,1,1,0) = 1478
f(2,1,0,1,1,1,0) = 1942
f(2,2,0,1,1,1,0) = 1749
f(2,0,1,1,1,1,0) = 1730
f(2,1,1,1,1,1,0) = 1666
f(2,2,1,1,1,1,0) = 1665
f(2,0,2,1,1,1,0) = 1857
f(2,1,2,1,1,1,0) = 1346
f(2,2,2,1,1,1,0) = 1345
f(2,0,0,2,1,1,0) = 1419
f(2,1,0,2,1,1,0) = 1793
f(2,2,0,2,1,1,0) = 1218
f(2,0,1,2,1,1,0) = 647
f(2,1,1,2,1,1,0) = 1638
f(2,2,1,2,1,1,0) = 1445
f(2,0,2,2,1,1,0) = 327
f(2,1,2,2,1,1,0) = 1430
f(2,2,2,2,1,1,0) = 1621
f(2,0,0,0,2,1,0) = 1921
f(2,1,0,0,2,1,0) = 915
f(2,2,0,0,2,1,0) = 723
f(2,0,1,0,2,1,0) = 643
f(2,1,1,0,2,1,0) = 422
f(2,2,1,0,2,1,0) = 613

f(0,2,2,1,0,1,0) = 1472
f(0,0,0,2,0,1,0) = 1435
f(0,1,0,2,0,1,0) = 1630
f(0,2,0,2,0,1,0) = 1246
f(0,0,1,2,0,1,0) = 949
f(0,1,1,2,0,1,0) = 1929
f(0,2,1,2,0,1,0) = 691
f(0,0,2,2,0,1,0) = 1344
f(0,1,2,2,0,1,0) = 1865
f(0,2,2,2,0,1,0) = 371
f(0,0,0,0,1,1,0) = 2038
f(0,1,0,0,1,1,0) = 780
f(0,2,0,0,1,1,0) = 1013
f(0,0,1,0,1,1,0) = 2029
f(0,1,1,0,1,1,0) = 1572
f(0,2,1,0,1,1,0) = 1197
f(0,0,2,0,1,1,0) = 2018
f(0,1,2,0,1,1,0) = 1316
f(0,2,2,0,1,1,0) = 1133
f(0,0,0,1,1,1,0) = 1569
f(0,1,0,1,1,1,0) = 1006
f(0,2,0,1,1,1,0) = 228
f(0,0,1,1,1,1,0) = 1673
f(0,1,1,1,1,1,0) = 906
f(0,2,1,1,1,1,0) = 1188
f(0,0,2,1,1,1,0) = 1654
f(0,1,2,1,1,1,0) = 842
f(0,2,2,1,1,1,0) = 1124
f(0,0,0,2,1,1,0) = 1691
f(0,1,0,2,1,1,0) = 990
f(0,2,0,2,1,1,0) = 212
f(0,0,1,2,1,1,0) = 1408
f(0,1,1,2,1,1,0) = 1581
f(0,2,1,2,1,1,0) = 755
f(0,0,2,2,1,1,0) = 1387
f(0,1,2,2,1,1,0) = 1325
f(0,2,2,2,1,1,0) = 499
f(0,0,0,0,2,1,0) = 1078
f(0,1,0,0,2,1,0) = 778
f(0,2,0,0,2,1,0) = 1011
f(0,0,1,0,2,1,0) = 2004
f(0,1,1,0,2,1,0) = 1570
f(0,2,1,0,2,1,0) = 1195
f(0,0,2,0,2,1,0) = 978

f(1,2,2,1,0,1,0) = 1529
f(1,0,0,2,0,1,0) = 1080
f(1,1,0,2,0,1,0) = 1458
f(1,2,0,2,0,1,0) = 1649
f(1,0,1,2,0,1,0) = 608
f(1,1,1,2,0,1,0) = 556
f(1,2,1,2,0,1,0) = 172
f(1,0,2,2,0,1,0) = 592
f(1,1,2,2,0,1,0) = 284
f(1,2,2,2,0,1,0) = 92
f(1,0,0,0,1,1,0) = 1522
f(1,1,0,0,1,1,0) = 1552
f(1,2,0,0,1,1,0) = 1168
f(1,0,1,0,1,1,0) = 1004
f(1,1,1,0,1,1,0) = 1972
f(1,2,1,0,1,1,0) = 1780
f(1,0,2,0,1,1,0) = 988
f(1,1,2,0,1,1,0) = 1908
f(1,2,2,0,1,1,0) = 1524
f(1,0,0,1,1,1,0) = 1656
f(1,1,0,1,1,1,0) = 538
f(1,2,0,1,1,1,0) = 153
f(1,0,1,1,1,1,0) = 1000
f(1,1,1,1,1,1,0) = 1584
f(1,2,1,1,1,1,0) = 1200
f(1,0,2,1,1,1,0) = 984
f(1,1,2,1,1,1,0) = 1328
f(1,2,2,1,1,1,0) = 1136
f(1,0,0,2,1,1,0) = 1530
f(1,1,0,2,1,1,0) = 1840
f(1,2,0,2,1,1,0) = 1264
f(1,0,1,2,1,1,0) = 1640
f(1,1,1,2,1,1,0) = 1836
f(1,2,1,2,1,1,0) = 1260
f(1,0,2,2,1,1,0) = 1624
f(1,1,2,2,1,1,0) = 1820
f(1,2,2,2,1,1,0) = 1244
f(1,0,0,0,2,1,0) = 1969
f(1,1,0,0,2,1,0) = 1018
f(1,2,0,0,2,1,0) = 1017
f(1,0,1,0,2,1,0) = 1778
f(1,1,1,0,2,1,0) = 810
f(1,2,1,0,2,1,0) = 233
f(1,0,2,0,2,1,0) = 1905

f(2,0,2,0,2,1,0) = 323
f(2,1,2,0,2,1,0) = 598
f(2,2,2,0,2,1,0) = 405
f(2,0,0,1,2,1,0) = 1094
f(2,1,0,1,2,1,0) = 1038
f(2,2,0,1,2,1,0) = 1037
f(2,0,1,1,2,1,0) = 486
f(2,1,1,1,2,1,0) = 1591
f(2,2,1,1,2,1,0) = 1207
f(2,0,2,1,2,1,0) = 917
f(2,1,2,1,2,1,0) = 1335
f(2,2,2,1,2,1,0) = 1143
f(2,0,0,2,2,1,0) = 1611
f(2,1,0,2,2,1,0) = 1986
f(2,2,0,2,2,1,0) = 1985
f(2,0,1,2,2,1,0) = 1687
f(2,1,1,2,2,1,0) = 1678
f(2,2,1,2,2,1,0) = 1677
f(2,0,2,2,2,1,0) = 1383
f(2,1,2,2,2,1,0) = 1358
f(2,2,2,2,2,1,0) = 1357
f(2,0,0,0,0,2,0) = 591
f(2,1,0,0,0,2,0) = 1807
f(2,2,0,0,0,2,0) = 1231
f(2,0,1,0,0,2,0) = 651
f(2,1,1,0,0,2,0) = 1967
f(2,2,1,0,0,2,0) = 1775
f(2,0,2,0,0,2,0) = 331
f(2,1,2,0,0,2,0) = 1887
f(2,2,2,0,0,2,0) = 1503
f(2,0,0,1,0,2,0) = 967
f(2,1,0,1,0,2,0) = 398
f(2,2,0,1,0,2,0) = 589
f(2,0,1,1,0,2,0) = 1455
f(2,1,1,1,0,2,0) = 1955
f(2,2,1,1,0,2,0) = 1763
f(2,0,2,1,0,2,0) = 1439
f(2,1,2,1,0,2,0) = 1875
f(2,2,2,1,0,2,0) = 1491
f(2,0,0,2,0,2,0) = 15
f(2,1,0,2,0,2,0) = 1026
f(2,2,0,2,0,2,0) = 1025
f(2,0,1,2,0,2,0) = 39
f(2,1,1,2,0,2,0) = 518

f(0,1,2,0,2,1,0) = 1314
f(0,2,2,0,2,1,0) = 1131
f(0,0,0,1,2,1,0) = 1313
f(0,1,0,1,2,1,0) = 1438
f(0,2,0,1,2,1,0) = 226
f(0,0,1,1,2,1,0) = 1164
f(0,1,1,1,2,1,0) = 522
f(0,2,1,1,2,1,0) = 1186
f(0,0,2,1,2,1,0) = 1609
f(0,1,2,1,2,1,0) = 266
f(0,2,2,1,2,1,0) = 1122
f(0,0,0,2,2,1,0) = 1373
f(0,1,0,2,2,1,0) = 1646
f(0,2,0,2,2,1,0) = 210
f(0,0,1,2,2,1,0) = 1333
f(0,1,1,2,2,1,0) = 1579
f(0,2,1,2,2,1,0) = 179
f(0,0,2,2,2,1,0) = 1371
f(0,1,2,2,2,1,0) = 1323
f(0,2,2,2,2,1,0) = 115
f(0,0,0,0,0,2,0) = 54
f(0,1,0,0,0,2,0) = 1801
f(0,2,0,0,0,2,0) = 1225
f(0,0,1,0,0,2,0) = 1044
f(0,1,1,0,0,2,0) = 566
f(0,2,1,0,0,2,0) = 191
f(0,0,2,0,0,2,0) = 1042
f(0,1,2,0,0,2,0) = 310
f(0,2,2,0,0,2,0) = 127
f(0,0,0,1,0,2,0) = 612
f(0,1,0,1,0,2,0) = 801
f(0,2,0,1,0,2,0) = 417
f(0,0,1,1,0,2,0) = 703
f(0,1,1,1,0,2,0) = 1676
f(0,2,1,1,0,2,0) = 182
f(0,0,2,1,0,2,0) = 1098
f(0,1,2,1,0,2,0) = 1356
f(0,2,2,1,0,2,0) = 118
f(0,0,0,2,0,2,0) = 596
f(0,1,0,2,0,2,0) = 785
f(0,2,0,2,0,2,0) = 593
f(0,0,1,2,0,2,0) = 1471
f(0,1,1,2,0,2,0) = 575
f(0,2,1,2,0,2,0) = 1717

f(1,1,2,0,2,1,0) = 794
f(1,2,2,0,2,1,0) = 217
f(1,0,0,1,2,1,0) = 1992
f(1,1,0,1,2,1,0) = 924
f(1,2,0,1,2,1,0) = 732
f(1,0,1,1,2,1,0) = 1448
f(1,1,1,1,2,1,0) = 1824
f(1,2,1,1,2,1,0) = 1248
f(1,0,2,1,2,1,0) = 1432
f(1,1,2,1,2,1,0) = 1808
f(1,2,2,1,2,1,0) = 1232
f(1,0,0,2,2,1,0) = 1146
f(1,1,0,2,2,1,0) = 1896
f(1,2,0,2,2,1,0) = 1512
f(1,0,1,2,2,1,0) = 553
f(1,1,1,2,2,1,0) = 956
f(1,2,1,2,2,1,0) = 764
f(1,0,2,2,2,1,0) = 90
f(1,1,2,2,2,1,0) = 892
f(1,2,2,2,2,1,0) = 508
f(1,0,0,0,0,2,0) = 60
f(1,1,0,0,0,2,0) = 1658
f(1,2,0,0,0,2,0) = 1465
f(1,0,1,0,0,2,0) = 490
f(1,1,1,0,0,2,0) = 1982
f(1,2,1,0,0,2,0) = 1789
f(1,0,2,0,0,2,0) = 921
f(1,1,2,0,0,2,0) = 1918
f(1,2,2,0,0,2,0) = 1533
f(1,0,0,1,0,2,0) = 636
f(1,1,0,1,0,2,0) = 540
f(1,2,0,1,0,2,0) = 156
f(1,0,1,1,0,2,0) = 1688
f(1,1,1,1,0,2,0) = 571
f(1,2,1,1,0,2,0) = 187
f(1,0,2,1,0,2,0) = 1384
f(1,1,2,1,0,2,0) = 315
f(1,2,2,1,0,2,0) = 123
f(1,0,0,2,0,2,0) = 443
f(1,1,0,2,0,2,0) = 288
f(1,2,0,2,0,2,0) = 96
f(1,0,1,2,0,2,0) = 170
f(1,1,1,2,0,2,0) = 936
f(1,2,1,2,0,2,0) = 744

| | | |
|---|---|---|
| f(2,2,1,2,0,2,0) = 133 | f(0,0,2,2,0,2,0) = 849 | f(1,0,2,2,0,2,0) = 281 |
| f(2,0,2,2,0,2,0) = 23 | f(0,1,2,2,0,2,0) = 319 | f(1,1,2,2,0,2,0) = 856 |
| f(2,1,2,2,0,2,0) = 262 | f(0,2,2,2,0,2,0) = 1397 | f(1,2,2,2,0,2,0) = 472 |
| f(2,2,2,2,0,2,0) = 69 | f(0,0,0,0,1,2,0) = 969 | f(1,0,0,0,1,2,0) = 126 |
| f(2,0,0,0,1,2,0) = 78 | f(0,1,0,0,1,2,0) = 1036 | f(1,1,0,0,1,2,0) = 1324 |
| f(2,1,0,0,1,2,0) = 1030 | f(0,2,0,0,1,2,0) = 1269 | f(1,2,0,0,1,2,0) = 1132 |
| f(2,2,0,0,1,2,0) = 1029 | f(0,0,1,0,1,2,0) = 1069 | f(1,0,1,0,1,2,0) = 1724 |
| f(2,0,1,0,1,2,0) = 142 | f(0,1,1,0,1,2,0) = 916 | f(1,1,1,0,1,2,0) = 1642 |
| f(2,1,1,0,1,2,0) = 1830 | f(0,2,1,0,1,2,0) = 733 | f(1,2,1,0,1,2,0) = 1449 |
| f(2,2,1,0,1,2,0) = 1253 | f(0,0,2,0,1,2,0) = 43 | f(1,0,2,0,1,2,0) = 1404 |
| f(2,0,2,0,1,2,0) = 269 | f(0,1,2,0,1,2,0) = 852 | f(1,1,2,0,1,2,0) = 1434 |
| f(2,1,2,0,1,2,0) = 1814 | f(0,2,2,0,1,2,0) = 477 | f(1,2,2,0,1,2,0) = 1625 |
| f(2,2,2,0,1,2,0) = 1237 | f(0,0,0,1,1,2,0) = 674 | f(1,0,0,1,1,2,0) = 436 |
| f(2,0,0,1,1,2,0) = 901 | f(0,1,0,1,1,2,0) = 1837 | f(1,1,0,1,1,2,0) = 62 |
| f(2,1,0,1,1,2,0) = 535 | f(0,2,0,1,1,2,0) = 401 | f(1,2,0,1,1,2,0) = 61 |
| f(2,2,0,1,1,2,0) = 151 | f(0,0,1,1,1,2,0) = 676 | f(1,0,1,1,1,2,0) = 664 |
| f(2,0,1,1,1,2,0) = 1957 | f(0,1,1,1,1,2,0) = 1932 | f(1,1,1,1,1,2,0) = 690 |
| f(2,1,1,1,1,2,0) = 1539 | f(0,2,1,1,1,2,0) = 724 | f(1,2,1,1,1,2,0) = 689 |
| f(2,2,1,1,1,2,0) = 1155 | f(0,0,2,1,1,2,0) = 714 | f(1,0,2,1,1,2,0) = 360 |
| f(2,0,2,1,1,2,0) = 1494 | f(0,1,2,1,1,2,0) = 1868 | f(1,1,2,1,1,2,0) = 370 |
| f(2,1,2,1,1,2,0) = 1283 | f(0,2,2,1,1,2,0) = 468 | f(1,2,2,1,1,2,0) = 369 |
| f(2,2,2,1,1,2,0) = 1091 | f(0,0,0,2,1,2,0) = 734 | f(1,0,0,2,1,2,0) = 953 |
| f(2,0,0,2,1,2,0) = 55 | f(0,1,0,2,1,2,0) = 1821 | f(1,1,0,2,1,2,0) = 1010 |
| f(2,1,0,2,1,2,0) = 1315 | f(0,2,0,2,1,2,0) = 609 | f(1,2,0,2,1,2,0) = 1009 |
| f(2,2,0,2,1,2,0) = 1123 | f(0,0,1,2,1,2,0) = 438 | f(1,0,1,2,1,2,0) = 1130 |
| f(2,0,1,2,1,2,0) = 615 | f(0,1,1,2,1,2,0) = 925 | f(1,1,1,2,1,2,0) = 904 |
| f(2,1,1,2,1,2,0) = 815 | f(0,2,1,2,1,2,0) = 1781 | f(1,2,1,2,1,2,0) = 712 |
| f(2,2,1,2,1,2,0) = 239 | f(0,0,2,2,1,2,0) = 883 | f(1,0,2,2,1,2,0) = 1561 |
| f(2,0,2,2,1,2,0) = 599 | f(0,1,2,2,1,2,0) = 861 | f(1,1,2,2,1,2,0) = 840 |
| f(2,1,2,2,1,2,0) = 799 | f(0,2,2,2,1,2,0) = 1525 | f(1,2,2,2,1,2,0) = 456 |
| f(2,2,2,2,1,2,0) = 223 | f(0,0,0,0,2,2,0) = 9 | f(1,0,0,0,2,2,0) = 573 |
| f(2,0,0,0,2,2,0) = 525 | f(0,1,0,0,2,2,0) = 1034 | f(1,1,0,0,2,2,0) = 50 |
| f(2,1,0,0,2,2,0) = 879 | f(0,2,0,0,2,2,0) = 1267 | f(1,2,0,0,2,2,0) = 49 |
| f(2,2,0,0,2,2,0) = 495 | f(0,0,1,0,2,2,0) = 29 | f(1,0,1,0,2,2,0) = 699 |
| f(2,0,1,0,2,2,0) = 1059 | f(0,1,1,0,2,2,0) = 914 | f(1,1,1,0,2,2,0) = 698 |
| f(2,1,1,0,2,2,0) = 523 | f(0,2,1,0,2,2,0) = 731 | f(1,2,1,0,2,2,0) = 697 |
| f(2,2,1,0,2,2,0) = 139 | f(0,0,2,0,2,2,0) = 18 | f(1,0,2,0,2,2,0) = 379 |
| f(2,0,2,0,2,2,0) = 1043 | f(0,1,2,0,2,2,0) = 850 | f(1,1,2,0,2,2,0) = 378 |
| f(2,1,2,0,2,2,0) = 267 | f(0,2,2,0,2,2,0) = 475 | f(1,2,2,0,2,2,0) = 377 |
| f(2,2,2,0,2,2,0) = 75 | f(0,0,0,1,2,2,0) = 356 | f(1,0,0,1,2,2,0) = 628 |
| f(2,0,0,1,2,2,0) = 517 | f(0,1,0,1,2,2,0) = 1835 | f(1,1,0,1,2,2,0) = 829 |
| f(2,1,0,1,2,2,0) = 783 | f(0,2,0,1,2,2,0) = 1057 | f(1,2,0,1,2,2,0) = 254 |
| f(2,2,0,1,2,2,0) = 207 | f(0,0,1,1,2,2,0) = 660 | f(1,0,1,1,2,2,0) = 1720 |
| f(2,0,1,1,2,2,0) = 423 | f(0,1,1,1,2,2,0) = 1548 | f(1,1,1,1,2,2,0) = 426 |

f(2,1,1,1,2,2,0) = 803
f(2,2,1,1,2,2,0) = 227
f(2,0,2,1,2,2,0) = 407
f(2,1,2,1,2,2,0) = 787
f(2,2,2,1,2,2,0) = 211
f(2,0,0,2,2,2,0) = 391
f(2,1,0,2,2,2,0) = 1894
f(2,2,0,2,2,2,0) = 1509
f(2,0,1,2,2,2,0) = 1063
f(2,1,1,2,2,2,0) = 911
f(2,2,1,2,2,2,0) = 719
f(2,0,2,2,2,2,0) = 1047
f(2,1,2,2,2,2,0) = 847
f(2,2,2,2,2,2,0) = 463
f(1,1,0,0,0,0,1) = 1470
f(2,2,0,0,0,0,1) = 1271
f(1,1,1,0,0,0,1) = 928
f(2,2,1,0,0,0,1) = 705
f(1,1,2,0,0,0,1) = 848
f(2,2,2,0,0,0,1) = 449
f(2,0,0,1,0,0,1) = 2039
f(0,2,0,1,0,0,1) = 2017
f(0,0,1,1,0,0,1) = 1740
f(0,1,1,1,0,0,1) = 1710
f(1,2,1,1,0,0,1) = 1184
f(1,0,2,1,0,0,1) = 368
f(1,1,2,1,0,0,1) = 1296
f(2,2,2,1,0,0,1) = 1239
f(1,1,0,2,0,0,1) = 2042
f(1,2,0,2,0,0,1) = 2041
f(1,0,1,2,0,0,1) = 1068
f(1,1,1,2,0,0,1) = 1980
f(1,2,1,2,0,0,1) = 1788
f(1,0,2,2,0,0,1) = 1052
f(1,1,2,2,0,0,1) = 1916
f(1,2,2,2,0,0,1) = 1532
f(2,0,0,0,1,0,1) = 1486
f(2,1,0,0,1,0,1) = 1798
f(2,2,0,0,1,0,1) = 1221
f(0,1,1,0,1,0,1) = 1556
f(0,2,1,0,1,0,1) = 1181
f(1,0,2,0,1,0,1) = 1436
f(1,1,2,0,1,0,1) = 1880
f(1,2,2,0,1,0,1) = 1496

f(0,2,1,1,2,2,0) = 722
f(0,0,2,1,2,2,0) = 639
f(0,1,2,1,2,2,0) = 1292
f(0,2,2,1,2,2,0) = 466
f(0,0,0,2,2,2,0) = 478
f(0,1,0,2,2,2,0) = 1819
f(0,2,0,2,2,2,0) = 1041
f(0,0,1,2,2,2,0) = 393
f(0,1,1,2,2,2,0) = 923
f(0,2,1,2,2,2,0) = 1205
f(0,0,2,2,2,2,0) = 374
f(0,1,2,2,2,2,0) = 859
f(0,2,2,2,2,2,0) = 1141
f(2,0,0,0,0,0,1) = 1999
f(2,1,0,0,0,0,1) = 1847
f(2,0,1,0,0,0,1) = 1711
f(2,1,1,0,0,0,1) = 898
f(2,0,2,0,0,0,1) = 1375
f(2,1,2,0,0,0,1) = 834
f(0,0,0,1,0,0,1) = 418
f(1,1,0,1,0,0,1) = 876
f(1,2,0,1,0,0,1) = 492
f(1,0,1,1,0,0,1) = 688
f(1,1,1,1,0,0,1) = 1568
f(2,2,1,1,0,0,1) = 1255
f(2,0,2,1,0,0,1) = 1431
f(2,1,2,1,0,0,1) = 1815
f(0,0,0,2,0,0,1) = 402
f(2,1,0,2,0,0,1) = 771
f(2,2,0,2,0,0,1) = 195
f(2,0,1,2,0,0,1) = 695
f(2,1,1,2,0,0,1) = 516
f(2,2,1,2,0,0,1) = 132
f(2,0,2,2,0,0,1) = 375
f(2,1,2,2,0,0,1) = 260
f(2,2,2,2,0,0,1) = 68
f(0,1,0,0,1,0,1) = 1996
f(0,2,0,0,1,0,1) = 204
f(1,0,1,0,1,0,1) = 1452
f(1,1,1,0,1,0,1) = 1960
f(1,2,1,0,1,0,1) = 1768
f(2,0,2,0,1,0,1) = 1869
f(2,1,2,0,1,0,1) = 1878
f(2,2,2,0,1,0,1) = 1493

f(1,2,1,1,2,2,0) = 617
f(1,0,2,1,2,2,0) = 1400
f(1,1,2,1,2,2,0) = 602
f(1,2,2,1,2,2,0) = 409
f(1,0,0,2,2,2,0) = 569
f(1,1,0,2,2,2,0) = 298
f(1,2,0,2,2,2,0) = 105
f(1,0,1,2,2,2,0) = 190
f(1,1,1,2,2,2,0) = 702
f(1,2,1,2,2,2,0) = 701
f(1,0,2,2,2,2,0) = 317
f(1,1,2,2,2,2,0) = 382
f(1,2,2,2,2,2,0) = 381
f(0,1,0,0,0,0,1) = 768
f(1,2,0,0,0,0,1) = 1661
f(0,1,1,0,0,0,1) = 512
f(1,2,1,0,0,0,1) = 736
f(0,1,2,0,0,0,1) = 256
f(1,2,2,0,0,0,1) = 464
f(1,0,0,1,0,0,1) = 624
f(2,1,0,1,0,0,1) = 781
f(2,2,0,1,0,0,1) = 206
f(2,0,1,1,0,0,1) = 1447
f(2,1,1,1,0,0,1) = 1831
f(0,0,2,1,0,0,1) = 110
f(0,1,2,1,0,0,1) = 1390
f(1,2,2,1,0,0,1) = 1104
f(1,0,0,2,0,0,1) = 2043
f(0,2,0,2,0,0,1) = 2001
f(0,0,1,2,0,0,1) = 1953
f(0,1,1,2,0,0,1) = 521
f(0,2,1,2,0,0,1) = 673
f(0,0,2,2,0,0,1) = 329
f(0,1,2,2,0,0,1) = 265
f(0,2,2,2,0,0,1) = 353
f(0,0,0,0,1,0,1) = 1984
f(1,1,0,0,1,0,1) = 1848
f(1,2,0,0,1,0,1) = 1272
f(2,0,1,0,1,0,1) = 1742
f(2,1,1,0,1,0,1) = 1958
f(2,2,1,0,1,0,1) = 1765
f(0,1,2,0,1,0,1) = 1300
f(0,2,2,0,1,0,1) = 1117
f(0,0,0,1,1,0,1) = 1774

| | | |
|---|---|---|
| f(1,0,0,1,1,0,1) = 1016 | f(2,0,0,1,1,0,1) = 1988 | f(0,1,0,1,1,0,1) = 1420 |
| f(1,1,0,1,1,0,1) = 1564 | f(2,1,0,1,1,0,1) = 918 | f(0,2,0,1,1,0,1) = 237 |
| f(1,2,0,1,1,0,1) = 1180 | f(2,2,0,1,1,0,1) = 725 | f(0,0,1,1,1,0,1) = 640 |
| f(1,0,1,1,1,0,1) = 1692 | f(2,0,1,1,1,0,1) = 2019 | f(0,1,1,1,1,0,1) = 1966 |
| f(1,1,1,1,1,0,1) = 1970 | f(2,1,1,1,1,0,1) = 1923 | f(0,2,1,1,1,0,1) = 1172 |
| f(1,2,1,1,1,0,1) = 1777 | f(2,2,1,1,1,0,1) = 1731 | f(0,0,2,1,1,0,1) = 610 |
| f(1,0,2,1,1,0,1) = 1388 | f(2,0,2,1,1,0,1) = 2003 | f(0,1,2,1,1,0,1) = 1902 |
| f(1,1,2,1,1,0,1) = 1906 | f(2,1,2,1,1,0,1) = 1859 | f(0,2,2,1,1,0,1) = 1108 |
| f(1,2,2,1,1,0,1) = 1521 | f(2,2,2,1,1,0,1) = 1475 | f(0,0,0,2,1,0,1) = 658 |
| f(1,0,0,2,1,0,1) = 1464 | f(2,0,0,2,1,0,1) = 11 | f(0,1,0,2,1,0,1) = 1612 |
| f(1,1,0,2,1,0,1) = 864 | f(2,1,0,2,1,0,1) = 1903 | f(0,2,0,2,1,0,1) = 221 |
| f(1,2,0,2,1,0,1) = 480 | f(2,2,0,2,1,0,1) = 1519 | f(0,0,1,2,1,0,1) = 1444 |
| f(1,0,1,2,1,0,1) = 1056 | f(2,0,1,2,1,0,1) = 549 | f(0,1,1,2,1,0,1) = 1565 |
| f(1,1,1,2,1,0,1) = 1964 | f(2,1,1,2,1,0,1) = 998 | f(0,2,1,2,1,0,1) = 737 |
| f(1,2,1,2,1,0,1) = 1772 | f(2,2,1,2,1,0,1) = 997 | f(0,0,2,2,1,0,1) = 354 |
| f(1,0,2,2,1,0,1) = 1040 | f(2,0,2,2,1,0,1) = 86 | f(0,1,2,2,1,0,1) = 1309 |
| f(1,1,2,2,1,0,1) = 1884 | f(2,1,2,2,1,0,1) = 982 | f(0,2,2,2,1,0,1) = 481 |
| f(1,2,2,2,1,0,1) = 1500 | f(2,2,2,2,1,0,1) = 981 | f(0,0,0,0,2,0,1) = 1024 |
| f(2,0,0,0,2,0,1) = 1933 | f(0,1,0,0,2,0,1) = 1994 | f(1,1,0,0,2,0,1) = 872 |
| f(2,1,0,0,2,0,1) = 779 | f(0,2,0,0,2,0,1) = 202 | f(1,2,0,0,2,0,1) = 488 |
| f(2,2,0,0,2,0,1) = 203 | f(0,0,1,0,2,0,1) = 20 | f(2,0,1,0,2,0,1) = 1703 |
| f(0,1,1,0,2,0,1) = 1554 | f(1,1,1,0,2,0,1) = 954 | f(2,1,1,0,2,0,1) = 903 |
| f(0,2,1,0,2,0,1) = 1179 | f(1,2,1,0,2,0,1) = 761 | f(2,2,1,0,2,0,1) = 711 |
| f(0,0,2,0,2,0,1) = 2011 | f(2,0,2,0,2,0,1) = 1367 | f(0,1,2,0,2,0,1) = 1298 |
| f(1,1,2,0,2,0,1) = 890 | f(2,1,2,0,2,0,1) = 839 | f(0,2,2,0,2,0,1) = 1115 |
| f(1,2,2,0,2,0,1) = 505 | f(2,2,2,0,2,0,1) = 455 | f(0,0,0,1,2,0,1) = 1518 |
| f(1,0,0,1,2,0,1) = 584 | f(2,0,0,1,2,0,1) = 1541 | f(0,1,0,1,2,0,1) = 394 |
| f(1,1,0,1,2,0,1) = 828 | f(2,1,0,1,2,0,1) = 782 | f(0,2,0,1,2,0,1) = 1259 |
| f(1,2,0,1,2,0,1) = 252 | f(2,2,0,1,2,0,1) = 205 | f(0,0,1,1,2,0,1) = 158 |
| f(1,0,1,1,2,0,1) = 416 | f(2,0,1,1,2,0,1) = 166 | f(0,1,1,1,2,0,1) = 1582 |
| f(1,1,1,1,2,0,1) = 1962 | f(2,1,1,1,2,0,1) = 1446 | f(0,2,1,1,2,0,1) = 1170 |
| f(1,2,1,1,2,0,1) = 1769 | f(2,2,1,1,2,0,1) = 1637 | f(0,0,2,1,2,0,1) = 619 |
| f(1,0,2,1,2,0,1) = 400 | f(2,0,2,1,2,0,1) = 277 | f(0,1,2,1,2,0,1) = 1326 |
| f(1,1,2,1,2,0,1) = 1882 | f(2,1,2,1,2,0,1) = 1622 | f(0,2,2,1,2,0,1) = 1106 |
| f(1,2,2,1,2,0,1) = 1497 | f(2,2,2,1,2,0,1) = 1429 | f(0,0,0,2,2,0,1) = 340 |
| f(2,0,0,2,2,0,1) = 971 | f(0,1,0,2,2,0,1) = 586 | f(1,1,0,2,2,0,1) = 1312 |
| f(2,1,0,2,2,0,1) = 1602 | f(0,2,0,2,2,0,1) = 1243 | f(1,2,0,2,2,0,1) = 1120 |
| f(2,2,0,2,2,0,1) = 1409 | f(0,0,1,2,2,0,1) = 1453 | f(1,0,1,2,2,0,1) = 1790 |
| f(2,0,1,2,2,0,1) = 1174 | f(0,1,1,2,2,0,1) = 1563 | f(1,1,1,2,2,0,1) = 952 |
| f(2,1,1,2,2,0,1) = 910 | f(0,2,1,2,2,0,1) = 161 | f(1,2,1,2,2,0,1) = 760 |
| f(2,2,1,2,2,0,1) = 717 | f(1,0,2,2,2,0,1) = 1917 | f(2,0,2,2,2,0,1) = 1317 |
| f(0,1,2,2,2,0,1) = 1307 | f(1,1,2,2,2,0,1) = 888 | f(2,1,2,2,2,0,1) = 846 |
| f(0,2,2,2,2,0,1) = 97 | f(1,2,2,2,2,0,1) = 504 | f(2,2,2,2,2,0,1) = 461 |
| f(0,0,0,0,0,1,1) = 1100 | f(1,0,0,0,0,1,1) = 1072 | f(2,0,0,0,0,1,1) = 1603 |
| f(0,1,0,0,0,1,1) = 1846 | f(1,1,0,0,0,1,1) = 912 | f(2,1,0,0,0,1,1) = 966 |
| f(0,2,0,0,0,1,1) = 1270 | f(1,2,0,0,0,1,1) = 720 | f(2,2,0,0,0,1,1) = 965 |
| f(0,0,1,0,0,1,1) = 996 | f(1,0,1,0,0,1,1) = 1696 | f(2,0,1,0,0,1,1) = 1158 |
| f(0,1,1,0,0,1,1) = 1536 | f(1,1,1,0,0,1,1) = 800 | f(2,1,1,0,0,1,1) = 642 |
| f(0,2,1,0,0,1,1) = 1161 | f(1,2,1,0,0,1,1) = 224 | f(2,2,1,0,0,1,1) = 641 |
| f(0,0,2,0,0,1,1) = 994 | f(1,0,2,0,0,1,1) = 1360 | f(2,0,2,0,0,1,1) = 1285 |
| f(0,1,2,0,0,1,1) = 1280 | f(1,1,2,0,0,1,1) = 784 | f(2,1,2,0,0,1,1) = 322 |
| f(0,2,2,0,0,1,1) = 1097 | f(1,2,2,0,0,1,1) = 208 | f(2,2,2,0,0,1,1) = 321 |
| f(0,0,0,1,0,1,1) = 1442 | f(1,0,0,1,0,1,1) = 1648 | f(2,0,0,1,0,1,1) = 1028 |
| f(0,1,0,1,0,1,1) = 1070 | f(1,1,0,1,0,1,1) = 830 | f(2,1,0,1,0,1,1) = 1558 |
| f(0,2,0,1,0,1,1) = 1249 | f(1,2,0,1,0,1,1) = 253 | f(2,2,0,1,0,1,1) = 1173 |
| f(1,0,1,1,0,1,1) = 2016 | f(2,0,1,1,0,1,1) = 2023 | f(0,1,1,1,0,1,1) = 686 |
| f(1,1,1,1,0,1,1) = 1952 | f(2,1,1,1,0,1,1) = 1686 | f(0,2,1,1,0,1,1) = 1152 |
| f(1,2,1,1,0,1,1) = 1760 | f(2,2,1,1,0,1,1) = 1685 | f(0,0,2,1,0,1,1) = 1600 |
| f(1,0,2,1,0,1,1) = 2000 | f(2,0,2,1,0,1,1) = 2007 | f(0,1,2,1,0,1,1) = 366 |
| f(1,1,2,1,0,1,1) = 1872 | f(2,1,2,1,0,1,1) = 1382 | f(0,2,2,1,0,1,1) = 1088 |
| f(1,2,2,1,0,1,1) = 1488 | f(2,2,2,1,0,1,1) = 1381 | f(0,0,0,2,0,1,1) = 1426 |
| f(2,0,0,2,0,1,1) = 1027 | f(0,1,0,2,0,1,1) = 1054 | f(1,1,0,2,0,1,1) = 1946 |
| f(2,1,0,2,0,1,1) = 1895 | f(0,2,0,2,0,1,1) = 1233 | f(1,2,0,2,0,1,1) = 1753 |
| f(2,2,0,2,0,1,1) = 1511 | f(0,0,1,2,0,1,1) = 929 | f(1,0,1,2,0,1,1) = 1644 |
| f(2,0,1,2,0,1,1) = 1750 | f(0,1,1,2,0,1,1) = 1545 | f(1,1,1,2,0,1,1) = 1596 |
| f(2,1,1,2,0,1,1) = 1540 | f(0,2,1,2,0,1,1) = 1715 | f(1,2,1,2,0,1,1) = 1212 |
| f(2,2,1,2,0,1,1) = 1156 | f(0,0,2,2,0,1,1) = 1353 | f(1,0,2,2,0,1,1) = 1628 |
| f(2,0,2,2,0,1,1) = 1893 | f(0,1,2,2,0,1,1) = 1289 | f(1,1,2,2,0,1,1) = 1340 |
| f(2,1,2,2,0,1,1) = 1284 | f(0,2,2,2,0,1,1) = 1395 | f(1,2,2,2,0,1,1) = 1148 |
| f(2,2,2,2,0,1,1) = 1092 | f(1,0,0,0,1,1,1) = 1138 | f(2,0,0,0,1,1,1) = 1090 |
| f(0,1,0,0,1,1,1) = 972 | f(1,1,0,0,1,1,1) = 1944 | f(2,1,0,0,1,1,1) = 1796 |
| f(0,2,0,0,1,1,1) = 2037 | f(1,2,0,0,1,1,1) = 1752 | f(2,2,0,0,1,1,1) = 1220 |
| f(0,0,1,0,1,1,1) = 2020 | f(1,0,1,0,1,1,1) = 2028 | f(2,0,1,0,1,1,1) = 644 |
| f(0,1,1,0,1,1,1) = 1956 | f(1,1,1,0,1,1,1) = 1976 | f(0,2,1,0,1,1,1) = 1773 |
| f(1,2,1,0,1,1,1) = 1784 | f(0,0,2,0,1,1,1) = 1482 | f(1,0,2,0,1,1,1) = 2012 |
| f(2,0,2,0,1,1,1) = 324 | f(0,1,2,0,1,1,1) = 1892 | f(1,1,2,0,1,1,1) = 1912 |
| f(0,2,2,0,1,1,1) = 1517 | f(1,2,2,0,1,1,1) = 1528 | f(1,0,0,1,1,1,1) = 2040 |
| f(2,0,0,1,1,1,1) = 1412 | f(0,1,0,1,1,1,1) = 2030 | f(1,1,0,1,1,1,1) = 1560 |
| f(2,1,0,1,1,1,1) = 772 | f(0,2,0,1,1,1,1) = 1252 | f(1,2,0,1,1,1,1) = 1176 |
| f(2,2,0,1,1,1,1) = 196 | f(0,0,1,1,1,1,1) = 1664 | f(1,0,1,1,1,1,1) = 992 |
| f(2,0,1,1,1,1,1) = 1154 | f(0,1,1,1,1,1,1) = 942 | f(1,1,1,1,1,1,1) = 1968 |
| f(2,1,1,1,1,1,1) = 1922 | f(0,2,1,1,1,1,1) = 1764 | f(1,2,1,1,1,1,1) = 1776 |
| f(2,2,1,1,1,1,1) = 1729 | f(0,0,2,1,1,1,1) = 1634 | f(1,0,2,1,1,1,1) = 976 |
| f(2,0,2,1,1,1,1) = 1281 | f(0,1,2,1,1,1,1) = 878 | f(1,1,2,1,1,1,1) = 1904 |
| f(2,1,2,1,1,1,1) = 1858 | f(0,2,2,1,1,1,1) = 1508 | f(1,2,2,1,1,1,1) = 1520 |
| f(2,2,2,1,1,1,1) = 1473 | f(0,0,0,2,1,1,1) = 1682 | f(1,0,0,2,1,1,1) = 1416 |
| f(2,0,0,2,1,1,1) = 1035 | f(0,1,0,2,1,1,1) = 2014 | f(1,1,0,2,1,1,1) = 1936 |
| f(2,1,0,2,1,1,1) = 1797 | f(0,2,0,2,1,1,1) = 1236 | f(1,2,0,2,1,1,1) = 1744 |
| f(2,2,0,2,1,1,1) = 1222 | f(1,0,1,2,1,1,1) = 1632 | f(2,0,1,2,1,1,1) = 1573 |
| f(0,1,1,2,1,1,1) = 1965 | f(1,1,1,2,1,1,1) = 1544 | f(2,1,1,2,1,1,1) = 2022 |
| f(0,2,1,2,1,1,1) = 1779 | f(1,2,1,2,1,1,1) = 1160 | f(2,2,1,2,1,1,1) = 2021 |
| f(0,0,2,2,1,1,1) = 1378 | f(1,0,2,2,1,1,1) = 1616 | f(2,0,2,2,1,1,1) = 1110 |
| f(0,1,2,2,1,1,1) = 1901 | f(1,1,2,2,1,1,1) = 1288 | f(2,1,2,2,1,1,1) = 2006 |
| f(0,2,2,2,1,1,1) = 1523 | f(1,2,2,2,1,1,1) = 1096 | f(2,2,2,2,1,1,1) = 2005 |
| f(0,0,0,0,2,1,1) = 1909 | f(1,0,0,0,2,1,1) = 1585 | f(2,0,0,0,2,1,1) = 1537 |
| f(0,1,0,0,2,1,1) = 970 | f(1,1,0,0,2,1,1) = 922 | f(2,1,0,0,2,1,1) = 773 |
| f(0,2,0,0,2,1,1) = 2035 | f(1,2,0,0,2,1,1) = 729 | f(2,2,0,0,2,1,1) = 198 |
| f(1,0,1,0,2,1,1) = 1202 | f(2,0,1,0,2,1,1) = 1667 | f(0,1,1,0,2,1,1) = 1954 |
| f(1,1,1,0,2,1,1) = 938 | f(2,1,1,0,2,1,1) = 935 | f(0,2,1,0,2,1,1) = 1771 |
| f(1,2,1,0,2,1,1) = 745 | f(2,2,1,0,2,1,1) = 743 | f(0,0,2,0,2,1,1) = 2002 |
| f(1,0,2,0,2,1,1) = 1329 | f(2,0,2,0,2,1,1) = 1347 | f(0,1,2,0,2,1,1) = 1890 |
| f(1,1,2,0,2,1,1) = 858 | f(2,1,2,0,2,1,1) = 855 | f(0,2,2,0,2,1,1) = 1515 |
| f(1,2,2,0,2,1,1) = 473 | f(2,2,2,0,2,1,1) = 471 | f(1,0,0,1,2,1,1) = 1608 |
| f(2,0,0,1,2,1,1) = 70 | f(0,1,0,1,2,1,1) = 1418 | f(1,1,0,1,2,1,1) = 920 |
| f(2,1,0,1,2,1,1) = 1422 | f(0,2,0,1,2,1,1) = 1250 | f(1,2,0,1,2,1,1) = 728 |
| f(2,2,0,1,2,1,1) = 1613 | f(0,0,1,1,2,1,1) = 1182 | f(1,0,1,1,2,1,1) = 1440 |
| f(2,0,1,1,2,1,1) = 1510 | f(0,1,1,1,2,1,1) = 558 | f(1,1,1,1,2,1,1) = 1578 |
| f(2,1,1,1,2,1,1) = 1062 | f(0,2,1,1,2,1,1) = 1762 | f(1,2,1,1,2,1,1) = 1193 |
| f(2,2,1,1,2,1,1) = 1061 | f(0,0,2,1,2,1,1) = 1643 | f(1,0,2,1,2,1,1) = 1424 |
| f(2,0,2,1,2,1,1) = 1941 | f(0,1,2,1,2,1,1) = 302 | f(1,1,2,1,2,1,1) = 1306 |
| f(2,1,2,1,2,1,1) = 1046 | f(0,2,2,1,2,1,1) = 1506 | f(1,2,2,1,2,1,1) = 1113 |
| f(2,2,2,1,2,1,1) = 1045 | f(0,0,0,2,2,1,1) = 1364 | f(1,0,0,2,2,1,1) = 1083 |
| f(2,0,0,2,2,1,1) = 1995 | f(0,1,0,2,2,1,1) = 1610 | f(1,1,0,2,2,1,1) = 1800 |
| f(2,1,0,2,2,1,1) = 962 | f(0,2,0,2,2,1,1) = 1234 | f(1,2,0,2,2,1,1) = 1224 |
| f(2,2,0,2,2,1,1) = 961 | f(1,0,1,2,2,1,1) = 1708 | f(2,0,1,2,2,1,1) = 1683 |
| f(0,1,1,2,2,1,1) = 1963 | f(1,1,1,2,2,1,1) = 948 | f(2,1,1,2,2,1,1) = 1934 |
| f(0,2,1,2,2,1,1) = 1203 | f(1,2,1,2,2,1,1) = 756 | f(2,2,1,2,2,1,1) = 1741 |
| f(0,0,2,2,2,1,1) = 1362 | f(1,0,2,2,2,1,1) = 1372 | f(2,0,2,2,2,1,1) = 1379 |
| f(0,1,2,2,2,1,1) = 1899 | f(1,1,2,2,2,1,1) = 884 | f(2,1,2,2,2,1,1) = 1870 |
| f(0,2,2,2,2,1,1) = 1139 | f(1,2,2,2,2,1,1) = 500 | f(2,2,2,2,2,1,1) = 1485 |
| f(0,0,0,0,0,2,1) = 885 | f(1,0,0,0,0,2,1) = 1084 | f(2,0,0,0,0,2,1) = 1615 |
| f(0,1,0,0,0,2,1) = 1792 | f(2,1,0,0,0,2,1) = 1805 | f(0,2,0,0,0,2,1) = 1216 |
| f(2,2,0,0,0,2,1) = 1230 | f(0,0,1,0,0,2,1) = 1589 | f(1,0,1,0,0,2,1) = 1514 |
| f(2,0,1,0,0,2,1) = 1675 | f(0,1,1,0,0,2,1) = 1590 | f(1,1,1,0,0,2,1) = 1690 |
| f(2,1,1,0,0,2,1) = 1935 | f(0,2,1,0,0,2,1) = 1215 | f(1,2,1,0,0,2,1) = 1689 |
| f(2,2,1,0,0,2,1) = 1743 | f(0,0,2,0,0,2,1) = 1587 | f(1,0,2,0,0,2,1) = 1945 |
| f(2,0,2,0,0,2,1) = 1355 | f(0,1,2,0,0,2,1) = 1334 | f(1,1,2,0,0,2,1) = 1386 |
| f(2,1,2,0,0,2,1) = 1871 | f(0,2,2,0,0,2,1) = 1151 | f(1,2,2,0,0,2,1) = 1385 |
| f(2,2,2,0,0,2,1) = 1487 | f(0,0,0,1,0,2,1) = 1636 | f(1,0,0,1,0,2,1) = 1660 |
| f(2,0,0,1,0,2,1) = 1015 | f(0,1,0,1,0,2,1) = 1825 | f(1,1,0,1,0,2,1) = 1948 |
| f(2,1,0,1,0,2,1) = 974 | f(0,2,0,1,0,2,1) = 1441 | f(1,2,0,1,0,2,1) = 1756 |
| f(2,2,0,1,0,2,1) = 973 | f(0,0,1,1,0,2,1) = 716 | f(1,0,1,1,0,2,1) = 1680 |
| f(2,0,1,1,0,2,1) = 933 | f(1,1,1,1,0,2,1) = 1595 | f(2,1,1,1,0,2,1) = 1959 |
| f(0,2,1,1,0,2,1) = 1206 | f(1,2,1,1,0,2,1) = 1211 | f(2,2,1,1,0,2,1) = 1767 |
| f(0,0,2,1,0,2,1) = 1134 | f(1,0,2,1,0,2,1) = 1376 | f(2,0,2,1,0,2,1) = 470 |
| f(1,1,2,1,0,2,1) = 1339 | f(2,1,2,1,0,2,1) = 1879 | f(0,2,2,1,0,2,1) = 1142 |
| f(1,2,2,1,0,2,1) = 1147 | f(2,2,2,1,0,2,1) = 1495 | f(0,0,0,2,0,2,1) = 1620 |
| f(1,0,0,2,0,2,1) = 1467 | f(2,0,0,2,0,2,1) = 1039 | f(0,1,0,2,0,2,1) = 1809 |
| f(1,1,0,2,0,2,1) = 1888 | f(2,1,0,2,0,2,1) = 1410 | f(0,2,0,2,0,2,1) = 1617 |
| f(1,2,0,2,0,2,1) = 1504 | f(2,2,0,2,0,2,1) = 1601 | f(1,0,1,2,0,2,1) = 746 |
| f(0,1,1,2,0,2,1) = 1599 | f(1,1,1,2,0,2,1) = 940 | f(2,1,1,2,0,2,1) = 678 |
| f(0,2,1,2,0,2,1) = 1697 | f(1,2,1,2,0,2,1) = 748 | f(2,2,1,2,0,2,1) = 677 |
| f(0,0,2,2,0,2,1) = 1873 | f(1,0,2,2,0,2,1) = 857 | f(0,1,2,2,0,2,1) = 1343 |
| f(1,1,2,2,0,2,1) = 860 | f(2,1,2,2,0,2,1) = 342 | f(0,2,2,2,0,2,1) = 1377 |
| f(1,2,2,2,0,2,1) = 476 | f(2,2,2,2,0,2,1) = 341 | f(0,0,0,0,1,2,1) = 76 |
| f(1,0,0,0,1,2,1) = 1150 | f(2,0,0,0,1,2,1) = 1102 | f(1,1,0,0,1,2,1) = 1900 |
| f(2,1,0,0,1,2,1) = 1414 | f(0,2,0,0,1,2,1) = 1228 | f(1,2,0,0,1,2,1) = 1516 |
| f(2,2,0,0,1,2,1) = 1605 | f(1,0,1,0,1,2,1) = 1704 | f(2,0,1,0,1,2,1) = 1166 |
| f(0,1,1,0,1,2,1) = 1940 | f(1,1,1,0,1,2,1) = 2026 | f(2,1,1,0,1,2,1) = 1829 |
| f(0,2,1,0,1,2,1) = 1757 | f(1,2,1,0,1,2,1) = 2025 | f(2,2,1,0,1,2,1) = 1254 |
| f(0,0,2,0,1,2,1) = 1067 | f(1,0,2,0,1,2,1) = 1368 | f(2,0,2,0,1,2,1) = 1293 |
| f(0,1,2,0,1,2,1) = 1876 | f(1,1,2,0,1,2,1) = 2010 | f(2,1,2,0,1,2,1) = 1813 |
| f(0,2,2,0,1,2,1) = 1501 | f(1,2,2,0,1,2,1) = 2009 | f(2,2,2,0,1,2,1) = 1238 |
| f(0,0,0,1,1,2,1) = 750 | f(1,0,0,1,1,2,1) = 1460 | f(2,0,0,1,1,2,1) = 1925 |
| f(0,1,0,1,1,2,1) = 1828 | f(1,1,0,1,1,2,1) = 58 | f(2,1,0,1,1,2,1) = 534 |
| f(0,2,0,1,1,2,1) = 1425 | f(1,2,0,1,1,2,1) = 57 | f(2,2,0,1,1,2,1) = 149 |
| f(0,0,1,1,1,2,1) = 1700 | f(1,0,1,1,1,2,1) = 656 | f(2,0,1,1,1,2,1) = 1443 |
| f(1,1,1,1,1,2,1) = 946 | f(2,1,1,1,1,2,1) = 1931 | f(0,2,1,1,1,2,1) = 1748 |
| f(1,2,1,1,1,2,1) = 753 | f(2,2,1,1,1,2,1) = 1739 | f(0,0,2,1,1,2,1) = 1738 |
| f(1,0,2,1,1,2,1) = 352 | f(2,0,2,1,1,2,1) = 1427 | f(1,1,2,1,1,2,1) = 882 |
| f(2,1,2,1,1,2,1) = 1867 | f(0,2,2,1,1,2,1) = 1492 | f(1,2,2,1,1,2,1) = 497 |
| f(2,2,2,1,1,2,1) = 1483 | f(0,0,0,2,1,2,1) = 1758 | f(1,0,0,2,1,2,1) = 440 |
| f(2,0,0,2,1,2,1) = 1079 | f(0,1,0,2,1,2,1) = 1812 | f(1,1,0,2,1,2,1) = 818 |
| f(2,1,0,2,1,2,1) = 1803 | f(0,2,0,2,1,2,1) = 1633 | f(1,2,0,2,1,2,1) = 241 |
| f(2,2,0,2,1,2,1) = 1227 | f(0,0,1,2,1,2,1) = 420 | f(1,0,1,2,1,2,1) = 1194 |
| f(2,0,1,2,1,2,1) = 726 | f(0,1,1,2,1,2,1) = 1949 | f(1,1,1,2,1,2,1) = 1928 |
| f(2,1,1,2,1,2,1) = 805 | f(0,2,1,2,1,2,1) = 1761 | f(1,2,1,2,1,2,1) = 1736 |
| f(2,2,1,2,1,2,1) = 230 | f(0,0,2,2,1,2,1) = 1907 | f(1,0,2,2,1,2,1) = 1305 |
| f(2,0,2,2,1,2,1) = 869 | f(0,1,2,2,1,2,1) = 1885 | f(1,1,2,2,1,2,1) = 1864 |
| f(2,1,2,2,1,2,1) = 789 | f(0,2,2,2,1,2,1) = 1505 | f(1,2,2,2,1,2,1) = 1480 |
| f(2,2,2,2,1,2,1) = 214 | f(0,0,0,0,2,2,1) = 1033 | f(1,0,0,0,2,2,1) = 1597 |
| f(2,0,0,0,2,2,1) = 1549 | f(1,1,0,0,2,2,1) = 626 | f(2,1,0,0,2,2,1) = 870 |
| f(0,2,0,0,2,2,1) = 1226 | f(1,2,0,0,2,2,1) = 433 | f(2,2,0,0,2,2,1) = 485 |
| f(0,0,1,0,2,2,1) = 989 | f(1,0,1,0,2,2,1) = 186 | f(2,0,1,0,2,2,1) = 1635 |
| f(0,1,1,0,2,2,1) = 1938 | f(1,1,1,0,2,2,1) = 570 | f(2,1,1,0,2,2,1) = 907 |
| f(0,2,1,0,2,2,1) = 1755 | f(1,2,1,0,2,2,1) = 185 | f(2,2,1,0,2,2,1) = 715 |
| f(0,0,2,0,2,2,1) = 563 | f(1,0,2,0,2,2,1) = 313 | f(2,0,2,0,2,2,1) = 1619 |
| f(0,1,2,0,2,2,1) = 1874 | f(1,1,2,0,2,2,1) = 314 | f(2,1,2,0,2,2,1) = 843 |
| f(0,2,2,0,2,2,1) = 1499 | f(1,2,2,0,2,2,1) = 121 | f(2,2,2,0,2,2,1) = 459 |
| f(0,0,0,1,2,2,1) = 494 | f(1,0,0,1,2,2,1) = 1652 | f(0,1,0,1,2,2,1) = 802 |
| f(1,1,0,1,2,2,1) = 1853 | f(2,1,0,1,2,2,1) = 590 | f(0,2,0,1,2,2,1) = 1459 |
| f(1,2,0,1,2,2,1) = 1278 | f(2,2,0,1,2,2,1) = 397 | f(0,0,1,1,2,2,1) = 1684 |
| f(1,0,1,1,2,2,1) = 1712 | f(1,1,1,1,2,2,1) = 1450 | f(2,1,1,1,2,2,1) = 931 |
| f(0,2,1,1,2,2,1) = 1746 | f(1,2,1,1,2,2,1) = 1641 | f(2,2,1,1,2,2,1) = 739 |
| f(0,0,2,1,2,2,1) = 1663 | f(1,0,2,1,2,2,1) = 1392 | f(1,1,2,1,2,2,1) = 1626 |
| f(2,1,2,1,2,2,1) = 851 | f(0,2,2,1,2,2,1) = 1490 | f(1,2,2,1,2,2,1) = 1433 |
| f(2,2,2,1,2,2,1) = 467 | f(0,0,0,2,2,2,1) = 1502 | f(1,0,0,2,2,2,1) = 1593 |
| f(2,0,0,2,2,2,1) = 1415 | f(0,1,0,2,2,2,1) = 786 | f(1,1,0,2,2,2,1) = 1466 |
| f(2,1,0,2,2,2,1) = 1606 | f(0,2,0,2,2,2,1) = 1651 | f(1,2,0,2,2,2,1) = 1657 |
| f(2,2,0,2,2,2,1) = 1413 | f(0,0,1,2,2,2,1) = 429 | f(1,0,1,2,2,2,1) = 1214 |
| f(2,0,1,2,2,2,1) = 150 | f(0,1,1,2,2,2,1) = 1947 | f(1,1,1,2,2,2,1) = 568 |
| f(2,1,1,2,2,2,1) = 646 | f(0,2,1,2,2,2,1) = 1185 | f(1,2,1,2,2,2,1) = 184 |
| f(2,2,1,2,2,2,1) = 645 | f(0,0,2,2,2,2,1) = 1398 | f(1,0,2,2,2,2,1) = 1341 |
| f(2,0,2,2,2,2,1) = 293 | f(0,1,2,2,2,2,1) = 1883 | f(1,1,2,2,2,2,1) = 312 |
| f(2,1,2,2,2,2,1) = 326 | f(0,2,2,2,2,2,1) = 1121 | f(1,2,2,2,2,2,1) = 120 |
| f(2,2,2,2,2,2,1) = 325 | f(1,0,0,0,0,0,2) = 48 | f(1,1,0,0,0,0,2) = 776 |
| f(2,1,0,0,0,0,2) = 386 | f(0,2,0,0,0,0,2) = 1279 | f(1,2,0,0,0,0,2) = 200 |
| f(2,2,0,0,0,0,2) = 577 | f(1,0,1,0,0,0,2) = 672 | f(1,1,1,0,0,0,2) = 1598 |
| f(2,1,1,0,0,0,2) = 1583 | f(0,2,1,0,0,0,2) = 1791 | f(1,2,1,0,0,0,2) = 1213 |
| f(2,2,1,0,0,0,2) = 1199 | f(1,0,2,0,0,0,2) = 336 | f(1,1,2,0,0,0,2) = 1342 |
| f(2,1,2,0,0,0,2) = 1311 | f(0,2,2,0,0,0,2) = 1535 | f(1,2,2,0,0,0,2) = 1149 |
| f(2,2,2,0,0,0,2) = 1119 | f(0,0,0,1,0,0,2) = 1645 | f(2,0,0,1,0,0,2) = 4 |
| f(0,1,0,1,0,0,2) = 46 | f(1,1,0,1,0,0,2) = 1852 | f(2,1,0,1,0,0,2) = 6 |
| f(1,2,0,1,0,0,2) = 1276 | f(2,2,0,1,0,0,2) = 5 | f(0,0,1,1,0,0,2) = 1718 |
| f(1,0,1,1,0,0,2) = 1672 | f(2,0,1,1,0,0,2) = 995 | f(0,1,1,1,0,0,2) = 1694 |
| f(1,1,1,1,0,0,2) = 1979 | f(2,1,1,1,0,0,2) = 515 | f(0,2,1,1,0,0,2) = 1782 |
| f(1,2,1,1,0,0,2) = 1787 | f(2,2,1,1,0,0,2) = 131 | f(0,0,2,1,0,0,2) = 94 |
| f(1,0,2,1,0,0,2) = 1352 | f(2,0,2,1,0,0,2) = 979 | f(0,1,2,1,0,0,2) = 1374 |
| f(1,1,2,1,0,0,2) = 1915 | f(2,1,2,1,0,0,2) = 259 | f(0,2,2,1,0,0,2) = 1526 |
| f(1,2,2,1,0,0,2) = 1531 | f(2,2,2,1,0,0,2) = 67 | f(0,0,0,2,0,0,2) = 1629 |
| f(1,0,0,2,0,0,2) = 8 | f(2,0,0,2,0,0,2) = 1423 | f(0,1,0,2,0,0,2) = 30 |
| f(1,1,0,2,0,0,2) = 1841 | f(2,1,0,2,0,0,2) = 1555 | f(1,2,0,2,0,0,2) = 1266 |
| f(2,2,0,2,0,0,2) = 1171 | f(0,0,1,2,0,0,2) = 1937 | f(1,0,1,2,0,0,2) = 616 |
| f(2,0,1,2,0,0,2) = 1679 | f(1,1,1,2,0,0,2) = 808 | f(2,1,1,2,0,0,2) = 943 |
| f(0,2,1,2,0,0,2) = 657 | f(1,2,1,2,0,0,2) = 232 | f(2,2,1,2,0,0,2) = 751 |
| f(0,0,2,2,0,0,2) = 307 | f(1,0,2,2,0,0,2) = 600 | f(2,0,2,2,0,0,2) = 1359 |
| f(1,1,2,2,0,0,2) = 792 | f(2,1,2,2,0,0,2) = 863 | f(0,2,2,2,0,0,2) = 337 |
| f(1,2,2,2,0,0,2) = 216 | f(2,2,2,2,0,0,2) = 479 | f(0,0,0,0,1,0,2) = 1023 |
| f(1,0,0,0,1,0,2) = 114 | f(0,1,0,0,1,0,2) = 1845 | f(1,1,0,0,1,0,2) = 1844 |
| f(2,1,0,0,1,0,2) = 1559 | f(0,2,0,0,1,0,2) = 53 | f(1,2,0,0,1,0,2) = 1268 |

f(2,2,0,0,1,0,2) = 1175
f(0,0,1,0,1,0,2) = 36
f(1,0,1,0,1,0,2) = 680
f(0,1,1,0,1,0,2) = 932
f(1,1,1,0,1,0,2) = 1592
f(2,1,1,0,1,0,2) = 1542
f(0,2,1,0,1,0,2) = 749
f(1,2,1,0,1,0,2) = 1208
f(2,2,1,0,1,0,2) = 1157
f(0,0,2,0,1,0,2) = 2027
f(1,0,2,0,1,0,2) = 344
f(0,1,2,0,1,0,2) = 868
f(1,1,2,0,1,0,2) = 1336
f(2,1,2,0,1,0,2) = 1286
f(0,2,2,0,1,0,2) = 493
f(1,2,2,0,1,0,2) = 1144
f(2,2,2,0,1,0,2) = 1093
f(0,0,0,1,1,0,2) = 1707
f(1,0,0,1,1,0,2) = 1076
f(0,1,0,1,1,0,2) = 804
f(1,1,0,1,1,0,2) = 638
f(2,1,0,1,1,0,2) = 927
f(0,2,0,1,1,0,2) = 1461
f(1,2,0,1,1,0,2) = 445
f(2,2,0,1,1,0,2) = 735
f(1,0,1,1,1,0,2) = 730
f(2,0,1,1,1,0,2) = 130
f(0,1,1,1,1,0,2) = 1950
f(1,1,1,1,1,0,2) = 1586
f(2,1,1,1,1,0,2) = 1543
f(0,2,1,1,1,0,2) = 740
f(1,2,1,1,1,0,2) = 1201
f(2,2,1,1,1,0,2) = 1159
f(0,0,2,1,1,0,2) = 594
f(1,0,2,1,1,0,2) = 873
f(2,0,2,1,1,0,2) = 257
f(0,1,2,1,1,0,2) = 1886
f(1,1,2,1,1,0,2) = 1330
f(2,1,2,1,1,0,2) = 1287
f(0,2,2,1,1,0,2) = 484
f(1,2,2,1,1,0,2) = 1137
f(2,2,2,1,1,0,2) = 1095
f(0,0,0,2,1,0,2) = 529
f(1,0,0,2,1,0,2) = 506
f(2,0,0,2,1,0,2) = 1463
f(0,1,0,2,1,0,2) = 788
f(1,1,0,2,1,0,2) = 1842
f(2,1,0,2,1,0,2) = 1795
f(0,2,0,2,1,0,2) = 1653
f(1,2,0,2,1,0,2) = 1265
f(2,2,0,2,1,0,2) = 1219
f(0,0,1,2,1,0,2) = 1428
f(1,0,1,2,1,0,2) = 1770
f(2,0,1,2,1,0,2) = 1647
f(0,1,1,2,1,0,2) = 941
f(1,1,1,2,1,0,2) = 618
f(2,1,1,2,1,0,2) = 550
f(0,2,1,2,1,0,2) = 721
f(1,2,1,2,1,0,2) = 425
f(2,2,1,2,1,0,2) = 165
f(0,0,2,2,1,0,2) = 1889
f(1,0,2,2,1,0,2) = 1881
f(2,0,2,2,1,0,2) = 1631
f(0,1,2,2,1,0,2) = 877
f(1,1,2,2,1,0,2) = 410
f(2,1,2,2,1,0,2) = 278
f(0,2,2,2,1,0,2) = 465
f(1,2,2,2,1,0,2) = 601
f(2,2,2,2,1,0,2) = 85
f(0,0,0,0,2,0,2) = 63
f(1,0,0,0,2,0,2) = 561
f(0,1,0,0,2,0,2) = 1843
f(1,1,0,0,2,0,2) = 826
f(2,1,0,0,2,0,2) = 775
f(0,2,0,0,2,0,2) = 51
f(1,2,0,0,2,0,2) = 249
f(2,2,0,0,2,0,2) = 199
f(1,0,1,0,2,0,2) = 178
f(2,0,1,0,2,0,2) = 611
f(0,1,1,0,2,0,2) = 930
f(1,1,1,0,2,0,2) = 554
f(2,1,1,0,2,0,2) = 551
f(0,2,1,0,2,0,2) = 747
f(1,2,1,0,2,0,2) = 169
f(2,2,1,0,2,0,2) = 167
f(1,0,2,0,2,0,2) = 305
f(2,0,2,0,2,0,2) = 595
f(0,1,2,0,2,0,2) = 866
f(1,1,2,0,2,0,2) = 282
f(2,1,2,0,2,0,2) = 279
f(0,2,2,0,2,0,2) = 491
f(1,2,2,0,2,0,2) = 89
f(2,2,2,0,2,0,2) = 87
f(0,0,0,1,2,0,2) = 1389
f(1,0,0,1,2,0,2) = 2036
f(2,0,0,1,2,0,2) = 583
f(0,1,0,1,2,0,2) = 1826
f(1,1,0,1,2,0,2) = 528
f(2,1,0,1,2,0,2) = 1567
f(0,2,0,1,2,0,2) = 435
f(1,2,0,1,2,0,2) = 144
f(2,2,0,1,2,0,2) = 1183
f(0,0,1,1,2,0,2) = 1693
f(1,0,1,1,2,0,2) = 1961
f(2,0,1,1,2,0,2) = 1007
f(0,1,1,1,2,0,2) = 1566
f(1,1,1,1,2,0,2) = 1066
f(2,1,1,1,2,0,2) = 547
f(0,2,1,1,2,0,2) = 738
f(1,2,1,1,2,0,2) = 1065
f(2,2,1,1,2,0,2) = 163
f(0,0,2,1,2,0,2) = 603
f(1,0,2,1,2,0,2) = 1498
f(2,0,2,1,2,0,2) = 991
f(0,1,2,1,2,0,2) = 1310
f(1,1,2,1,2,0,2) = 1050
f(2,1,2,1,2,0,2) = 275
f(0,2,2,1,2,0,2) = 482
f(1,2,2,1,2,0,2) = 1049
f(2,2,2,1,2,0,2) = 83
f(0,0,0,2,2,0,2) = 273
f(1,0,0,2,2,0,2) = 59
f(2,0,0,2,2,0,2) = 1031
f(0,1,0,2,2,0,2) = 1810
f(1,1,0,2,2,0,2) = 1322
f(2,1,0,2,2,0,2) = 867
f(0,2,0,2,2,0,2) = 627
f(1,2,0,2,2,0,2) = 1129
f(2,2,0,2,2,0,2) = 483
f(0,0,1,2,2,0,2) = 1437
f(1,0,1,2,2,0,2) = 44
f(2,0,1,2,2,0,2) = 659

f(0,1,1,2,2,0,2) = 939        f(1,1,1,2,2,0,2) = 572        f(2,1,1,2,2,0,2) = 526
f(0,2,1,2,2,0,2) = 145        f(1,2,1,2,2,0,2) = 188        f(2,2,1,2,2,0,2) = 141
f(0,0,2,2,2,0,2) = 1407       f(1,0,2,2,2,0,2) = 28         f(2,0,2,2,2,0,2) = 355
f(0,1,2,2,2,0,2) = 875        f(1,1,2,2,2,0,2) = 316        f(2,1,2,2,2,0,2) = 270
f(0,2,2,2,2,0,2) = 81         f(1,2,2,2,2,0,2) = 124        f(2,2,2,2,2,0,2) = 77
f(0,0,0,0,0,1,2) = 1162       f(1,0,0,0,0,1,2) = 432        f(2,0,0,0,0,1,2) = 963
f(0,1,0,0,0,1,2) = 831        f(1,1,0,0,0,1,2) = 817        f(0,2,0,0,0,1,2) = 255
f(1,2,0,0,0,1,2) = 242        f(0,0,1,0,0,1,2) = 460        f(1,0,1,0,0,1,2) = 692
f(2,0,1,0,0,1,2) = 102        f(0,1,1,0,0,1,2) = 896        f(1,1,1,0,0,1,2) = 560
f(2,1,1,0,0,1,2) = 662        f(0,2,1,0,0,1,2) = 713        f(1,2,1,0,0,1,2) = 176
f(2,2,1,0,0,1,2) = 661        f(0,0,2,0,0,1,2) = 458        f(1,0,2,0,0,1,2) = 372
f(2,0,2,0,0,1,2) = 533        f(0,1,2,0,0,1,2) = 832        f(1,1,2,0,0,1,2) = 304
f(2,1,2,0,0,1,2) = 358        f(0,2,2,0,0,1,2) = 457        f(1,2,2,0,0,1,2) = 112
f(2,2,2,0,0,1,2) = 357        f(0,0,0,1,0,1,2) = 427        f(1,0,0,1,0,1,2) = 1008
f(2,0,0,1,0,1,2) = 580        f(0,1,0,1,0,1,2) = 430        f(1,1,0,1,0,1,2) = 446
f(2,1,0,1,0,1,2) = 543        f(0,2,0,1,0,1,2) = 238        f(1,2,0,1,0,1,2) = 637
f(2,2,0,1,0,1,2) = 159        f(0,0,1,1,0,1,2) = 174        f(2,0,1,1,0,1,2) = 1190
f(0,1,1,1,0,1,2) = 670        f(1,1,1,1,0,1,2) = 1706       f(2,1,1,1,0,1,2) = 1571
f(0,2,1,1,0,1,2) = 704        f(1,2,1,1,0,1,2) = 1705       f(2,2,1,1,0,1,2) = 1187
f(2,0,2,1,0,1,2) = 1301       f(0,1,2,1,0,1,2) = 350        f(1,1,2,1,0,1,2) = 1370
f(2,1,2,1,0,1,2) = 1299       f(0,2,2,1,0,1,2) = 448        f(1,2,2,1,0,1,2) = 1369
f(2,2,2,1,0,1,2) = 1107       f(0,0,0,2,0,1,2) = 411        f(1,0,0,2,0,1,2) = 1032
f(2,0,0,2,0,1,2) = 387        f(0,1,0,2,0,1,2) = 606        f(1,1,0,2,0,1,2) = 1074
f(2,1,0,2,0,1,2) = 291        f(0,2,0,2,0,1,2) = 222        f(1,2,0,2,0,1,2) = 1073
f(2,2,0,2,0,1,2) = 99         f(0,0,1,2,0,1,2) = 913        f(1,0,1,2,0,1,2) = 1577
f(2,0,1,2,0,1,2) = 671        f(0,1,1,2,0,1,2) = 905        f(1,1,1,2,0,1,2) = 552
f(2,1,1,2,0,1,2) = 900        f(1,2,1,2,0,1,2) = 168        f(2,2,1,2,0,1,2) = 708
f(0,0,2,2,0,1,2) = 1331       f(1,0,2,2,0,1,2) = 1114       f(2,0,2,2,0,1,2) = 367
f(0,1,2,2,0,1,2) = 841        f(1,1,2,2,0,1,2) = 280        f(2,1,2,2,0,1,2) = 836
f(1,2,2,2,0,1,2) = 88         f(2,2,2,2,0,1,2) = 452        f(0,0,0,0,1,1,2) = 1014
f(1,0,0,0,1,1,2) = 498        f(2,0,0,0,1,1,2) = 450        f(0,1,0,0,1,1,2) = 821
f(1,1,0,0,1,1,2) = 1562       f(2,1,0,0,1,1,2) = 1614       f(1,2,0,0,1,1,2) = 1177
f(2,2,0,0,1,1,2) = 1421       f(0,0,1,0,1,1,2) = 1484       f(1,0,1,0,1,1,2) = 428
f(2,0,1,0,1,1,2) = 1734       f(0,1,1,0,1,1,2) = 548        f(1,1,1,0,1,1,2) = 1588
f(2,1,1,0,1,1,2) = 1926       f(0,2,1,0,1,1,2) = 173        f(1,2,1,0,1,1,2) = 1204
f(2,2,1,0,1,1,2) = 1733       f(0,0,2,0,1,1,2) = 1058       f(1,0,2,0,1,1,2) = 412
f(2,0,2,0,1,1,2) = 1861       f(0,1,2,0,1,1,2) = 292        f(1,1,2,0,1,1,2) = 1332
f(2,1,2,0,1,1,2) = 1862       f(0,2,2,0,1,1,2) = 109        f(1,2,2,0,1,1,2) = 1140
f(2,2,2,0,1,1,2) = 1477       f(0,0,0,1,1,1,2) = 545        f(1,0,0,1,1,1,2) = 632
f(2,0,0,1,1,1,2) = 454        f(0,1,0,1,1,1,2) = 396        f(1,1,0,1,1,1,2) = 634
f(2,1,0,1,1,1,2) = 390        f(0,2,0,1,1,1,2) = 1261       f(1,2,0,1,1,1,2) = 441
f(2,2,0,1,1,1,2) = 581        f(0,0,1,1,1,1,2) = 649        f(1,0,1,1,1,1,2) = 1754
f(2,0,1,1,1,1,2) = 706        f(0,1,1,1,1,1,2) = 926        f(1,1,1,1,1,1,2) = 1722
f(2,1,1,1,1,1,2) = 1927       f(0,2,1,1,1,1,2) = 164        f(1,2,1,1,1,1,2) = 1721

f(2,2,1,1,1,2) = 1735
f(2,0,2,1,1,2) = 833
f(2,1,2,1,1,2) = 1863
f(2,2,2,1,1,2) = 1479
f(0,1,0,2,1,1,2) = 588
f(0,2,0,2,1,1,2) = 1245
f(0,0,1,2,1,1,2) = 384
f(1,1,1,2,1,1,2) = 1580
f(2,2,1,2,1,1,2) = 421
f(0,1,2,2,1,1,2) = 301
f(1,2,2,2,1,1,2) = 1116
f(1,0,0,0,2,1,2) = 945
f(1,1,0,0,2,1,2) = 442
f(2,2,0,0,2,1,2) = 147
f(2,0,1,0,2,1,2) = 679
f(2,1,1,0,2,1,2) = 38
f(2,2,1,0,2,1,2) = 37
f(0,1,2,0,2,1,2) = 290
f(0,2,2,0,2,1,2) = 107
f(0,0,0,1,2,1,2) = 289
f(0,1,0,1,2,1,2) = 414
f(0,2,0,1,2,1,2) = 235
f(0,0,1,1,2,1,2) = 140
f(0,1,1,1,2,1,2) = 542
f(0,2,1,1,2,1,2) = 162
f(0,0,2,1,2,1,2) = 1627
f(0,1,2,1,2,1,2) = 286
f(0,2,2,1,2,1,2) = 98
f(0,0,0,2,2,1,2) = 1297
f(0,1,0,2,2,1,2) = 622
f(0,2,0,2,2,1,2) = 219
f(0,0,1,2,2,1,2) = 309
f(0,1,1,2,2,1,2) = 555
f(1,2,1,2,2,1,2) = 180
f(1,0,2,2,2,1,2) = 604
f(1,1,2,2,2,1,2) = 308
f(2,2,2,2,2,1,2) = 1101
f(2,0,0,0,0,2,2) = 975
f(2,1,0,0,0,2,2) = 1327
f(2,2,0,0,0,2,2) = 1135
f(2,0,1,0,0,2,2) = 687
f(2,1,1,0,0,2,2) = 1839
f(2,2,1,0,0,2,2) = 1263
f(2,0,2,0,0,2,2) = 351

f(0,0,2,1,1,1,2) = 1618
f(0,1,2,1,1,1,2) = 862
f(0,2,2,1,1,1,2) = 100
f(0,0,0,2,1,1,2) = 1553
f(1,1,0,2,1,1,2) = 1650
f(1,2,0,2,1,1,2) = 1457
f(2,0,1,2,1,1,2) = 655
f(2,1,1,2,1,1,2) = 614
f(0,0,2,2,1,1,2) = 363
f(1,1,2,2,1,1,2) = 1308
f(2,2,2,2,1,1,2) = 597
f(2,0,0,0,2,1,2) = 897
f(2,1,0,0,2,1,2) = 531
f(0,0,1,0,2,1,2) = 980
f(0,1,1,0,2,1,2) = 546
f(0,2,1,0,2,1,2) = 171
f(1,0,2,0,2,1,2) = 881
f(1,1,2,0,2,1,2) = 793
f(1,2,2,0,2,1,2) = 218
f(1,0,0,1,2,1,2) = 968
f(1,1,0,1,2,1,2) = 820
f(1,2,0,1,2,1,2) = 244
f(1,0,1,1,2,1,2) = 1178
f(1,1,1,1,2,1,2) = 1833
f(1,2,1,1,2,1,2) = 1258
f(1,0,2,1,2,1,2) = 1321
f(1,1,2,1,2,1,2) = 1817
f(1,2,2,1,2,1,2) = 1242
f(1,0,0,2,2,1,2) = 122
f(1,1,0,2,2,1,2) = 1898
f(1,2,0,2,2,1,2) = 1513
f(1,0,1,2,2,1,2) = 620
f(1,1,1,2,2,1,2) = 564
f(2,2,1,2,2,1,2) = 1165
f(2,0,2,2,2,1,2) = 1391
f(2,1,2,2,2,1,2) = 1294
f(0,0,0,0,2,2) = 947
f(0,1,0,0,0,2,2) = 777
f(0,2,0,0,0,2,2) = 201
f(0,0,1,0,0,2,2) = 1053
f(0,1,1,0,0,2,2) = 950
f(0,2,1,0,0,2,2) = 767
f(0,0,2,0,0,2,2) = 1051
f(0,1,2,0,0,2,2) = 886

f(1,0,2,1,1,1,2) = 1897
f(1,1,2,1,1,1,2) = 1402
f(1,2,2,1,1,1,2) = 1401
f(2,0,0,2,1,1,2) = 395
f(2,1,0,2,1,1,2) = 769
f(2,2,0,2,1,1,2) = 194
f(0,1,1,2,1,1,2) = 557
f(1,2,1,2,1,1,2) = 1196
f(2,0,2,2,1,1,2) = 335
f(2,1,2,2,1,1,2) = 406
f(0,0,0,0,2,1,2) = 1971
f(0,1,0,0,2,1,2) = 819
f(1,2,0,0,2,1,2) = 633
f(1,0,1,0,2,1,2) = 754
f(1,1,1,0,2,1,2) = 809
f(1,2,1,0,2,1,2) = 234
f(2,0,2,0,2,1,2) = 343
f(2,1,2,0,2,1,2) = 22
f(2,2,2,0,2,1,2) = 21
f(2,0,0,1,2,1,2) = 1607
f(2,1,0,1,2,1,2) = 1806
f(2,2,0,1,2,1,2) = 1229
f(2,0,1,1,2,1,2) = 742
f(2,1,1,1,2,1,2) = 567
f(2,2,1,1,2,1,2) = 183
f(2,0,2,1,2,1,2) = 853
f(2,1,2,1,2,1,2) = 311
f(2,2,2,1,2,1,2) = 119
f(2,0,0,2,2,1,2) = 587
f(2,1,0,2,2,1,2) = 1990
f(2,2,0,2,2,1,2) = 1989
f(2,0,1,2,2,1,2) = 1695
f(2,1,1,2,2,1,2) = 1550
f(0,0,2,2,2,1,2) = 347
f(0,1,2,2,2,1,2) = 299
f(1,2,2,2,2,1,2) = 116
f(1,0,0,0,0,2,2) = 444
f(1,1,0,0,0,2,2) = 1082
f(1,2,0,0,0,2,2) = 1081
f(1,0,1,0,0,2,2) = 762
f(1,1,1,0,0,2,2) = 1726
f(1,2,1,0,0,2,2) = 1725
f(1,0,2,0,0,2,2) = 889
f(1,1,2,0,0,2,2) = 1406

| | | |
|---|---|---|
| f(2,1,2,0,0,2,2) = 1823 | f(0,2,2,0,0,2,2) = 511 | f(1,2,2,0,0,2,2) = 1405 |
| f(2,2,2,0,0,2,2) = 1247 | f(0,0,0,1,0,2,2) = 621 | f(1,0,0,1,0,2,2) = 1020 |
| f(0,1,0,1,0,2,2) = 814 | f(1,1,0,1,0,2,2) = 536 | f(2,1,0,1,0,2,2) = 294 |
| f(0,2,0,1,0,2,2) = 993 | f(1,2,0,1,0,2,2) = 152 | f(2,2,0,1,0,2,2) = 101 |
| f(0,0,1,1,0,2,2) = 694 | f(1,0,1,1,0,2,2) = 154 | f(2,0,1,1,0,2,2) = 419 |
| f(0,1,1,1,0,2,2) = 652 | f(1,1,1,1,0,2,2) = 955 | f(2,1,1,1,0,2,2) = 899 |
| f(0,2,1,1,0,2,2) = 758 | f(1,2,1,1,0,2,2) = 763 | f(2,2,1,1,0,2,2) = 707 |
| f(0,0,2,1,0,2,2) = 1118 | f(1,0,2,1,0,2,2) = 297 | f(2,0,2,1,0,2,2) = 403 |
| f(0,1,2,1,0,2,2) = 332 | f(1,1,2,1,0,2,2) = 891 | f(2,1,2,1,0,2,2) = 835 |
| f(0,2,2,1,0,2,2) = 502 | f(1,2,2,1,0,2,2) = 507 | f(2,2,2,1,0,2,2) = 451 |
| f(0,0,0,2,0,2,2) = 605 | f(1,0,0,2,0,2,2) = 1019 | f(2,0,0,2,0,2,2) = 399 |
| f(0,1,0,2,0,2,2) = 798 | f(1,1,0,2,0,2,2) = 874 | f(2,1,0,2,0,2,2) = 1794 |
| f(0,2,0,2,0,2,2) = 977 | f(1,2,0,2,0,2,2) = 489 | f(2,2,0,2,0,2,2) = 1217 |
| f(0,0,1,2,0,2,2) = 447 | f(1,0,1,2,0,2,2) = 40 | f(2,0,1,2,0,2,2) = 47 |
| f(0,1,1,2,0,2,2) = 959 | f(1,1,1,2,0,2,2) = 666 | f(2,1,1,2,0,2,2) = 559 |
| f(0,2,1,2,0,2,2) = 1681 | f(1,2,1,2,0,2,2) = 665 | f(2,2,1,2,0,2,2) = 175 |
| f(1,0,2,2,0,2,2) = 24 | f(2,0,2,2,0,2,2) = 31 | f(0,1,2,2,0,2,2) = 895 |
| f(1,1,2,2,0,2,2) = 362 | f(2,1,2,2,0,2,2) = 287 | f(0,2,2,2,0,2,2) = 1361 |
| f(1,2,2,2,0,2,2) = 361 | f(2,2,2,2,0,2,2) = 95 | f(0,0,0,0,1,2,2) = 138 |
| f(1,0,0,0,1,2,2) = 510 | f(2,0,0,0,1,2,2) = 462 | f(0,1,0,0,1,2,2) = 12 |
| f(1,1,0,0,1,2,2) = 1849 | f(2,1,0,0,1,2,2) = 1318 | f(0,2,0,0,1,2,2) = 1077 |
| f(1,2,0,0,1,2,2) = 1274 | f(2,2,0,0,1,2,2) = 1125 | f(0,0,1,0,1,2,2) = 45 |
| f(1,0,1,0,1,2,2) = 700 | f(2,0,1,0,1,2,2) = 718 | f(0,1,1,0,1,2,2) = 532 |
| f(1,1,1,0,1,2,2) = 1576 | f(2,1,1,0,1,2,2) = 1574 | f(0,2,1,0,1,2,2) = 157 |
| f(1,2,1,0,1,2,2) = 1192 | f(2,2,1,0,1,2,2) = 1189 | f(1,0,2,0,1,2,2) = 380 |
| f(2,0,2,0,1,2,2) = 845 | f(0,1,2,0,1,2,2) = 276 | f(1,1,2,0,1,2,2) = 1304 |
| f(2,1,2,0,1,2,2) = 1302 | f(0,2,2,0,1,2,2) = 93 | f(1,2,2,0,1,2,2) = 1112 |
| f(2,2,2,0,1,2,2) = 1109 | f(0,0,0,1,1,2,2) = 683 | f(1,0,0,1,1,2,2) = 52 |
| f(2,0,0,1,1,2,2) = 964 | f(0,1,0,1,1,2,2) = 813 | f(1,1,0,1,1,2,2) = 1086 |
| f(2,1,0,1,1,2,2) = 823 | f(0,2,0,1,1,2,2) = 437 | f(1,2,0,1,1,2,2) = 1085 |
| f(2,2,0,1,1,2,2) = 247 | f(0,0,1,1,1,2,2) = 685 | f(1,0,1,1,1,2,2) = 668 |
| f(2,0,1,1,1,2,2) = 675 | f(0,1,1,1,1,2,2) = 908 | f(1,1,1,1,1,2,2) = 562 |
| f(2,1,1,1,1,2,2) = 1547 | f(0,2,1,1,1,2,2) = 148 | f(1,2,1,1,1,2,2) = 177 |
| f(2,2,1,1,1,2,2) = 1163 | f(1,0,2,1,1,2,2) = 364 | f(2,0,2,1,1,2,2) = 339 |
| f(0,1,2,1,1,2,2) = 844 | f(1,1,2,1,1,2,2) = 306 | f(2,1,2,1,1,2,2) = 1291 |
| f(0,2,2,1,1,2,2) = 84 | f(1,2,2,1,1,2,2) = 113 | f(2,2,2,1,1,2,2) = 1099 |
| f(1,0,0,2,1,2,2) = 1977 | f(2,0,0,2,1,2,2) = 439 | f(0,1,0,2,1,2,2) = 797 |
| f(1,1,0,2,1,2,2) = 434 | f(2,1,0,2,1,2,2) = 1319 | f(0,2,0,2,1,2,2) = 629 |
| f(1,2,0,2,1,2,2) = 625 | f(2,2,0,2,1,2,2) = 1127 | f(0,0,1,2,1,2,2) = 404 |
| f(1,0,1,2,1,2,2) = 106 | f(2,0,1,2,1,2,2) = 623 | f(0,1,1,2,1,2,2) = 541 |
| f(1,1,1,2,1,2,2) = 1002 | f(2,1,1,2,1,2,2) = 934 | f(0,2,1,2,1,2,2) = 1745 |
| f(1,2,1,2,1,2,2) = 1001 | f(2,2,1,2,1,2,2) = 741 | f(0,0,2,2,1,2,2) = 865 |
| f(1,0,2,2,1,2,2) = 537 | f(2,0,2,2,1,2,2) = 607 | f(0,1,2,2,1,2,2) = 285 |
| f(1,1,2,2,1,2,2) = 986 | f(2,1,2,2,1,2,2) = 854 | f(0,2,2,2,1,2,2) = 1489 |
| f(1,2,2,2,1,2,2) = 985 | f(2,2,2,2,1,2,2) = 469 | f(1,0,0,0,2,2,2) = 957 |
| f(2,0,0,0,2,2,2) = 909 | f(0,1,0,0,2,2,2) = 10 | f(1,1,0,0,2,2,2) = 827 |
| f(2,1,0,0,2,2,2) = 295 | f(0,2,0,0,2,2,2) = 1075 | f(1,2,0,0,2,2,2) = 251 |
| f(2,2,0,0,2,2,2) = 103 | f(0,0,1,0,2,2,2) = 565 | f(1,0,1,0,2,2,2) = 1723 |
| f(2,0,1,0,2,2,2) = 35 | f(0,1,1,0,2,2,2) = 530 | f(2,1,1,0,2,2,2) = 519 |
| f(0,2,1,0,2,2,2) = 155 | f(2,2,1,0,2,2,2) = 135 | f(0,0,2,0,2,2,2) = 27 |
| f(1,0,2,0,2,2,2) = 1403 | f(2,0,2,0,2,2,2) = 19 | f(0,1,2,0,2,2,2) = 274 |
| f(2,1,2,0,2,2,2) = 263 | f(0,2,2,0,2,2,2) = 91 | f(2,2,2,0,2,2,2) = 71 |
| f(0,0,0,1,2,2,2) = 365 | f(1,0,0,1,2,2,2) = 1012 | f(2,0,0,1,2,2,2) = 631 |
| f(0,1,0,1,2,2,2) = 811 | f(1,1,0,1,2,2,2) = 825 | f(2,1,0,1,2,2,2) = 303 |
| f(0,2,0,1,2,2,2) = 33 | f(1,2,0,1,2,2,2) = 250 | f(2,2,0,1,2,2,2) = 111 |
| f(0,0,1,1,2,2,2) = 669 | f(1,0,1,1,2,2,2) = 937 | f(2,0,1,1,2,2,2) = 431 |
| f(0,1,1,1,2,2,2) = 524 | f(1,1,1,1,2,2,2) = 42 | f(2,1,1,1,2,2,2) = 951 |
| f(0,2,1,1,2,2,2) = 146 | f(1,2,1,1,2,2,2) = 41 | f(2,2,1,1,2,2,2) = 759 |
| f(1,0,2,1,2,2,2) = 474 | f(2,0,2,1,2,2,2) = 415 | f(0,1,2,1,2,2,2) = 268 |
| f(1,1,2,1,2,2,2) = 26 | f(2,1,2,1,2,2,2) = 887 | f(0,2,2,1,2,2,2) = 82 |
| f(1,2,2,1,2,2,2) = 25 | f(2,2,2,1,2,2,2) = 503 | f(1,0,0,2,2,2,2) = 635 |
| f(2,0,0,2,2,2,2) = 7 | f(0,1,0,2,2,2,2) = 795 | f(1,1,0,2,2,2,2) = 1851 |
| f(2,1,0,2,2,2,2) = 871 | f(0,2,0,2,2,2,2) = 17 | f(1,2,0,2,2,2,2) = 1275 |
| f(2,2,0,2,2,2,2) = 487 | f(0,0,1,2,2,2,2) = 413 | f(1,0,1,2,2,2,2) = 766 |
| f(2,0,1,2,2,2,2) = 1071 | f(0,1,1,2,2,2,2) = 539 | f(1,1,1,2,2,2,2) = 574 |
| f(2,1,1,2,2,2,2) = 527 | f(0,2,1,2,2,2,2) = 1169 | f(1,2,1,2,2,2,2) = 189 |
| f(2,2,1,2,2,2,2) = 143 | f(0,0,2,2,2,2,2) = 383 | f(1,0,2,2,2,2,2) = 893 |
| f(2,0,2,2,2,2,2) = 1055 | f(0,1,2,2,2,2,2) = 283 | f(1,1,2,2,2,2,2) = 318 |
| f(2,1,2,2,2,2,2) = 271 | f(0,2,2,2,2,2,2) = 1105 | f(1,2,2,2,2,2,2) = 125 |
| f(2,2,2,2,2,2,2) = 79 | | |

List 4

The completed list of weighings of the third algorithm to sort 11 coins

*The 1-th weighing*
W( ) = {1,2,3}:{4,5,6}

*The 2-th weighing*
w(0) = {1,7,8}:{4,9,10}         w(1) = {1,7,8}:{2,9,10}         w(2) = {1,7,8}:{2,9,10}

*The 3-th weighing*
w(0,0) = {1,2,7}:{3,4,8}        w(1,0) = {1,5,7,9}:{2,6,8,10}   w(2,0) = {1,5,7,9}:{2,6,8,10}
w(0,1) = {7,9}:{8,10}           w(1,1) = {5,7,9}:{6,8,10}       w(2,1) = {5,7,9}:{6,8,10}
w(0,2) = {7,9}:{8,10}           w(1,2) = {5,7,9}:{6,8,10}       w(2,2) = {5,7,9}:{6,8,10}

*The 4-th weighing*
w(0,0,0) = {1,3}:{2,4}          w(1,0,0) = {1,2,4}:{3,7,10}     w(2,0,0) = {1,2,4}:{3,7,10}
w(0,1,0) = {5,7,10}:{6,8,9}     w(1,1,0) = {2,9,11}:{3,4,5}     w(2,1,0) = {1,7,11}:{3,4,5}
w(0,2,0) = {5,7,10}:{6,8,9}     w(1,2,0) = {1,7,11}:{3,4,5}     w(2,2,0) = {2,9,11}:{3,4,5}
w(0,0,1) = {5,9}:{6,10}         w(1,0,1) = {2,3,6}:{5,8,9}      w(2,0,1) = {1,3,5}:{6,8,9}
w(0,1,1) = {1,5,6}:{2,3,7}      w(1,1,1) = {3,4,7,10}:{1,2,6,11} w(2,1,1) = {3,4,7,10}:{1,2,5,11}
w(0,2,1) = {1,5,6}:{2,3,8}      w(1,2,1) = {3,4,8,9}:{1,2,6,11}  w(2,2,1) = {3,4,8,9}:{1,2,5,11}
w(0,0,2) = {5,9}:{6,10}         w(1,0,2) = {1,3,5}:{6,8,9}       w(2,0,2) = {2,3,6}:{5,8,9}
w(0,1,2) = {1,5,6}:{2,3,8}      w(1,1,2) = {3,4,8,9}:{1,2,5,11}  w(2,1,2) = {3,4,8,9}:{1,2,6,11}
w(0,2,2) = {1,5,6}:{2,3,7}      w(1,2,2) = {3,4,7,10}:{1,2,5,11} w(2,2,2) = {3,4,7,10}:{1,2,6,11}

*The 5-th weighing*
w(0,0,0,0) = {1,2}:{10,11}      w(1,0,0,0) = {1}:{2}             w(2,0,0,0) = {1}:{2}
w(0,1,0,0) = {2}:{3}            w(1,1,0,0) = {2,6,8}:{3,5,11}    w(2,1,0,0) = {1,6,10}:{3,5,11}
w(0,2,0,0) = {2}:{3}            w(1,2,0,0) = {1,6,10}:{3,5,11}   w(2,2,0,0) = {2,6,8}:{3,5,11}
w(0,0,1,0) = {1}:{2}            w(1,0,1,0) = {2,5,9}:{3,4,7}     w(2,0,1,0) = {1,6,8}:{3,4,10}
w(0,1,1,0) = {2}:{3}            w(1,1,1,0) = {1,2}:{7,11}        w(2,1,1,0) = {1,2}:{10,11}
w(0,2,1,0) = {2}:{3}            w(1,2,1,0) = {1,2}:{9,11}        w(2,2,1,0) = {1,2}:{8,11}
w(0,0,2,0) = {1}:{2}            w(1,0,2,0) = {1,6,8}:{3,4,10}    w(2,0,2,0) = {2,5,9}:{3,4,7}
w(0,1,2,0) = {2}:{3}            w(1,1,2,0) = {1,2}:{8,11}        w(2,1,2,0) = {1,2}:{9,11}
w(0,2,2,0) = {2}:{3}            w(1,2,2,0) = {1,2}:{10,11}       w(2,2,2,0) = {1,2}:{7,11}
w(0,0,0,1) = {1}:{7}            w(1,0,0,1) = {3,7,10}:{4,5,6}    w(2,0,0,1) = {1,2}:{8,9}
w(0,1,0,1) = {2,5}:{3,7}        w(1,1,0,1) = {1,5,9}:{2,4,11}    w(2,1,0,1) = {2,4,11}:{3,5,8}
w(0,2,0,1) = {2,6}:{3,8}        w(1,2,0,1) = {2,5,7}:{1,4,11}    w(2,2,0,1) = {1,4,11}:{3,5,10}
w(0,0,1,1) = {1,4,10}:{3,7,11}  w(1,0,1,1) = {1,6,11}:{3,7,10}   w(2,0,1,1) = {3,6,7}:{8,10,11}
w(0,1,1,1) = {1,10}:{8,9}       w(1,1,1,1) = {4,9}:{5,8}         w(2,1,1,1) = {2,5,6}:{4,10,11}
w(0,2,1,1) = {2}:{3}            w(1,2,1,1) = {4,7}:{5,10}        w(2,2,1,1) = {1,5,6}:{4,8,11}
w(0,0,2,1) = {1,4,9}:{3,7,11}   w(1,0,2,1) = {2,5,11}:{3,7,10}   w(2,0,2,1) = {3,5,10}:{7,9,11}
w(0,1,2,1) = {1,9}:{7,10}       w(1,1,2,1) = {4,10}:{6,7}        w(2,1,2,1) = {2,5,6}:{4,9,11}
w(0,2,2,1) = {2}:{3}            w(1,2,2,1) = {4,8}:{6,9}         w(2,2,2,1) = {1,5,6}:{4,7,11}
w(0,0,0,2) = {1}:{7}            w(1,0,0,2) = {1,2}:{8,9}         w(2,0,0,2) = {3,7,10}:{4,5,6}
w(0,1,0,2) = {2,6}:{3,8}        w(1,1,0,2) = {1,4,11}:{3,5,10}   w(2,1,0,2) = {2,5,7}:{1,4,11}
w(0,2,0,2) = {2,5}:{3,7}        w(1,2,0,2) = {2,4,11}:{3,5,8}    w(2,2,0,2) = {1,5,9}:{2,4,11}
w(0,0,1,2) = {1,4,9}:{3,7,11}   w(1,0,1,2) = {3,5,10}:{7,9,11}   w(2,0,1,2) = {2,5,11}:{3,7,10}

w(0,1,1,2) = {2}:{3}   w(1,1,1,2) = {1,5,6}:{4,7,11}   w(2,1,1,2) = {4,8}:{6,9}
w(0,2,1,2) = {1,9}:{7,10}   w(1,2,1,2) = {2,5,6}:{4,9,11}   w(2,2,1,2) = {4,10}:{6,7}
w(0,0,2,2) = {1,4,10}:{3,7,11}   w(1,0,2,2) = {3,6,7}:{8,10,11}   w(2,0,2,2) = {1,6,11}:{3,7,10}
w(0,1,2,2) = {2}:{3}   w(1,1,2,2) = {1,5,6}:{4,8,11}   w(2,1,2,2) = {4,7}:{5,10}
w(0,2,2,2) = {1,10}:{8,9}   w(1,2,2,2) = {2,5,6}:{4,10,11}   w(2,2,2,2) = {4,9}:{5,8}

*The 6-th weighing*

w(0,0,0,0,0) = {1}:{2}   w(1,0,0,0,0) = {3,4}:{9,11}   w(2,0,0,0,0) = {3,4}:{9,11}
w(0,1,0,0,0) = {1,11}:{2,3}   w(1,1,0,0,0) = {4}:{9}   w(2,1,0,0,0) = {4}:{7}
w(0,2,0,0,0) = {1,11}:{2,3}   w(1,2,0,0,0) = {4}:{7}   w(2,2,0,0,0) = {4}:{9}
w(0,0,1,0,0) = {4,7}:{10,11}   w(1,0,1,0,0) = {2,9}:{3,11}   w(2,0,1,0,0) = {1,8}:{3,11}
w(0,1,1,0,0) = {5,6}:{9,11}   w(1,1,1,0,0) = {2,8}:{7,10}   w(2,1,1,0,0) = {1,9}:{7,10}
w(0,2,1,0,0) = {5,6}:{10,11}   w(1,2,1,0,0) = {1,10}:{8,9}   w(2,2,1,0,0) = {2,7}:{8,9}
w(0,0,2,0,0) = {4,7}:{10,11}   w(1,0,2,0,0) = {1,8}:{3,11}   w(2,0,2,0,0) = {2,9}:{3,11}
w(0,1,2,0,0) = {5,6}:{10,11}   w(1,1,2,0,0) = {2,7}:{8,9}   w(2,1,2,0,0) = {1,10}:{8,9}
w(0,2,2,0,0) = {5,6}:{9,11}   w(1,2,2,0,0) = {1,9}:{7,10}   w(2,2,2,0,0) = {2,8}:{7,10}
w(0,0,0,1,0) = {5,11}:{6,9}   w(1,0,0,1,0) = {3,4}:{9,11}   w(2,0,0,1,0) = {1,2}:{10,11}
w(0,1,0,1,0) = {1,10}:{8,11}   w(1,1,0,1,0) = {1,4}:{2,6}   w(2,1,0,1,0) = {7,9}:{10,11}
w(0,2,0,1,0) = {1,9}:{7,11}   w(1,2,0,1,0) = {2,4}:{1,6}   w(2,2,0,1,0) = {7,9}:{8,11}
w(0,0,1,1,0) = {1}:{2}   w(1,0,1,1,0) = {5}:{6}   w(2,0,1,1,0) = {1,4}:{7,11}
w(0,1,1,1,0) = {1,3,11}:{2,5,6}   w(1,1,1,1,0) = {1,6}:{8,11}   w(2,1,1,1,0) = {8,9}:{5,11}
w(0,2,1,1,0) = {5,6}:{10,11}   w(1,2,1,1,0) = {2,6}:{10,11}   w(2,2,1,1,0) = {7,10}:{5,11}
w(0,0,2,1,0) = {1}:{2}   w(1,0,2,1,0) = {6}:{5}   w(2,0,2,1,0) = {2,4}:{10,11}
w(0,1,2,1,0) = {1,3,11}:{2,5,6}   w(1,1,2,1,0) = {1,5}:{7,11}   w(2,1,2,1,0) = {7,10}:{6,11}
w(0,2,2,1,0) = {5,6}:{9,11}   w(1,2,2,1,0) = {2,5}:{9,11}   w(2,2,2,1,0) = {8,9}:{6,11}
w(0,0,0,2,0) = {5,11}:{6,9}   w(1,0,0,2,0) = {1,2}:{10,11}   w(2,0,0,2,0) = {3,4}:{9,11}
w(0,1,0,2,0) = {1,9}:{7,11}   w(1,1,0,2,0) = {7,9}:{8,11}   w(2,1,0,2,0) = {2,4}:{1,6}
w(0,2,0,2,0) = {1,10}:{8,11}   w(1,2,0,2,0) = {7,9}:{10,11}   w(2,2,0,2,0) = {1,4}:{2,6}
w(0,0,1,2,0) = {1}:{2}   w(1,0,1,2,0) = {2,4}:{10,11}   w(2,0,1,2,0) = {6}:{5}
w(0,1,1,2,0) = {5,6}:{9,11}   w(1,1,1,2,0) = {8,9}:{6,11}   w(2,1,1,2,0) = {2,5}:{9,11}
w(0,2,1,2,0) = {1,3,11}:{2,5,6}   w(1,2,1,2,0) = {7,10}:{6,11}   w(2,2,1,2,0) = {1,5}:{7,11}
w(0,0,2,2,0) = {1}:{2}   w(1,0,2,2,0) = {1,4}:{7,11}   w(2,0,2,2,0) = {5}:{6}
w(0,1,2,2,0) = {5,6}:{10,11}   w(1,1,2,2,0) = {7,10}:{5,11}   w(2,1,2,2,0) = {2,6}:{10,11}
w(0,2,2,2,0) = {1,3,11}:{2,5,6}   w(1,2,2,2,0) = {8,9}:{5,11}   w(2,2,2,2,0) = {1,6}:{8,11}
w(0,0,0,0,1) = {10}:{11}   w(1,0,0,0,1) = {3,4}:{9,11}   w(2,0,0,0,1) = {3,4}:{9,11}
w(0,1,0,0,1) = {3,11}:{9,10}   w(1,1,0,0,1) = {1,6}:{5,8}   w(2,1,0,0,1) = {2,6}:{7,10}
w(0,2,0,0,1) = {2,11}:{9,10}   w(1,2,0,0,1) = {2,6}:{5,10}   w(2,2,0,0,1) = {1,6}:{8,9}
w(0,0,1,0,1) = {5}:{6}   w(1,0,1,0,1) = {5,11}:{7,9}   w(2,0,1,0,1) = {2,4}:{10,11}
w(0,1,1,0,1) = {5,9}:{6,11}   w(1,1,1,0,1) = {2,4}:{5,9}   w(2,1,1,0,1) = {6}:{8}
w(0,2,1,0,1) = {5,10}:{6,11}   w(1,2,1,0,1) = {1,4}:{5,7}   w(2,2,1,0,1) = {6}:{10}
w(0,0,2,0,1) = {5}:{6}   w(1,0,2,0,1) = {6,11}:{8,10}   w(2,0,2,0,1) = {1,4}:{7,11}
w(0,1,2,0,1) = {5,10}:{6,11}   w(1,1,2,0,1) = {2,4}:{6,10}   w(2,1,2,0,1) = {5}:{7}
w(0,2,2,0,1) = {5,9}:{6,11}   w(1,2,2,0,1) = {1,4}:{6,8}   w(2,2,2,0,1) = {5}:{9}
w(0,0,0,1,1) = {3,5}:{1,11}   w(1,0,0,1,1) = {3,8}:{4,11}   w(2,0,0,1,1) = {1,10}:{2,11}
w(0,1,0,1,1) = {1,11}:{2,8}   w(1,1,0,1,1) = {3,6}:{7,10}   w(2,1,0,1,1) = {5,6}:{9,11}

w(0,2,0,1,1) = {1,5}:{6,11}    w(1,2,0,1,1) = {3,6}:{8,9}     w(2,2,0,1,1) = {5,6}:{7,11}
w(0,0,1,1,1) = {3,9}:{4,11}    w(1,0,1,1,1) = {1,5}:{7,11}    w(2,0,1,1,1) = {1,9}:{7,11}
w(0,1,1,1,1) = {1,3,11}:{2,5,6}  w(1,1,1,1,1) = {1,2}:{4,11}  w(2,1,1,1,1) = {1,4}:{3,8}
w(0,2,1,1,1) = {5,10}:{6,11}   w(1,2,1,1,1) = {1,2}:{4,11}    w(2,2,1,1,1) = {2,4}:{3,10}
w(0,0,2,1,1) = {3,11}:{7,9}    w(1,0,2,1,1) = {2,6}:{10,11}   w(2,0,2,1,1) = {2,8}:{10,11}
w(0,1,2,1,1) = {1,3,11}:{2,5,6}  w(1,1,2,1,1) = {1,2}:{4,11}  w(2,1,2,1,1) = {1,4}:{3,7}
w(0,2,2,1,1) = {5,9}:{6,11}    w(1,2,2,1,1) = {1,2}:{4,11}    w(2,2,2,1,1) = {2,4}:{3,9}
w(0,0,0,2,1) = {5,11}:{6,9}    w(1,0,0,2,1) = {1,5}:{7,11}    w(2,0,0,2,1) = {3,5}:{9,11}
w(0,1,0,2,1) = {1,11}:{2,7}    w(1,1,0,2,1) = {2,9}:{10,11}   w(2,1,0,2,1) = {2,4}:{7,10}
w(0,2,0,2,1) = {1,6}:{5,11}    w(1,2,0,2,1) = {1,7}:{8,11}    w(2,2,0,2,1) = {1,4}:{8,9}
w(0,0,1,2,1) = {3,10}:{4,11}   w(1,0,1,2,1) = {2}:{10}        w(2,0,1,2,1) = {6,7}:{1,11}
w(0,1,1,2,1) = {5,9}:{6,11}    w(1,1,1,2,1) = {8,9}:{3,11}    w(2,1,1,2,1) = {1,8}:{7,11}
w(0,2,1,2,1) = {1,3,11}:{2,5,6}  w(1,2,1,2,1) = {7,10}:{3,11} w(2,2,1,2,1) = {2,10}:{9,11}
w(0,0,2,2,1) = {3,11}:{7,10}   w(1,0,2,2,1) = {1}:{7}         w(2,0,2,2,1) = {5,10}:{2,11}
w(0,1,2,2,1) = {5,10}:{6,11}   w(1,1,2,2,1) = {7,10}:{3,11}   w(2,1,2,2,1) = {1,7}:{8,11}
w(0,2,2,2,1) = {1,3,11}:{2,5,6}  w(1,2,2,2,1) = {8,9}:{3,11}  w(2,2,2,2,1) = {2,9}:{10,11}
w(0,0,0,0,2) = {10}:{11}       w(1,0,0,0,2) = {3,4}:{9,11}    w(2,0,0,0,2) = {3,4}:{9,11}
w(0,1,0,0,2) = {2,11}:{9,10}   w(1,1,0,0,2) = {1,6}:{8,9}     w(2,1,0,0,2) = {2,6}:{5,10}
w(0,2,0,0,2) = {3,11}:{9,10}   w(1,2,0,0,2) = {2,6}:{7,10}    w(2,2,0,0,2) = {1,6}:{5,8}
w(0,0,1,0,2) = {5}:{6}         w(1,0,1,0,2) = {1,4}:{7,11}    w(2,0,1,0,2) = {6,11}:{8,10}
w(0,1,1,0,2) = {5,9}:{6,11}    w(1,1,1,0,2) = {5}:{9}         w(2,1,1,0,2) = {1,4}:{6,8}
w(0,2,1,0,2) = {5,10}:{6,11}   w(1,2,1,0,2) = {5}:{7}         w(2,2,1,0,2) = {2,4}:{6,10}
w(0,0,2,0,2) = {5}:{6}         w(1,0,2,0,2) = {2,4}:{10,11}   w(2,0,2,0,2) = {5,11}:{7,9}
w(0,1,2,0,2) = {5,10}:{6,11}   w(1,1,2,0,2) = {6}:{10}        w(2,1,2,0,2) = {1,4}:{5,7}
w(0,2,2,0,2) = {5,9}:{6,11}    w(1,2,2,0,2) = {6}:{8}         w(2,2,2,0,2) = {2,4}:{5,9}
w(0,0,0,1,2) = {5,11}:{6,9}    w(1,0,0,1,2) = {3,5}:{9,11}    w(2,0,0,1,2) = {1,5}:{7,11}
w(0,1,0,1,2) = {1,6}:{5,11}    w(1,1,0,1,2) = {1,4}:{8,9}     w(2,1,0,1,2) = {1,7}:{8,11}
w(0,2,0,1,2) = {1,11}:{2,7}    w(1,2,0,1,2) = {2,4}:{7,10}    w(2,2,0,1,2) = {2,9}:{10,11}
w(0,0,1,1,2) = {3,11}:{7,10}   w(1,0,1,1,2) = {5,10}:{2,11}   w(2,0,1,1,2) = {1}:{7}
w(0,1,1,1,2) = {1,3,11}:{2,5,6}  w(1,1,1,1,2) = {2,9}:{10,11} w(2,1,1,1,2) = {8,9}:{3,11}
w(0,2,1,1,2) = {5,10}:{6,11}   w(1,2,1,1,2) = {1,7}:{8,11}    w(2,2,1,1,2) = {7,10}:{3,11}
w(0,0,2,1,2) = {3,10}:{4,11}   w(1,0,2,1,2) = {6,7}:{1,11}    w(2,0,2,1,2) = {2}:{10}
w(0,1,2,1,2) = {1,3,11}:{2,5,6}  w(1,1,2,1,2) = {2,10}:{9,11} w(2,1,2,1,2) = {7,10}:{3,11}
w(0,2,2,1,2) = {5,9}:{6,11}    w(1,2,2,1,2) = {1,8}:{7,11}    w(2,2,2,1,2) = {8,9}:{3,11}
w(0,0,0,2,2) = {3,5}:{1,11}    w(1,0,0,2,2) = {1,10}:{2,11}   w(2,0,0,2,2) = {3,8}:{4,11}
w(0,1,0,2,2) = {1,5}:{6,11}    w(1,1,0,2,2) = {5,6}:{7,11}    w(2,1,0,2,2) = {3,6}:{8,9}
w(0,2,0,2,2) = {1,11}:{2,8}    w(1,2,0,2,2) = {5,6}:{9,11}    w(2,2,0,2,2) = {3,6}:{7,10}
w(0,0,1,2,2) = {3,11}:{7,9}    w(1,0,1,2,2) = {2,8}:{10,11}   w(2,0,1,2,2) = {2,6}:{10,11}
w(0,1,1,2,2) = {5,9}:{6,11}    w(1,1,1,2,2) = {2,4}:{3,9}     w(2,1,1,2,2) = {1,2}:{4,11}
w(0,2,1,2,2) = {1,3,11}:{2,5,6}  w(1,2,1,2,2) = {1,4}:{3,7}   w(2,2,1,2,2) = {1,2}:{4,11}
w(0,0,2,2,2) = {3,9}:{4,11}    w(1,0,2,2,2) = {1,9}:{7,11}    w(2,0,2,2,2) = {1,5}:{7,11}
w(0,1,2,2,2) = {5,10}:{6,11}   w(1,1,2,2,2) = {2,4}:{3,10}    w(2,1,2,2,2) = {1,2}:{4,11}
w(0,2,2,2,2) = {1,3,11}:{2,5,6}  w(1,2,2,2,2) = {1,4}:{3,8}   w(2,2,2,2,2) = {1,2}:{4,11}

*The 7-th weighing*

w(0,0,0,0,0) = { }:{ }         w(1,0,0,0,0) = {6}:{8}          w(2,0,0,0,0) = {6}:{8}
w(0,1,0,0,0) = {6}:{9}         w(1,1,0,0,0) = {10}:{11}        w(2,1,0,0,0) = {8}:{11}
w(0,2,0,0,0) = {6}:{9}         w(1,2,0,0,0) = {8}:{11}         w(2,2,0,0,0) = {10}:{11}
w(0,0,1,0,0) = {4}:{6}         w(1,0,1,0,0) = {10}:{11}        w(2,0,1,0,0) = {7}:{11}
w(0,1,1,0,0) = {9}:{10}        w(1,1,1,0,0) = {5}:{9}          w(2,1,1,0,0) = {6}:{8}
w(0,2,1,0,0) = {10}:{9}        w(1,2,1,0,0) = {5}:{7}          w(2,2,1,0,0) = {6}:{10}
w(0,0,2,0,0) = {4}:{6}         w(1,0,2,0,0) = {7}:{11}         w(2,0,2,0,0) = {10}:{11}
w(0,1,2,0,0) = {10}:{9}        w(1,1,2,0,0) = {6}:{10}         w(2,1,2,0,0) = {5}:{7}
w(0,2,2,0,0) = {9}:{10}        w(1,2,2,0,0) = {6}:{8}          w(2,2,2,0,0) = {5}:{9}
w(0,0,0,1,0) = {9}:{10}        w(1,0,0,1,0) = {3}:{5}          w(2,0,0,1,0) = {3,4}:{6,7}
w(0,1,0,1,0) = {6}:{9}         w(1,1,0,1,0) = {5}:{7}          w(2,1,0,1,0) = {2}:{9}
w(0,2,0,1,0) = {5}:{10}        w(1,2,0,1,0) = {5}:{9}          w(2,2,0,1,0) = {1}:{7}
w(0,0,1,1,0) = {6}:{7}         w(1,0,1,1,0) = {4}:{6}          w(2,0,1,1,0) = {7}:{11}
w(0,1,1,1,0) = {5}:{6}         w(1,1,1,1,0) = {1,2}:{6,7}      w(2,1,1,1,0) = {1,4}:{3,6}
w(0,2,1,1,0) = {10}:{9}        w(1,2,1,1,0) = {1,2}:{6,9}      w(2,2,1,1,0) = {2,4}:{3,6}
w(0,0,2,1,0) = {5}:{7}         w(1,0,2,1,0) = {4}:{5}          w(2,0,2,1,0) = {10}:{11}
w(0,1,2,1,0) = {5}:{6}         w(1,1,2,1,0) = {1,2}:{5,8}      w(2,1,2,1,0) = {1,4}:{3,5}
w(0,2,2,1,0) = {9}:{10}        w(1,2,2,1,0) = {1,2}:{5,10}     w(2,2,2,1,0) = {2,4}:{3,5}
w(0,0,0,2,0) = {9}:{10}        w(1,0,0,2,0) = {3,4}:{6,7}      w(2,0,0,2,0) = {3}:{5}
w(0,1,0,2,0) = {5}:{10}        w(1,1,0,2,0) = {1}:{7}          w(2,1,0,2,0) = {5}:{9}
w(0,2,0,2,0) = {6}:{9}         w(1,2,0,2,0) = {2}:{9}          w(2,2,0,2,0) = {5}:{7}
w(0,0,1,2,0) = {5}:{7}         w(1,0,1,2,0) = {10}:{11}        w(2,0,1,2,0) = {4}:{5}
w(0,1,1,2,0) = {9}:{10}        w(1,1,1,2,0) = {2,4}:{3,5}      w(2,1,1,2,0) = {1,2}:{5,10}
w(0,2,1,2,0) = {5}:{6}         w(1,2,1,2,0) = {1,4}:{3,5}      w(2,2,1,2,0) = {1,2}:{5,8}
w(0,0,2,2,0) = {6}:{7}         w(1,0,2,2,0) = {7}:{11}         w(2,0,2,2,0) = {4}:{6}
w(0,1,2,2,0) = {10}:{9}        w(1,1,2,2,0) = {2,4}:{3,6}      w(2,1,2,2,0) = {1,2}:{6,9}
w(0,2,2,2,0) = {5}:{6}         w(1,2,2,2,0) = {1,4}:{3,6}      w(2,2,2,2,0) = {1,2}:{6,7}
w(0,0,0,0,1) = {1}:{2}         w(1,0,0,0,1) = {6}:{8}          w(2,0,0,0,1) = {10}:{11}
w(0,1,0,0,1) = {1}:{7}         w(1,1,0,0,1) = {3}:{4}          w(2,1,0,0,1) = {1}:{9}
w(0,2,0,0,1) = {1}:{7}         w(1,2,0,0,1) = {3}:{4}          w(2,2,0,0,1) = {2}:{7}
w(0,0,1,0,1) = {8}:{11}        w(1,0,1,0,1) = {10}:{11}        w(2,0,1,0,1) = {1}:{4}
w(0,1,1,0,1) = {9}:{11}        w(1,1,1,0,1) = {5}:{8}          w(2,1,1,0,1) = {5,10}:{8,9}
w(0,2,1,0,1) = {10}:{11}       w(1,2,1,0,1) = {5}:{10}         w(2,2,1,0,1) = {5,8}:{7,10}
w(0,0,2,0,1) = {8}:{11}        w(1,0,2,0,1) = {7}:{11}         w(2,0,2,0,1) = {2}:{4}
w(0,1,2,0,1) = {10}:{11}       w(1,1,2,0,1) = {6}:{7}          w(2,1,2,0,1) = {6,9}:{7,10}
w(0,2,2,0,1) = {9}:{11}        w(1,2,2,0,1) = {6}:{9}          w(2,2,2,0,1) = {6,7}:{8,9}
w(0,0,0,1,1) = {9}:{10}        w(1,0,0,1,1) = {3}:{4}          w(2,0,0,1,1) = {2}:{11}
w(0,1,0,1,1) = {6}:{8}         w(1,1,0,1,1) = {2}:{11}         w(2,1,0,1,1) = {1,4}:{5,6}
w(0,2,0,1,1) = {5}:{7}         w(1,2,0,1,1) = {1}:{11}         w(2,2,0,1,1) = {2,4}:{5,6}
w(0,0,1,1,1) = {1,3}:{2,8}     w(1,0,1,1,1) = {2}:{3}          w(2,0,1,1,1) = {2,3}:{9,10}
w(0,1,1,1,1) = {5}:{6}         w(1,1,1,1,1) = {6,7}:{8,9}      w(2,1,1,1,1) = {3}:{8}
w(0,2,1,1,1) = {10}:{11}       w(1,2,1,1,1) = {6,9}:{7,10}     w(2,2,1,1,1) = {3}:{10}
w(0,0,2,1,1) = {2}:{11}        w(1,0,2,1,1) = {1}:{3}          w(2,0,2,1,1) = {1,3}:{7,8}
w(0,1,2,1,1) = {5}:{6}         w(1,1,2,1,1) = {5,8}:{7,10}     w(2,1,2,1,1) = {3}:{7}

| | | |
|---|---|---|
| w(0,2,2,1,1,0) = {9}:{11} | w(1,2,2,1,1,0) = {5,10}:{8,9} | w(2,2,2,1,1,0) = {3}:{9} |
| w(0,0,0,2,1,0) = {5}:{6} | w(1,0,0,2,1,0) = {2}:{11} | w(2,0,0,2,1,0) = {3}:{4} |
| w(0,1,0,2,1,0) = {5}:{7} | w(1,1,0,2,1,0) = {2}:{7} | w(2,1,0,2,1,0) = {1,5}:{3,4} |
| w(0,2,0,2,1,0) = {6}:{8} | w(1,2,0,2,1,0) = {1}:{9} | w(2,2,0,2,1,0) = {2,5}:{3,4} |
| w(0,0,1,2,1,0) = {1,3}:{2,8} | w(1,0,1,2,1,0) = {2,4}:{5,6} | w(2,0,1,2,1,0) = {2,4}:{5,6} |
| w(0,1,1,2,1,0) = {9}:{11} | w(1,1,1,2,1,0) = {2}:{3} | w(2,1,1,2,1,0) = {5,6}:{8,9} |
| w(0,2,1,2,1,0) = {5}:{6} | w(1,2,1,2,1,0) = {1}:{3} | w(2,2,1,2,1,0) = {5,6}:{7,10} |
| w(0,0,2,2,1,0) = {2}:{11} | w(1,0,2,2,1,0) = {1,4}:{5,6} | w(2,0,2,2,1,0) = {1,4}:{5,6} |
| w(0,1,2,2,1,0) = {10}:{11} | w(1,1,2,2,1,0) = {2}:{3} | w(2,1,2,2,1,0) = {5,6}:{7,10} |
| w(0,2,2,2,1,0) = {5}:{6} | w(1,2,2,2,1,0) = {1}:{3} | w(2,2,2,2,1,0) = {5,6}:{8,9} |
| w(0,0,0,0,2,0) = {1}:{2} | w(1,0,0,0,2,0) = {10}:{11} | w(2,0,0,0,2,0) = {6}:{8} |
| w(0,1,0,0,2,0) = {1}:{7} | w(1,1,0,0,2,0) = {2}:{7} | w(2,1,0,0,2,0) = {3}:{4} |
| w(0,2,0,0,2,0) = {1}:{7} | w(1,2,0,0,2,0) = {1}:{9} | w(2,2,0,0,2,0) = {3}:{4} |
| w(0,0,1,0,2,0) = {8}:{11} | w(1,0,1,0,2,0) = {2}:{4} | w(2,0,1,0,2,0) = {7}:{11} |
| w(0,1,1,0,2,0) = {9}:{11} | w(1,1,1,0,2,0) = {6,7}:{8,9} | w(2,1,1,0,2,0) = {6}:{9} |
| w(0,2,1,0,2,0) = {10}:{11} | w(1,2,1,0,2,0) = {6,9}:{7,10} | w(2,2,1,0,2,0) = {6}:{7} |
| w(0,0,2,0,2,0) = {8}:{11} | w(1,0,2,0,2,0) = {1}:{4} | w(2,0,2,0,2,0) = {10}:{11} |
| w(0,1,2,0,2,0) = {10}:{11} | w(1,1,2,0,2,0) = {5,8}:{7,10} | w(2,1,2,0,2,0) = {5}:{10} |
| w(0,2,2,0,2,0) = {9}:{11} | w(1,2,2,0,2,0) = {5,10}:{8,9} | w(2,2,2,0,2,0) = {5}:{8} |
| w(0,0,0,1,2,0) = {5}:{6} | w(1,0,0,1,2,0) = {3}:{4} | w(2,0,0,1,2,0) = {2}:{11} |
| w(0,1,0,1,2,0) = {6}:{8} | w(1,1,0,1,2,0) = {2,5}:{3,4} | w(2,1,0,1,2,0) = {1}:{9} |
| w(0,2,0,1,2,0) = {5}:{7} | w(1,2,0,1,2,0) = {1,5}:{3,4} | w(2,2,0,1,2,0) = {2}:{7} |
| w(0,0,1,1,2,0) = {2}:{11} | w(1,0,1,1,2,0) = {1,4}:{5,6} | w(2,0,1,1,2,0) = {1,4}:{5,6} |
| w(0,1,1,1,2,0) = {5}:{6} | w(1,1,1,1,2,0) = {5,6}:{8,9} | w(2,1,1,1,2,0) = {1}:{3} |
| w(0,2,1,1,2,0) = {10}:{11} | w(1,2,1,1,2,0) = {5,6}:{7,10} | w(2,2,1,1,2,0) = {2}:{3} |
| w(0,0,2,1,2,0) = {1,3}:{2,8} | w(1,0,2,1,2,0) = {2,4}:{5,6} | w(2,0,2,1,2,0) = {2,4}:{5,6} |
| w(0,1,2,1,2,0) = {5}:{6} | w(1,1,2,1,2,0) = {5,6}:{7,10} | w(2,1,2,1,2,0) = {1}:{3} |
| w(0,2,2,1,2,0) = {9}:{11} | w(1,2,2,1,2,0) = {5,6}:{8,9} | w(2,2,2,1,2,0) = {2}:{3} |
| w(0,0,0,2,2,0) = {9}:{10} | w(1,0,0,2,2,0) = {2}:{11} | w(2,0,0,2,2,0) = {3}:{4} |
| w(0,1,0,2,2,0) = {5}:{7} | w(1,1,0,2,2,0) = {2,4}:{5,6} | w(2,1,0,2,2,0) = {1}:{11} |
| w(0,2,0,2,2,0) = {6}:{8} | w(1,2,0,2,2,0) = {1,4}:{5,6} | w(2,2,0,2,2,0) = {2}:{11} |
| w(0,0,1,2,2,0) = {2}:{11} | w(1,0,1,2,2,0) = {1,3}:{7,8} | w(2,0,1,2,2,0) = {1}:{3} |
| w(0,1,1,2,2,0) = {9}:{11} | w(1,1,1,2,2,0) = {3}:{9} | w(2,1,1,2,2,0) = {5,10}:{8,9} |
| w(0,2,1,2,2,0) = {5}:{6} | w(1,2,1,2,2,0) = {3}:{7} | w(2,2,1,2,2,0) = {5,8}:{7,10} |
| w(0,0,2,2,2,0) = {1,3}:{2,8} | w(1,0,2,2,2,0) = {2,3}:{9,10} | w(2,0,2,2,2,0) = {2}:{3} |
| w(0,1,2,2,2,0) = {10}:{11} | w(1,1,2,2,2,0) = {3}:{10} | w(2,1,2,2,2,0) = {6,9}:{7,10} |
| w(0,2,2,2,2,0) = {5}:{6} | w(1,2,2,2,2,0) = {3}:{8} | w(2,2,2,2,2,0) = {6,7}:{8,9} |
| w(0,0,0,0,0,1) = {10}:{11} | w(1,0,0,0,0,1) = {8}:{11} | w(2,0,0,0,0,1) = {8}:{11} |
| w(0,1,0,0,0,1) = {1}:{11} | w(1,1,0,0,0,1) = {1}:{8} | w(2,1,0,0,0,1) = {1}:{9} |
| w(0,2,0,0,0,1) = {1}:{11} | w(1,2,0,0,0,1) = {2}:{10} | w(2,2,0,0,0,1) = {2}:{7} |
| w(0,0,1,0,0,1) = {10}:{11} | w(1,0,1,0,0,1) = {3}:{11} | w(2,0,1,0,0,1) = {7}:{11} |
| w(0,1,1,0,0,1) = {9}:{11} | w(1,1,1,0,0,1) = {5}:{9} | w(2,1,1,0,0,1) = {3,5}:{6,8} |
| w(0,2,1,0,0,1) = {10}:{11} | w(1,2,1,0,0,1) = {5}:{7} | w(2,2,1,0,0,1) = {3,5}:{6,10} |
| w(0,0,2,0,0,1) = {8}:{11} | w(1,0,2,0,0,1) = {3}:{11} | w(2,0,2,0,0,1) = {10}:{11} |
| w(0,1,2,0,0,1) = {10}:{11} | w(1,1,2,0,0,1) = {6}:{10} | w(2,1,2,0,0,1) = {3,6}:{5,7} |
| w(0,2,2,0,0,1) = {9}:{11} | w(1,2,2,0,0,1) = {6}:{8} | w(2,2,2,0,0,1) = {3,6}:{5,9} |
| w(0,0,0,1,0,1) = {10}:{11} | w(1,0,0,1,0,1) = {9}:{11} | w(2,0,0,1,0,1) = {10}:{11} |
| w(0,1,0,1,0,1) = {8}:{11} | w(1,1,0,1,0,1) = {8}:{10} | w(2,1,0,1,0,1) = {1}:{11} |
| w(0,2,0,1,0,1) = {2}:{11} | w(1,2,0,1,0,1) = {10}:{8} | w(2,2,0,1,0,1) = {2}:{11} |
| w(0,0,1,1,0,1) = {5}:{7} | w(1,0,1,1,0,1) = {4}:{6} | w(2,0,1,1,0,1) = {3,4}:{9,10} |
| w(0,1,1,1,0,1) = {5}:{6} | w(1,1,1,1,0,1) = {2}:{9} | w(2,1,1,1,0,1) = {1}:{3} |
| w(0,2,1,1,0,1) = {10}:{11} | w(1,2,1,1,0,1) = {1}:{7} | w(2,2,1,1,0,1) = {2}:{3} |
| w(0,0,2,1,0,1) = {5}:{7} | w(1,0,2,1,0,1) = {4}:{5} | w(2,0,2,1,0,1) = {3,4}:{7,8} |
| w(0,1,2,1,0,1) = {5}:{6} | w(1,1,2,1,0,1) = {2}:{10} | w(2,1,2,1,0,1) = {1}:{3} |
| w(0,2,2,1,0,1) = {9}:{11} | w(1,2,2,1,0,1) = {1}:{8} | w(2,2,2,1,0,1) = {2}:{3} |
| w(0,0,0,2,0,1) = {10}:{11} | w(1,0,0,2,0,1) = {4}:{5} | w(2,0,0,2,0,1) = {9}:{11} |
| w(0,1,0,2,0,1) = {7}:{11} | w(1,1,0,2,0,1) = {6}:{8} | w(2,1,0,2,0,1) = {1,4}:{3,10} |
| w(0,2,0,2,0,1) = {2}:{11} | w(1,2,0,2,0,1) = {6}:{10} | w(2,2,0,2,0,1) = {2,4}:{3,8} |
| w(0,0,1,2,0,1) = {6}:{7} | w(1,0,1,2,0,1) = {1}:{3} | w(2,0,1,2,0,1) = {4}:{11} |
| w(0,1,1,2,0,1) = {9}:{11} | w(1,1,1,2,0,1) = {2,4}:{3,5} | w(2,1,1,2,0,1) = {9}:{11} |
| w(0,2,1,2,0,1) = {9}:{11} | w(1,2,1,2,0,1) = {1,4}:{3,5} | w(2,2,1,2,0,1) = {7}:{11} |
| w(0,0,2,2,0,1) = {6}:{7} | w(1,0,2,2,0,1) = {2}:{3} | w(2,0,2,2,0,1) = {4}:{11} |
| w(0,1,2,2,0,1) = {10}:{11} | w(1,1,2,2,0,1) = {2,4}:{3,6} | w(2,1,2,2,0,1) = {10}:{11} |
| w(0,2,2,2,0,1) = {10}:{11} | w(1,2,2,2,0,1) = {1,4}:{3,6} | w(2,2,2,2,0,1) = {8}:{11} |
| w(0,0,0,0,1,1) = { }:{ } | w(1,0,0,0,1,1) = {8}:{11} | w(2,0,0,0,1,1) = {8}:{11} |
| w(0,1,0,0,1,1) = {1}:{7} | w(1,1,0,0,1,1) = {2}:{8} | w(2,1,0,0,1,1) = {2}:{9} |
| w(0,2,0,0,1,1) = {2}:{11} | w(1,2,0,0,1,1) = {1}:{10} | w(2,2,0,0,1,1) = {1}:{7} |
| w(0,0,1,0,1,1) = {9}:{11} | w(1,0,1,0,1,1) = {5,9}:{7,11} | w(2,0,1,0,1,1) = {2,11}:{4,10} |
| w(0,1,1,0,1,1) = {2,3}:{9,11} | w(1,1,1,0,1,1) = {3}:{8} | w(2,1,1,0,1,1) = {9}:{7} |
| w(0,2,1,0,1,1) = {2,3}:{10,11} | w(1,2,1,0,1,1) = {3}:{10} | w(2,2,1,0,1,1) = {7}:{9} |
| w(0,0,2,0,1,1) = {10}:{11} | w(1,0,2,0,1,1) = {6,8}:{10,11} | w(2,0,2,0,1,1) = {1,11}:{4,7} |
| w(0,1,2,0,1,1) = {2,3}:{10,11} | w(1,1,2,0,1,1) = {3}:{7} | w(2,1,2,0,1,1) = {10}:{8} |
| w(0,2,2,0,1,1) = {2,3}:{9,11} | w(1,2,2,0,1,1) = {3}:{9} | w(2,2,2,0,1,1) = {8}:{10} |
| w(0,0,0,1,1,1) = {9}:{10} | w(1,0,0,1,1,1) = {3,4}:{9,11} | w(2,0,0,1,1,1) = {2}:{11} |
| w(0,1,0,1,1,1) = {4,5}:{10,11} | w(1,1,0,1,1,1) = {6}:{11} | w(2,1,0,1,1,1) = {2}:{10} |
| w(0,2,0,1,1,1) = {1,4}:{9,11} | w(1,2,0,1,1,1) = {6}:{11} | w(2,2,0,1,1,1) = {1}:{8} |
| w(0,0,1,1,1,1) = {6}:{7} | w(1,0,1,1,1,1) = {4}:{6} | w(2,0,1,1,1,1) = {7}:{11} |
| w(0,1,1,1,1,1) = {5}:{6} | w(1,1,1,1,1,1) = {2}:{9} | w(2,1,1,1,1,1) = {1,3}:{4,9} |
| w(0,2,1,1,1,1) = {2,3}:{10,11} | w(1,2,1,1,1,1) = {1}:{7} | w(2,2,1,1,1,1) = {2,3}:{4,7} |
| w(0,0,2,1,1,1) = {5,6}:{8,9} | w(1,0,2,1,1,1) = {4}:{5} | w(2,0,2,1,1,1) = {10}:{11} |
| w(0,1,2,1,1,1) = {5}:{6} | w(1,1,2,1,1,1) = {2}:{10} | w(2,1,2,1,1,1) = {1,3}:{4,10} |
| w(0,2,2,1,1,1) = {2,3}:{9,11} | w(1,2,2,1,1,1) = {1}:{8} | w(2,2,2,1,1,1) = {2,3}:{4,8} |
| w(0,0,0,2,1,1) = {10}:{11} | w(1,0,0,2,1,1) = {2}:{4} | w(2,0,0,2,1,1) = {9}:{11} |
| w(0,1,0,2,1,1) = {4,6}:{9,11} | w(1,1,0,2,1,1) = {7}:{8} | w(2,1,0,2,1,1) = {1,2}:{3,11} |
| w(0,2,0,2,1,1) = {1,4}:{10,11} | w(1,2,0,2,1,1) = {9}:{10} | w(2,2,0,2,1,1) = {1,2}:{3,11} |
| w(0,0,1,2,1,1) = {5}:{7} | w(1,0,1,2,1,1) = {4}:{6} | w(2,0,1,2,1,1) = {4}:{6} |
| w(0,1,1,2,1,1) = {2,3}:{9,11} | w(1,1,1,2,1,1) = {3,6}:{4,9} | w(2,1,1,2,1,1) = {5,10}:{8,11} |
| w(0,2,1,2,1,1) = {9}:{11} | w(1,2,1,2,1,1) = {3,6}:{4,7} | w(2,2,1,2,1,1) = {5,8}:{10,11} |
| w(0,0,2,2,1,1) = {6,7}:{5,8} | w(1,0,2,2,1,1) = {4}:{5} | w(2,0,2,2,1,1) = {4}:{5} |
| w(0,1,2,2,1,1) = {2,3}:{10,11} | w(1,1,2,2,1,1) = {3,5}:{4,10} | w(2,1,2,2,1,1) = {6,9}:{7,11} |
| w(0,2,2,2,1,1) = {10}:{11} | w(1,2,2,2,1,1) = {3,5}:{4,8} | w(2,2,2,2,1,1) = {6,7}:{9,11} |
| w(0,0,0,0,2,1) = { }:{ } | w(1,0,0,0,2,1) = {8}:{11} | w(2,0,0,0,2,1) = {8}:{11} |
| w(0,1,0,0,2,1) = {1}:{7} | w(1,1,0,0,2,1) = {6}:{10} | w(2,1,0,0,2,1) = {5}:{9} |
| w(0,2,0,0,2,1) = {3}:{11} | w(1,2,0,0,2,1) = {6}:{8} | w(2,2,0,0,2,1) = {5}:{7} |
| w(0,0,1,0,2,1) = {9}:{11} | w(1,0,1,0,2,1) = {7}:{11} | w(2,0,1,0,2,1) = {3}:{11} |
| w(0,1,1,0,2,1) = {2,3}:{9,11} | w(1,1,1,0,2,1) = {1}:{8} | w(2,1,1,0,2,1) = {5,6}:{9,10} |
| w(0,2,1,0,2,1) = {2,3}:{10,11} | w(1,2,1,0,2,1) = {2}:{10} | w(2,2,1,0,2,1) = {5,6}:{7,8} |
| w(0,0,2,0,2,1) = {10}:{11} | w(1,0,2,0,2,1) = {10}:{11} | w(2,0,2,0,2,1) = {3}:{11} |
| w(0,1,2,0,2,1) = {2,3}:{10,11} | w(1,1,2,0,2,1) = {1}:{7} | w(2,1,2,0,2,1) = {5,6}:{9,10} |
| w(0,2,2,0,2,1) = {2,3}:{9,11} | w(1,2,2,0,2,1) = {2}:{9} | w(2,2,2,0,2,1) = {5,6}:{7,8} |
| w(0,0,0,1,2,1) = {10}:{11} | w(1,0,0,1,2,1) = {9}:{11} | w(2,0,0,1,2,1) = {1,11}:{3,4} |
| w(0,1,0,1,2,1) = {5}:{11} | w(1,1,0,1,2,1) = {3}:{4} | w(2,1,0,1,2,1) = {10}:{8} |
| w(0,2,0,1,2,1) = {1}:{11} | w(1,2,0,1,2,1) = {3}:{4} | w(2,2,0,1,2,1) = {8}:{10} |
| w(0,0,1,1,2,1) = {4,5}:{6,7} | w(1,0,1,1,2,1) = {4}:{6} | w(2,0,1,1,2,1) = {10}:{11} |
| w(0,1,1,1,2,1) = {5}:{6} | w(1,1,1,1,2,1) = {1}:{2} | w(2,1,1,1,2,1) = {1}:{11} |
| w(0,2,1,1,2,1) = {2,3}:{10,11} | w(1,2,1,1,2,1) = {2}:{1} | w(2,2,1,1,2,1) = {2}:{11} |
| w(0,0,2,1,2,1) = {5,6}:{10,11} | w(1,0,2,1,2,1) = {4}:{5} | w(2,0,2,1,2,1) = {7}:{11} |
| w(0,1,2,1,2,1) = {5}:{6} | w(1,1,2,1,2,1) = {1}:{2} | w(2,1,2,1,2,1) = {1}:{11} |
| w(0,2,2,1,2,1) = {2,3}:{9,11} | w(1,2,2,1,2,1) = {2}:{1} | w(2,2,2,1,2,1) = {2}:{11} |
| w(0,0,0,2,2,1) = {6}:{11} | w(1,0,0,2,2,1) = {7}:{11} | w(2,0,0,2,2,1) = {3,4}:{9,11} |
| w(0,1,0,2,2,1) = {6}:{11} | w(1,1,0,2,2,1) = {2,6}:{4,5} | w(2,1,0,2,2,1) = {3,11}:{4,8} |
| w(0,2,0,2,2,1) = {1}:{11} | w(1,2,0,2,2,1) = {1,6}:{4,5} | w(2,2,0,2,2,1) = {3,11}:{4,10} |
| w(0,0,1,2,2,1) = {4,5}:{6,8} | w(1,0,1,2,2,1) = {7}:{8} | w(2,0,1,2,2,1) = {3,4}:{5,6} |
| w(0,1,1,2,2,1) = {2,3}:{9,11} | w(1,1,1,2,2,1) = {3,4}:{8,9} | w(2,1,1,2,2,1) = {5,10}:{8,9} |
| w(0,2,1,2,2,1) = {9}:{11} | w(1,2,1,2,2,1) = {3,4}:{7,10} | w(2,2,1,2,2,1) = {5,8}:{7,10} |
| w(0,0,2,2,2,1) = {5,6}:{10,11} | w(1,0,2,2,2,1) = {10}:{9} | w(2,0,2,2,2,1) = {3,4}:{5,6} |
| w(0,1,2,2,2,1) = {2,3}:{10,11} | w(1,1,2,2,2,1) = {3,4}:{7,10} | w(2,1,2,2,2,1) = {6,9}:{7,10} |
| w(0,2,2,2,2,1) = {10}:{11} | w(1,2,2,2,2,1) = {3,4}:{8,9} | w(2,2,2,2,2,1) = {6,7}:{8,9} |
| w(0,0,0,0,0,2) = {10}:{11} | w(1,0,0,0,0,2) = {8}:{11} | w(2,0,0,0,0,2) = {8}:{11} |
| w(0,1,0,0,0,2) = {1}:{11} | w(1,1,0,0,0,2) = {2}:{7} | w(2,1,0,0,0,2) = {2}:{10} |
| w(0,2,0,0,0,2) = {1}:{11} | w(1,2,0,0,0,2) = {1}:{9} | w(2,2,0,0,0,2) = {1}:{8} |
| w(0,0,1,0,0,2) = {8}:{11} | w(1,0,1,0,0,2) = {10}:{11} | w(2,0,1,0,0,2) = {3}:{11} |
| w(0,1,1,0,0,2) = {9}:{11} | w(1,1,1,0,0,2) = {3,6}:{5,9} | w(2,1,1,0,0,2) = {6}:{8} |
| w(0,2,1,0,0,2) = {10}:{11} | w(1,2,1,0,0,2) = {3,6}:{5,7} | w(2,2,1,0,0,2) = {6}:{10} |
| w(0,0,2,0,0,2) = {10}:{11} | w(1,0,2,0,0,2) = {7}:{11} | w(2,0,2,0,0,2) = {3}:{11} |
| w(0,1,2,0,0,2) = {10}:{11} | w(1,1,2,0,0,2) = {3,5}:{6,10} | w(2,1,2,0,0,2) = {5}:{7} |
| w(0,2,2,0,0,2) = {9}:{11} | w(1,2,2,0,0,2) = {3,5}:{6,8} | w(2,2,2,0,0,2) = {5}:{9} |
| w(0,0,0,1,0,2) = {10}:{11} | w(1,0,0,1,0,2) = {9}:{11} | w(2,0,0,1,0,2) = {4}:{5} |
| w(0,1,0,1,0,2) = {2}:{11} | w(1,1,0,1,0,2) = {2,4}:{3,8} | w(2,1,0,1,0,2) = {6}:{10} |
| w(0,2,0,1,0,2) = {7}:{11} | w(1,2,0,1,0,2) = {1,4}:{3,10} | w(2,2,0,1,0,2) = {6}:{8} |
| w(0,0,1,1,0,2) = {6}:{7} | w(1,0,1,1,0,2) = {4}:{11} | w(2,0,1,1,0,2) = {2}:{3} |
| w(0,1,1,1,0,2) = {10}:{11} | w(1,1,1,1,0,2) = {8}:{11} | w(2,1,1,1,0,2) = {1,4}:{3,6} |
| w(0,2,1,1,0,2) = {10}:{11} | w(1,2,1,1,0,2) = {10}:{11} | w(2,2,1,1,0,2) = {2,4}:{3,6} |
| w(0,0,2,1,0,2) = {6}:{7} | w(1,0,2,1,0,2) = {4}:{11} | w(2,0,2,1,0,2) = {1}:{3} |
| w(0,1,2,1,0,2) = {9}:{11} | w(1,1,2,1,0,2) = {7}:{11} | w(2,1,2,1,0,2) = {1,4}:{3,5} |
| w(0,2,2,1,0,2) = {9}:{11} | w(1,2,2,1,0,2) = {9}:{11} | w(2,2,2,1,0,2) = {2,4}:{3,5} |
| w(0,0,0,2,0,2) = {10}:{11} | w(1,0,0,2,0,2) = {10}:{11} | w(2,0,0,2,0,2) = {9}:{11} |
| w(0,1,0,2,0,2) = {2}:{11} | w(1,1,0,2,0,2) = {2}:{11} | w(2,1,0,2,0,2) = {10}:{8} |
| w(0,2,0,2,0,2) = {8}:{11} | w(1,2,0,2,0,2) = {1}:{11} | w(2,2,0,2,0,2) = {8}:{10} |
| w(0,0,1,2,0,2) = {5}:{7} | w(1,0,1,2,0,2) = {3,4}:{7,8} | w(2,0,1,2,0,2) = {4}:{5} |
| w(0,1,1,2,0,2) = {9}:{11} | w(1,1,1,2,0,2) = {2}:{3} | w(2,1,1,2,0,2) = {1}:{8} |
| w(0,2,1,2,0,2) = {5}:{6} | w(1,2,1,2,0,2) = {1}:{3} | w(2,2,1,2,0,2) = {2}:{10} |
| w(0,0,2,2,0,2) = {5}:{7} | w(1,0,2,2,0,2) = {3,4}:{9,10} | w(2,0,2,2,0,2) = {4}:{6} |
| w(0,1,2,2,0,2) = {10}:{11} | w(1,1,2,2,0,2) = {2}:{3} | w(2,1,2,2,0,2) = {1}:{7} |
| w(0,2,2,2,0,2) = {5}:{6} | w(1,2,2,2,0,2) = {1}:{3} | w(2,2,2,2,0,2) = {2}:{9} |
| w(0,0,0,0,1,2) = { }:{ } | w(1,0,0,0,1,2) = {8}:{11} | w(2,0,0,0,1,2) = {8}:{11} |
| w(0,1,0,0,1,2) = {3}:{11} | w(1,1,0,0,1,2) = {5}:{7} | w(2,1,0,0,1,2) = {6}:{8} |
| w(0,2,0,0,1,2) = {1}:{7} | w(1,2,0,0,1,2) = {5}:{9} | w(2,2,0,0,1,2) = {6}:{10} |
| w(0,0,1,0,1,2) = {10}:{11} | w(1,0,1,0,1,2) = {3}:{11} | w(2,0,1,0,1,2) = {10}:{11} |
| w(0,1,1,0,1,2) = {2,3}:{9,11} | w(1,1,1,0,1,2) = {5,6}:{7,8} | w(2,1,1,0,1,2) = {2}:{9} |
| w(0,2,1,0,1,2) = {2,3}:{10,11} | w(1,2,1,0,1,2) = {5,6}:{9,10} | w(2,2,1,0,1,2) = {1}:{7} |
| w(0,0,2,0,1,2) = {9}:{11} | w(1,0,2,0,1,2) = {3}:{11} | w(2,0,2,0,1,2) = {7}:{11} |
| w(0,1,2,0,1,2) = {2,3}:{10,11} | w(1,1,2,0,1,2) = {5,6}:{7,8} | w(2,1,2,0,1,2) = {2}:{10} |
| w(0,2,2,0,1,2) = {2,3}:{9,11} | w(1,2,2,0,1,2) = {5,6}:{9,10} | w(2,2,2,0,1,2) = {1}:{8} |
| w(0,0,0,1,1,2) = {6}:{11} | w(1,0,0,1,1,2) = {3,4}:{9,11} | w(2,0,0,1,1,2) = {7}:{11} |
| w(0,1,0,1,1,2) = {1}:{11} | w(1,1,0,1,1,2) = {3,11}:{4,10} | w(2,1,0,1,1,2) = {1,6}:{4,5} |
| w(0,2,0,1,1,2) = {6}:{11} | w(1,2,0,1,1,2) = {3,11}:{4,8} | w(2,2,0,1,1,2) = {2,6}:{4,5} |
| w(0,0,1,1,1,2) = {5,6}:{10,11} | w(1,0,1,1,1,2) = {3,4}:{5,6} | w(2,0,1,1,1,2) = {10}:{9} |
| w(0,1,1,1,1,2) = {10}:{11} | w(1,1,1,1,1,2) = {6,7}:{8,9} | w(2,1,1,1,1,2) = {3,4}:{8,9} |
| w(0,2,1,1,1,2) = {2,3}:{10,11} | w(1,2,1,1,1,2) = {6,9}:{7,10} | w(2,2,1,1,1,2) = {3,4}:{7,10} |
| w(0,0,2,1,1,2) = {4,5}:{6,8} | w(1,0,2,1,1,2) = {3,4}:{5,6} | w(2,0,2,1,1,2) = {7}:{8} |
| w(0,1,2,1,1,2) = {9}:{11} | w(1,1,2,1,1,2) = {5,8}:{7,10} | w(2,1,2,1,1,2) = {3,4}:{7,10} |
| w(0,2,2,1,1,2) = {2,3}:{9,11} | w(1,2,2,1,1,2) = {5,10}:{8,9} | w(2,2,2,1,1,2) = {3,4}:{8,9} |
| w(0,0,0,2,1,2) = {10}:{11} | w(1,0,0,2,1,2) = {1,11}:{3,4} | w(2,0,0,2,1,2) = {9}:{11} |
| w(0,1,0,2,1,2) = {1}:{11} | w(1,1,0,2,1,2) = {8}:{10} | w(2,1,0,2,1,2) = {3}:{4} |
| w(0,2,0,2,1,2) = {5}:{11} | w(1,2,0,2,1,2) = {10}:{8} | w(2,2,0,2,1,2) = {3}:{4} |
| w(0,0,1,2,1,2) = {5,6}:{10,11} | w(1,0,1,2,1,2) = {7}:{11} | w(2,0,1,2,1,2) = {4}:{5} |
| w(0,1,1,2,1,2) = {2,3}:{9,11} | w(1,1,1,2,1,2) = {2}:{11} | w(2,1,1,2,1,2) = {2}:{1} |
| w(0,2,1,2,1,2) = {5}:{6} | w(1,2,1,2,1,2) = {1}:{11} | w(2,2,1,2,1,2) = {1}:{2} |
| w(0,0,2,2,1,2) = {4,5}:{6,7} | w(1,0,2,2,1,2) = {10}:{11} | w(2,0,2,2,1,2) = {4}:{6} |
| w(0,1,2,2,1,2) = {2,3}:{10,11} | w(1,1,2,2,1,2) = {2}:{11} | w(2,1,2,2,1,2) = {2}:{1} |
| w(0,2,2,2,1,2) = {5}:{6} | w(1,2,2,2,1,2) = {1}:{11} | w(2,2,2,2,1,2) = {1}:{2} |
| w(0,0,0,0,2,2) = { }:{ } | w(1,0,0,0,2,2) = {8}:{11} | w(2,0,0,0,2,2) = {8}:{11} |
| w(0,1,0,0,2,2) = {2}:{11} | w(1,1,0,0,2,2) = {1}:{7} | w(2,1,0,0,2,2) = {1}:{10} |
| w(0,2,0,0,2,2) = {1}:{7} | w(1,2,0,0,2,2) = {2}:{9} | w(2,2,0,0,2,2) = {2}:{8} |
| w(0,0,1,0,2,2) = {10}:{11} | w(1,0,1,0,2,2) = {1,11}:{4,7} | w(2,0,1,0,2,2) = {6,8}:{10,11} |

w(0,1,1,0,2,2) = {2,3}:{9,11}    w(1,1,1,0,2,2) = {8}:{10}         w(2,1,1,0,2,2) = {3}:{9}
w(0,2,1,0,2,2) = {2,3}:{10,11}   w(1,2,1,0,2,2) = {10}:{8}         w(2,2,1,0,2,2) = {3}:{7}
w(0,0,2,0,2,2) = {9}:{11}        w(1,0,2,0,2,2) = {2,11}:{4,10}    w(2,0,2,0,2,2) = {5,9}:{7,11}
w(0,1,2,0,2,2) = {2,3}:{10,11}   w(1,1,2,0,2,2) = {7}:{9}          w(2,1,2,0,2,2) = {3}:{10}
w(0,2,2,0,2,2) = {2,3}:{9,11}    w(1,2,2,0,2,2) = {9}:{7}          w(2,2,2,0,2,2) = {3}:{8}
w(0,0,0,1,2,2) = {10}:{11}       w(1,0,0,1,2,2) = {9}:{11}         w(2,0,0,1,2,2) = {2}:{4}
w(0,1,0,1,2,2) = {1,4}:{10,11}   w(1,1,0,1,2,2) = {1,2}:{3,11}     w(2,1,0,1,2,2) = {9}:{10}
w(0,2,0,1,2,2) = {4,6}:{9,11}    w(1,2,0,1,2,2) = {1,2}:{3,11}     w(2,2,0,1,2,2) = {7}:{8}
w(0,0,1,1,2,2) = {6,7}:{5,8}     w(1,0,1,1,2,2) = {4}:{5}          w(2,0,1,1,2,2) = {4}:{5}
w(0,1,1,1,2,2) = {10}:{11}       w(1,1,1,1,2,2) = {6,7}:{9,11}     w(2,1,1,1,2,2) = {3,5}:{4,8}
w(0,2,1,1,2,2) = {2,3}:{10,11}   w(1,2,1,1,2,2) = {6,9}:{7,11}     w(2,2,1,1,2,2) = {3,5}:{4,10}
w(0,0,2,1,2,2) = {5}:{7}         w(1,0,2,1,2,2) = {4}:{6}          w(2,0,2,1,2,2) = {4}:{6}
w(0,1,2,1,2,2) = {9}:{11}        w(1,1,2,1,2,2) = {5,8}:{10,11}    w(2,1,2,1,2,2) = {3,6}:{4,7}
w(0,2,2,1,2,2) = {2,3}:{9,11}    w(1,2,2,1,2,2) = {5,10}:{8,11}    w(2,2,2,1,2,2) = {3,6}:{4,9}
w(0,0,0,2,2,2) = {9}:{10}        w(1,0,0,2,2,2) = {2}:{11}         w(2,0,0,2,2,2) = {3,4}:{9,11}
w(0,1,0,2,2,2) = {1,4}:{9,11}    w(1,1,0,2,2,2) = {1}:{8}          w(2,1,0,2,2,2) = {6}:{11}
w(0,2,0,2,2,2) = {4,5}:{10,11}   w(1,2,0,2,2,2) = {2}:{10}         w(2,2,0,2,2,2) = {6}:{11}
w(0,0,1,2,2,2) = {5,6}:{8,9}     w(1,0,1,2,2,2) = {10}:{11}        w(2,0,1,2,2,2) = {4}:{5}
w(0,1,1,2,2,2) = {2,3}:{9,11}    w(1,1,1,2,2,2) = {2,3}:{4,8}      w(2,1,1,2,2,2) = {1}:{8}
w(0,2,1,2,2,2) = {5}:{6}         w(1,2,1,2,2,2) = {1,3}:{4,10}     w(2,2,1,2,2,2) = {2}:{10}
w(0,0,2,2,2,2) = {6}:{7}         w(1,0,2,2,2,2) = {7}:{11}         w(2,0,2,2,2,2) = {4}:{6}
w(0,1,2,2,2,2) = {2,3}:{10,11}   w(1,1,2,2,2,2) = {2,3}:{4,7}      w(2,1,2,2,2,2) = {1}:{7}
w(0,2,2,2,2,2) = {5}:{6}         w(1,2,2,2,2,2) = {1,3}:{4,9}      w(2,2,2,2,2,2) = {2}:{9}

List 5

The completed map of the third algorithm to sort 11 coins

| | | |
|---|---|---|
| f(0,0,0,0,0,0) = 2047 | f(1,0,0,0,0,0) = 1468 | f(2,0,0,0,0,0) = 579 |
| f(0,1,0,0,0,0) = 1855 | f(1,1,0,0,0,0) = 2046 | f(2,1,0,0,0,0) = 2 |
| f(0,2,0,0,0,0) = 192 | f(1,2,0,0,0,0) = 2045 | f(2,2,0,0,0,0) = 1 |
| f(1,0,1,0,0,0) = 1786 | f(2,0,1,0,0,0) = 134 | f(0,1,1,0,0,0) = 1974 |
| f(1,1,1,0,0,0) = 944 | f(2,1,1,0,0,0) = 1551 | f(0,2,1,0,0,0) = 137 |
| f(1,2,1,0,0,0) = 752 | f(2,2,1,0,0,0) = 1167 | f(1,0,2,0,0,0) = 1913 |
| f(2,0,2,0,0,0) = 261 | f(0,1,2,0,0,0) = 1910 | f(1,1,2,0,0,0) = 880 |
| f(2,1,2,0,0,0) = 1295 | f(0,2,2,0,0,0) = 73 | f(1,2,2,0,0,0) = 496 |
| f(2,2,2,0,0,0) = 1103 | f(1,0,0,1,0,0) = 2044 | f(2,0,0,1,0,0) = 1991 |
| f(0,1,0,1,0,0) = 1838 | f(1,1,0,1,0,0) = 300 | f(2,1,0,1,0,0) = 774 |
| f(0,2,0,1,0,0) = 225 | f(1,2,0,1,0,0) = 108 | f(2,2,0,1,0,0) = 197 |
| f(0,0,1,1,0,0) = 1664 | f(1,0,1,1,0,0) = 696 | f(2,0,1,1,0,0) = 2031 |
| f(0,1,1,1,0,0) = 1674 | f(1,1,1,1,0,0) = 1594 | f(2,1,1,1,0,0) = 1827 |
| f(0,2,1,1,0,0) = 128 | f(1,2,1,1,0,0) = 1209 | f(2,2,1,1,0,0) = 1251 |
| f(0,0,2,1,0,0) = 639 | f(1,0,2,1,0,0) = 376 | f(2,0,2,1,0,0) = 2015 |
| f(0,1,2,1,0,0) = 1354 | f(1,1,2,1,0,0) = 1338 | f(2,1,2,1,0,0) = 1811 |
| f(0,2,2,1,0,0) = 64 | f(1,2,2,1,0,0) = 1145 | f(2,2,2,1,0,0) = 1235 |
| f(1,0,0,2,0,0) = 56 | f(2,0,0,2,0,0) = 3 | f(0,1,0,2,0,0) = 1822 |
| f(1,1,0,2,0,0) = 1850 | f(2,1,0,2,0,0) = 1939 | f(0,2,0,2,0,0) = 209 |
| f(1,2,0,2,0,0) = 1273 | f(2,2,0,2,0,0) = 1747 | f(0,0,1,2,0,0) = 1408 |
| f(1,0,1,2,0,0) = 32 | f(2,0,1,2,0,0) = 1671 | f(0,1,1,2,0,0) = 1983 |
| f(1,1,1,2,0,0) = 812 | f(2,1,1,2,0,0) = 902 | f(0,2,1,2,0,0) = 693 |
| f(1,2,1,2,0,0) = 236 | f(2,2,1,2,0,0) = 709 | f(0,0,2,2,0,0) = 383 |
| f(1,0,2,2,0,0) = 16 | f(2,0,2,2,0,0) = 1351 | f(0,1,2,2,0,0) = 1919 |
| f(1,1,2,2,0,0) = 796 | f(2,1,2,2,0,0) = 838 | f(0,2,2,2,0,0) = 373 |
| f(1,2,2,2,0,0) = 220 | f(2,2,2,2,0,0) = 453 | f(0,0,0,0,1,0) = 1984 |
| f(1,0,0,0,1,0) = 1534 | f(2,0,0,0,1,0) = 66 | f(0,1,0,0,1,0) = 1804 |
| f(1,1,0,0,1,0) = 1854 | f(2,1,0,0,1,0) = 1943 | f(0,2,0,0,1,0) = 245 |
| f(1,2,0,0,1,0) = 1277 | f(2,2,0,0,1,0) = 1751 | f(0,0,1,0,1,0) = 1162 |
| f(2,0,1,0,1,0) = 710 | f(1,1,1,0,1,0) = 1832 | f(2,1,1,0,1,0) = 1702 |
| f(1,2,1,0,1,0) = 1256 | f(2,2,1,0,1,0) = 1701 | f(0,0,2,0,1,0) = 74 |
| f(2,0,2,0,1,0) = 837 | f(1,1,2,0,1,0) = 1816 | f(2,1,2,0,1,0) = 1366 |
| f(1,2,2,0,1,0) = 1240 | f(2,2,2,0,1,0) = 1365 | f(0,0,0,1,1,0) = 1134 |
| f(2,0,0,1,1,0) = 388 | f(1,1,0,1,1,0) = 1662 | f(2,1,0,1,1,0) = 919 |
| f(1,2,0,1,1,0) = 1469 | f(2,2,0,1,1,0) = 727 | f(0,0,1,1,1,0) = 1700 |
| f(1,0,1,1,1,0) = 648 | f(2,0,1,1,1,0) = 1699 | f(0,1,1,1,1,0) = 1930 |
| f(1,1,1,1,1,0) = 1714 | f(2,1,1,1,1,0) = 1538 | f(1,2,1,1,1,0) = 1713 |
| f(2,2,1,1,1,0) = 1153 | f(0,0,2,1,1,0) = 1634 | f(1,0,2,1,1,0) = 328 |
| f(2,0,2,1,1,0) = 1363 | f(0,1,2,1,1,0) = 1866 | f(1,1,2,1,1,0) = 1394 |
| f(2,1,2,1,1,0) = 1282 | f(1,2,2,1,1,0) = 1393 | f(2,2,2,1,1,0) = 1089 |
| f(1,0,0,2,1,0) = 392 | f(1,1,0,2,1,0) = 2034 | f(2,1,0,2,1,0) = 1799 |
| f(1,2,0,2,1,0,0) = 2033 | f(2,2,0,2,1,0,0) = 1223 | f(0,0,1,2,1,0,0) = 1428 |
| f(1,0,1,2,1,0,0) = 1064 | f(2,0,1,2,1,0,0) = 1639 | f(1,1,1,2,1,0,0) = 520 |
| f(2,1,1,2,1,0,0) = 806 | f(0,2,1,2,1,0,0) = 757 | f(1,2,1,2,1,0,0) = 136 |
| f(2,2,1,2,1,0,0) = 229 | f(0,0,2,2,1,0,0) = 1362 | f(1,0,2,2,1,0,0) = 1048 |
| f(2,0,2,2,1,0,0) = 1623 | f(1,1,2,2,1,0,0) = 264 | f(2,1,2,2,1,0,0) = 790 |
| f(0,2,2,2,1,0,0) = 501 | f(1,2,2,2,1,0,0) = 72 | f(2,2,2,2,1,0,0) = 213 |
| f(0,0,0,0,2,0,0) = 63 | f(1,0,0,0,2,0,0) = 1981 | f(2,0,0,0,2,0,0) = 513 |
| f(0,1,0,0,2,0,0) = 1802 | f(1,1,0,0,2,0,0) = 296 | f(2,1,0,0,2,0,0) = 770 |
| f(0,2,0,0,2,0,0) = 243 | f(1,2,0,0,2,0,0) = 104 | f(2,2,0,0,2,0,0) = 193 |
| f(0,0,1,0,2,0,0) = 1973 | f(1,0,1,0,2,0,0) = 1210 | f(1,1,1,0,2,0,0) = 682 |
| f(2,1,1,0,2,0,0) = 807 | f(1,2,1,0,2,0,0) = 681 | f(2,2,1,0,2,0,0) = 231 |
| f(0,0,2,0,2,0,0) = 885 | f(1,0,2,0,2,0,0) = 1337 | f(1,1,2,0,2,0,0) = 346 |
| f(2,1,2,0,2,0,0) = 791 | f(1,2,2,0,2,0,0) = 345 | f(2,2,2,0,2,0,0) = 215 |
| f(2,0,0,1,2,0,0) = 1655 | f(1,1,0,1,2,0,0) = 824 | f(2,1,0,1,2,0,0) = 14 |
| f(1,2,0,1,2,0,0) = 248 | f(2,2,0,1,2,0,0) = 13 | f(0,0,1,1,2,0,0) = 685 |
| f(1,0,1,1,2,0,0) = 424 | f(2,0,1,1,2,0,0) = 999 | f(0,1,1,1,2,0,0) = 1546 |
| f(1,1,1,1,2,0,0) = 1834 | f(2,1,1,1,2,0,0) = 1975 | f(1,2,1,1,2,0,0) = 1257 |
| f(2,2,1,1,2,0,0) = 1783 | f(0,0,2,1,2,0,0) = 619 | f(1,0,2,1,2,0,0) = 408 |
| f(2,0,2,1,2,0,0) = 983 | f(0,1,2,1,2,0,0) = 1290 | f(1,1,2,1,2,0,0) = 1818 |
| f(2,1,2,1,2,0,0) = 1911 | f(1,2,2,1,2,0,0) = 1241 | f(2,2,2,1,2,0,0) = 1527 |
| f(0,0,0,2,2,0,0) = 913 | f(1,0,0,2,2,0,0) = 1659 | f(1,1,0,2,2,0,0) = 1320 |
| f(2,1,0,2,2,0,0) = 578 | f(1,2,0,2,2,0,0) = 1128 | f(2,2,0,2,2,0,0) = 385 |
| f(0,0,1,2,2,0,0) = 413 | f(1,0,1,2,2,0,0) = 684 | f(2,0,1,2,2,0,0) = 1719 |
| f(1,1,1,2,2,0,0) = 958 | f(2,1,1,2,2,0,0) = 654 | f(0,2,1,2,2,0,0) = 181 |
| f(1,2,1,2,2,0,0) = 765 | f(2,2,1,2,2,0,0) = 653 | f(0,0,2,2,2,0,0) = 347 |
| f(1,0,2,2,2,0,0) = 348 | f(2,0,2,2,2,0,0) = 1399 | f(1,1,2,2,2,0,0) = 894 |
| f(2,1,2,2,2,0,0) = 334 | f(0,2,2,2,2,0,0) = 117 | f(1,2,2,2,2,0,0) = 509 |
| f(2,2,2,2,2,0,0) = 333 | f(1,0,0,0,0,1,0) = 1456 | f(2,0,0,0,0,1,0) = 1987 |
| f(0,1,0,0,0,1,0) = 822 | f(1,1,0,0,0,1,0) = 816 | f(2,1,0,0,0,1,0) = 582 |
| f(0,2,0,0,0,1,0) = 246 | f(1,2,0,0,0,1,0) = 240 | f(2,2,0,0,0,1,0) = 389 |
| f(0,0,1,0,0,1,0) = 1684 | f(1,0,1,0,0,1,0) = 1716 | f(2,0,1,0,0,1,0) = 1126 |
| f(0,1,1,0,0,1,0) = 1920 | f(1,1,1,0,0,1,0) = 544 | f(2,1,1,0,0,1,0) = 514 |
| f(0,2,1,0,0,1,0) = 1737 | f(1,2,1,0,0,1,0) = 160 | f(2,2,1,0,0,1,0) = 129 |
| f(0,0,2,0,0,1,0) = 1971 | f(1,0,2,0,0,1,0) = 1396 | f(2,0,2,0,0,1,0) = 1557 |
| f(0,1,2,0,0,1,0) = 1856 | f(1,1,2,0,0,1,0) = 272 | f(2,1,2,0,0,1,0) = 258 |
| f(0,2,2,0,0,1,0) = 1481 | f(1,2,2,0,0,1,0) = 80 | f(2,2,2,0,0,1,0) = 65 |
| f(0,0,0,1,0,1,0) = 418 | f(1,0,0,1,0,1,0) = 2032 | f(2,0,0,1,0,1,0) = 1604 |
| f(0,1,0,1,0,1,0) = 1454 | f(1,1,0,1,0,1,0) = 1022 | f(2,1,0,1,0,1,0) = 1951 |
| f(0,2,0,1,0,1,0) = 1262 | f(1,2,0,1,0,1,0) = 1021 | f(2,2,0,1,0,1,0) = 1759 |
| f(0,0,1,1,0,1,0) = 174 | f(1,0,1,1,0,1,0) = 2024 | f(2,0,1,1,0,1,0) = 1766 |
| f(0,1,1,1,0,1,0) = 650 | f(1,1,1,1,0,1,0) = 1978 | f(2,1,1,1,0,1,0) = 1575 |
| f(0,2,1,1,0,1,0) = 1728 | f(1,2,1,1,0,1,0) = 1785 | f(2,2,1,1,0,1,0) = 1191 |
| f(0,0,2,1,0,1,0) = 34 | f(1,0,2,1,0,1,0) = 2008 | f(2,0,2,1,0,1,0) = 1877 |
| f(0,1,2,1,0,1,0) = 330 | f(1,1,2,1,0,1,0) = 1914 | f(2,1,2,1,0,1,0) = 1303 |

f(0,2,2,1,0,1,0) = 1472   f(1,2,2,1,0,1,0) = 1529   f(2,2,2,1,0,1,0) = 1111
f(0,0,0,2,0,1,0) = 365    f(1,0,0,2,0,1,0) = 1080   f(2,0,0,2,0,1,0) = 1411
f(0,1,0,2,0,1,0) = 1630   f(1,1,0,2,0,1,0) = 1458   f(2,1,0,2,0,1,0) = 1891
f(0,2,0,2,0,1,0) = 1246   f(1,2,0,2,0,1,0) = 1649   f(2,2,0,2,0,1,0) = 1507
f(0,0,1,2,0,1,0) = 158    f(1,0,1,2,0,1,0) = 608    f(2,0,1,2,0,1,0) = 663
f(0,1,1,2,0,1,0) = 1929   f(1,1,1,2,0,1,0) = 556    f(2,1,1,2,0,1,0) = 1924
f(0,2,1,2,0,1,0) = 691    f(1,2,1,2,0,1,0) = 172    f(2,2,1,2,0,1,0) = 1732
f(0,0,2,2,0,1,0) = 18     f(1,0,2,2,0,1,0) = 592    f(2,0,2,2,0,1,0) = 359
f(0,1,2,2,0,1,0) = 1865   f(1,1,2,2,0,1,0) = 284    f(2,1,2,2,0,1,0) = 1860
f(0,2,2,2,0,1,0) = 371    f(1,2,2,2,0,1,0) = 92     f(2,2,2,2,0,1,0) = 1476
f(0,0,0,0,1,1) = 1024     f(1,0,0,0,1,1,0) = 1522   f(2,0,0,0,1,1,0) = 1474
f(0,1,0,0,1,1,0) = 780    f(1,1,0,0,1,1,0) = 1552   f(2,1,0,0,1,1,0) = 1998
f(0,2,0,0,1,1,0) = 1013   f(1,2,0,0,1,1,0) = 1168   f(2,2,0,0,1,1,0) = 1997
f(0,0,1,0,1,1,0) = 1518   f(1,0,1,0,1,1,0) = 1004   f(2,0,1,0,1,1,0) = 1668
f(0,1,1,0,1,1,0) = 1572   f(1,1,1,0,1,1,0) = 1972   f(2,1,1,0,1,1,0) = 1670
f(0,2,1,0,1,1,0) = 1197   f(1,2,1,0,1,1,0) = 1780   f(2,2,1,0,1,1,0) = 1669
f(0,0,2,0,1,1,0) = 354    f(1,0,2,0,1,1,0) = 988    f(2,0,2,0,1,1,0) = 1348
f(0,1,2,0,1,1,0) = 1316   f(1,1,2,0,1,1,0) = 1908   f(2,1,2,0,1,1,0) = 1350
f(0,2,2,0,1,1,0) = 1133   f(1,2,2,0,1,1,0) = 1524   f(2,2,2,0,1,1,0) = 1349
f(1,0,0,1,1,1,0) = 1656   f(2,0,0,1,1,1,0) = 1478   f(0,1,0,1,1,1,0) = 1006
f(1,1,0,1,1,1,0) = 538    f(2,1,0,1,1,1,0) = 1942   f(0,2,0,1,1,1,0) = 228
f(1,2,0,1,1,1,0) = 153    f(2,2,0,1,1,1,0) = 1749   f(0,0,1,1,1,1,0) = 1774
f(1,0,1,1,1,1,0) = 1000   f(2,0,1,1,1,1,0) = 1730   f(0,1,1,1,1,1,0) = 906
f(1,1,1,1,1,1,0) = 1584   f(2,1,1,1,1,1,0) = 1666   f(0,2,1,1,1,1,0) = 1188
f(1,2,1,1,1,1,0) = 1200   f(2,2,1,1,1,1,0) = 1665   f(0,0,2,1,1,1,0) = 576
f(1,0,2,1,1,1,0) = 984    f(2,0,2,1,1,1,0) = 1857   f(0,1,2,1,1,1,0) = 842
f(1,1,2,1,1,1,0) = 1328   f(2,1,2,1,1,1,0) = 1346   f(0,2,2,1,1,1,0) = 1124
f(1,2,2,1,1,1,0) = 1136   f(2,2,2,1,1,1,0) = 1345   f(0,0,0,2,1,1,0) = 356
f(1,0,0,2,1,1,0) = 1530   f(2,0,0,2,1,1,0) = 1419   f(0,1,0,2,1,1,0) = 990
f(1,1,0,2,1,1,0) = 1840   f(2,1,0,2,1,1,0) = 1793   f(0,2,0,2,1,1,0) = 212
f(1,2,0,2,1,1,0) = 1264   f(2,2,0,2,1,1,0) = 1218   f(0,0,1,2,1,1,0) = 1502
f(1,0,1,2,1,1,0) = 1640   f(2,0,1,2,1,1,0) = 647    f(0,1,1,2,1,1,0) = 1581
f(1,1,1,2,1,1,0) = 1836   f(2,1,1,2,1,1,0) = 1638   f(0,2,1,2,1,1,0) = 755
f(1,2,1,2,1,1,0) = 1260   f(2,2,1,2,1,1,0) = 1445   f(0,0,2,2,1,1,0) = 338
f(1,0,2,2,1,1,0) = 1624   f(2,0,2,2,1,1,0) = 327    f(0,1,2,2,1,1,0) = 1325
f(1,1,2,2,1,1,0) = 1820   f(2,1,2,2,1,1,0) = 1430   f(0,2,2,2,1,1,0) = 499
f(1,2,2,2,1,1,0) = 1244   f(2,2,2,2,1,1,0) = 1621   f(0,0,0,0,2,1) = 1087
f(1,0,0,0,2,1,0) = 1969   f(2,0,0,0,2,1,0) = 1921   f(0,1,0,0,2,1,0) = 778
f(1,1,0,0,2,1,0) = 1018   f(2,1,0,0,2,1,0) = 915    f(0,2,0,0,2,1,0) = 1011
f(1,2,0,0,2,1,0) = 1017   f(2,2,0,0,2,1,0) = 723    f(0,0,1,0,2,1,0) = 1453
f(1,0,1,0,2,1,0) = 1778   f(2,0,1,0,2,1,0) = 643    f(0,1,1,0,2,1,0) = 1570
f(1,1,1,0,2,1,0) = 810    f(2,1,1,0,2,1,0) = 422    f(0,2,1,0,2,1,0) = 1195
f(1,2,1,0,2,1,0) = 233    f(2,2,1,0,2,1,0) = 613    f(0,0,2,0,2,1,0) = 289
f(1,0,2,0,2,1,0) = 1905   f(2,0,2,0,2,1,0) = 323    f(0,1,2,0,2,1,0) = 1314

| | | |
|---|---|---|
| f(1,1,2,0,2,1,0) = 794 | f(2,1,2,0,2,1,0) = 598 | f(0,2,2,0,2,1,0) = 1131 |
| f(1,2,2,0,2,1,0) = 217 | f(2,2,2,0,2,1,0) = 405 | f(0,0,0,1,2,1,0) = 427 |
| f(1,0,0,1,2,1,0) = 1992 | f(2,0,0,1,2,1,0) = 1094 | f(0,1,0,1,2,1,0) = 1438 |
| f(1,1,0,1,2,1,0) = 924 | f(2,1,0,1,2,1,0) = 1038 | f(0,2,0,1,2,1,0) = 226 |
| f(1,2,0,1,2,1,0) = 732 | f(2,2,0,1,2,1,0) = 1037 | f(0,0,1,1,2,1,0) = 640 |
| f(1,0,1,1,2,1,0) = 1448 | f(2,0,1,1,2,1,0) = 486 | f(0,1,1,1,2,1,0) = 522 |
| f(1,1,1,1,2,1,0) = 1824 | f(2,1,1,1,2,1,0) = 1591 | f(0,2,1,1,2,1,0) = 1186 |
| f(1,2,1,1,2,1,0) = 1248 | f(2,2,1,1,2,1,0) = 1207 | f(0,0,2,1,2,1,0) = 1067 |
| f(1,0,2,1,2,1,0) = 1432 | f(2,0,2,1,2,1,0) = 917 | f(0,1,2,1,2,1,0) = 266 |
| f(1,1,2,1,2,1,0) = 1808 | f(2,1,2,1,2,1,0) = 1335 | f(0,2,2,1,2,1,0) = 1122 |
| f(1,2,2,1,2,1,0) = 1232 | f(2,2,2,1,2,1,0) = 1143 | f(0,0,0,2,2,1,0) = 1953 |
| f(1,0,0,2,2,1,0) = 1146 | f(2,0,0,2,2,1,0) = 1611 | f(0,1,0,2,2,1,0) = 1646 |
| f(1,1,0,2,2,1,0) = 1896 | f(2,1,0,2,2,1,0) = 1986 | f(0,2,0,2,2,1,0) = 210 |
| f(1,2,0,2,2,1,0) = 1512 | f(2,2,0,2,2,1,0) = 1985 | f(0,0,1,2,2,1,0) = 393 |
| f(1,0,1,2,2,1,0) = 553 | f(2,0,1,2,2,1,0) = 1687 | f(0,1,1,2,2,1,0) = 1579 |
| f(1,1,1,2,2,1,0) = 956 | f(2,1,1,2,2,1,0) = 1678 | f(0,2,1,2,2,1,0) = 179 |
| f(1,2,1,2,2,1,0) = 764 | f(2,2,1,2,2,1,0) = 1677 | f(0,0,2,2,2,1,0) = 1051 |
| f(1,0,2,2,2,1,0) = 90 | f(2,0,2,2,2,1,0) = 1383 | f(0,1,2,2,2,1,0) = 1323 |
| f(1,1,2,2,2,1,0) = 892 | f(2,1,2,2,2,1,0) = 1358 | f(0,2,2,2,2,1,0) = 115 |
| f(1,2,2,2,2,1,0) = 508 | f(2,2,2,2,2,1,0) = 1357 | f(1,0,0,0,0,2,0) = 60 |
| f(2,0,0,0,0,2,0) = 591 | f(0,1,0,0,0,2,0) = 1801 | f(1,1,0,0,0,2,0) = 1658 |
| f(2,1,0,0,0,2,0) = 1807 | f(0,2,0,0,0,2,0) = 1225 | f(1,2,0,0,0,2,0) = 1465 |
| f(2,2,0,0,0,2,0) = 1231 | f(0,0,1,0,0,2,0) = 76 | f(1,0,1,0,0,2,0) = 490 |
| f(2,0,1,0,0,2,0) = 651 | f(0,1,1,0,0,2,0) = 566 | f(1,1,1,0,0,2,0) = 1982 |
| f(2,1,1,0,0,2,0) = 1967 | f(0,2,1,0,0,2,0) = 191 | f(1,2,1,0,0,2,0) = 1789 |
| f(2,2,1,0,0,2,0) = 1775 | f(0,0,2,0,0,2,0) = 363 | f(1,0,2,0,0,2,0) = 921 |
| f(2,0,2,0,0,2,0) = 331 | f(0,1,2,0,0,2,0) = 310 | f(1,1,2,0,0,2,0) = 1918 |
| f(2,1,2,0,0,2,0) = 1887 | f(0,2,2,0,0,2,0) = 127 | f(1,2,2,0,0,2,0) = 1533 |
| f(2,2,2,0,0,2,0) = 1503 | f(0,0,0,1,0,2,0) = 1682 | f(1,0,0,1,0,2,0) = 636 |
| f(2,0,0,1,0,2,0) = 967 | f(0,1,0,1,0,2,0) = 801 | f(1,1,0,1,0,2,0) = 540 |
| f(2,1,0,1,0,2,0) = 398 | f(0,2,0,1,0,2,0) = 417 | f(1,2,0,1,0,2,0) = 156 |
| f(2,2,0,1,0,2,0) = 589 | f(0,0,1,1,0,2,0) = 2029 | f(1,0,1,1,0,2,0) = 1688 |
| f(2,0,1,1,0,2,0) = 1455 | f(0,1,1,1,0,2,0) = 1676 | f(1,1,1,1,0,2,0) = 571 |
| f(2,1,1,1,0,2,0) = 1955 | f(0,2,1,1,0,2,0) = 182 | f(1,2,1,1,0,2,0) = 187 |
| f(2,2,1,1,0,2,0) = 1763 | f(0,0,2,1,0,2,0) = 1889 | f(1,0,2,1,0,2,0) = 1384 |
| f(2,0,2,1,0,2,0) = 1439 | f(0,1,2,1,0,2,0) = 1356 | f(1,1,2,1,0,2,0) = 315 |
| f(2,1,2,1,0,2,0) = 1875 | f(0,2,2,1,0,2,0) = 118 | f(1,2,2,1,0,2,0) = 123 |
| f(2,2,2,1,0,2,0) = 1491 | f(0,0,0,2,0,2,0) = 1629 | f(1,0,0,2,0,2,0) = 443 |
| f(2,0,0,2,0,2,0) = 15 | f(0,1,0,2,0,2,0) = 785 | f(1,1,0,2,0,2,0) = 288 |
| f(2,1,0,2,0,2,0) = 1026 | f(0,2,0,2,0,2,0) = 593 | f(1,2,0,2,0,2,0) = 96 |
| f(2,2,0,2,0,2,0) = 1025 | f(0,0,1,2,0,2,0) = 2013 | f(1,0,1,2,0,2,0) = 170 |
| f(2,0,1,2,0,2,0) = 39 | f(0,1,1,2,0,2,0) = 575 | f(1,1,1,2,0,2,0) = 936 |
| f(2,1,1,2,0,2,0) = 518 | f(0,2,1,2,0,2,0) = 1717 | f(1,2,1,2,0,2,0) = 744 |
| f(2,2,1,2,0,2,0) = 133 | f(0,0,2,2,0,2,0) = 1873 | f(1,0,2,2,0,2,0) = 281 |
| f(2,0,2,2,0,2,0) = 23 | f(0,1,2,2,0,2,0) = 319 | f(1,1,2,2,0,2,0) = 856 |
| f(2,1,2,2,0,2,0) = 262 | f(0,2,2,2,0,2,0) = 1397 | f(1,2,2,2,0,2,0) = 472 |
| f(2,2,2,2,0,2,0) = 69 | f(0,0,0,0,1,2) = 960 | f(1,0,0,0,1,2,0) = 126 |
| f(2,0,0,0,1,2,0) = 78 | f(0,1,0,0,1,2,0) = 1036 | f(1,1,0,0,1,2,0) = 1324 |
| f(2,1,0,0,1,2,0) = 1030 | f(0,2,0,0,1,2,0) = 1269 | f(1,2,0,0,1,2,0) = 1132 |
| f(2,2,0,0,1,2,0) = 1029 | f(0,0,1,0,1,2,0) = 1758 | f(1,0,1,0,1,2,0) = 1724 |
| f(2,0,1,0,1,2,0) = 142 | f(0,1,1,0,1,2,0) = 916 | f(1,1,1,0,1,2,0) = 1642 |
| f(2,1,1,0,1,2,0) = 1830 | f(0,2,1,0,1,2,0) = 733 | f(1,2,1,0,1,2,0) = 1449 |
| f(2,2,1,0,1,2,0) = 1253 | f(0,0,2,0,1,2,0) = 594 | f(1,0,2,0,1,2,0) = 1404 |
| f(2,0,2,0,1,2,0) = 269 | f(0,1,2,0,1,2,0) = 852 | f(1,1,2,0,1,2,0) = 1434 |
| f(2,1,2,0,1,2,0) = 1814 | f(0,2,2,0,1,2,0) = 477 | f(1,2,2,0,1,2,0) = 1625 |
| f(2,2,2,0,1,2,0) = 1237 | f(0,0,0,1,1,2,0) = 94 | f(1,0,0,1,1,2,0) = 436 |
| f(2,0,0,1,1,2,0) = 901 | f(0,1,0,1,1,2,0) = 1837 | f(1,1,0,1,1,2,0) = 62 |
| f(2,1,0,1,1,2,0) = 535 | f(0,2,0,1,1,2,0) = 401 | f(1,2,0,1,1,2,0) = 61 |
| f(2,2,0,1,1,2,0) = 151 | f(0,0,1,1,1,2,0) = 996 | f(1,0,1,1,1,2,0) = 664 |
| f(2,0,1,1,1,2,0) = 1957 | f(0,1,1,1,1,2,0) = 1932 | f(1,1,1,1,1,2,0) = 690 |
| f(2,1,1,1,1,2,0) = 1539 | f(0,2,1,1,1,2,0) = 724 | f(1,2,1,1,1,2,0) = 689 |
| f(2,2,1,1,1,2,0) = 1155 | f(0,0,2,1,1,2,0) = 1654 | f(1,0,2,1,1,2,0) = 360 |
| f(2,0,2,1,1,2,0) = 1494 | f(0,1,2,1,1,2,0) = 1868 | f(1,1,2,1,1,2,0) = 370 |
| f(2,1,2,1,1,2,0) = 1283 | f(0,2,2,1,1,2,0) = 468 | f(1,2,2,1,1,2,0) = 369 |
| f(2,2,2,1,1,2,0) = 1091 | f(0,0,0,2,1,2,0) = 1620 | f(1,0,0,2,1,2,0) = 953 |
| f(2,0,0,2,1,2,0) = 55 | f(0,1,0,2,1,2,0) = 1821 | f(1,1,0,2,1,2,0) = 1010 |
| f(2,1,0,2,1,2,0) = 1315 | f(0,2,0,2,1,2,0) = 609 | f(1,2,0,2,1,2,0) = 1009 |
| f(2,2,0,2,1,2,0) = 1123 | f(0,0,1,2,1,2,0) = 980 | f(1,0,1,2,1,2,0) = 1130 |
| f(2,0,1,2,1,2,0) = 615 | f(0,1,1,2,1,2,0) = 925 | f(1,1,1,2,1,2,0) = 904 |
| f(2,1,1,2,1,2,0) = 815 | f(0,2,1,2,1,2,0) = 1781 | f(1,2,1,2,1,2,0) = 712 |
| f(2,2,1,2,1,2,0) = 239 | f(0,0,2,2,1,2,0) = 1407 | f(1,0,2,2,1,2,0) = 1561 |
| f(2,0,2,2,1,2,0) = 599 | f(0,1,2,2,1,2,0) = 861 | f(1,1,2,2,1,2,0) = 840 |
| f(2,1,2,2,1,2,0) = 799 | f(0,2,2,2,1,2,0) = 1525 | f(1,2,2,2,1,2,0) = 456 |
| f(2,2,2,2,1,2,0) = 223 | f(0,0,0,0,2,2) = 1023 | f(1,0,0,0,2,2,0) = 573 |
| f(2,0,0,0,2,2,0) = 525 | f(0,1,0,0,2,2,0) = 1034 | f(1,1,0,0,2,2,0) = 50 |
| f(2,1,0,0,2,2,0) = 879 | f(0,2,0,0,2,2,0) = 1267 | f(1,2,0,0,2,2,0) = 49 |
| f(2,2,0,0,2,2,0) = 495 | f(0,0,1,0,2,2,0) = 1693 | f(1,0,1,0,2,2,0) = 699 |
| f(2,0,1,0,2,2,0) = 1059 | f(0,1,1,0,2,2,0) = 914 | f(1,1,1,0,2,2,0) = 698 |
| f(2,1,1,0,2,2,0) = 523 | f(0,2,1,0,2,2,0) = 731 | f(1,2,1,0,2,2,0) = 697 |
| f(2,2,1,0,2,2,0) = 139 | f(0,0,2,0,2,2,0) = 529 | f(1,0,2,0,2,2,0) = 379 |
| f(2,0,2,0,2,2,0) = 1043 | f(0,1,2,0,2,2,0) = 850 | f(1,1,2,0,2,2,0) = 378 |
| f(2,1,2,0,2,2,0) = 267 | f(0,2,2,0,2,2,0) = 475 | f(1,2,2,0,2,2,0) = 377 |
| f(2,2,2,0,2,2,0) = 75 | f(0,0,0,1,2,2,0) = 1691 | f(1,0,0,1,2,2,0) = 628 |
| f(2,0,0,1,2,2,0) = 517 | f(0,1,0,1,2,2,0) = 1835 | f(1,1,0,1,2,2,0) = 829 |
| f(2,1,0,1,2,2,0) = 783 | f(0,2,0,1,2,2,0) = 1057 | f(1,2,0,1,2,2,0) = 254 |
| f(2,2,0,1,2,2,0) = 207 | f(0,0,1,1,2,2,0) = 1709 | f(1,0,1,1,2,2,0) = 1720 |
| f(2,0,1,1,2,2,0) = 423 | f(0,1,1,1,2,2,0) = 1548 | f(1,1,1,1,2,2,0) = 426 |
| f(2,1,1,1,2,2,0) = 803 | f(0,2,1,1,2,2,0) = 722 | f(1,2,1,1,2,2,0) = 617 |
| f(2,2,1,1,2,2,0) = 227 | f(0,0,2,1,2,2,0) = 545 | f(1,0,2,1,2,2,0) = 1400 |
| f(2,0,2,1,2,2,0) = 407 | f(0,1,2,1,2,2,0) = 1292 | f(1,1,2,1,2,2,0) = 602 |
| f(2,1,2,1,2,2,0) = 787 | f(0,2,2,1,2,2,0) = 466 | f(1,2,2,1,2,2,0) = 409 |
| f(2,2,2,1,2,2,0) = 211 | f(1,0,0,2,2,2,0) = 569 | f(2,0,0,2,2,2,0) = 391 |
| f(0,1,0,2,2,2,0) = 1819 | f(1,1,0,2,2,2,0) = 298 | f(2,1,0,2,2,2,0) = 1894 |
| f(0,2,0,2,2,2,0) = 1041 | f(1,2,0,2,2,2,0) = 105 | f(2,2,0,2,2,2,0) = 1509 |
| f(0,0,1,2,2,2,0) = 1471 | f(1,0,1,2,2,2,0) = 190 | f(2,0,1,2,2,2,0) = 1063 |
| f(0,1,1,2,2,2,0) = 923 | f(1,1,1,2,2,2,0) = 702 | f(2,1,1,2,2,2,0) = 911 |
| f(0,2,1,2,2,2,0) = 1205 | f(1,2,1,2,2,2,0) = 701 | f(2,2,1,2,2,2,0) = 719 |
| f(0,0,2,2,2,2,0) = 273 | f(1,0,2,2,2,2,0) = 317 | f(2,0,2,2,2,2,0) = 1047 |
| f(0,1,2,2,2,2,0) = 859 | f(1,1,2,2,2,2,0) = 382 | f(2,1,2,2,2,2,0) = 847 |
| f(0,2,2,2,2,2,0) = 1141 | f(1,2,2,2,2,2,0) = 381 | f(2,2,2,2,2,2,0) = 463 |
| f(2,0,0,0,0,0,1) = 1999 | f(0,1,0,0,0,0,1) = 768 | f(1,1,0,0,0,0,1) = 1470 |
| f(2,1,0,0,0,0,1) = 1847 | f(1,2,0,0,0,0,1) = 1661 | f(2,2,0,0,0,0,1) = 1271 |
| f(0,0,1,0,0,0,1) = 420 | f(2,0,1,0,0,0,1) = 1711 | f(0,1,1,0,0,0,1) = 512 |
| f(1,1,1,0,0,0,1) = 928 | f(2,1,1,0,0,0,1) = 898 | f(1,2,1,0,0,0,1) = 736 |
| f(2,2,1,0,0,0,1) = 705 | f(0,0,2,0,0,0,1) = 883 | f(2,0,2,0,0,0,1) = 1375 |
| f(0,1,2,0,0,0,1) = 256 | f(1,1,2,0,0,0,1) = 848 | f(2,1,2,0,0,0,1) = 834 |
| f(1,2,2,0,0,0,1) = 464 | f(2,2,2,0,0,0,1) = 449 | f(0,0,0,1,0,0,1) = 1698 |
| f(1,0,0,1,0,0,1) = 624 | f(2,0,0,1,0,0,1) = 2039 | f(1,1,0,1,0,0,1) = 876 |
| f(2,1,0,1,0,0,1) = 781 | f(0,2,0,1,0,0,1) = 2017 | f(1,2,0,1,0,0,1) = 492 |
| f(2,2,0,1,0,0,1) = 206 | f(0,0,1,1,0,0,1) = 716 | f(1,0,1,1,0,0,1) = 688 |
| f(2,0,1,1,0,0,1) = 1447 | f(0,1,1,1,0,0,1) = 1710 | f(1,1,1,1,0,0,1) = 1568 |
| f(2,1,1,1,0,0,1) = 1831 | f(1,2,1,1,0,0,1) = 1184 | f(2,2,1,1,0,0,1) = 1255 |
| f(0,0,2,1,0,0,1) = 1643 | f(1,0,2,1,0,0,1) = 368 | f(2,0,2,1,0,0,1) = 1431 |
| f(0,1,2,1,0,0,1) = 1390 | f(1,1,2,1,0,0,1) = 1296 | f(2,1,2,1,0,0,1) = 1815 |
| f(1,2,2,1,0,0,1) = 1104 | f(2,2,2,1,0,0,1) = 1239 | f(0,0,0,2,0,0,1) = 1645 |
| f(1,0,0,2,0,0,1) = 2043 | f(1,1,0,2,0,0,1) = 2042 | f(2,1,0,2,0,0,1) = 771 |
| f(0,2,0,2,0,0,1) = 2001 | f(1,2,0,2,0,0,1) = 2041 | f(2,2,0,2,0,0,1) = 195 |
| f(0,0,1,2,0,0,1) = 460 | f(1,0,1,2,0,0,1) = 1068 | f(2,0,1,2,0,0,1) = 695 |
| f(0,1,1,2,0,0,1) = 521 | f(1,1,1,2,0,0,1) = 1980 | f(2,1,1,2,0,0,1) = 516 |
| f(0,2,1,2,0,0,1) = 673 | f(1,2,1,2,0,0,1) = 1788 | f(2,2,1,2,0,0,1) = 132 |
| f(0,0,2,2,0,0,1) = 1371 | f(1,0,2,2,0,0,1) = 1052 | f(2,0,2,2,0,0,1) = 375 |
| f(0,1,2,2,0,0,1) = 265 | f(1,1,2,2,0,0,1) = 1916 | f(2,1,2,2,0,0,1) = 260 |
| f(0,2,2,2,0,0,1) = 353 | f(1,2,2,2,0,0,1) = 1532 | f(2,2,2,2,0,0,1) = 68 |
| f(0,0,0,0,1,0,1) = 2038 | f(2,0,0,0,1,0,1) = 1486 | f(0,1,0,0,1,0,1) = 1996 |
| f(1,1,0,0,1,0,1) = 1848 | f(2,1,0,0,1,0,1) = 1798 | f(0,2,0,0,1,0,1) = 204 |
| f(1,2,0,0,1,0,1) = 1272 | f(2,2,0,0,1,0,1) = 1221 | f(1,0,1,0,1,0,1) = 1452 |
| f(2,0,1,0,1,0,1) = 1742 | f(0,1,1,0,1,0,1) = 1556 | f(1,1,1,0,1,0,1) = 1960 |
| f(2,1,1,0,1,0,1) = 1958 | f(0,2,1,0,1,0,1) = 1181 | f(1,2,1,0,1,0,1) = 1768 |
| f(2,2,1,0,1,0,1) = 1765 | f(0,0,2,0,1,0,1) = 1098 | f(1,0,2,0,1,0,1) = 1436 |
| f(2,0,2,0,1,0,1) = 1869 | f(0,1,2,0,1,0,1) = 1300 | f(1,1,2,0,1,0,1) = 1880 |
| f(2,1,2,0,1,0,1) = 1878 | f(0,2,2,0,1,0,1) = 1117 | f(1,2,2,0,1,0,1) = 1496 |
| f(2,2,2,0,1,0,1) = 1493 | f(0,0,0,1,1,0,1) = 714 | f(1,0,0,1,1,0,1) = 1016 |
| f(2,0,0,1,1,0,1) = 1988 | f(0,1,0,1,1,0,1) = 1420 | f(1,1,0,1,1,0,1) = 1564 |
| f(2,1,0,1,1,0,1) = 918 | f(0,2,0,1,1,0,1) = 237 | f(1,2,0,1,1,0,1) = 1180 |
| f(2,2,0,1,1,0,1) = 725 | f(0,0,1,1,1,0,1) = 1718 | f(1,0,1,1,1,0,1) = 1692 |
| f(2,0,1,1,1,0,1) = 2019 | f(0,1,1,1,1,0,1) = 1966 | f(1,1,1,1,1,0,1) = 1970 |
| f(2,1,1,1,1,0,1) = 1923 | f(0,2,1,1,1,0,1) = 1172 | f(1,2,1,1,1,0,1) = 1777 |
| f(2,2,1,1,1,0,1) = 1731 | f(0,0,2,1,1,0,1) = 1600 | f(1,0,2,1,1,0,1) = 1388 |
| f(2,0,2,1,1,0,1) = 2003 | f(0,1,2,1,1,0,1) = 1902 | f(1,1,2,1,1,0,1) = 1906 |
| f(2,1,2,1,1,0,1) = 1859 | f(0,2,2,1,1,0,1) = 1108 | f(1,2,2,1,1,0,1) = 1521 |
| f(2,2,2,1,1,0,1) = 1475 | f(0,0,0,2,1,0,1) = 1636 | f(1,0,0,2,1,0,1) = 1464 |
| f(2,0,0,2,1,0,1) = 11 | f(0,1,0,2,1,0,1) = 1612 | f(1,1,0,2,1,0,1) = 864 |
| f(2,1,0,2,1,0,1) = 1903 | f(0,2,0,2,1,0,1) = 221 | f(1,2,0,2,1,0,1) = 480 |
| f(2,2,0,2,1,0,1) = 1519 | f(0,0,1,2,1,0,1) = 1462 | f(1,0,1,2,1,0,1) = 1056 |
| f(2,0,1,2,1,0,1) = 549 | f(0,1,1,2,1,0,1) = 1565 | f(1,1,1,2,1,0,1) = 1964 |
| f(2,1,1,2,1,0,1) = 998 | f(0,2,1,2,1,0,1) = 737 | f(1,2,1,2,1,0,1) = 1772 |
| f(2,2,1,2,1,0,1) = 997 | f(0,0,2,2,1,0,1) = 1344 | f(1,0,2,2,1,0,1) = 1040 |
| f(2,0,2,2,1,0,1) = 86 | f(0,1,2,2,1,0,1) = 1309 | f(1,1,2,2,1,0,1) = 1884 |
| f(2,1,2,2,1,0,1) = 982 | f(0,2,2,2,1,0,1) = 481 | f(1,2,2,2,1,0,1) = 1500 |
| f(2,2,2,2,1,0,1) = 981 | f(0,0,0,0,2,0,1) = 54 | f(2,0,0,0,2,0,1) = 1933 |
| f(0,1,0,0,2,0,1) = 1994 | f(1,1,0,0,2,0,1) = 872 | f(2,1,0,0,2,0,1) = 779 |
| f(0,2,0,0,2,0,1) = 202 | f(1,2,0,0,2,0,1) = 488 | f(2,2,0,0,2,0,1) = 203 |
| f(2,0,1,0,2,0,1) = 1703 | f(0,1,1,0,2,0,1) = 1554 | f(1,1,1,0,2,0,1) = 954 |
| f(2,1,1,0,2,0,1) = 903 | f(0,2,1,0,2,0,1) = 1179 | f(1,2,1,0,2,0,1) = 761 |
| f(2,2,1,0,2,0,1) = 711 | f(0,0,2,0,2,0,1) = 1909 | f(2,0,2,0,2,0,1) = 1367 |
| f(0,1,2,0,2,0,1) = 1298 | f(1,1,2,0,2,0,1) = 890 | f(2,1,2,0,2,0,1) = 839 |
| f(0,2,2,0,2,0,1) = 1115 | f(1,2,2,0,2,0,1) = 505 | f(2,2,2,0,2,0,1) = 455 |
| f(0,0,0,1,2,0,1) = 1707 | f(1,0,0,1,2,0,1) = 584 | f(2,0,0,1,2,0,1) = 1541 |
| f(0,1,0,1,2,0,1) = 394 | f(1,1,0,1,2,0,1) = 828 | f(2,1,0,1,2,0,1) = 782 |
| f(0,2,0,1,2,0,1) = 1259 | f(1,2,0,1,2,0,1) = 252 | f(2,2,0,1,2,0,1) = 205 |
| f(0,0,1,1,2,0,1) = 1673 | f(1,0,1,1,2,0,1) = 416 | f(2,0,1,1,2,0,1) = 166 |
| f(0,1,1,1,2,0,1) = 1582 | f(1,1,1,1,2,0,1) = 1962 | f(2,1,1,1,2,0,1) = 1446 |
| f(0,2,1,1,2,0,1) = 1170 | f(1,2,1,1,2,0,1) = 1769 | f(2,2,1,1,2,0,1) = 1637 |
| f(0,0,2,1,2,0,1) = 1003 | f(1,0,2,1,2,0,1) = 400 | f(2,0,2,1,2,0,1) = 277 |
| f(0,1,2,1,2,0,1) = 1326 | f(1,1,2,1,2,0,1) = 1882 | f(2,1,2,1,2,0,1) = 1622 |
| f(0,2,2,1,2,0,1) = 1106 | f(1,2,2,1,2,0,1) = 1497 | f(2,2,2,1,2,0,1) = 1429 |
| f(0,0,0,2,2,0,1) = 1589 | f(2,0,0,2,2,0,1) = 971 | f(0,1,0,2,2,0,1) = 586 |
| f(1,1,0,2,2,0,1) = 1312 | f(2,1,0,2,2,0,1) = 1602 | f(0,2,0,2,2,0,1) = 1243 |
| f(1,2,0,2,2,0,1) = 1120 | f(2,2,0,2,2,0,1) = 1409 | f(0,0,1,2,2,0,1) = 1417 |
| f(1,0,1,2,2,0,1) = 1790 | f(2,0,1,2,2,0,1) = 1174 | f(0,1,1,2,2,0,1) = 1563 |
| f(1,1,1,2,2,0,1) = 952 | f(2,1,1,2,2,0,1) = 910 | f(0,2,1,2,2,0,1) = 161 |
| f(1,2,1,2,2,0,1) = 760 | f(2,2,1,2,2,0,1) = 717 | f(0,0,2,2,2,0,1) = 987 |
| f(1,0,2,2,2,0,1) = 1917 | f(2,0,2,2,2,0,1) = 1317 | f(0,1,2,2,2,0,1) = 1307 |
| f(1,1,2,2,2,0,1) = 888 | f(2,1,2,2,2,0,1) = 846 | f(0,2,2,2,2,0,1) = 97 |
| f(1,2,2,2,2,0,1) = 504 | f(2,2,2,2,2,0,1) = 461 | f(0,0,0,0,0,1,1) = 1078 |
| f(1,0,0,0,0,1,1) = 1072 | f(2,0,0,0,0,1,1) = 1603 | f(0,1,0,0,0,1,1) = 1846 |
| f(1,1,0,0,0,1,1) = 912 | f(2,1,0,0,0,1,1) = 966 | f(0,2,0,0,0,1,1) = 1270 |
| f(1,2,0,0,0,1,1) = 720 | f(2,2,0,0,0,1,1) = 965 | f(0,0,1,0,0,1,1) = 1444 |
| f(1,0,1,0,0,1,1) = 1696 | f(2,0,1,0,0,1,1) = 1158 | f(0,1,1,0,0,1,1) = 1536 |
| f(1,1,1,0,0,1,1) = 800 | f(2,1,1,0,0,1,1) = 642 | f(0,2,1,0,0,1,1) = 1161 |
| f(1,2,1,0,0,1,1) = 224 | f(2,2,1,0,0,1,1) = 641 | f(0,0,2,0,0,1,1) = 1907 |
| f(1,0,2,0,0,1,1) = 1360 | f(2,0,2,0,0,1,1) = 1285 | f(0,1,2,0,0,1,1) = 1280 |
| f(1,1,2,0,0,1,1) = 784 | f(2,1,2,0,0,1,1) = 322 | f(0,2,2,0,0,1,1) = 1097 |
| f(1,2,2,0,0,1,1) = 208 | f(2,2,2,0,0,1,1) = 321 | f(0,0,0,1,0,1,1) = 1442 |
| f(1,0,0,1,0,1,1) = 1648 | f(2,0,0,1,0,1,1) = 1028 | f(0,1,0,1,0,1,1) = 1070 |
| f(1,1,0,1,0,1,1) = 830 | f(2,1,0,1,0,1,1) = 1558 | f(0,2,0,1,0,1,1) = 1249 |
| f(1,2,0,1,0,1,1) = 253 | f(2,2,0,1,0,1,1) = 1173 | f(0,0,1,1,0,1,1) = 750 |
| f(1,0,1,1,0,1,1) = 2016 | f(2,0,1,1,0,1,1) = 2023 | f(0,1,1,1,0,1,1) = 686 |
| f(1,1,1,1,0,1,1) = 1952 | f(2,1,1,1,0,1,1) = 1686 | f(0,2,1,1,0,1,1) = 1152 |
| f(1,2,1,1,0,1,1) = 1760 | f(2,2,1,1,0,1,1) = 1685 | f(0,0,2,1,0,1,1) = 994 |
| f(1,0,2,1,0,1,1) = 2000 | f(2,0,2,1,0,1,1) = 2007 | f(0,1,2,1,0,1,1) = 366 |
| f(1,1,2,1,0,1,1) = 1872 | f(2,1,2,1,0,1,1) = 1382 | f(0,2,2,1,0,1,1) = 1088 |
| f(1,2,2,1,0,1,1) = 1488 | f(2,2,2,1,0,1,1) = 1381 | f(0,0,0,2,0,1,1) = 1389 |
| f(2,0,0,2,0,1,1) = 1027 | f(0,1,0,2,0,1,1) = 1054 | f(1,1,0,2,0,1,1) = 1946 |
| f(2,1,0,2,0,1,1) = 1895 | f(0,2,0,2,0,1,1) = 1233 | f(1,2,0,2,0,1,1) = 1753 |
| f(2,2,0,2,0,1,1) = 1511 | f(0,0,1,2,0,1,1) = 478 | f(1,0,1,2,0,1,1) = 1644 |
| f(2,0,1,2,0,1,1) = 1750 | f(0,1,1,2,0,1,1) = 1545 | f(1,1,1,2,0,1,1) = 1596 |
| f(2,1,1,2,0,1,1) = 1540 | f(0,2,1,2,0,1,1) = 1715 | f(1,2,1,2,0,1,1) = 1212 |
| f(2,2,1,2,0,1,1) = 1156 | f(0,0,2,2,0,1,1) = 978 | f(1,0,2,2,0,1,1) = 1628 |
| f(2,0,2,2,0,1,1) = 1893 | f(0,1,2,2,0,1,1) = 1289 | f(1,1,2,2,0,1,1) = 1340 |
| f(2,1,2,2,0,1,1) = 1284 | f(0,2,2,2,0,1,1) = 1395 | f(1,2,2,2,0,1,1) = 1148 |
| f(2,2,2,2,0,1,1) = 1092 | f(1,0,0,0,1,1,1) = 1138 | f(2,0,0,0,1,1,1) = 1090 |
| f(0,1,0,0,1,1,1) = 972 | f(1,1,0,0,1,1,1) = 1944 | f(2,1,0,0,1,1,1) = 1796 |
| f(0,2,0,0,1,1,1) = 2037 | f(1,2,0,0,1,1,1) = 1752 | f(2,2,0,0,1,1,1) = 1220 |
| f(1,0,1,0,1,1,1) = 2028 | f(2,0,1,0,1,1,1) = 644 | f(0,1,1,0,1,1,1) = 1956 |
| f(1,1,1,0,1,1,1) = 1976 | f(0,2,1,0,1,1,1) = 1773 | f(1,2,1,0,1,1,1) = 1784 |
| f(0,0,2,0,1,1,1) = 1378 | f(1,0,2,0,1,1,1) = 2012 | f(2,0,2,0,1,1,1) = 324 |
| f(0,1,2,0,1,1,1) = 1892 | f(1,1,2,0,1,1,1) = 1912 | f(0,2,2,0,1,1,1) = 1517 |
| f(1,2,2,0,1,1,1) = 1528 | f(0,0,0,1,1,1,1) = 1738 | f(1,0,0,1,1,1,1) = 2040 |
| f(2,0,0,1,1,1,1) = 1412 | f(0,1,0,1,1,1,1) = 2030 | f(1,1,0,1,1,1,1) = 1560 |
| f(2,1,0,1,1,1,1) = 772 | f(0,2,0,1,1,1,1) = 1252 | f(1,2,0,1,1,1,1) = 1176 |
| f(2,2,0,1,1,1,1) = 196 | f(0,0,1,1,1,1,1) = 1740 | f(1,0,1,1,1,1,1) = 992 |
| f(2,0,1,1,1,1,1) = 1154 | f(0,1,1,1,1,1,1) = 942 | f(1,1,1,1,1,1,1) = 1968 |
| f(2,1,1,1,1,1,1) = 1922 | f(0,2,1,1,1,1,1) = 1764 | f(1,2,1,1,1,1,1) = 1776 |
| f(2,2,1,1,1,1,1) = 1729 | f(0,0,2,1,1,1,1) = 2018 | f(1,0,2,1,1,1,1) = 976 |
| f(2,0,2,1,1,1,1) = 1281 | f(0,1,2,1,1,1,1) = 878 | f(1,1,2,1,1,1,1) = 1904 |
| f(2,1,2,1,1,1,1) = 1858 | f(0,2,2,1,1,1,1) = 1508 | f(1,2,2,1,1,1,1) = 1520 |
| f(2,2,2,1,1,1,1) = 1473 | f(0,0,0,2,1,1,1) = 1380 | f(1,0,0,2,1,1,1) = 1416 |
| f(2,0,0,2,1,1,1) = 1035 | f(0,1,0,2,1,1,1) = 2014 | f(1,1,0,2,1,1,1) = 1936 |
| f(2,1,0,2,1,1,1) = 1797 | f(0,2,0,2,1,1,1) = 1236 | f(1,2,0,2,1,1,1) = 1744 |
| f(2,2,0,2,1,1,1) = 1222 | f(0,0,1,2,1,1,1) = 1484 | f(1,0,1,2,1,1,1) = 1632 |
| f(2,0,1,2,1,1,1) = 1573 | f(0,1,1,2,1,1,1) = 1965 | f(1,1,1,2,1,1,1) = 1544 |
| f(2,1,1,2,1,1,1) = 2022 | f(0,2,1,2,1,1,1) = 1779 | f(1,2,1,2,1,1,1) = 1160 |
| f(2,2,1,2,1,1,1) = 2021 | f(0,0,2,2,1,1,1) = 2002 | f(1,0,2,2,1,1,1) = 1616 |
| f(2,0,2,2,1,1,1) = 1110 | f(0,1,2,2,1,1,1) = 1901 | f(1,1,2,2,1,1,1) = 1288 |
| f(2,1,2,2,1,1,1) = 2006 | f(0,2,2,2,1,1,1) = 1523 | f(1,2,2,2,1,1,1) = 1096 |
| f(2,2,2,2,1,1,1) = 2005 | f(1,0,0,0,2,1,1) = 1585 | f(2,0,0,0,2,1,1) = 1537 |
| f(0,1,0,0,2,1,1) = 970 | f(1,1,0,0,2,1,1) = 922 | f(2,1,0,0,2,1,1) = 773 |
| f(0,2,0,0,2,1,1) = 2035 | f(1,2,0,0,2,1,1) = 729 | f(2,2,0,0,2,1,1) = 198 |
| f(1,0,1,0,2,1,1) = 1202 | f(2,0,1,0,2,1,1) = 1667 | f(0,1,1,0,2,1,1) = 1954 |
| f(1,1,1,0,2,1,1) = 938 | f(2,1,1,0,2,1,1) = 935 | f(0,2,1,0,2,1,1) = 1771 |
| f(1,2,1,0,2,1,1) = 745 | f(2,2,1,0,2,1,1) = 743 | f(0,0,2,0,2,1,1) = 1313 |
| f(1,0,2,0,2,1,1) = 1329 | f(2,0,2,0,2,1,1) = 1347 | f(0,1,2,0,2,1,1) = 1890 |
| f(1,1,2,0,2,1,1) = 858 | f(2,1,2,0,2,1,1) = 855 | f(0,2,2,0,2,1,1) = 1515 |
| f(1,2,2,0,2,1,1) = 473 | f(2,2,2,0,2,1,1) = 471 | f(0,0,0,1,2,1,1) = 1451 |
| f(1,0,0,1,2,1,1) = 1608 | f(2,0,0,1,2,1,1) = 70 | f(0,1,0,1,2,1,1) = 1418 |
| f(1,1,0,1,2,1,1) = 920 | f(2,1,0,1,2,1,1) = 1422 | f(0,2,0,1,2,1,1) = 1250 |
| f(1,2,0,1,2,1,1) = 728 | f(2,2,0,1,2,1,1) = 1613 | f(0,0,1,1,2,1,1) = 1005 |
| f(1,0,1,1,2,1,1) = 1440 | f(2,0,1,1,2,1,1) = 1510 | f(0,1,1,1,2,1,1) = 558 |
| f(1,1,1,1,2,1,1) = 1578 | f(2,1,1,1,2,1,1) = 1062 | f(0,2,1,1,2,1,1) = 1762 |
| f(1,2,1,1,2,1,1) = 1193 | f(2,2,1,1,2,1,1) = 1061 | f(0,0,2,1,2,1,1) = 2027 |
| f(1,0,2,1,2,1,1) = 1424 | f(2,0,2,1,2,1,1) = 1941 | f(0,1,2,1,2,1,1) = 302 |
| f(1,1,2,1,2,1,1) = 1306 | f(2,1,2,1,2,1,1) = 1046 | f(0,2,2,1,2,1,1) = 1506 |
| f(1,2,2,1,2,1,1) = 1113 | f(2,2,2,1,2,1,1) = 1045 | f(0,0,0,2,2,1,1) = 1937 |
| f(1,0,0,2,2,1,1) = 1083 | f(2,0,0,2,2,1,1) = 1995 | f(0,1,0,2,2,1,1) = 1610 |
| f(1,1,0,2,2,1,1) = 1800 | f(2,1,0,2,2,1,1) = 962 | f(0,2,0,2,2,1,1) = 1234 |
| f(1,2,0,2,2,1,1) = 1224 | f(2,2,0,2,2,1,1) = 961 | f(0,0,1,2,2,1,1) = 384 |
| f(1,0,1,2,2,1,1) = 1708 | f(2,0,1,2,2,1,1) = 1683 | f(0,1,1,2,2,1,1) = 1963 |
| f(1,1,1,2,2,1,1) = 948 | f(2,1,1,2,2,1,1) = 1934 | f(0,2,1,2,2,1,1) = 1203 |
| f(1,2,1,2,2,1,1) = 756 | f(2,2,1,2,2,1,1) = 1741 | f(0,0,2,2,2,1,1) = 2011 |
| f(1,0,2,2,2,1,1) = 1372 | f(2,0,2,2,2,1,1) = 1379 | f(0,1,2,2,2,1,1) = 1899 |
| f(1,1,2,2,2,1,1) = 884 | f(2,1,2,2,2,1,1) = 1870 | f(0,2,2,2,2,1,1) = 1139 |
| f(1,2,2,2,2,1,1) = 500 | f(2,2,2,2,2,1,1) = 1485 | f(0,0,0,0,0,2,1) = 1033 |
| f(1,0,0,0,0,2,1) = 1084 | f(2,0,0,0,0,2,1) = 1615 | f(0,1,0,0,0,2,1) = 1792 |
| f(2,1,0,0,0,2,1) = 1805 | f(0,2,0,0,0,2,1) = 1216 | f(2,2,0,0,0,2,1) = 1230 |
| f(0,0,1,0,0,2,1) = 1100 | f(1,0,1,0,0,2,1) = 1514 | f(2,0,1,0,0,2,1) = 1675 |
| f(0,1,1,0,0,2,1) = 1590 | f(1,1,1,0,0,2,1) = 1690 | f(2,1,1,0,0,2,1) = 1935 |
| f(0,2,1,0,0,2,1) = 1215 | f(1,2,1,0,0,2,1) = 1689 | f(2,2,1,0,0,2,1) = 1743 |
| f(0,0,2,0,0,2,1) = 1387 | f(1,0,2,0,0,2,1) = 1945 | f(2,0,2,0,0,2,1) = 1355 |
| f(0,1,2,0,0,2,1) = 1334 | f(1,1,2,0,0,2,1) = 1386 | f(2,1,2,0,0,2,1) = 1871 |
| f(0,2,2,0,0,2,1) = 1151 | f(1,2,2,0,0,2,1) = 1385 | f(2,2,2,0,0,2,1) = 1487 |
| f(0,0,0,1,0,2,1) = 1426 | f(1,0,0,1,0,2,1) = 1660 | f(2,0,0,1,0,2,1) = 1015 |
| f(0,1,0,1,0,2,1) = 1825 | f(1,1,0,1,0,2,1) = 1948 | f(2,1,0,1,0,2,1) = 974 |
| f(0,2,0,1,0,2,1) = 1441 | f(1,2,0,1,0,2,1) = 1756 | f(2,2,0,1,0,2,1) = 973 |
| f(1,0,1,1,0,2,1) = 1680 | f(2,0,1,1,0,2,1) = 933 | f(1,1,1,1,0,2,1) = 1595 |
| f(2,1,1,1,0,2,1) = 1959 | f(0,2,1,1,0,2,1) = 1206 | f(1,2,1,1,0,2,1) = 1211 |
| f(2,2,1,1,0,2,1) = 1767 | f(0,0,2,1,0,2,1) = 1609 | f(1,0,2,1,0,2,1) = 1376 |
| f(2,0,2,1,0,2,1) = 470 | f(1,1,2,1,0,2,1) = 1339 | f(2,1,2,1,0,2,1) = 1879 |
| f(0,2,2,1,0,2,1) = 1142 | f(1,2,2,1,0,2,1) = 1147 | f(2,2,2,1,0,2,1) = 1495 |
| f(0,0,0,2,0,2,1) = 1373 | f(1,0,0,2,0,2,1) = 1467 | f(2,0,0,2,0,2,1) = 1039 |
| f(0,1,0,2,0,2,1) = 1809 | f(1,1,0,2,0,2,1) = 1888 | f(2,1,0,2,0,2,1) = 1410 |
| f(0,2,0,2,0,2,1) = 1617 | f(1,2,0,2,0,2,1) = 1504 | f(2,2,0,2,0,2,1) = 1601 |
| f(1,0,1,2,0,2,1) = 746 | f(0,1,1,2,0,2,1) = 1599 | f(1,1,1,2,0,2,1) = 940 |
| f(2,1,1,2,0,2,1) = 678 | f(0,2,1,2,0,2,1) = 1697 | f(1,2,1,2,0,2,1) = 748 |
| f(2,2,1,2,0,2,1) = 677 | f(0,0,2,2,0,2,1) = 1353 | f(1,0,2,2,0,2,1) = 857 |
| f(0,1,2,2,0,2,1) = 1343 | f(1,1,2,2,0,2,1) = 860 | f(2,1,2,2,0,2,1) = 342 |
| f(0,2,2,2,0,2,1) = 1377 | f(1,2,2,2,0,2,1) = 476 | f(2,2,2,2,0,2,1) = 341 |
| f(1,0,0,0,1,2,1) = 1150 | f(2,0,0,0,1,2,1) = 1102 | f(1,1,0,0,1,2,1) = 1900 |
| f(2,1,0,0,1,2,1) = 1414 | f(0,2,0,0,1,2,1) = 1228 | f(1,2,0,0,1,2,1) = 1516 |
| f(2,2,0,0,1,2,1) = 1605 | f(1,0,1,0,1,2,1) = 1704 | f(2,0,1,0,1,2,1) = 1166 |
| f(0,1,1,0,1,2,1) = 1940 | f(1,1,1,0,1,2,1) = 2026 | f(2,1,1,0,1,2,1) = 1829 |
| f(0,2,1,0,1,2,1) = 1757 | f(1,2,1,0,1,2,1) = 2025 | f(2,2,1,0,1,2,1) = 1254 |
| f(0,0,2,0,1,2,1) = 1618 | f(1,0,2,0,1,2,1) = 1368 | f(2,0,2,0,1,2,1) = 1293 |
| f(0,1,2,0,1,2,1) = 1876 | f(1,1,2,0,1,2,1) = 2010 | f(2,1,2,0,1,2,1) = 1813 |
| f(0,2,2,0,1,2,1) = 1501 | f(1,2,2,0,1,2,1) = 2009 | f(2,2,2,0,1,2,1) = 1238 |
| f(0,0,0,1,1,2,1) = 1118 | f(1,0,0,1,1,2,1) = 1460 | f(2,0,0,1,1,2,1) = 1925 |
| f(0,1,0,1,1,2,1) = 1828 | f(1,1,0,1,1,2,1) = 58 | f(2,1,0,1,1,2,1) = 534 |
| f(0,2,0,1,1,2,1) = 1425 | f(1,2,0,1,1,2,1) = 57 | f(2,2,0,1,1,2,1) = 149 |
| f(0,0,1,1,1,2,1) = 2020 | f(1,0,1,1,1,2,1) = 656 | f(2,0,1,1,1,2,1) = 1443 |
| f(1,1,1,1,1,2,1) = 946 | f(2,1,1,1,1,2,1) = 1931 | f(0,2,1,1,1,2,1) = 1748 |
| f(1,2,1,1,1,2,1) = 753 | f(2,2,1,1,1,2,1) = 1739 | f(0,0,2,1,1,2,1) = 1058 |
| f(1,0,2,1,1,2,1) = 352 | f(2,0,2,1,1,2,1) = 1427 | f(1,1,2,1,1,2,1) = 882 |
| f(2,1,2,1,1,2,1) = 1867 | f(0,2,2,1,1,2,1) = 1492 | f(1,2,2,1,1,2,1) = 497 |
| f(2,2,2,1,1,2,1) = 1483 | f(0,0,0,2,1,2,1) = 1364 | f(1,0,0,2,1,2,1) = 440 |
| f(2,0,0,2,1,2,1) = 1079 | f(0,1,0,2,1,2,1) = 1812 | f(1,1,0,2,1,2,1) = 818 |
| f(2,1,0,2,1,2,1) = 1803 | f(0,2,0,2,1,2,1) = 1633 | f(1,2,0,2,1,2,1) = 241 |
| f(2,2,0,2,1,2,1) = 1227 | f(0,0,1,2,1,2,1) = 2004 | f(1,0,1,2,1,2,1) = 1194 |
| f(2,0,1,2,1,2,1) = 726 | f(0,1,1,2,1,2,1) = 1949 | f(1,1,1,2,1,2,1) = 1928 |
| f(2,1,1,2,1,2,1) = 805 | f(0,2,1,2,1,2,1) = 1761 | f(1,2,1,2,1,2,1) = 1736 |
| f(2,2,1,2,1,2,1) = 230 | f(0,0,2,2,1,2,1) = 1398 | f(1,0,2,2,1,2,1) = 1305 |
| f(2,0,2,2,1,2,1) = 869 | f(0,1,2,2,1,2,1) = 1885 | f(1,1,2,2,1,2,1) = 1864 |
| f(2,1,2,2,1,2,1) = 789 | f(0,2,2,2,1,2,1) = 1505 | f(1,2,2,2,1,2,1) = 1480 |
| f(2,2,2,2,1,2,1) = 214 | f(1,0,0,0,2,2,1) = 1597 | f(2,0,0,0,2,2,1) = 1549 |
| f(1,1,0,0,2,2,1) = 626 | f(2,1,0,0,2,2,1) = 870 | f(0,2,0,0,2,2,1) = 1226 |
| f(1,2,0,0,2,2,1) = 433 | f(2,2,0,0,2,2,1) = 485 | f(1,0,1,0,2,2,1) = 186 |
| f(2,0,1,0,2,2,1) = 1635 | f(0,1,1,0,2,2,1) = 1938 | f(1,1,1,0,2,2,1) = 570 |
| f(2,1,1,0,2,2,1) = 907 | f(0,2,1,0,2,2,1) = 1755 | f(1,2,1,0,2,2,1) = 185 |
| f(2,2,1,0,2,2,1) = 715 | f(0,0,2,0,2,2,1) = 1553 | f(1,0,2,0,2,2,1) = 313 |

f(2,0,2,0,2,2,1) = 1619
f(2,1,2,0,2,2,1) = 843
f(2,2,2,0,2,2,1) = 459
f(0,1,0,1,2,2,1) = 802
f(0,2,0,1,2,2,1) = 1459
f(0,0,1,1,2,2,1) = 1727
f(2,1,1,1,2,2,1) = 931
f(2,2,1,1,2,2,1) = 739
f(1,1,2,1,2,2,1) = 1626
f(1,2,2,1,2,2,1) = 1433
f(1,0,0,2,2,2,1) = 1593
f(1,1,0,2,2,2,1) = 1466
f(1,2,0,2,2,2,1) = 1657
f(1,0,1,2,2,2,1) = 1214
f(1,1,1,2,2,2,1) = 568
f(1,2,1,2,2,2,1) = 184
f(1,0,2,2,2,2,1) = 1341
f(1,1,2,2,2,2,1) = 312
f(1,2,2,2,2,2,1) = 120
f(1,1,0,0,0,0,2) = 776
f(1,2,0,0,0,0,2) = 200
f(1,0,1,0,0,0,2) = 672
f(0,2,1,0,0,0,2) = 1791
f(0,0,2,0,0,0,2) = 1627
f(2,1,2,0,0,0,2) = 1311
f(2,2,2,0,0,0,2) = 1119
f(0,1,0,1,0,0,2) = 46
f(1,2,0,1,0,0,2) = 1276
f(1,0,1,1,0,0,2) = 1672
f(1,1,1,1,0,0,2) = 1979
f(1,2,1,1,0,0,2) = 1787
f(1,0,2,1,0,0,2) = 1352
f(1,1,2,1,0,0,2) = 1915
f(1,2,2,1,0,0,2) = 1531
f(1,0,0,2,0,0,2) = 8
f(1,1,0,2,0,0,2) = 1841
f(2,2,0,2,0,0,2) = 1171
f(2,0,1,2,0,0,2) = 1679
f(0,2,1,2,0,0,2) = 657
f(0,0,2,2,0,0,2) = 1331
f(1,1,2,2,0,0,2) = 792
f(1,2,2,2,0,0,2) = 216
f(1,0,0,0,1,0,2) = 114
f(2,1,0,0,1,0,2) = 1559

f(0,1,2,0,2,2,1) = 1874
f(0,2,2,0,2,2,1) = 1499
f(0,0,0,1,2,2,1) = 1435
f(1,1,0,1,2,2,1) = 1853
f(1,2,0,1,2,2,1) = 1278
f(1,0,1,1,2,2,1) = 1712
f(0,2,1,1,2,2,1) = 1746
f(0,0,2,1,2,2,1) = 865
f(2,1,2,1,2,2,1) = 851
f(2,2,2,1,2,2,1) = 467
f(2,0,0,2,2,2,1) = 1415
f(2,1,0,2,2,2,1) = 1606
f(2,2,0,2,2,2,1) = 1413
f(2,0,1,2,2,2,1) = 150
f(2,1,1,2,2,2,1) = 646
f(2,2,1,2,2,2,1) = 645
f(2,0,2,2,2,2,1) = 293
f(2,1,2,2,2,2,1) = 326
f(2,2,2,2,2,2,1) = 325
f(2,1,0,0,0,0,2) = 386
f(2,2,0,0,0,0,2) = 577
f(1,1,1,0,0,0,2) = 1598
f(1,2,1,0,0,0,2) = 1213
f(1,0,2,0,0,0,2) = 336
f(0,2,2,0,0,0,2) = 1535
f(0,0,0,1,0,0,2) = 402
f(1,1,0,1,0,0,2) = 1852
f(2,2,0,1,0,0,2) = 5
f(2,0,1,1,0,0,2) = 995
f(2,1,1,1,0,0,2) = 515
f(2,2,1,1,0,0,2) = 131
f(2,0,2,1,0,0,2) = 979
f(2,1,2,1,0,0,2) = 259
f(2,2,2,1,0,0,2) = 67
f(2,0,0,2,0,0,2) = 1423
f(2,1,0,2,0,0,2) = 1555
f(0,0,1,2,0,0,2) = 404
f(1,1,1,2,0,0,2) = 808
f(1,2,1,2,0,0,2) = 232
f(1,0,2,2,0,0,2) = 600
f(2,1,2,2,0,0,2) = 863
f(2,2,2,2,0,0,2) = 479
f(0,1,0,0,1,0,2) = 1845
f(0,2,0,0,1,0,2) = 53

f(1,1,2,0,2,2,1) = 314
f(1,2,2,0,2,2,1) = 121
f(1,0,0,1,2,2,1) = 1652
f(2,1,0,1,2,2,1) = 590
f(2,2,0,1,2,2,1) = 397
f(1,1,1,1,2,2,1) = 1450
f(1,2,1,1,2,2,1) = 1641
f(1,0,2,1,2,2,1) = 1392
f(0,2,2,1,2,2,1) = 1490
f(0,0,0,2,2,2,1) = 565
f(0,1,0,2,2,2,1) = 786
f(0,2,0,2,2,2,1) = 1651
f(0,0,1,2,2,2,1) = 1437
f(0,1,1,2,2,2,1) = 1947
f(0,2,1,2,2,2,1) = 1185
f(0,0,2,2,2,2,1) = 849
f(0,1,2,2,2,2,1) = 1883
f(0,2,2,2,2,2,1) = 1121
f(1,0,0,0,0,0,2) = 48
f(0,2,0,0,0,0,2) = 1279
f(0,0,1,0,0,0,2) = 1164
f(2,1,1,0,0,0,2) = 1583
f(2,2,1,0,0,0,2) = 1199
f(1,1,2,0,0,0,2) = 1342
f(1,2,2,0,0,0,2) = 1149
f(2,0,0,1,0,0,2) = 4
f(2,1,0,1,0,0,2) = 6
f(0,0,1,1,0,0,2) = 676
f(0,1,1,1,0,0,2) = 1694
f(0,2,1,1,0,0,2) = 1782
f(0,0,2,1,0,0,2) = 1587
f(0,1,2,1,0,0,2) = 1374
f(0,2,2,1,0,0,2) = 1526
f(0,0,0,2,0,0,2) = 349
f(0,1,0,2,0,0,2) = 30
f(1,2,0,2,0,0,2) = 1266
f(1,0,1,2,0,0,2) = 616
f(2,1,1,2,0,0,2) = 943
f(2,2,1,2,0,0,2) = 751
f(2,0,2,2,0,0,2) = 1359
f(0,2,2,2,0,0,2) = 337
f(0,0,0,0,1,0,2) = 1993
f(1,1,0,0,1,0,2) = 1844
f(1,2,0,0,1,0,2) = 1268

| | | |
|---|---|---|
| f(2,2,0,0,1,0,2) = 1175 | f(0,0,1,0,1,0,2) = 138 | f(1,0,1,0,1,0,2) = 680 |
| f(0,1,1,0,1,0,2) = 932 | f(1,1,1,0,1,0,2) = 1592 | f(2,1,1,0,1,0,2) = 1542 |
| f(0,2,1,0,1,0,2) = 749 | f(1,2,1,0,1,0,2) = 1208 | f(2,2,1,0,1,0,2) = 1157 |
| f(1,0,2,0,1,0,2) = 344 | f(0,1,2,0,1,0,2) = 868 | f(1,1,2,0,1,0,2) = 1336 |
| f(2,1,2,0,1,0,2) = 1286 | f(0,2,2,0,1,0,2) = 493 | f(1,2,2,0,1,0,2) = 1144 |
| f(2,2,2,0,1,0,2) = 1093 | f(0,0,0,1,1,0,2) = 458 | f(1,0,0,1,1,0,2) = 1076 |
| f(0,1,0,1,1,0,2) = 804 | f(1,1,0,1,1,0,2) = 638 | f(2,1,0,1,1,0,2) = 927 |
| f(0,2,0,1,1,0,2) = 1461 | f(1,2,0,1,1,0,2) = 445 | f(2,2,0,1,1,0,2) = 735 |
| f(0,0,1,1,1,0,2) = 1060 | f(1,0,1,1,1,0,2) = 730 | f(2,0,1,1,1,0,2) = 130 |
| f(0,1,1,1,1,0,2) = 1950 | f(1,1,1,1,1,0,2) = 1586 | f(2,1,1,1,1,0,2) = 1543 |
| f(0,2,1,1,1,0,2) = 740 | f(1,2,1,1,1,0,2) = 1201 | f(2,2,1,1,1,0,2) = 1159 |
| f(0,0,2,1,1,0,2) = 630 | f(1,0,2,1,1,0,2) = 873 | f(2,0,2,1,1,0,2) = 257 |
| f(0,1,2,1,1,0,2) = 1886 | f(1,1,2,1,1,0,2) = 1330 | f(2,1,2,1,1,0,2) = 1287 |
| f(0,2,2,1,1,0,2) = 484 | f(1,2,2,1,1,0,2) = 1137 | f(2,2,2,1,1,0,2) = 1095 |
| f(0,0,0,2,1,0,2) = 340 | f(1,0,0,2,1,0,2) = 506 | f(2,0,0,2,1,0,2) = 1463 |
| f(0,1,0,2,1,0,2) = 788 | f(1,1,0,2,1,0,2) = 1842 | f(2,1,0,2,1,0,2) = 1795 |
| f(0,2,0,2,1,0,2) = 1653 | f(1,2,0,2,1,0,2) = 1265 | f(2,2,0,2,1,0,2) = 1219 |
| f(0,0,1,2,1,0,2) = 1044 | f(1,0,1,2,1,0,2) = 1770 | f(2,0,1,2,1,0,2) = 1647 |
| f(0,1,1,2,1,0,2) = 941 | f(1,1,1,2,1,0,2) = 618 | f(2,1,1,2,1,0,2) = 550 |
| f(0,2,1,2,1,0,2) = 721 | f(1,2,1,2,1,0,2) = 425 | f(2,2,1,2,1,0,2) = 165 |
| f(0,0,2,2,1,0,2) = 374 | f(1,0,2,2,1,0,2) = 1881 | f(2,0,2,2,1,0,2) = 1631 |
| f(0,1,2,2,1,0,2) = 877 | f(1,1,2,2,1,0,2) = 410 | f(2,1,2,2,1,0,2) = 278 |
| f(0,2,2,2,1,0,2) = 465 | f(1,2,2,2,1,0,2) = 601 | f(2,2,2,2,1,0,2) = 85 |
| f(0,0,0,0,2,0,2) = 9 | f(1,0,0,0,2,0,2) = 561 | f(0,1,0,0,2,0,2) = 1843 |
| f(1,1,0,0,2,0,2) = 826 | f(2,1,0,0,2,0,2) = 775 | f(0,2,0,0,2,0,2) = 51 |
| f(1,2,0,0,2,0,2) = 249 | f(2,2,0,0,2,0,2) = 199 | f(0,0,1,0,2,0,2) = 949 |
| f(1,0,1,0,2,0,2) = 178 | f(2,0,1,0,2,0,2) = 611 | f(0,1,1,0,2,0,2) = 930 |
| f(1,1,1,0,2,0,2) = 554 | f(2,1,1,0,2,0,2) = 551 | f(0,2,1,0,2,0,2) = 747 |
| f(1,2,1,0,2,0,2) = 169 | f(2,2,1,0,2,0,2) = 167 | f(1,0,2,0,2,0,2) = 305 |
| f(2,0,2,0,2,0,2) = 595 | f(0,1,2,0,2,0,2) = 866 | f(1,1,2,0,2,0,2) = 282 |
| f(2,1,2,0,2,0,2) = 279 | f(0,2,2,0,2,0,2) = 491 | f(1,2,2,0,2,0,2) = 89 |
| f(2,2,2,0,2,0,2) = 87 | f(0,0,0,1,2,0,2) = 411 | f(1,0,0,1,2,0,2) = 2036 |
| f(2,0,0,1,2,0,2) = 583 | f(0,1,0,1,2,0,2) = 1826 | f(1,1,0,1,2,0,2) = 528 |
| f(2,1,0,1,2,0,2) = 1567 | f(0,2,0,1,2,0,2) = 435 | f(1,2,0,1,2,0,2) = 144 |
| f(2,2,0,1,2,0,2) = 1183 | f(0,0,1,1,2,0,2) = 703 | f(1,0,1,1,2,0,2) = 1961 |
| f(2,0,1,1,2,0,2) = 1007 | f(0,1,1,1,2,0,2) = 1566 | f(1,1,1,1,2,0,2) = 1066 |
| f(2,1,1,1,2,0,2) = 547 | f(0,2,1,1,2,0,2) = 738 | f(1,2,1,1,2,0,2) = 1065 |
| f(2,2,1,1,2,0,2) = 163 | f(0,0,2,1,2,0,2) = 585 | f(1,0,2,1,2,0,2) = 1498 |
| f(2,0,2,1,2,0,2) = 991 | f(0,1,2,1,2,0,2) = 1310 | f(1,1,2,1,2,0,2) = 1050 |
| f(2,1,2,1,2,0,2) = 275 | f(0,2,2,1,2,0,2) = 482 | f(1,2,2,1,2,0,2) = 1049 |
| f(2,2,2,1,2,0,2) = 83 | f(0,0,0,2,2,0,2) = 1333 | f(1,0,0,2,2,0,2) = 59 |
| f(2,0,0,2,2,0,2) = 1031 | f(0,1,0,2,2,0,2) = 1810 | f(1,1,0,2,2,0,2) = 1322 |
| f(2,1,0,2,2,0,2) = 867 | f(0,2,0,2,2,0,2) = 627 | f(1,2,0,2,2,0,2) = 1129 |
| f(2,2,0,2,2,0,2) = 483 | f(0,0,1,2,2,0,2) = 447 | f(1,0,1,2,2,0,2) = 44 |
| f(2,0,1,2,2,0,2) = 659 | f(0,1,1,2,2,0,2) = 939 | f(1,1,1,2,2,0,2) = 572 |
| f(2,1,1,2,2,0,2) = 526 | f(0,2,1,2,2,0,2) = 145 | f(1,2,1,2,2,0,2) = 188 |
| f(2,2,1,2,2,0,2) = 141 | f(0,0,2,2,2,0,2) = 329 | f(1,0,2,2,2,0,2) = 28 |
| f(2,0,2,2,2,0,2) = 355 | f(0,1,2,2,2,0,2) = 875 | f(1,1,2,2,2,0,2) = 316 |
| f(2,1,2,2,2,0,2) = 270 | f(0,2,2,2,2,0,2) = 81 | f(1,2,2,2,2,0,2) = 124 |
| f(2,2,2,2,2,0,2) = 77 | f(0,0,0,0,0,1,2) = 1014 | f(1,0,0,0,0,1,2) = 432 |
| f(2,0,0,0,0,1,2) = 963 | f(0,1,0,0,0,1,2) = 831 | f(1,1,0,0,0,1,2) = 817 |
| f(0,2,0,0,0,1,2) = 255 | f(1,2,0,0,0,1,2) = 242 | f(0,0,1,0,0,1,2) = 660 |
| f(1,0,1,0,0,1,2) = 692 | f(2,0,1,0,0,1,2) = 102 | f(0,1,1,0,0,1,2) = 896 |
| f(1,1,1,0,0,1,2) = 560 | f(2,1,1,0,0,1,2) = 662 | f(0,2,1,0,0,1,2) = 713 |
| f(1,2,1,0,0,1,2) = 176 | f(2,2,1,0,0,1,2) = 661 | f(0,0,2,0,0,1,2) = 947 |
| f(1,0,2,0,0,1,2) = 372 | f(2,0,2,0,0,1,2) = 533 | f(0,1,2,0,0,1,2) = 832 |
| f(1,1,2,0,0,1,2) = 304 | f(2,1,2,0,0,1,2) = 358 | f(0,2,2,0,0,1,2) = 457 |
| f(1,2,2,0,0,1,2) = 112 | f(2,2,2,0,0,1,2) = 357 | f(0,0,0,1,0,1,2) = 674 |
| f(1,0,0,1,0,1,2) = 1008 | f(2,0,0,1,0,1,2) = 580 | f(0,1,0,1,0,1,2) = 430 |
| f(1,1,0,1,0,1,2) = 446 | f(2,1,0,1,0,1,2) = 543 | f(0,2,0,1,0,1,2) = 238 |
| f(1,2,0,1,0,1,2) = 637 | f(2,2,0,1,0,1,2) = 159 | f(0,0,1,1,0,1,2) = 694 |
| f(2,0,1,1,0,1,2) = 1190 | f(0,1,1,1,0,1,2) = 670 | f(1,1,1,1,0,1,2) = 1706 |
| f(2,1,1,1,0,1,2) = 1571 | f(0,2,1,1,0,1,2) = 704 | f(1,2,1,1,0,1,2) = 1705 |
| f(2,2,1,1,0,1,2) = 1187 | f(2,0,2,1,0,1,2) = 1301 | f(0,1,2,1,0,1,2) = 350 |
| f(1,1,2,1,0,1,2) = 1370 | f(2,1,2,1,0,1,2) = 1299 | f(0,2,2,1,0,1,2) = 448 |
| f(1,2,2,1,0,1,2) = 1369 | f(2,2,2,1,0,1,2) = 1107 | f(0,0,0,2,0,1,2) = 621 |
| f(1,0,0,2,0,1,2) = 1032 | f(2,0,0,2,0,1,2) = 387 | f(0,1,0,2,0,1,2) = 606 |
| f(1,1,0,2,0,1,2) = 1074 | f(2,1,0,2,0,1,2) = 291 | f(0,2,0,2,0,1,2) = 222 |
| f(1,2,0,2,0,1,2) = 1073 | f(2,2,0,2,0,1,2) = 99 | f(0,0,1,2,0,1,2) = 438 |
| f(1,0,1,2,0,1,2) = 1577 | f(2,0,1,2,0,1,2) = 671 | f(0,1,1,2,0,1,2) = 905 |
| f(1,1,1,2,0,1,2) = 552 | f(2,1,1,2,0,1,2) = 900 | f(1,2,1,2,0,1,2) = 168 |
| f(2,2,1,2,0,1,2) = 708 | f(1,0,2,2,0,1,2) = 1114 | f(2,0,2,2,0,1,2) = 367 |
| f(0,1,2,2,0,1,2) = 841 | f(1,1,2,2,0,1,2) = 280 | f(2,1,2,2,0,1,2) = 836 |
| f(1,2,2,2,0,1,2) = 88 | f(2,2,2,2,0,1,2) = 452 | f(1,0,0,0,1,1,2) = 498 |
| f(2,0,0,0,1,1,2) = 450 | f(0,1,0,0,1,1,2) = 821 | f(1,1,0,0,1,1,2) = 1562 |
| f(2,1,0,0,1,1,2) = 1614 | f(1,2,0,0,1,1,2) = 1177 | f(2,2,0,0,1,1,2) = 1421 |
| f(0,0,1,0,1,1,2) = 494 | f(1,0,1,0,1,1,2) = 428 | f(2,0,1,0,1,1,2) = 1734 |
| f(0,1,1,0,1,1,2) = 548 | f(1,1,1,0,1,1,2) = 1588 | f(2,1,1,0,1,1,2) = 1926 |
| f(0,2,1,0,1,1,2) = 173 | f(1,2,1,0,1,1,2) = 1204 | f(2,2,1,0,1,1,2) = 1733 |
| f(1,0,2,0,1,1,2) = 412 | f(2,0,2,0,1,1,2) = 1861 | f(0,1,2,0,1,1,2) = 292 |
| f(1,1,2,0,1,1,2) = 1332 | f(2,1,2,0,1,1,2) = 1862 | f(0,2,2,0,1,1,2) = 109 |
| f(1,2,2,0,1,1,2) = 1140 | f(2,2,2,0,1,1,2) = 1477 | f(0,0,0,1,1,1,2) = 1482 |
| f(1,0,0,1,1,1,2) = 632 | f(2,0,0,1,1,1,2) = 454 | f(0,1,0,1,1,1,2) = 396 |
| f(1,1,0,1,1,1,2) = 634 | f(2,1,0,1,1,1,2) = 390 | f(0,2,0,1,1,1,2) = 1261 |
| f(1,2,0,1,1,1,2) = 441 | f(2,2,0,1,1,1,2) = 581 | f(0,0,1,1,1,1,2) = 1198 |
| f(1,0,1,1,1,1,2) = 1754 | f(2,0,1,1,1,1,2) = 706 | f(0,1,1,1,1,1,2) = 926 |
| f(1,1,1,1,1,1,2) = 1722 | f(2,1,1,1,1,1,2) = 1927 | f(0,2,1,1,1,1,2) = 164 |
| f(1,2,1,1,1,1,2) = 1721 | f(2,2,1,1,1,1,2) = 1735 | f(0,0,2,1,1,1,2) = 610 |
| f(1,0,2,1,1,1,2) = 1897 | f(2,0,2,1,1,1,2) = 833 | f(0,1,2,1,1,1,2) = 862 |
| f(1,1,2,1,1,1,2) = 1402 | f(2,1,2,1,1,1,2) = 1863 | f(0,2,2,1,1,1,2) = 100 |
| f(1,2,2,1,1,1,2) = 1401 | f(2,2,2,1,1,1,2) = 1479 | f(0,0,0,2,1,1,2) = 612 |
| f(2,0,0,2,1,1,2) = 395 | f(0,1,0,2,1,1,2) = 588 | f(1,1,0,2,1,1,2) = 1650 |
| f(2,1,0,2,1,1,2) = 769 | f(0,2,0,2,1,1,2) = 1245 | f(1,2,0,2,1,1,2) = 1457 |
| f(2,2,0,2,1,1,2) = 194 | f(0,0,1,2,1,1,2) = 1182 | f(2,0,1,2,1,1,2) = 655 |
| f(0,1,1,2,1,1,2) = 557 | f(1,1,1,2,1,1,2) = 1580 | f(2,1,1,2,1,1,2) = 614 |
| f(1,2,1,2,1,1,2) = 1196 | f(2,2,1,2,1,1,2) = 421 | f(0,0,2,2,1,1,2) = 320 |
| f(2,0,2,2,1,1,2) = 335 | f(0,1,2,2,1,1,2) = 301 | f(1,1,2,2,1,1,2) = 1308 |
| f(2,1,2,2,1,1,2) = 406 | f(1,2,2,2,1,1,2) = 1116 | f(2,2,2,2,1,1,2) = 597 |
| f(1,0,0,0,2,1,2) = 945 | f(2,0,0,0,2,1,2) = 897 | f(0,1,0,0,2,1,2) = 819 |
| f(1,1,0,0,2,1,2) = 442 | f(2,1,0,0,2,1,2) = 531 | f(1,2,0,0,2,1,2) = 633 |
| f(2,2,0,0,2,1,2) = 147 | f(0,0,1,0,2,1,2) = 429 | f(1,0,1,0,2,1,2) = 754 |
| f(2,0,1,0,2,1,2) = 679 | f(0,1,1,0,2,1,2) = 546 | f(1,1,1,0,2,1,2) = 809 |
| f(2,1,1,0,2,1,2) = 38 | f(0,2,1,0,2,1,2) = 171 | f(1,2,1,0,2,1,2) = 234 |
| f(2,2,1,0,2,1,2) = 37 | f(1,0,2,0,2,1,2) = 881 | f(2,0,2,0,2,1,2) = 343 |
| f(0,1,2,0,2,1,2) = 290 | f(1,1,2,0,2,1,2) = 793 | f(2,1,2,0,2,1,2) = 22 |
| f(0,2,2,0,2,1,2) = 107 | f(1,2,2,0,2,1,2) = 218 | f(2,2,2,0,2,1,2) = 21 |
| f(0,0,0,1,2,1,2) = 683 | f(1,0,0,1,2,1,2) = 968 | f(2,0,0,1,2,1,2) = 1607 |
| f(0,1,0,1,2,1,2) = 414 | f(1,1,0,1,2,1,2) = 820 | f(2,1,0,1,2,1,2) = 1806 |
| f(0,2,0,1,2,1,2) = 235 | f(1,2,0,1,2,1,2) = 244 | f(2,2,0,1,2,1,2) = 1229 |
| f(0,0,1,1,2,1,2) = 649 | f(1,0,1,1,2,1,2) = 1178 | f(2,0,1,1,2,1,2) = 742 |
| f(0,1,1,1,2,1,2) = 542 | f(1,1,1,1,2,1,2) = 1833 | f(2,1,1,1,2,1,2) = 567 |
| f(0,2,1,1,2,1,2) = 162 | f(1,2,1,1,2,1,2) = 1258 | f(2,2,1,1,2,1,2) = 183 |
| f(0,0,2,1,2,1,2) = 43 | f(1,0,2,1,2,1,2) = 1321 | f(2,0,2,1,2,1,2) = 853 |
| f(0,1,2,1,2,1,2) = 286 | f(1,1,2,1,2,1,2) = 1817 | f(2,1,2,1,2,1,2) = 311 |
| f(0,2,2,1,2,1,2) = 98 | f(1,2,2,1,2,1,2) = 1242 | f(2,2,2,1,2,1,2) = 119 |
| f(0,0,0,2,2,1,2) = 929 | f(1,0,0,2,2,1,2) = 122 | f(2,0,0,2,2,1,2) = 587 |
| f(0,1,0,2,2,1,2) = 622 | f(1,1,0,2,2,1,2) = 1898 | f(2,1,0,2,2,1,2) = 1990 |
| f(0,2,0,2,2,1,2) = 219 | f(1,2,0,2,2,1,2) = 1513 | f(2,2,0,2,2,1,2) = 1989 |
| f(0,0,1,2,2,1,2) = 989 | f(1,0,1,2,2,1,2) = 620 | f(2,0,1,2,2,1,2) = 1695 |
| f(0,1,1,2,2,1,2) = 555 | f(1,1,1,2,2,1,2) = 564 | f(2,1,1,2,2,1,2) = 1550 |
| f(1,2,1,2,2,1,2) = 180 | f(2,2,1,2,2,1,2) = 1165 | f(0,0,2,2,2,1,2) = 27 |
| f(1,0,2,2,2,1,2) = 604 | f(2,0,2,2,2,1,2) = 1391 | f(0,1,2,2,2,1,2) = 299 |
| f(1,1,2,2,2,1,2) = 308 | f(2,1,2,2,2,1,2) = 1294 | f(1,2,2,2,2,1,2) = 116 |
| f(2,2,2,2,2,1,2) = 1101 | f(0,0,0,0,2,2) = 969 | f(1,0,0,0,2,2) = 444 |
| f(2,0,0,0,2,2) = 975 | f(0,1,0,0,2,2) = 777 | f(1,1,0,0,2,2) = 1082 |
| f(2,1,0,0,2,2) = 1327 | f(0,2,0,0,2,2) = 201 | f(1,2,0,0,2,2) = 1081 |
| f(2,2,0,0,2,2) = 1135 | f(0,0,1,0,2,2) = 140 | f(1,0,1,0,2,2) = 762 |
| f(2,0,1,0,2,2) = 687 | f(0,1,1,0,2,2) = 950 | f(1,1,1,0,2,2) = 1726 |
| f(2,1,1,0,2,2) = 1839 | f(0,2,1,0,2,2) = 767 | f(1,2,1,0,2,2) = 1725 |
| f(2,2,1,0,2,2) = 1263 | f(0,0,2,0,2,2) = 603 | f(1,0,2,0,2,2) = 889 |
| f(2,0,2,0,2,2) = 351 | f(0,1,2,0,2,2) = 886 | f(1,1,2,0,2,2) = 1406 |
| f(2,1,2,0,2,2) = 1823 | f(0,2,2,0,2,2) = 511 | f(1,2,2,0,2,2) = 1405 |
| f(2,2,2,0,0,2,2) = 1247 | f(0,0,0,1,0,2,2) = 658 | f(1,0,0,1,0,2,2) = 1020 |
| f(0,1,0,1,0,2,2) = 814 | f(1,1,0,1,0,2,2) = 536 | f(2,1,0,1,0,2,2) = 294 |
| f(0,2,0,1,0,2,2) = 993 | f(1,2,0,1,0,2,2) = 152 | f(2,2,0,1,0,2,2) = 101 |
| f(0,0,1,1,0,2,2) = 1069 | f(1,0,1,1,0,2,2) = 154 | f(2,0,1,1,0,2,2) = 419 |
| f(0,1,1,1,0,2,2) = 652 | f(1,1,1,1,0,2,2) = 955 | f(2,1,1,1,0,2,2) = 899 |
| f(0,2,1,1,0,2,2) = 758 | f(1,2,1,1,0,2,2) = 763 | f(2,2,1,1,0,2,2) = 707 |
| f(0,0,2,1,0,2,2) = 1569 | f(1,0,2,1,0,2,2) = 297 | f(2,0,2,1,0,2,2) = 403 |
| f(0,1,2,1,0,2,2) = 332 | f(1,1,2,1,0,2,2) = 891 | f(2,1,2,1,0,2,2) = 835 |
| f(0,2,2,1,0,2,2) = 502 | f(1,2,2,1,0,2,2) = 507 | f(2,2,2,1,0,2,2) = 451 |
| f(0,0,0,2,0,2,2) = 605 | f(1,0,0,2,0,2,2) = 1019 | f(2,0,0,2,0,2,2) = 399 |
| f(0,1,0,2,0,2,2) = 798 | f(1,1,0,2,0,2,2) = 874 | f(2,1,0,2,0,2,2) = 1794 |
| f(0,2,0,2,0,2,2) = 977 | f(1,2,0,2,0,2,2) = 489 | f(2,2,0,2,0,2,2) = 1217 |
| f(0,0,1,2,0,2,2) = 1053 | f(1,0,1,2,0,2,2) = 40 | f(2,0,1,2,0,2,2) = 47 |
| f(0,1,1,2,0,2,2) = 959 | f(1,1,1,2,0,2,2) = 666 | f(2,1,1,2,0,2,2) = 559 |
| f(0,2,1,2,0,2,2) = 1681 | f(1,2,1,2,0,2,2) = 665 | f(2,2,1,2,0,2,2) = 175 |
| f(0,0,2,2,0,2,2) = 1297 | f(1,0,2,2,0,2,2) = 24 | f(2,0,2,2,0,2,2) = 31 |
| f(0,1,2,2,0,2,2) = 895 | f(1,1,2,2,0,2,2) = 362 | f(2,1,2,2,0,2,2) = 287 |
| f(0,2,2,2,0,2,2) = 1361 | f(1,2,2,2,0,2,2) = 361 | f(2,2,2,2,0,2,2) = 95 |
| f(1,0,0,0,1,2,2) = 510 | f(2,0,0,0,1,2,2) = 462 | f(0,1,0,0,1,2,2) = 12 |
| f(1,1,0,0,1,2,2) = 1849 | f(2,1,0,0,1,2,2) = 1318 | f(0,2,0,0,1,2,2) = 1077 |
| f(1,2,0,0,1,2,2) = 1274 | f(2,2,0,0,1,2,2) = 1125 | f(0,0,1,0,1,2,2) = 734 |
| f(1,0,1,0,1,2,2) = 700 | f(2,0,1,0,1,2,2) = 718 | f(0,1,1,0,1,2,2) = 532 |
| f(1,1,1,0,1,2,2) = 1576 | f(2,1,1,0,1,2,2) = 1574 | f(0,2,1,0,1,2,2) = 157 |
| f(1,2,1,0,1,2,2) = 1192 | f(2,2,1,0,1,2,2) = 1189 | f(1,0,2,0,1,2,2) = 380 |
| f(2,0,2,0,1,2,2) = 845 | f(0,1,2,0,1,2,2) = 276 | f(1,1,2,0,1,2,2) = 1304 |
| f(2,1,2,0,1,2,2) = 1302 | f(0,2,2,0,1,2,2) = 93 | f(1,2,2,0,1,2,2) = 1112 |
| f(2,2,2,0,1,2,2) = 1109 | f(0,0,0,1,1,2,2) = 110 | f(1,0,0,1,1,2,2) = 52 |
| f(2,0,0,1,1,2,2) = 964 | f(0,1,0,1,1,2,2) = 813 | f(1,1,0,1,1,2,2) = 1086 |
| f(2,1,0,1,1,2,2) = 823 | f(0,2,0,1,1,2,2) = 437 | f(1,2,0,1,1,2,2) = 1085 |
| f(2,2,0,1,1,2,2) = 247 | f(0,0,1,1,1,2,2) = 36 | f(1,0,1,1,1,2,2) = 668 |
| f(2,0,1,1,1,2,2) = 675 | f(0,1,1,1,1,2,2) = 908 | f(1,1,1,1,1,2,2) = 562 |
| f(2,1,1,1,1,2,2) = 1547 | f(0,2,1,1,1,2,2) = 148 | f(1,2,1,1,1,2,2) = 177 |
| f(2,2,1,1,1,2,2) = 1163 | f(0,0,2,1,1,2,2) = 1663 | f(1,0,2,1,1,2,2) = 364 |
| f(2,0,2,1,1,2,2) = 339 | f(0,1,2,1,1,2,2) = 844 | f(1,1,2,1,1,2,2) = 306 |
| f(2,1,2,1,1,2,2) = 1291 | f(0,2,2,1,1,2,2) = 84 | f(1,2,2,1,1,2,2) = 113 |
| f(2,2,2,1,1,2,2) = 1099 | f(0,0,0,2,1,2,2) = 596 | f(1,0,0,2,1,2,2) = 1977 |
| f(2,0,0,2,1,2,2) = 439 | f(0,1,0,2,1,2,2) = 797 | f(1,1,0,2,1,2,2) = 434 |
| f(2,1,0,2,1,2,2) = 1319 | f(0,2,0,2,1,2,2) = 629 | f(1,2,0,2,1,2,2) = 625 |
| f(2,2,0,2,1,2,2) = 1127 | f(0,0,1,2,1,2,2) = 20 | f(1,0,1,2,1,2,2) = 106 |
| f(2,0,1,2,1,2,2) = 623 | f(0,1,1,2,1,2,2) = 541 | f(1,1,1,2,1,2,2) = 1002 |
| f(2,1,1,2,1,2,2) = 934 | f(0,2,1,2,1,2,2) = 1745 | f(1,2,1,2,1,2,2) = 1001 |
| f(2,2,1,2,1,2,2) = 741 | f(0,0,2,2,1,2,2) = 1042 | f(1,0,2,2,1,2,2) = 537 |
| f(2,0,2,2,1,2,2) = 607 | f(0,1,2,2,1,2,2) = 285 | f(1,1,2,2,1,2,2) = 986 |
| f(2,1,2,2,1,2,2) = 854 | f(0,2,2,2,1,2,2) = 1489 | f(1,2,2,2,1,2,2) = 985 |

f(2,2,2,2,1,2,2) = 469
f(0,1,0,0,2,2,2) = 10
f(0,2,0,0,2,2,2) = 1075
f(0,0,1,0,2,2,2) = 669
f(0,1,1,0,2,2,2) = 530
f(2,2,1,0,2,2,2) = 135
f(0,1,2,0,2,2,2) = 274
f(2,2,2,0,2,2,2) = 71
f(2,0,0,1,2,2,2) = 631
f(2,1,0,1,2,2,2) = 303
f(2,2,0,1,2,2,2) = 111
f(2,0,1,1,2,2,2) = 431
f(2,1,1,1,2,2,2) = 951
f(2,2,1,1,2,2,2) = 759
f(2,0,2,1,2,2,2) = 415
f(2,1,2,1,2,2,2) = 887
f(2,2,2,1,2,2,2) = 503
f(2,0,0,2,2,2,2) = 7
f(2,1,0,2,2,2,2) = 871
f(2,2,0,2,2,2,2) = 487
f(2,0,1,2,2,2,2) = 1071
f(2,1,1,2,2,2,2) = 527
f(2,2,1,2,2,2,2) = 143
f(2,0,2,2,2,2,2) = 1055
f(2,1,2,2,2,2,2) = 271
f(2,2,2,2,2,2,2) = 79

f(1,0,0,0,2,2,2) = 957
f(1,1,0,0,2,2,2) = 827
f(1,2,0,0,2,2,2) = 251
f(1,0,1,0,2,2,2) = 1723
f(2,1,1,0,2,2,2) = 519
f(1,0,2,0,2,2,2) = 1403
f(2,1,2,0,2,2,2) = 263
f(0,0,0,1,2,2,2) = 667
f(0,1,0,1,2,2,2) = 811
f(0,2,0,1,2,2,2) = 33
f(0,0,1,1,2,2,2) = 45
f(0,1,1,1,2,2,2) = 524
f(0,2,1,1,2,2,2) = 146
f(0,0,2,1,2,2,2) = 563
f(0,1,2,1,2,2,2) = 268
f(0,2,2,1,2,2,2) = 82
f(0,0,0,2,2,2,2) = 309
f(0,1,0,2,2,2,2) = 795
f(0,2,0,2,2,2,2) = 17
f(0,0,1,2,2,2,2) = 29
f(0,1,1,2,2,2,2) = 539
f(0,2,1,2,2,2,2) = 1169
f(0,0,2,2,2,2,2) = 307
f(0,1,2,2,2,2,2) = 283
f(0,2,2,2,2,2,2) = 1105

f(2,0,0,0,2,2,2) = 909
f(2,1,0,0,2,2,2) = 295
f(2,2,0,0,2,2,2) = 103
f(2,0,1,0,2,2,2) = 35
f(0,2,1,0,2,2,2) = 155
f(2,0,2,0,2,2,2) = 19
f(0,2,2,0,2,2,2) = 91
f(1,0,0,1,2,2,2) = 1012
f(1,1,0,1,2,2,2) = 825
f(1,2,0,1,2,2,2) = 250
f(1,0,1,1,2,2,2) = 937
f(1,1,1,1,2,2,2) = 42
f(1,2,1,1,2,2,2) = 41
f(1,0,2,1,2,2,2) = 474
f(1,1,2,1,2,2,2) = 26
f(1,2,2,1,2,2,2) = 25
f(1,0,0,2,2,2,2) = 635
f(1,1,0,2,2,2,2) = 1851
f(1,2,0,2,2,2,2) = 1275
f(1,0,1,2,2,2,2) = 766
f(1,1,1,2,2,2,2) = 574
f(1,2,1,2,2,2,2) = 189
f(1,0,2,2,2,2,2) = 893
f(1,1,2,2,2,2,2) = 318
f(1,2,2,2,2,2,2) = 125

List 6